\documentclass[10pt,a4paper,landscape,twocolumn]{article}\hoffset = -20pt\evensidemargin = 3pt\oddsidemargin = 0pt\textwidth = 260mm\marginparsep = 0pt\marginparwidth = 0pt\voffset = -54pt\textheight = 187mm\topmargin = -20pt\headheight = 0pt\headsep = 25pt\footskip = 15pt\columnsep = 80pt\columnseprule = 1pt

\usepackage{srcltx,enumerate} 
\newcommand{\lb}[1]{\label{#1}}\newcommand{\reff}{\ref}\newcommand{\rf}[1]{(\reff{#1})}\newcommand{\cit}[1]{\cite{#1}}\newcommand{\bibi}[1]{\bibitem{#1}}

\usepackage{pstricks,pst-plot}
\newcommand{\figu}[3]{ \begin{figure}[ht]\begin{center} #2 \end{center}\caption{\lb{#1}#3.}\end{figure}}
\usepackage{amsmath,theorem}
\usepackage{epsfig}
\usepackage{amssymb}
\usepackage[latin1]{inputenc}
\usepackage[T1]{fontenc}
\usepackage[francais]{babel}
\newcommand{\fp}{\fontfamily{ppl}\selectfont\small}
\newcommand{\ft}{\fontfamily{cmr}\selectfont\normalsize}
\numberwithin{equation}{section}
\numberwithin{figure}{section}
\newtheorem{lm}{{\bf Lemme}}[section]
\newtheorem{theor}[lm]{{\bf Théorème}}
\newtheorem{deff}[lm]{{\bf Définition}}
\newtheorem{deffs}[lm]{{\bf Définitions}}
\newtheorem{cj}[lm]{{\bf Conjecture}}
\newtheorem{rem}[lm]{{\bf Remarque}}
\newtheorem{cor}[lm]{{\bf Corollaire}}
\newtheorem{exemple}[lm]{{\bf Exemple}}
\newtheorem{prop}[lm]{{\bf Proposition}}
\newtheorem{propodef}[lm]{{\bf Proposition-Définition}}
{\theorembodyfont{\rmfamily}

}
\newcommand{\sep}{$\!${\rm.}\ }
\renewcommand{\rq}{
\fontfamily{cmr}
\selectfont\normalsize{\noindent\sc Remarque} \sep}
\newcommand{\rqs}{
\fontfamily{cmr}
\selectfont\normalsize{\noindent\sc Remarques} \sep}

\newcommand{\df}[2]{\begin{deff}{\lb{#1}}\sep{\rm #2}\end{deff}}

\newcommand{\lem}[2]{\begin{lm}{\lb{#1}}\sep{\sl #2}\end{lm}}
\newcommand{\theo}[2]{\begin{theor}{\lb{#1}}\sep{\sl 
#2}\end{theor}}

\newcommand{\coro}[2]{\begin{cor}{\lb{#1}}\sep{\sl #2}\end{cor}}

\newcommand{\propo}[2]{\begin{prop}{\lb{#1}}\sep{\sl #2}\end{prop}}

\renewcommand{\sec}[2]{\section{#2.}{\lb{#1}}}
\newcommand{\sub}[2]{\subsection{#2.}{\lb{#1}}}
\newcommand{\ssub}[2]{\subsubsection{#2.}{\lb{#1}}}
\newcommand{\eq}[2]{\begin{equation}#2\ \ \  \ \lb{#1}\end{equation}}
\newcommand{\rff}[2]{{(\reff{#1},\reff{#2})}}
\newcommand{\rft}[2]{{(\reff{#1}--\reff{#2})}}
\newcommand{\bul}{\medskip\noindent$\bullet$ \ }
\newcommand{\pr}[1]{\fontfamily{ppl}\fp{\selectfont\normalsize\noindent {\sl#1\;\sep }}}
\newcommand{\ep}{\ft\hfill \framebox[2mm]{\ } \medskip}
\newcommand{\tq}{\ ;\ }

\newcommand{\be}{\begin{enumerate}}
\newcommand{\ee}{\end{enumerate}}
\newcommand{\noi}{\noindent}
\newcommand{\med}{\medskip}
\newcommand{\ts}{\textstyle}
\newcommand{\ds}{\displaystyle}
\newcommand{\lp}{\left(}
\newcommand{\rp}{\right)}

\newcommand{\norm}[1]{\left|#1\right|}
\newcommand{\dnorm}[1]{\left|\left|#1\right|\right|}
\newcommand{\gk}[1]{\left(#1\right)}

\renewcommand{\~}{\widetilde}
\renewcommand{\hat}{\widehat}

\renewcommand{\a}{\alpha}

\renewcommand{\b}{\beta}
\renewcommand{\d}{\delta}
\newcommand{\D}{{\cal D}}
\newcommand{\De}{\Delta}
\newcommand{\eps}{\varepsilon}
\newcommand{\e}{\eta}
\newcommand{\f}{\varphi}
\newcommand{\g}{\gamma}
\newcommand{\G}{{\cal G}}
\renewcommand{\H}{{\cal H}}
\newcommand{\J}{{\cal J}}
\renewcommand{\L}{{\cal L}}
\renewcommand{\l}{\lambda}
\newcommand{\M}{{\cal M}}
\newcommand{\E}{{\cal E}}

\newcommand{\Pc}{{\widehat P}}
\newcommand{\Phic}{\widehat{\Phi}}
\renewcommand{\P}{\widehat{\bf Q}}
\renewcommand{\O}{{\cal O}}

\newcommand{\s}{\sigma}
\renewcommand{\S}{{\bf S}}

\newcommand{\T}{{\bf T}}
\newcommand{\U}{{\cal U}}
\newcommand{\vu}{{\bf v}}
\newcommand{\V}{{\cal V}}
\newcommand{\w}{\omega}
\newcommand{\W}{\Omega}
\newcommand{\N}{\mathbb{N}}
\newcommand{\Z}{\mathbb{Z}}
\newcommand{\R}{\mathbb{R}}
\newcommand{\C}{\mathbb{C}}

\newcommand{\uc}{{\widehat u}}
\newcommand{\vc}{{\widehat v}}
\newcommand{\ac}{{\widehat a}}

\newcommand{\yc}{{\widehat y}}
\newcommand{\zc}{{\widehat z}}

\newcommand{\cc}{{\cal C}}
\newcommand{\cch}{{\widehat{\cal C}}}

\newcommand{\cl}{{\rm cl}}
\newcommand{\re}{{\rm Re}\,}
\newcommand{\im}{{\rm Im}\,}

\newcommand{\val}{{\rm val}}

\newcommand{\apriori}{{\sl a priori}}
\newcommand{\Apriori}{{\sl A priori}}
\newcommand{\cf}{{\sl c.f. }}
\newcommand{\ie}{{\sl i.e. }}
\newcommand{\usp}{\frac1p}

\newcommand{\ed}{\end{document}}

\newcommand{\dac}{{\sc dac}}
\newcommand{\Dacs}{\dac}
\newcommand{\dacs}{\dac}

\title{Développements asymptotiques combinés et\\
points tournants d'équations différentielles\\singulièrement perturbées}
\author{A. {\sc Fruchard} et R. {\sc Schäfke}}
\begin{document}
\begin{center}
{\Large Développements asymptotiques combinés et\\
points tournants d'équations différentielles\\
singulièrement perturbées\\\ \\}
A. {\sc Fruchard} et R. {\sc Schäfke}, \today
\end{center}
%
%
%
\sec{0.}{Résumé}
L'objet de ce  mémoire est d'élaborer une théorie de développements asymptotiques pour des fonctions de deux variables, qui combinent à la fois des fonctions d'une des deux variables et  des fonctions du quotient de ces deux variables. Ces {\sl développements asymptotiques combinés} (\dacs{} en abrégé) se révèlent particulièrement bien adaptés à la description des solutions d'équations  différentielles ordinaires  singulièrement perturbées au voisinage de points tournants. Décrivons le contexte en quelques mots.
Etant donnée une équation de la forme
\eq{fb}{
\eps y'=\Phi(x,y,\eps)
}
où $\Phi$ est de classe $\cc^\infty$, $x$ et $y$ sont des variables réelles ou complexes et $\eps$ est un petit paramètre, réel positif ou dans un secteur du plan complexe, nous appelons {\sl point tournant} un point $(x^*,y^*)$ de la {\sl variété lente} $\L$ d'équation $\Phi(x^*,y^*,0)=0$, tel que $\frac{\partial \Phi}{\partial y}(x^*,y^*,0)=0$.

Au voisinage d'un {\sl point régulier} $(\~x,\~y)$ de $\L$, \ie tel que $\frac{\partial \Phi}{\partial y}(\~x,\~y,0)\neq0$,
on vérifie facilement que \rf{fb} a une unique solution formelle $\yc=\sum_{n\geq0}y_n(x)\eps^n$ et il est connu que les solutions de \rf{fb} admettent $\yc$ pour développement asymptotique dans des domaines adéquats.
La théorie {\sl classique} des développements 
combinés 
permet aussi de décrire la {\sl couche limite} (appelée aussi la couche intérieure) d'une solution ayant en $\~x$ une valeur initiale assez proche de $\~y$. Par exemple  dans le cadre réel, si le point $(\~x,\~y)$ est {\sl attractif}, \ie $\frac{\partial \Phi}{\partial y}(\~x,\~y,0)<0$, on peut donner une approximation d'une solution $y=y(x,\eps)$, uniforme sur un intervalle $[\~x,\~x+\d]$, sous la forme 
$\ts\sum_{n\geq0}\Big(y_n(x)+z_n\big(\ts\frac{x-\~x}\eps\big)\Big)\eps^n$, comportant d'une part la solution formelle $\yc$ et d'autre part des fonctions $z_n$ à  décroissance exponentielle à  l'infini.

En un point tournant $(x^*,y^*)$,  en général les coefficients de la solution formelle présentent des singularités de type pôle ou ramification et cette méthode des développements combinés classiques n'est plus applicable. La méthode la plus répandue pour obtenir une approximation des solutions est le recollement de deux développements dits {\sl intérieur} et {\sl extérieur}.
C'est ce qu'on appelle le {\sl matching} dans la littérature anglo-saxonne. Plus qu'une méthode,  le matching est une idée très générale, qui regroupe des méthodes diverses dans de nombreuses situations où apparaissent à  des équations fonctionnelles (différentielles ordinaires, aux dérivées partielles, etc.)

L'objet du mémoire est de présenter de nouveaux  \dacs, de la forme 
$
\sum_{n\geq0}$ $\Big(a_n(x)+$ $\!g_n\big(\ts\frac {x-x^*}\e\big)\Big)\e^n.
$
Dans le cadre des équations de la forme \rf{fb}, 
la variable $\e$ est une racine du petit paramètre $\eps$ et 
les fonctions $a_n$ (la partie dite {\sl lente} du \dac) sont liées aux fonctions $y_n$ de la solution formelle $\yc$ mais ne leur sont pas égales~: ce sont les parties régulières de ces fonctions $y_n$. 
De manière symétrique, la partie {\sl rapide} du \dac, constituée des fonctions $g_n$, correspond à la partie régulière à l'infini de la solution formelle de l'{\sl équation intérieure} obtenue par le changement de variable $x=\eta X$.
L'avantage des développements combinés est de donner une approximation uniforme des solutions dans une région contenant à  la fois des points loin du point tournant et des points dans un petit voisinage de ce point tournant.

Nous présentons deux méthodes pour obtenir ces \dacs. La première est une méthode directe selon le modèle classique~: étude de solution formelles, puis preuve d'existence d'une solution analytique ayant cette solution formelle pour \dac.
La deuxième méthode est plus indirecte et utilise un résultat de type Ramis-Sibuya~: 
nous montrons d'abord qu'il existe des solutions de \rf{fb}
pour $\eps$ et pour $x$ dans des secteurs formant des recouvrements de l'origine, 
puis que ces solutions satisfont des estimations exponentielles. 
Cette méthode complexe peut sembler bien sophistiquée à  
un lecteur habitué au contexte des équations 
dans le champ réel, mais elle est la plus simple à  notre connaissance pour la perturbation singulière. 
Par ailleurs, cette construction 
par des recouvrements a le mérite de fournir 
gratuitement des estimations de type Gevrey, non seulement pour le développement   mais aussi pour les développements asymptotiques 
des fonctions $g_n$. Le lien étroit entre la perturbation singulière et la théorie Gevrey est connu~; parmi les trois applications que nous présentons, deux utilisent de manière incontournable le Gevrey. 

En dépit du caractère naturel, presque familier, de ces développements, nous n'en avons trouvé presque aucune trace 
 dans la littérature existante. En particulier leur similitude apparente avec les \dacs{} classiques cache des différences profondes. La comparaison de nos résultats avec l'abondante littérature sur les développements combinés classiques et le matching sera mentionnée au fur et à mesure du mémoire,
ainsi qu'à la fin.
Dans ce mémoire, l'accent est mis sur les équations différentielles ordinaires, mais nous sommes convaincus que ces \dacs~peuvent se révéler très utiles pour d'autres types d'équations fonctionnelles.
\med\\
\noindent{\bf Mots-clés :} 
point tournant, développement asymptotique combiné, série Gevrey, canard, perturbation singulière, équation différentielle complexe.
\med\\
\noindent{\bf Classification AMS :} \ 
34E,	30E10, 41A60, 34M30, 34M60.
\tableofcontents
%
%
%
%
%
\sec{1.2}{Quatre exemples introductifs}
Nous présentons ici des exemples, les plus simples possibles, qui montrent que les 
solutions d'équations singulièrement perturbées ont naturellement des \dacs{} près des points tournants. 
Les quatre exemples sont des équations linéaires.
Le premier exemple est le plus simple présentant un point tournant. Le deuxième contient un paramètre de contrôle, permettant de \og chasser le canard\fg. Le troisième exemple contient aussi un paramètre de contrôle mais le point tournant n'est plus simple, ce qui a pour conséquence que les canards ne sont pas des solutions surstables
au sens de Guy Wallet \cit{w}. 
Enfin le quatrième exemple concerne des \og faux canards \fg~: la courbe lente est d'abord répulsive puis attractive. Dans cette situation, toute solution de condition initiale bornée au point tournant est définie et bornée sur un intervalle contenant ce point tournant, mais pour que cette solution ait un \dac, il faut au moins que cette condition initiale ait un développement en puissances de $\eta$. Nous verrons que cette condition nécessaire est aussi suffisante.
\sub{1.2.1}{Premier exemple}
Commençons par l'équation 
\eq{1b}{
\eps\frac{dy}{dx}=2xy+\eps g(x),
}
où $\eps>0$ est le petit paramètre et $x,y$ sont des variables réelles. 
On suppose que la fonction $g:\R\to\R$ est de classe $\cc^\infty$ et bornée, ainsi que toutes ses dérivées. 
Ces hypothèses sont là  pour simplifier la présentation, mais beaucoup des résultats qui suivent sont valables avec des hypothèses plus faibles. 
Par exemple l'existence de la solution $y^-$ ci-dessous n'utilise que la continuité de $g$ et pour montrer que cette solution a localement un développement asymptotique à  l'ordre de $N$, il suffit que $g$ soit de classe  $\cc^N$. 

\figu{f2.1}{\epsfxsize5cm\epsfbox{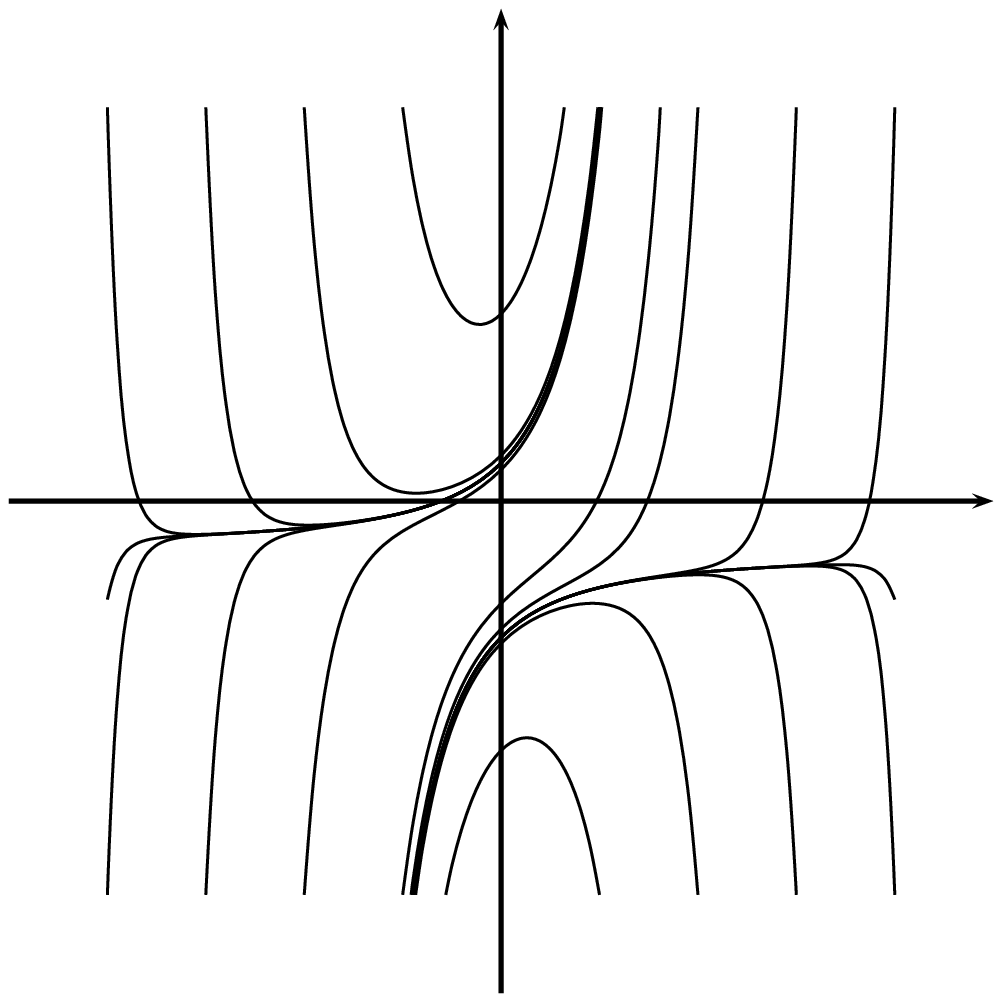}\hspace{1.5cm}\epsfxsize5cm\epsfbox{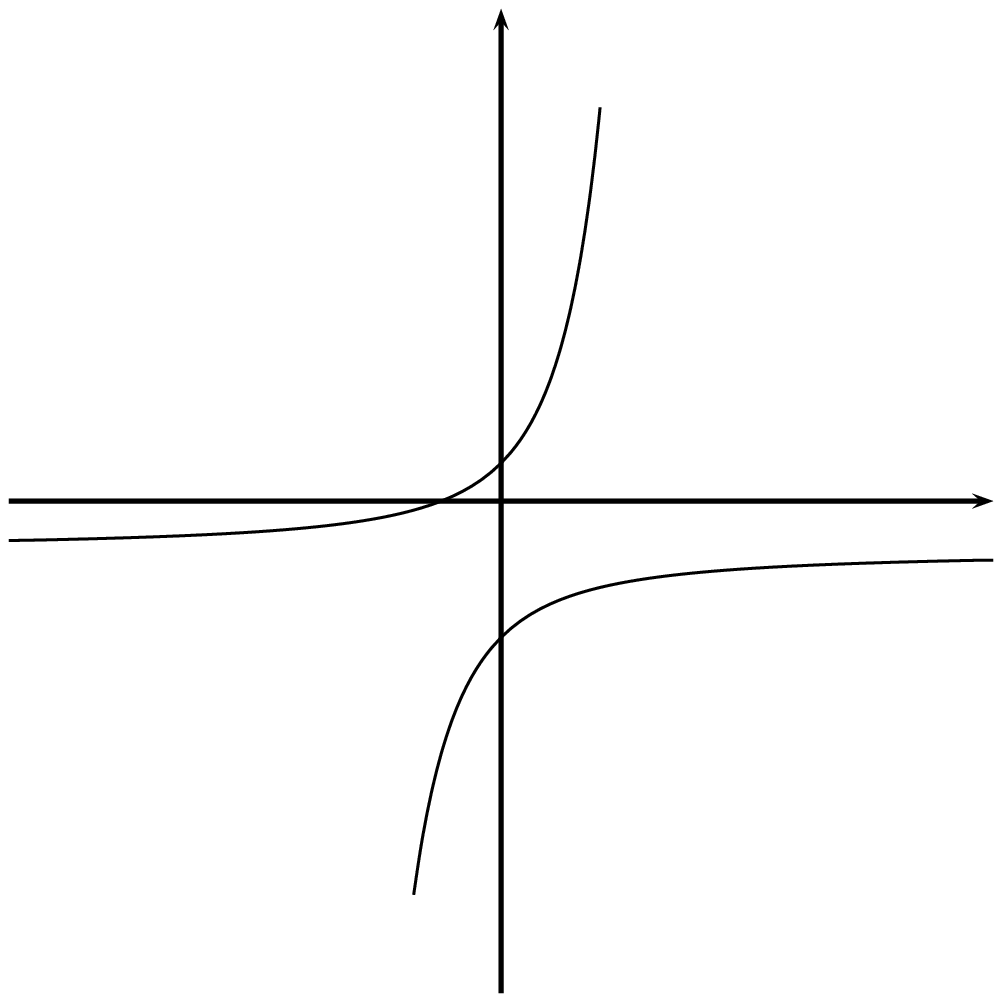}
\vspace{-8mm}
}
{Quelques solutions de \rf{1b} avec $\eps=1$, $|x|\leq4$, $|y|\leq4$ et $g(x)=x+1$.
À droite, les solutions $y^-$ et $y^+$}

Avec des modifications mineures, on peut aussi traiter le cas d'un intervalle fini ou de fonctions non bornées.
Par ailleurs, la présentation de cet exemple pourrait se faire dans le cadre complexe, ce qui correspondrait mieux à  l'esprit du mémoire.
Cependant, nous tenons à  présenter ces exemples dans le cadre réel pour bien illustrer le fait que ces \dacs{} sont aussi très utiles dans ce cadre. Les applications présentées dans la partie \reff{6.} ont d'ailleurs leur origine dans le réel.

Puisque l'équation \rf{1b} est linéaire, ses solutions sont définies sur tout $\R$. 
La courbe lente mentionnée dans l'introduction est ici $y=0$~; elle est attractive pour $x<0$ et répulsive pour $x>0$.
L'équation \rf{1b} présente un point tournant simple en $x=0$.
Pour chaque $\eps>0$ fixé, on peut voir qu'il existe une unique fonction $y^-(\cdot,\eps)$ bornée sur $\R^-$ qui est solution de \rf{1b} pour la valeur $\eps$ du petit paramètre. Cette solution est donnée par la formule de variation de la constante
\eq{3b}{
y^-(x,\eps)=e^{x^2/\eps}\int_{-\infty}^xe^{-t^2/\eps}g(t)dt.
}
Puisque dans ce mémoire $\eps$ est une variable, ce que nous appelons ``solution'' est en fait une famille de solutions dépendant de $\eps$. La formule \rf{3b} définit une famille de solutions de \rf{1b} qui sont non seulement bornées sur $\R$ prises isolément, mais de plus bornées uniformément par rapport à  $\eps$ (ou encore, c'est une fonction  des deux variables $x$ et $\eps$ qui est bornée sur $\R^-\times\,]0,\eps_0]$ pour $\eps_0>0$ fixé).
Dans toute la suite, nous utiliserons parfois l'expression ``bornée'' pour ``bornée uniformément par rapport à  $\eps$''.

Nous voulons avoir une description quantitative de cette solution, non seulement pour $x<0$ mais aussi pour $x$ proche de $0$ et pour $x>0$. Nous commençons avec le cas $x<0$. 

Par une succession d'intégrations par parties, on montre aisément que, pour tout $\d>0$ fixé, la solution $y^-$ admet aussi un développement asymptotique au sens de Poincaré uniforme sur $]-\infty,-\d[$, de la forme $\yc=\sum_{n\geq0}y_n(x)\eps^n$~: pour tout entier $N>0$, il existe une constante $C_N$ telle que pour tout $x\in\,]-\infty,-\d[$ et tout $\eps\in\,]0,\eps_0]$
$$
\bigg|\,y^-(x,\eps)-\sum_{n=0}^{N-1}y_n(x)\eps^n\bigg|\leq C_N\eps^N.
$$
De fait, ce développement est l'unique solution formelle de \rf{1b}~; elle est donnée par les premiers termes
\eq{6b}{
y_0(x)=0,\quad y_{1}(x)=-\frac1{2x}\,g(x),
}
puis récursivement par
\eq{6bis}{
y_{n+1}(x)=\frac1{2x}\,y'_{n}(x),\quad n\geq1.
}
Pour voir que $\yc$ est bien un développement asymptotique de $y^-$, on peut par exemple écrire $y=y^{(N)}+z\eps^N$ avec $y^{(N)}=\sum_{n=0}^{N-1}y_n\eps^n$ et vérifier que $z$ satisfait une équation du même genre que \rf{1b}, donc est bornée sur la demi-droite $]-\infty,-\d[$.

Si on remplace $-\infty$ par $+\infty$, la même formule \rf{3b} fournit aussi une unique solution $y^+$ bornée sur $\R^+$, qui admet un développement asymptotique sur $]\d,+\infty[$ pour tout $\d>0$. Puisque ce développement est l'unique solution formelle de \rf{1b}, c'est le même que celui de $y^-$.

\bul 
Dans le cas très particulier où $g$ est impaire, alors d'une part on a $y^-=y^+$ et d'autre part les formules \rf{6b} et \rf{6bis} 
impliquent que pour tout $n\in\N$ la fonction $y_n$ est paire et sans pôle en $x=0$~; ainsi la solution formelle $\yc$ reste définie en $x=0$. 
Il est donc naturel de se demander si le développement commun $\sum_{n\geq0}y_n(x)\eps^n$, valide pour $x$ loin de $0$, reste valide près de $0$. Dans cet exemple c'est le cas et on peut le montrer en utilisant $y^{(N)}$ comme précédemment. 

Pour des équations analogues, par exemple en changeant $2x$ par $4x^3$ dans \rf{1b}, \cf {2.2.2}, ce résultat n'est plus vrai. 
On a toujours $y^-=y^+$ et donc une solution bornée sur tout $\R$, mais les coefficients 
de la solution formelle 
admettent en général des pôles en 
$x=0$ et les sommes partielles $y^{(N)}$ de $\hat y$ 
ne peuvent plus être des approximations uniformes de $y$. 
L'équation \rf{1b} est l'un des exemples les plus simples où la théorie de la surstabilité peut s'appliquer. Nous ne poursuivons pas plus loin la discussion dans cette direction car nous voulons présenter les \dacs{} et non la surstabilité. 

\bul  
Lorsque $g$ n'est pas impaire, le développement de $y^+$ permet toutefois d'avoir aussi une approximation de $y^-$ sur $]\d,+\infty[$. En effet, on a $y^-(x,\eps)=y^+(x,\eps)+I(\eps)e^{x^2/\eps}$ avec $I(\eps)=\ds\int_{-\infty}^{+\infty}e^{-t^2/\eps}g(t)dt$. 
Si  $g$ n'est pas impaire, alors la fonction $I$ est non nulle. 
Pour voir ceci, on peut écrire $I(\eps)=\int_{0}^{+\infty}e^{-s/\eps}\big(g(\sqrt{s})+
g(-\sqrt{s})\big)\frac1{2\sqrt{s}}ds$ et utiliser l'injectivité de la transformation de 
Laplace.
Plus concrètement, 
si la partie paire de $g$ --- donnée par $g^+(x)=\frac12\big(g(x)+g(-x)\big)$ --- n'est pas plate, alors elle satisfait $g^+(x)\sim Cx^{2N}$ avec $C\neq0$ et $N\in\N$, et obtient $I(\eps)\sim C'\eps^{N+1/2}$.
Par conséquent, en un point fixé $x>0$, $y^-(x,\eps)$ prend une valeur exponentiellement grande par rapport à $\eps$~: il existe $c,a,\eps_0>0$ (qui dépendent de $x$) tels que pour tout $\eps\in\,]0,\eps_0[$, $|y^-(x,\eps)|\geq c\exp\big(\frac a\eps\big)$.
\med

Il est possible de décrire précisément en fonction de $N$ le domaine où $y^-$ reste borné et le domaine où $y^-$ tend vers l'infini, mais nous ne faisons pas une étude exhaustive ici~; nous allons décrire
$y^-(x,\eps)$ en des points $x$ de l'ordre de $\e=\sqrt\eps$, c'est-à-dire lorsque $x$ et $\eps$ tendent vers $0$ avec $\frac x\e$ borné.
Tout d'abord, on peut voir qu'en de tels points $y^-$ reste bornée, et même tend vers 
$0$. En effet le changement de variable $x=\e X$ (avec $\e=\sqrt\eps$) dans \rf{3b} donne
pour tout $K$ réel
\eq{0b}{
y^-(\e X,\eps)=\e\int_{-\infty}^Xe^{X^2-T^2}g(\e T)dT=\O(\e)}
quand $X\leq K$.

Nous allons voir que $y^-$ admet un développement en puissances de $\e$, mettant en jeu à  la fois des fonctions de la variable lente $x$ et de la variable rapide $X=\frac x\e$. Ceci peut se voir par une succession d'intégrations par parties. En effet, notons $Sg$ la fonction définie par 
\eq S{
g(x)=g(0)+xSg(x).
}
Puisque $g$ est $\cc^\infty$ et bornée sur $\R$, $Sg$ l'est aussi (et même tend vers 0 à  l'infini). Une première intégration par parties donne
$$
y^-(x,\eps)=e^{x^2/\eps}
\bigg(g(0)\int_{-\infty}^xe^{-t^2/\eps}dt+\int_{-\infty}^xe^{-t^2/\eps}Sg(t)tdt
\bigg)
$$
$$
=g(0)\e\,U^-\big(\ts\frac x\e\big)-\frac\eps2\,Sg(x)+
\ts\frac\eps2\,e^{x^2/\eps}\ds\int_{-\infty}^xe^{-t^2/\eps}(Sg)'(t)dt.
$$
avec $U^-(X)=e^{X^2}\int_{-\infty}^X e^{-T^2}\,dT.$
En appliquant \rf{0b} à  $(Sg)'$ au lieu de $g$, on a ainsi pour tout $K$ réel
$$
y^-(x,\eps)=g(0)\e U^-\big(\ts\frac x\e\big)-\frac\eps2Sg(x)+\O\big(\e^3\big)
$$
quand $\e\to0$ uniformément sur l'ensemble des $x$ avec $\frac x\e\leq K$.
En itérant l'intégration par parties, on obtient, avec l'opérateur $S$ donné par \rf S et l'opérateur $D=\frac d{dx}$
\begin{eqnarray}\lb{onze}
y^-(x,\eps)&=&\sum_{n=0}^{N-1}\Big(\big(\ts\frac12DS\big)^ng\Big)(0)\e^{2n+1} U^-\big(\ts\frac x\e\big)
-\\\nonumber&&\frac12\,\ds\sum_{n=0}^{N-1}S\Big(\big(\ts\frac12DS\big)^ng\Big)(x)\e^{2n+2}
+\O\big(\e^{2N+1}\big)
\end{eqnarray}
quand $\e\to0$ uniformément sur l'ensemble des $x$ avec $\frac x\e\leq K$.
Il s'agit d'un exemple de développement combiné, de la forme 
$\sum_{n\geq0}\Big(a_n(x)+g_n^-\big(\ts\frac x\e\big)\Big)\e^n$, avec ici 
\eq{8bis}{
a_0=0,\,a_{2n}=-\ts\frac12\,S\Big(\big(\ts\frac12DS\big)^{n-1}g\Big),\,a_{2n+1}=0
}
et 
$$g_n^-=c_nU^-\mbox{ avec }
c_{2n}=0,\,c_{2n+1}=\Big(\big(\ts\frac12DS\big)^ng\Big)(0).
$$
De plus, la fonction $U^-$ admet un développement asymptotique à  l'infini, donné par 
\begin{eqnarray}\lb{U}
U^-(X)&\sim&\sum_{n\geq0}(-1)^{n+1}\,1.3...(2n-1)2^{-n-1}X^{-2n-1}\\\nonumber&=&-\frac1{2X}+\frac1{4X^3}-\frac3{8X^5}+...,\quad X\to-\infty.
\end{eqnarray}
On a aussi une formule analogue pour la solution $y^+$ bornée sur $\R^+$~:
\begin{eqnarray*}
y^+(x,\eps)&=&\sum_{n=0}^{N-1}\Big(\big(\ts\frac12DS\big)^ng\Big)(0)\e^{2n+1} U^+\big(\ts\frac x\e\big)
-\\\nonumber&&\frac12\,\ds\sum_{n=0}^{N-1}S\Big(\big(\ts\frac12DS\big)^ng\Big)(x)\e^{2n+2}
+\O\big(\e^{2N+1}\big)
\end{eqnarray*}
avec $U^+(X)=e^{X^2}\ds\int_{+\infty}^X e^{-T^2}\,dT = -U^-(-X)$.
Par ailleurs, dans le cas où $g$ est impaire, on a $(DS)^ng$ impaire pour tout $n\in\N$, donc la première partie du développement \rf{onze} est identiquement nulle. On retrouve ainsi le fait que $y^-$ a un développement asymptotique classique en puissances de $\e^2=\eps$, avec des coefficients de la variable $x$ uniquement.

Nous sommes dans une situation proche de la définition \reff{d2.2}. La principale différence est qu'ici le cadre est réel, alors que la définition  \reff{d2.2} est dans le cadre complexe, mais des modifications mineures permettent de définir les solutions $y^-$, $y^+$ et les fonctions $U^-$, $U^+$ sur des secteurs contenant les demi-axes réels. 
Par exemple, si la  fonction $g$ est définie et analytique bornée dans une bande horizontale $S=\{x\in\C\tq|\im x|<d\}$, alors on peut montrer que la formule \rf{3b} avec pour chemin d'intégration la demi-droite à  gauche issue de $x$, définit une fonction analytique dans $S$ qui, de plus, est bornée dans $S^-=\big\{x\in S\tq|\arg x-\pi|<\frac{3\pi}4\big\}$. 
Ceci permet par exemple de comparer les \dacs{} de $y^-$ et de $y^+$ dans l'intersection $S^-\cap(-{S^-})=\{x\in\C\tq|\re x|<|\im x|<d\}$.\vspace*{-12mm}
\figu{f2.2}{ 
\epsfxsize10cm\epsfbox{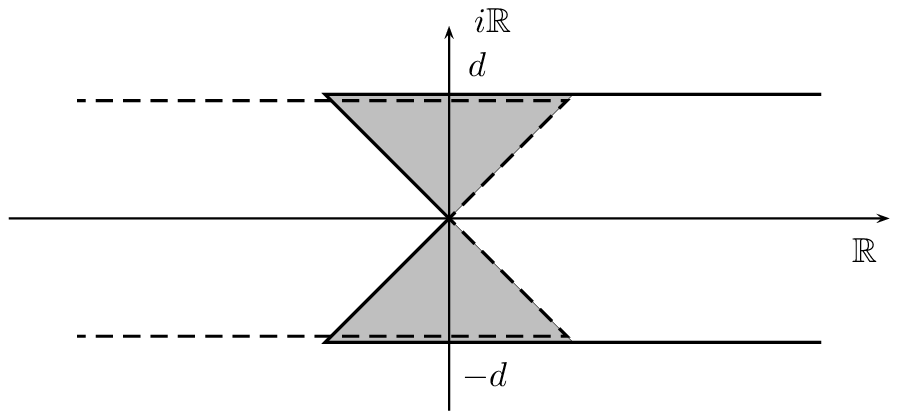}
\vspace{-16mm}}
{Le  bord de $S^-$ en pointillés, celui de $-S^-$ en trait plein,  l'intersection $S^-\cap(-S^-)$ en grisé}

À l'aide de normes de Nagumo \cit{crss}, on montre facilement que le développement \rf{onze} est Gevrey  d'ordre 1/2 en $\e$. 
Observons d'ailleurs que le développement \rf U est aussi un développement Gevrey d'ordre 1/2.

Notons au passage que les termes $a_n$ et $g_n^-$ sont nuls pour la moitié d'entre eux~; il serait donc possible de réécrire \rf{onze} sous forme de séries en puissances de $\eps$, mais cela est très particulier à  cet exemple. Dans les applications de la partie \reff{6.}, nous verrons que, bien souvent, les termes de la partie lente d'un \dac{} sont effectivement nuls sauf pour les puissances qui sont des multiples de $p$, mais que les termes de la partie rapide n'ont pas de raison \apriori{} de s'annuler.
\sub{1.2.2}{Extensions}
{\noi\bf2.2.1}.
Avant de présenter le deuxième exemple, nous voudrions explorer des généralisations et extensions du premier exemple. La première généralisation concerne la nature purement locale du résultat. Tout d'abord, toute autre solution de \rf{1b}, de condition initiale $y(x_0,\eps)$ bornée (pour $\eps\in\,]0,\eps_0]$) en un point fixé $x_0<0$, admet un \dac{} sur $[x_0+\d,0]$ pour tout $\d\in\,]0,|x_0|[$~; de plus ce \dac{} est le même que celui de $y^-$ puisque les deux solutions sont exponentiellement proches l'une de l'autre sur $[x_0+\d,0]$. 
Pour les mêmes raisons, le résultat reste valide avec une hypothèse seulement locale  sur $g$~: si $g$ est de classe $\cc^\infty$ sur un  intervalle $]-r,r[$, alors pour tout $x_0\in\,]\!-r,0[$, tout $\d\in\,]0,|x_0|[$ et toute fonction $c=c(\eps)$  bornée sur $]0,\eps_0]$, la solution de \rf{1b} de condition initiale $y(x_0,\eps)=c(\eps)$ admet un \dac{} de la forme \rf{onze} sur $[x_0+\d,0]$~; il en est de même pour les solutions à  droite, \ie avec $x_0\in\,]0,r[$. Ces \dacs{} ne dépendent pas des conditions initiales~; ils sont \apriori{} différents à  gauche et à  droite mais ils ont la même partie lente $\sum_{n\geq0}a_n(x)\e^n$ avec $a_n(x)$ données par \rf{8bis}.
\bigskip\\
{\noi\bf2.2.2}.
La deuxième généralisation est de remplacer le terme $2x$ dans \rf{1b} par $p\,x^{p-1}$, où $p$ est un entier pair%
\footnote{\ Ici la situation qui nous intéresse est avec $p$ pair, mais le cas $p$ impair a aussi son intérêt, \cf l'équation \rf{ujmod} dans la partie \reff{6.2}.}.
 Nous avons toujours une unique solution $y^-$ bornée sur $\R^-$ et  une unique solution $y^+$ bornée sur $\R^+$. Elles sont données à   présent par
$y^\pm(x,\eps)=e^{x^p/\eps}\ds\int_{\pm\infty}^xe^{-t^p/\eps}g(t)dt$. La condition $y^-=y^+$ est toujours équivalente à  $g$ impaire.
La recherche d'une solution formelle $\yc=\sum_{n\geq0}y_n(x)\eps^n$ aboutit à 
$$
y_0(x)=0,~~y_{1}(x)=-\frac1{p\,x^{p-1}}\,g(x),
~~\mbox{puis}~~y_{n+1}(x)=\frac1{p\,x^{p-1}}\,y'_{n}(x).
$$
En général cette solution formelle n'est pas définie en $x=0$, même lorsque $g$ est impaire. Dans le cas où $g$ est impaire, le développement de $y^-$ est valide aussi bien pour les $x$ positifs que pour les $x$ négatifs, mais ne peut pas être valide au voisinage de $0$. Cependant la même méthode d'intégrations par parties successives permet de montrer (que $g$ soit impaire on non) que $y^-$ possède un \dac, mêlant 
à la fois des fonctions de $x$ et des fonctions 
de la variable rapide $X=\frac x\e$, avec $\e=\eps^{1/p}$. Les calculs sont plus longs et compliqués sans être plus difficiles~; ils font apparaître les $p-1$ \og fonctions spéciales\fg \ suivantes
$$
U^-_{k}(X)=e^{X^p}\ds\int_{-\infty}^Xe^{-T^p}T^{k-1}dT,\quad k=1,...,p-1.
$$
On obtient finalement un \dac{} pour $y^-$ de la forme $\sum_{n\geq0}\Big(a_n(x)+$$g_n\big(\ts\frac x\e\big)\Big)\e^n$ avec $\e=\eps^{1/p}$ et  de $g_n$ de la forme $g_n=c_{n1}U^-_{1}+\dots+c_{n,p-1}U^-_{p-1}$. De même que précédemment pour $U^-$, ces fonctions $U^-_{kp}$ ont aussi un développement asymptotique lorsque $X$ tend vers $-\infty$.
\bigskip\\
{\noi\bf2.2.3}.
Une troisième extension concerne les équations où la fonction $g$ dépend de $\eps$. Si $g$ a un développement asymptotique en puissances de $\eps$, ainsi que toutes ses dérivées par rapport à  $x$, alors on peut montrer que les fonctions $y^\pm$ ont encore des \dacs{} dont les fonctions $g_n^\pm$ sont proportionnelles à $U^\pm$~; le facteur est le même pour les deux signes. Précisément ce \dac{} est donné comme avant par
\begin{eqnarray}\lb{onzea}
y^\pm(x,\eps)&=&
\sum_{n=0}^{N-1}A_{N-n,\,n}(0,\eps)\e^{2n+1} U^\pm\big(\ts\frac x\e\big)
-\\\nonumber&&\frac12\,\ds\sum_{n=0}^{N-1}B_{N-n,\,n}(x,\eps)\e^{2n+2}
+\O\big(\e^{2N+1}\big)
\end{eqnarray}
où $A_{mn}:(x,\eps)\mapsto\ds\sum_{k=0}^mA_{mnk}(x)\eps^k$ est le jet d'ordre $m$ par rapport à  $\eps$ de la fonction $\big(\ts\frac12DS\big)^ng$ et $B_{mn}$ est le jet d'ordre $m$ par rapport à  $\eps$ de la fonction $S\left(\big(\ts\frac12DS\big)^ng\right)$.

La seule modification à  apporter est la condition 
nécessaire et suffisante pour avoir $y^-=y^+$. A la place de $g$ impaire, cette condition devient 
$$
\int_{-\infty}^\infty e^{-t^2/\eps}\,g(t,\eps)\,dt=0.
$$
 Lorsque $y^-=y^+$, on obtient à  nouveau que la solution bornée sur $\R$ a un développement classique dans le cas $p=2$ car les facteurs de $U^\pm$ dans \rf{onzea} s'annulent. Par contre lorsque $p\geq4$, ce n'est plus forcément le cas.
\bigskip\\
{\noi\bf2.2.4}. 
Notre dernière extension est d'avoir $f(x)$ à  la place de $p\,x^{p-1}$, où $f$ est une fonction de classe $\cc^\infty$ vérifiant $xf(x)>0$ si $x\neq0$ et $f(x)\sim ax^{p-1},\,x\to0$ avec  $a\neq0$.
La solution $y^-$ s'écrit alors $y^-(x,\eps)=e^{F(x)/\eps}\ds\int_{-\infty}^xe^{-F(t)/\eps}g(t)dt$ avec $F(x)=\int_0^xf(t)dt$, si on ajoute l'hypothèse que
$\ds\int_{-\infty}^0 e^{-F(t)/\eps}\,dt$ converge pour $\eps>0$ assez petit.

Un difféomorphisme  $x=\f(\xi)$ permet de se ramener à  une équation de la forme
$$
\eps\frac{dz}{d\xi}=p\,\xi^{p-1}z+\eps h(\xi),
$$
ce qui permet d'obtenir pour $y^-$ un \dac{} de la forme 
$$
\sum_{n\geq0}\Big(a_n(x)+g_n\big(\ts\frac{ \f^{-1}(x)}\e\big)\Big)\e^n.
$$

Notre théorie générale des \dacs{} permet de montrer qu'un tel développement peut aussi se transformer en un \dac{} de la variable $X=\frac x\e$, \ie de la forme $\sum_{n\geq0}\gk{b_n(x)+h_n\big(\ts\frac x\e\big)}\e^n$, \cf proposition \reff{l2.4} (c).
Remarquons que, pour cette extension, on a besoin d'autres fonctions ``rapides''
que $U^-$ et $U^+$. 
\sub{1.2.3}{Deuxième exemple}
Considérons  à  présent une équation déjà  considérée sous une forme voisine par 
Claude Lobry dans le chapitre introductif \cit{l1}. Il s'agit de l'équation
\eq{e1.1}{
\eps\frac{dy}{dx}=2xy+\eps g(x)+\eps\a,
}
où $\eps>0$ est le petit paramètre, $\a\in\R$ un paramètre de contrôle, et où la fonction $g:\R\to\R$ est de classe $\cc^\infty$, et bornéee ainsi que toutes ses dérivées. La question est la suivante.
\med
 
{\sl Existe-t-il des valeurs de $\a$ pour lesquelles il existe une solution %
bornée%
\footnote{\ Rappelons que `` bornée '' signifie uniformément par rapport à  $\eps$ dans un intervalle $]0,\eps_0]$. Dans le présent contexte, il se trouve que, pour tout $\eps$ fixé, il existe une unique valeur $\a=\a(\eps)$ pour laquelle l'équation \rf{e1.1} a une solution $y=y(x,\eps)$ bornée sur $\R$ au sens classique, et que la fonction $y$ ainsi définie est aussi bornée sur $\R\times\,]0,\eps_0]$.}
sur tout $\R$~?}
\med

La réponse est \og{oui}\fg. Pour le voir, on procède ainsi~: étant donné $\a$ arbitraire, il existe toujours l'unique solution bornée sur $\R^-$, notée $y^-$ et donnée par
\eq{e1.3}{
y^-(x,\eps)=e^{x^2/\eps}\int_{-\infty}^xe^{-t^2/\eps}\big(\a+g(t)\big)dt.
}
En remplaçant $-\infty$ par $+\infty$, la même formule fournit aussi une unique solution $y^+$ bornée sur $\R^+$.
On en déduit qu'il existe une solution $y$ bornée sur tout $\R$ si et seulement si $y^+=y^-$, ce qui donne une équation pour le paramètre $\a$, dont la solution est
\eq{e1.4}{
\a(\eps)=-\bigg(\int_{-\infty}^{+\infty}e^{-t^2/\eps}g(t)dt\bigg)\bigg/\bigg(\int_{-\infty}^{+\infty}e^{-t^2/\eps}dt\bigg).
}
Comme nous avons supposé $g$ de classe $\cc^\infty$, 
on en déduit que $\a$ admet un développement asymptotique quand  $\eps$ tend vers 0.
Pour l'étude des solutions $y^\pm$  correspondant à  cette valeur de $\a$, 
on est dans la situation de la deuxième extension (avec pour fonction $g$ la fonction $(x,\eps)\mapsto g(x)+\a(\eps)$) et
on a $y^-=y^+$ par le choix de $\a$.
La solution $y$ admet donc aussi un développement asymptotique, dont les coefficients
sont des fonctions ${\cal C}^\infty$, y compris en $x=0$.

On calcule directement les premiers termes
\eq{e1.6}{
y_0(x)=0,\quad \a_0=-g(0),\quad y_{1}(x)=-\frac1{2x}\big(g(x)+\a_0\big)  .
}
Pour la suite, il suffit de remarquer que les $\a_n$ et les $y_n(x)$
sont déterminées uniquement par le fait que $y_n$ n'a pas de pôle en $x=0$.
On les calcule donc récursivement par
\eq{e1.6bis}{
y_{n+1}(x)=\frac1{2x}\big(y'_{n}(x)-\a_n\big),\quad \a_n=y'_{n}(0).
}

Dans le champ complexe, supposons par exemple $g$ analytique et vérifiant 
$|g(x)|\leq Me^{M|x|^2}$ dans une bande horizontale
$$
S=\{x\in\C\tq|\im x|<d\}
$$
pour un certain $d>0$. L'équation \rf{e1.1} a un relief constitué des courbes de niveau de la fonction $R:x\mapsto\re(x^2)$. On peut montrer que, pour toute valeur de $\a=\a(\eps)$  bornée et pour tout $\d\in\,]0,d[$ fixé, la solution $y^+$ est bornée  dans le domaine
$$
D_\d^+=\big\{x\in S\tq|\arg x|<\ts\frac{3\pi}4-\d,\,|x|>\d\big\},
$$
qui comprend la majeure partie de la montagne à  l'est et des deux vallées au nord et au sud dans la bande $S$, \cf figure \reff{f2.3}.
\figu{f2.3}{
\epsfxsize8.6cm\epsfbox{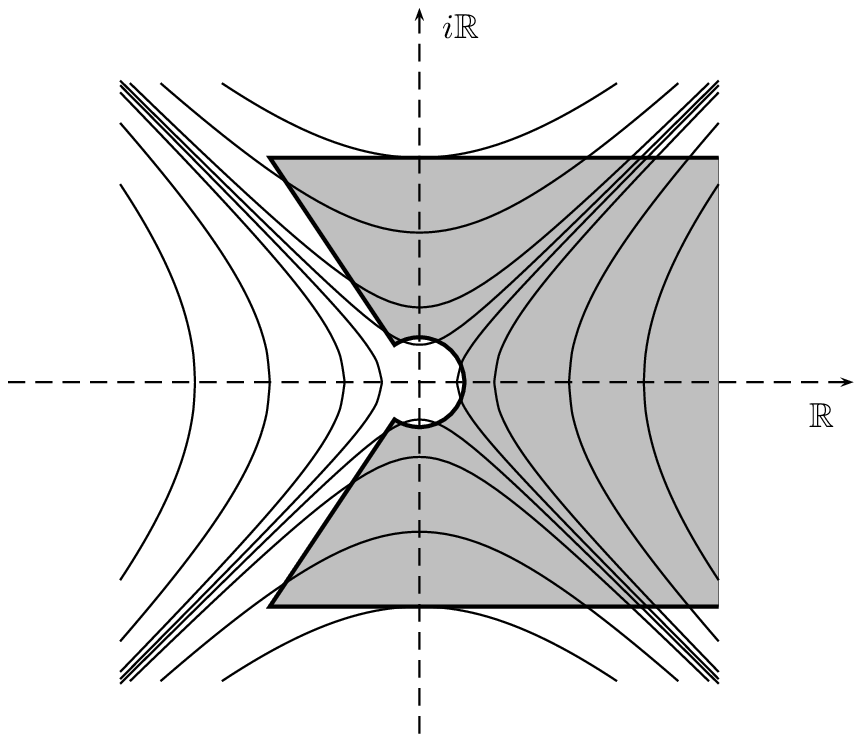}
}
{Les courbes de niveau $\re(x^2)=$ constante~; en gris le domaine $D_\d^+$}
Si de plus  $\a$ admet un développement asymptotique en puissances de $\eps$, alors $y^+$ aussi 
a un développement asymptotique $y^+(x,\eps)\sim\sum_{n\geq0}y_n^+(x)\eps^n$.
Précisément, pour tout $\d>0$ et tout $N\in\N$, il existe $M=M(\d,N)>0$ tel que 
\eq{e1.5}{
\Big|y^+(x,\eps)-\sum_{n=0}^{N-1}y_n^+(x)\eps^n\Big|\leq M\eps^N
}
pour tout $x\in D^+_\d$.
Il en est de même pour $y^-$ dans $D_\d^-=\{x\in S\tq|\arg x-\pi|<\frac{3\pi}2-\d,\,|x|>\d\}$.

Dans le cas où $\a$ est donné par \rf{e1.4}, c'est-à-dire lorsque $y^+=y^-$, ces deux développements coïncident~: en effet, par unicité du développement asymptotique, on a $y^+_n=y^-_n$ sur
$D_\d^-\cap D_\d^+$, donc sur $D_\d^-\cup D_\d^+=S\setminus D'(0,\d)$ 
par unicité du prolongement analytique,
où $D'(0,\d)$ est le disque fermé de centre $0$ et de rayon $\d$.

Ce développement est donné par \rf{e1.6} et \rf{e1.6bis}.

Puisque $y^+$ est analytique dans tout $S$ et que \rf{e1.5} est satisfait pour $x\in S\setminus D(0,\d)$, par le principe du maximum l'inégalité est satisfaite pour tout $x$ dans $S$. En résumé, du fait que le point tournant $x=0$ est simple, le relief ne présente que deux montagnes, si bien qu'une solution définie et bornée sur des parties de ces deux montagnes est automatiquement bornée sur tout un voisinage du point tournant. On dit d'une telle solution qu'elle est {\sl surstable}, suivant la terminologie adoptée par Guy Wallet \cit{w}.
\sub{1.2.4}{Troisième exemple}

Remplaçons le terme $2x$ par $4x^3$~; autrement dit, considérons l'équation
\eq{e1.7}{
\eps\frac{dy}{dx}=4x^3y+\eps g(x)+\eps \a,
}
Pour tout $a=a(\eps)$, il existe encore une unique solution $y^+$ bornée sur $\R^+$ et  une unique solution $y^-$ bornée sur $\R^-$. Il existe aussi une unique valeur de $\a$ pour laquelle ces deux solutions coïncident pour donner une solution $y$ bornée sur $\R$. Cette solution est donnée par
\begin{eqnarray}\lb{e1.8}
y(x,\eps)&=&e^{x^4/\eps}\int_{-\infty}^xe^{-t^4/\eps}\big(\a+g(t)\big)dt\quad\\\nonumber
\mbox{avec}\quad\a&=&\a(\eps)=-\frac{\int_{-\infty}^{+\infty}e^{-t^4/\eps}g(t)dt}{\int_{-\infty}^{+\infty}e^{-t^4/\eps}dt}.
\end{eqnarray}
\figu{f2.4}{ 
\epsfxsize54mm\epsfbox{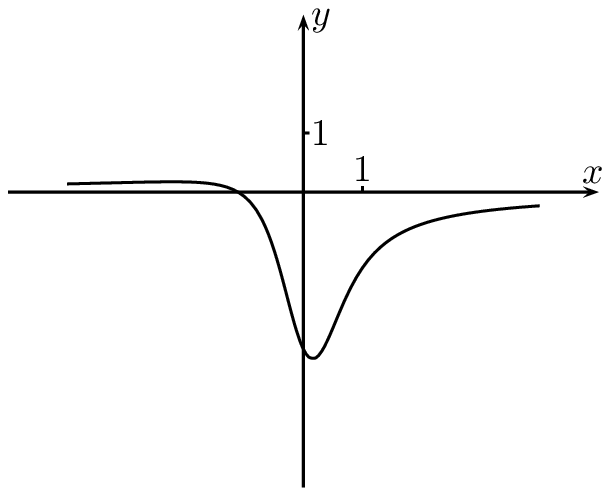}\hspace{7mm}\epsfxsize54mm\epsfbox{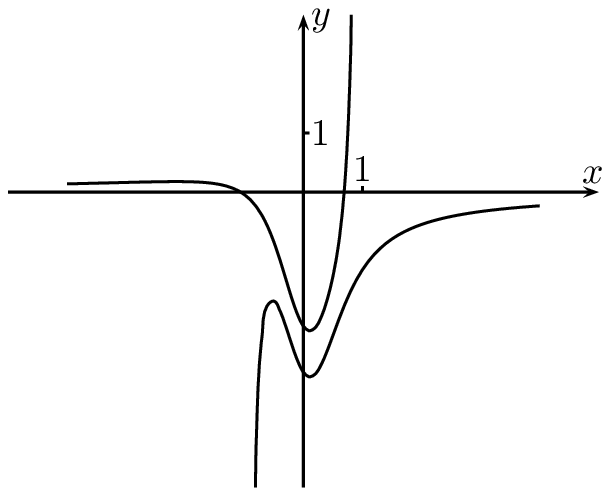}
\vspace{-8mm}}
{Les solutions $y^+$ et $y^-$ pour $g(x)=3x^2+3x$. 
À gauche pour $\a=-0.507$, à droite pour $\a=-0.35$.
Ici $|x|\leq5$, $-5\leq y \leq3$ et $\eps=\frac14$}
À nouveau on peut montrer que $\a(\eps)$ admet un développement asymptotique de la forme 
$\a(\eps)\sim\ds\sum_{n=0}^\infty \a_n \eps^{n/2}$.
En revanche, la solution $y$ n'a en général pas de développement asymptotique en puissances de $\eps^{1/2}$ dans un voisinage réel de $0$. En effet, un tel développement serait une solution formelle  $\yc=\sum_{n\geq0}y_n(x)\eps^{n/2},\ 
\ac=\sum_{n\geq0}\a_n\eps^{n/2}$
devant vérifier $y_0(x)=0$, $y_2(x)=-\frac1{4x^3}\big(g(x)+\a_0\big)$ et $y_{2n+2}(x)=\frac1{4x^3}\big(y'_{2n}(x)-\a_{2n}\big)$. Dès qu'un coefficient $y_{2n}$ a une dérivée dont le développement de Taylor contient un terme non nul en $x$ ou $x^2$, le coefficient suivant présente un pôle en $x=0$, quelque soit le choix de $\a_{2n}$. 
\med

Nous allons voir qu'il est cependant possible de donner une approximation de $y$ valide dans tout  un voisinage réel de $0$, à  l'aide de fonctions de $\frac x\e$, où $\e=\eps^{1/4}$. Cette idée
était à la base de la thèse de Th.\ Forget \cit{fo}.

En effet,  isolons les premiers termes du développement de $g$ en écrivant 
$g(x)=g_0+g_1x+g_2x^2+4x^3h(x)$. La première formule de \rf{e1.8} devient
\begin{eqnarray}\lb{e1.8bis}
y(x,\eps)&=&(\a(\eps)+g_0)U^-_0\big(\ts\frac x\e\big)+\e g_1U_1^-\big(\frac x\e\big)+
\\\nonumber
&&\e^2g_2U_2^-\big(\frac x\e\big)+e^{x^4/\eps}\ds\int_{-\infty}^xe^{-t^4/\eps}4t^3h(t)dt
\end{eqnarray}
où on a posé
$$
U^-_j(X)=-e^{X^4}\int_{-\infty}^Xe^{-T^4}T^jdT.
$$
Une intégration par parties donne 
$$
e^{x^4/\eps}\int_{-\infty}^xe^{-t^4/\eps}4t^3h(t)dt=\eps h(x)-\eps e^{x^4/\eps}\int_{-\infty}^xe^{-t^4/\eps}h'(t)dt
$$
qui est de la même forme que dans la formule \rf{e1.8bis}, ce qui permet de réitérer.
De la même manière que dans {2.2.2}, on obtient à  la fin un \dac{}.
Le fait que $y^+$ soit égal à  $y^-$ entraîne qu'il s'agit d'un \dac{} sur {\sl tout l'axe réel}.
Précisément, il existe des fonctions $a_n\in C^\infty(\R)$ et des combinaisons linéaires
$g_n$  des fonctions $U^-_0,U^-_1,U^-_2$
telles que $y^\pm(x,\eps)=\sum_{n=0}^{N-1}\Big(a_n(x)+g_n\big(\tfrac x\e\big)\Big)+
{\cal O}(\e^N)$ uniformément sur $\R$.

Dans le champ complexe, le relief associé à  \rf{e1.7}, qui est donné par la partie réelle de $x^4$, comprend quatre montagnes. La solution étudiée est proche de la courbe lente sur deux montagnes, à  l'est et à  l'ouest, mais \apriori{} pas sur les deux autres montagnes nord et sud~; elle n'est donc pas surstable.
\sub{1.2.5}{Quatrième exemple}
Remplaçons le terme $2x$ par $-2x$ dans le premier exemple \rf{1b}, \ie considérons 
\eq{1b-}{\eps\frac{dy}{dx}=-2xy+\eps g(x)}
avec $g:\R\to\R$ de classe $C^\infty$ et bornée ainsi que toutes ses dérivées.
Contrairement à \rf{1b}, la courbe lente $y=0$ est d'abord répulsive pour $x<0$,
puis attractive pour $x>0$. 

Fixons une  valeur initiale $c=c(\eps)$, bornée quand $0<\eps\to0$.
La solution de condition initiale $y(0,c(\eps),\eps)=c(\eps)$ est alors bornée sur tout $\R$ quand $\eps\to0$. Elle est donnée par la formule 
\eq{vc}{
y(x,c(\eps),\eps)=e^{-x^2/\eps}\int_0^xe^{t^2/\eps}g(t)\,dt + 
 c(\eps) e^{-x^2/\eps}.
}

Par une succession d'intégrations par parties, on montre comme dans la partie \reff{1.2.1} que
$y(x,c(\eps),\eps)$ admet un développement asymptotique $\ds\sum_{n\geq1}y_n(x)\eps^n$ pour $x$ loin de $0$. Précisément, pour tout $\d>0$, on a
$$
y(x,c(\eps),\eps)\sim\sum_{n\geq1}y_n(x)\eps^n
$$
quand $0<\eps\to0$, uniformément sur $]-\infty,-\d]\cup[\d,+\infty[$. 
Les coefficients $y_n(x)$ sont déterminés de manière analogue
à \rf{6b} et \rf{6bis} par $y_1(x)=\frac1{2x}g(x)$ et $y_{n+1}(x)=-\frac1{2x}y_n'(x)$.
Ce développement est indépendant de la condition initiale $c(\eps)$~; le terme exponentiel
dans \rf{vc} ne contribue pas en dehors de voisinages de $0$.

Comme dans les parties précédentes, il se pose la question du comportement de ces 
solutions au voisinage de $x=0$. On procède de nouveau par intégration par parties en
utilisant $g(x)=g(0)+x\,Sg(x)$~; on obtient
\begin{eqnarray*}
y(x,c(\eps),\eps)\!\!&=&g(0)\e U(X)+\tfrac\eps2Sg(x)+\big(c(\eps)-
   \\\nonumber
  && 
  \tfrac\eps2Sg(0)\big)e^{-x^2/\eps}-\tfrac\eps2e^{-x^2/\eps}\int_0^xe^{t^2/\eps} (Sg)'(t)\,dt
\end{eqnarray*}  
avec $\e=\sqrt\eps$ et $U(X)=\int_0^Xe^{-X^2+T^2}dT$.
En répétant plusieurs fois l'intégration par parties, on obtient comme dans la partie
\reff{1.2.1}
\begin{eqnarray}\lb{onze-}
y^-(x,\eps)\!\!\!&=&\!\!\ds\sum_{n=0}^{N-1}
   \Big(\big(\ts-\frac12DS\big)^ng\Big)(0)\e^{2n+1} U\big(\ts\frac x\e\big)+
   \\
  \nonumber&&\frac12\,\ds\sum_{n=0}^{N-1}S\Big(\big(-\ts\frac12DS\big)^ng\Big)(x)\e^{2n+2}
  +\\\nonumber&&
  \gk{c(\eps)-
     \tfrac12\,\ds\sum_{n=0}^{N-1}S\Big(\big(-\ts\frac12DS\big)^ng\Big)(0)\e^{2n+2}}
   e^{-x^2/\eps}
+\O\big(\e^{2N+1}\big).
\end{eqnarray}
Supposons maintenant que $c(\eps)$ admette un développement 
asymptotique en puissances de $\e$~: $c(\eps)\sim\sum_{n\geq0}c_n\e^n$.
Alors on obtient à nouveau un développement combiné 
$$
y(x,c(\eps),\eps)=\sum_{n\geq0}\gk{a_n(x)+g_n(\tfrac x\e)}\e^n
$$
 avec pour  $n=0,1,2,...$
$$
a_0=0,\quad  a_{2n+2}=\tfrac12S\gk{-\tfrac12DS}^{n}g,\quad a_{2n+1}=0
$$
ainsi que
$$
g_0(X)=c_0e^{-X^2}~\mbox{ et pour }~n\geq1\quad g_{2n}(X)=d_ne^{-X^2}
$$
avec
$$
d_n={c_{2n}-\tfrac12S\Big(\big(-\ts\frac12DS\big)^{n-1}g\Big)(0)}
$$
et pour $n\geq0$
$$
g_{2n+1}(X)=b_{n}U(X)+c_{2n+1}e^{-X^2}~\mbox{ avec}~\quad b_{n}=\Big(\big(\ts-\frac12DS\big)^ng\Big)(0).  
$$
Contrairement à la partie \reff{1.2.1}, ce développement dépend de la condition
initiale, mais il est le même pour $x$ positif ou négatif.
On peut généraliser et étendre cet exemple, mais cette exploration, analogue
à la partie \reff{1.2.2}, est laissée au lecteur. 
\sub{1.3}{Les développements combinés}
Si l'on veut généraliser la méthode de la partie précédente à  des équations non linéaires, cela nécessite d'élargir la famille de fonctions dans lesquelles s'écrivent les solutions. En particulier, il est nécessaire de prendre en compte les produits de fonctions $U^\pm_j\,U^\pm_k$,
les solutions d'équations différentielles $\eps y'=p x^{p-1} y + U^\pm\big(\tfrac x\e\big)$,
ainsi que des produits de fonctions de $x$ et de fonctions de $\frac x\e$.
Une stratégie
est de construire une algèbre contenant les fonctions $x\mapsto x^n$ et $(x,\eps)\mapsto U_j^\pm\big(\frac x\e\big)$ et stable par les opérateurs
 $\J^\pm$ suivants. Pour chaque signe $+$ et $-$, $\J^\pm$ associe, à  une fonction $v$ à  croissance polynomiale (\ie vérifiant $\exists N\in\N,C>0\;\forall x\in\R,\; 
|v(x)|\leq C\,|x|^N$),
l'unique solution à  croissance polynomiale sur $\R^\pm$ de l'équation 
$
\ds\frac{dU}{dX}=4X^3U+v(X)
$ 
i.e. $\J^\pm$ est donné par 
$$
\J^\pm v(X)=e^{X^4}\ds\int_{\pm\infty}^Xe^{-T^4}v(T)dT.
$$
La construction de la plus petite algèbre avec ces propriétés
conduit à  définir un grand nombre de ``fonctions spéciales''~; ceci était la stratégie adoptée par Th.\ Forget dans sa thèse 
pour l'approximation de solutions canard comme 
dans \reff{1.2.4}. Une complication
additionnelle provient du caractère non unique de l'écriture. 
Par exemple, le terme $x$ peut être considéré aussi bien comme un terme fonction de $x$ d'ordre $0$ en $\e$ qu'un terme fonction de $\frac x\e$ d'ordre $1$ en $\e$, s'il est écrit $\frac x\e\e$.

La  stratégie que nous avons adoptée est de considérer d'emblée une algèbre plus grosse. Un avantage de cette stratégie est la simplicité~; un inconvénient est de donner moins de renseignements sur les coefficients.
%
%
\sec{2.}{Développements asymptotiques combinés : \\ étude générale}
Nous présentons dans cette partie la théorie générale des \dacs~:
leur définition et leur comportement vis-à-vis des opérations élémentaires d'addition, multiplication, division, dérivation, intégration, composition et prolongement. Dans la section \reff{2.3}, nous faisons aussi le lien  entre nos \dacs{} et les développements intérieurs et extérieurs de la classique méthode de matching. Ces derniers sont par ailleurs un bon procédé pour déterminer les coefficients d'un développement combiné en pratique, à  condition d'avoir pu montrer par ailleurs l'existence de ce développement. 

Beaucoup des problèmes résolus par l'utilisation des \dacs{} ont leur origine dans un cadre réel, si bien qu'une présentation purement réelle semble suffire. 
Cependant un élément essentiel dans la résolution de ces problèmes est l'aspect Gevrey de ces \dacs, qui sera développé dans la partie \reff{2bis}. Pour obtenir ces propriétés Gevrey,
plusieurs méthodes sont possibles, certaines dans le cadre réel, mais à  nos yeux la méthode la plus efficace est d'appliquer 
notre ``théorème-clé'' de type Ramis-Sibuya \reff{t3.1}, pour lequel le cadre complexe est indispensable.
C'est pourquoi la présentation faite ici est exclusivement complexe~; on trouvera dans \cit{fs2} une présentation très voisine dans le cadre réel. 
\sub{2.1}{Notations}
Le disque ouvert de centre $0$ et de rayon $r$ est noté $D(0,r)$ et le disque fermé $D'(0,r)$.
On note $\cl(E)$ l'adhérence d'une partie $E$ de $\C$.
Étant donnés $\a<\b\leq\a+2\pi$ et $0<r\leq\infty$, $S(\a,\b,r)$ est le secteur 
$$
S(\a,\b,r)=\{x\in\C\tq0<|x|<r,\,\a<\arg x<\b\}.
$$
Il est d'usage de considérer un secteur comme une partie de la surface de Riemann du logarithme $\~{\C^*}$. Cependant, puisque nos secteurs auront toujours une ouverture inférieure à  $2\pi$, nous les considérerons comme des sous-ensembles de $\C^*$.
Étant donné un secteur $S=S(\a,\b,r)$ et $\mu>0$, on note $V(\a,\b,r,\mu)$ l'union du secteur $S$ et du disque $D(0,r)$~:
\eq{V1}{
V(\a,\b,r,\mu)=\{x\in\C \tq (|x|<r\mbox{ et }\a<\arg z<\b)\mbox{ ou bien }\norm x<\mu\}.
}
Nous serons aussi amenés à  définir des \dacs{} dans des secteurs privés d'un disque fermé, \ie des arcs de couronnes. On définit donc, pour $\mu<0$
\eq{V2}{
V(\a,\b,r,\mu)=\{x\in \C\tq-\mu<|x|<r\mbox{ et }\a<\arg x<\b\}.
}
\figu{f3.1}{
\vspace{-0.5cm}
\epsfxsize9cm\epsfbox{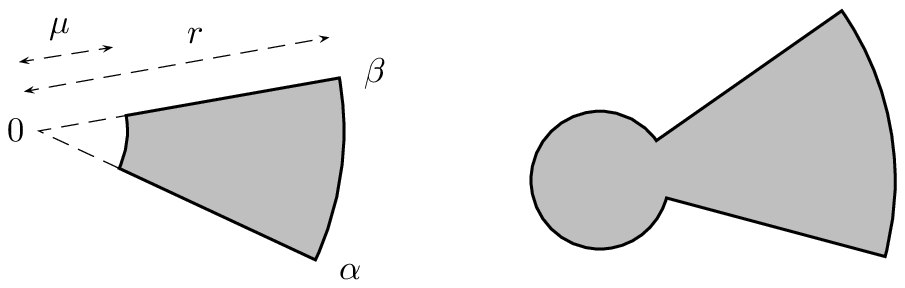}
\vspace{-7mm}
}
{Deux exemples de quasi-secteurs, à gauche pour $\mu<0$, à droite pour $\mu>0$}
Pour faire court, nous appellerons {\sl quasi-secteurs} ces ensembles, que $\mu$ soit positif ou négatif. Souvent pour simplifier, nous ne considérerons que le cas $\mu>0$~; les modifications à  apporter dans le cas $\mu<0$ sont mineures et seront indiquées au fur et à  mesure.

Étant donné un quasi-secteur infini $V=V(\a,\b,\infty,\mu)$, 
$\G(V)$ désigne l'espace vectoriel des fonctions $g$ holomorphes 
bornées dans $V$ et ayant un développement asymptotique au sens 
de Poincaré à  l'infini sans terme 
constant $g(X)\sim\sum_{\nu\geq1}g_\nu X^{-\nu},\,V\ni X\to\infty$, \ie
$$
\forall N\in\N\quad\exists C_N>0\quad\forall X\in V,\qquad\Big|g(X)-\ds\sum_{\nu=1}^{ N-1}g_\nu X^{-\nu}\Big|\leq C_N|X|^{-N}.
$$
Soit $\T:\G(V)\to\G(V)$ l'opérateur qui, à  une fonction $g$, associe la fonction $\T g$ donnée par
\eq4{
\T g(X)=Xg(X)-g_1
}
où $g_1$ est le premier terme du développement asymptotique de $g$ à  l'infini.
\med\\
Par abus, la notation $V(0,2\pi,\infty,-r)$ désigne toute la couronne des $x\in\C$
avec $\norm x>r\geq 0$, et pas cette couronne privée de la demi-droite $]r,+\infty[$.
L'espace de Banach $\G(V)$ est alors l'espace des fonctions
holomorphes bornées sur $V$ tendant vers 0 quand $X\to\infty$.
\med

Pour un nombre $r_0>0$, on note $\H(r_0)$ l'espace vectoriel des  fonctions $a$ holomorphes bornées dans le disque $D(0,r_0)$ centré en $0$ et de rayon $r_0$.
De manière analogue à  $\T$, soit $\S:\H(r_0)\to\H(r_0)$ l'opérateur qui, à  une fonction $a$, associe la fonction $\S a$ donnée par
\eq5{
\S a(x)=\frac{a(x)-a(0)}x.
}
Sur les développements, les opérateurs $\S$ et $\T$ ont pour action de décaler vers la gauche et d'effacer le premier terme~:
si $g(X)\sim\ds\sum_{\nu\geq1}g_\nu X^{-\nu}$, alors $\T g(X)\sim\ds\sum_{\nu\geq1}g_{\nu+1} X^{-\nu}$ et si $a(x)=\ds\sum_{\nu=0}^\infty a_\nu x^{\nu}$, alors $\S a(x)=\ds\sum_{\nu=0}^\infty a_{\nu+1} x^{\nu}$.
\sub{2.2}{Séries formelles combinées}
\df{d2.1}{
 Soit $V=V(\a,\b,\infty,\mu)$ un quasi-secteur infini (avec $\mu$ positif ou négatif) et $r_0>0$. 
Une {\sl série formelle combinée} associée à  $V$ et à  $D(0,r_0)$ est une expression de la forme
\eq2{
\yc(x,\e)=\sum_{n\geq0}\lp a_n(x)+g_n\big(\ts\frac x\e\big)\rp\e^n,
}
où $a_n\in\H(r_0)$ et $g_n\in\G(V)$.
\med

Les fonctions $a_n$ forment la {\sl partie lente} de la série formelle combinée, et les $g_n$ la {\sl partie rapide}.}
\rqs1.\ Une série formelle combinée est précisément un élément de $(\H(r_0)\times\G(V))^\N$.
De manière analogue à  ce qui se fait pour les séries formelles classiques, nous pourrions donc représenter une série formelle combinée sous la forme $\sum_{n\geq0}\big( a_n(x)+g_n(X)\big)\e^n$ 
--- ou sous la forme $\sum_{n\geq0}\big( a_n(x),g_n(X)\big)\e^n$ --- 
à  l'aide de trois variables. Cependant, nous n'aurons pas à  considérer de fonctions de trois variables asymptotiques à  une série formelle combinée~: dans l'usage que nous en ferons les trois variables seront liées par la relation $x=\e X$. 
\med\\
2.\ \Apriori{}, pour l'unicité du développement, il suffirait de demander aux fonctions $g_n$ de tendre vers $0$ à  l'infini. Cependant, la formule de multiplication \rf{prodcomb} montre qu'il est nécessaire que les fonctions de $\frac x\e$ aient un développement asymptotique complet à  l'infini si l'on veut pouvoir écrire les développements combinés par exemple de produits $x^N\,g\big(\ts\frac x\e\big)$. Nous demandons aux $g_n$ d'être analytiques dans des quasi-secteurs, et pas simplement des secteurs, car ceci sera nécessaire dans les applications.

On note $\cch(r_0,V)$ l'ensemble des séries formelles combinées associées à  $V$ et à  $D(0,r_0)$ et on le munit 
de l'addition et de la multiplication scalaire canoniques ainsi que de la distance 
ultramétrique
\eq d{
d(\yc_1,\yc_2)=2^{-\val(\yc_1-\yc_2)},\mbox{ où }\val(\yc)=\min\{n\geq0\tq a_n\mbox{ ou } g_n\not\equiv0\}
}
et de la topologie induite par cette distance. 
On note $I$ l'inclusion canonique de $\H(r_0)$ dans $\cch(r_0,V)$ et, pour ne pas multiplier les notations, par la même lettre celle de $\G(V)$ dans $\cch(r_0,V)$. Le symbole $\e$ désigne à  la fois le nombre réel, la fonction de $(x,\e)$ à  valeur $\e$ et la série formelle combinée comprenant un seul terme $a_1=1$.

\

{\noi\bf L'algèbre ``différentielle'' $\cch(r_0,V)$} \sep
\med\\
Grâce au fait que les $g_n(X)$ admettent
un développement asymptotique quand $X\to\infty$, les opérateurs $\S$ et $\T$ permettent de munir $\cch(r_0,V)$ d'une structure d'algèbre de la manière suivante.

On commence par effectuer le produit de deux séries formelles combinées $\yc$ et $\zc$ en développant le produit terme à  terme~: si $\yc(x,\e)=\sum_{n\geq0}\Big( a_n(x)+$ $g_n\big(\ts\frac x\e\big)\Big)\e^n$ et $\zc(x,\e)=\sum_{n\geq0}\lp b_n(x)+h_n\big(\ts\frac x\e\big)\rp\e^n$, alors on pose
\begin{eqnarray*}
\yc.\zc(x,\e)&=&\sum_{n\geq0}\Bigg(\sum_{\nu=0}^n\Big( I(a_\nu)(x,\e)+
  I(g_\nu)( x,\e)\Big)\cdot
\\
\nonumber
&&
\ \ \ \ \ \ \ \ \ \ \ \ \ \cdot\Big( I(b_{n-\nu})(x,\e)+I(h_{n-\nu})(x,\e)\Big)\Bigg)\e^n\ \ ;
\end{eqnarray*}
ceci est une série convergente pour la topologie de $\cch(r_0,V)$ et il reste donc à  définir
les produits d'images par $I$.
Les ensembles $I(\H(r_0))$ et $I(\G(V))$ étant naturellement munis d'une structure d'algèbre, 
il nous suffit donc de définir le produit d'un élément $I(a)(x,\e)$, $a\in\H(r_0)$ par un élément $I(g)(x,\e)$,
$g\in\G(V)$.
Pour cela, en notant
$
a(x)=\ds\sum_{\nu=0}^\infty a_\nu x^{\nu}\mbox{  et  }g(X)\sim\ds\sum_{\nu>0}g_\nu X^{-\nu},
$
on observe que $a(x)g\big(\frac x\e\big)=\big( a_0+x\S a(x)\big)\,g\big(\frac x\e\big)$ et $x\,g\big(\frac x\e\big)=\lp g_1+\T g\big(\frac x\e\big)\rp\e$.
Autrement dit, un développement asymptotique du produit de fonctions 
$a(x)\,g\gk{\frac x\e}$ par rapport à  $\e$ peut être obtenu par
\eq{prodcomb}{
a(x)\,g\big(\ts\frac x\e\big)=a_0\,g\big(\ts\frac x\e\big)+g_1\,\S a(x)\e+
\S a(x)\,\T g\big(\ts\frac x\e\big)\e.
}
En itérant la formule correspondante pour les images par $I$, 
on définit ainsi, avec la convention $g_0=0$
\eq{dprod}{
I(a)(x,\e)\,I(g)(x,\e)=\ds\sum_{\nu\geq0}\lp g_\nu\,(\S^\nu a)(x) + a_\nu\,(\T^\nu g)\big(\ts\frac x\e\big)\rp\e^\nu.
}
\rqs
 1.\ Les séries combinées classiques  \cit{vb,bef} entrent dans le cadre des nôtres~: il s'agit du cas où les fonctions $g_n$ sont à  décroissance exponentielle. Rappelons qu'une fonction $g:J=]\mu,+\infty[\to\R$ est {\sl à  décroissance exponentielle} s'il existe $C,A>0$ vérifiant 
$$
\forall X\in J,\ |g(X)|\leq C\exp(-AX).
$$
Une fonction $g$ à  décroissance exponentielle satisfait en particulier 
$g(X)=\O(X^{-N}),$ $X\to+\infty$ pour tout entier $N$, donc est {\sl plate}~: elle admet la série nulle pour développement asymptotique. 
\med\\
2.\ Dans le cas des séries combinées classiques, la partie  lente d'un produit ne dépend que des parties lentes des facteurs. Ceci peut se voir par exemple grâce à  \rf{dprod}~: lorsque tous les $g_\nu$ sont nuls, le produit d'un terme lent et d'un terme rapide n'engendre que des termes rapides; \cf aussi 
la remarque (2), page 11 de \cit{bef}. 
En revanche, dans le cas des séries combinées du présent article, la formule \rf{dprod} montre que le produit d'un terme lent avec un terme rapide fait aussi apparaître des termes lents, si bien que tout est imbriqué.
\med\\
Les séries combinées sont aussi compatibles avec la composition à  gauche.

\lem{lm2.1}{
Soit $\yc\in\e\cch(r_0,V)$ une série formelle combinée sans terme constant, \ie avec $a_0(x)\equiv0$ et $g_0(X)\equiv0$.
Soit $\Pc\in\cch(r_0,V)[\!\,[y,\e]\!\,]$ 
une série formelle à  deux variables dont les coefficients sont des séries formelles combinées : $\Pc(x,y,\e)=\ds\sum_{j,k\geq0}p_{j,k}(x,\e)y^j\e^k$ avec $p_{j,k}\in\cch(r_0,V)$. 
Alors l'expression
$\P(\yc)(x,\e)=\Pc(x,\yc(x,\e),\e):=\ds\sum_{j,k\geq0}p_{j,k}(x,\e)\yc(x,\e)^j\e^k$
définit une série formelle combinée.
De plus, l'application $\P:\e\cch(r_0,V)\to\cch(r_0,V)$ ainsi définie est 
$1$-lipschitzienne.
}

La preuve est immédiate, grâce à la convergence de la série définissant $\P(\yc)(x,\e)$.
\med

Nous voulons à  présent définir la dérivée d'une série formelle combinée.
Du fait que les dérivées de fonctions de $\H(r_0)$ ne sont plus nécessairement bornées
et que celles de fonctions de $\G(V)$ n'ont plus nécessairement un développement asymptotique,
cette dérivation est un peu plus délicate à  traiter, bien que la formule soit simple. En particulier, cela nécessite de réduire un peu les domaines de définition des fonctions.
Dans le cadre réel, il n'est pas toujours possible de définir la dérivée d'une série formelle combinée.
\df{deriv}{
Si $\yc\in\cch(r_0,V(\a,\b,\infty,\mu))$  est une série formelle combinée donnée par 
$$
\yc(x,\e)=\sum_{n\geq0}\lp a_n(x)+g_n\big(\ts\frac x\e\big)\rp\e^n,
$$
telle que le premier terme rapide $g_0$ est identiquement nul, alors sa {\sl dérivée par rapport à  $x$}, $\frac{d\yc}{dx}$, est donnée par
$$
\frac{d\yc}{dx}(x,\e)=\sum_{n\geq0}\lp a_n'(x)+g_{n+1}'\big(\ts\frac x\e\big)\rp\e^n\ \ ;
$$
Ceci définit un élément de $\cch(\~ r_0,V(\~\a,\~\b,\infty,\~\mu))$ pour tout
$\~ r_0\in\,]0,r_0[$, $\a<\~\a<\~\b<\b$  et tout 
$\~\mu<\mu$.}
\rq La dérivée d'une série formelle combinée n'est définie \apriori{} que si le premier terme rapide $g_0$ est identiquement nul. En revanche, l'opérateur $\e\,\frac d{dx}=\frac d{dx}(\e\;\cdot\ )$ est défini sans condition sur $g_0$.
\sub{2.2bis}{Développements combinés}
Jusqu'à  présent, les objets que nous avons considérés étaient des expressions formelles.
Nous voulons maintenant définir le développement combiné d'une fonction de deux variables $x$ et $\e$. Le plus simple et le plus naturel serait de considérer des fonctions définies sur un produit de secteurs en $x$ et en $\e$. 
Cependant, pour certaines applications, 
il sera commode que le domaine par rapport à  $x$ contienne un voisinage de $0$ de taille proportionnelle à  $\norm\e$. Pour d'autres applications, au contraire, il nous sera nécessaire d'ôter un voisinage de $0$. C'est la raison pour laquelle nous avons introduit les quasi-secteurs \rf{V1} et \rf{V2}. 
\df{d2.2}{Soit $V=V(\a,\b,\infty,\mu)$ un quasi-secteur infini, soit
$S_2=S(\a_2,\b_2,\e_0)$ un secteur fini et soit $\a_1<\b_1$ tels que
$\a\leq\a_1-\b_2<\b_1-\a_2\leq\b$. Soit
$y(x,\e)$
une fonction holomorphe définie pour $\e\in S_2$ et $x\in V(\a_1,\b_1,r_0,\mu\norm\e)$.
Enfin, soit 
$\yc(x,\e)=\ds\sum_{n\geq0}\gk{a_n(x)+g_n\big(\ts\frac x\e\big)}\e^n\in\cch(r_0,V)$.
Nous disons que {\sl $y$ admet $\yc$ comme \dac{}}
et nous écrivons {\sl $y(x,\e)\sim\yc(x,\e)$, $S_2\ni\e\to0$, $x\in V(\a_1,\b_1,r_0,\mu\norm\e)$},
si, pour tout entier $N$, il existe une constante $K_N$ telle que, pour tout
$\e\in S_2$ et tout $x\in V(\a_1,\b_1,r_0,\mu\norm\e)$
\eq{defcomb}{\norm{y(x,\e)-\sum_{n=0}^{N-1}\gk{a_n(x)+g_n\big(\ts\frac x\e\big)}\e^n}\leq K_N \norm\e^N  .}}
Ici encore, les fonctions $a_n$ forment la  {\sl partie lente} du \dac{} et les $g_n$ la {\sl partie rapide}.
Les conditions sur les angles $\a_j,\b_j$ assurent l'implication~: si $\e\in S_2$ et $x\in V(\a_1,\b_1,r_0,\mu\norm\e)$ alors $x/\e\in V$.
\med\\
\rqs 1.\ Dans le cas d'une couronne $V=V(0,2\pi,\infty,-r)$, $r>0$, 
il n'y a pas besoin de
conditions sur les angles. Un développement combiné 
$y(x,\e)\sim\sum_{n\geq0}(a_n(x)+g_n(\tfrac x\e))\e^n$ est alors une autre forme
d'un développement monomial introduit par \cit{cms}.  Quand on pose $u=\e/x$,
alors $\e=xu$ et la fonction $z(x,u)=y(x,xu)$ est définie sur un
{\sl secteur en $xu$} défini dans \cit{cms} après la définition 3.4, \ie l'ensemble des 
$(x,u)$ tels que $\norm x<r_0, \norm u<\min\big(\frac1r,\frac{\e_0}{r_0}\big)$ et 
$\arg(xu)\in]\a_1,\b_1[$, et y admet le développement 
$z(x,u)\sim\sum_{n\geq0}\big(a_n(x)+b_n(u)\big)(xu)^n$ 
défini dans \cit{cms}, définition 3.6, avec $b_n(u)=g_n(1/u)$.
\med\\
2.\
Dans un but de simplicité, nous demandons aux fonctions $a_n$ d'être holomorphes dans tout le disque $D(0,r_0)$, alors que la fonction $y$ elle-même n'est définie que pour $x\in V(\a_1,\b_1,r_0,\mu\norm\e)$. Dans la partie \reff{2.3} nous serons amenés à  généraliser la définition de \dac{} au cas où les fonctions $a_n$ sont holomorphes dans un domaine contenant $0$, \cf la remarque après la proposition \reff{p3.11bis}.
\med\\
3.\
Une fonction $y(x,\e)$  ne peut pas avoir deux \dacs{} différents quand $S_2\ni\e\to0$ et $x\in V(\a_1,\b_1,r_0,\mu\norm\e)$. En effet, on a $\ds\lim_{\e\to0}y(x,\e)=a_0(x)$ pour 
$x\in S(\a_1,\b_1,r_0)$ donc la fonction holomorphe $a_0\in \H(r_0)$ est déterminée de 
manière unique, donc aussi $a_0(0)$. On continue avec
$\ds\lim_{\e\to0}y(\e X,\e)=a_0(0)+g_0(X)$, et ainsi de suite. 
Il convient de noter que, pour démontrer cette unicité, seule la propriété que les $g_n(X)$ tendent vers $0$ quand $X\to\infty$ a été utilisée~; ceci sera utile dans la 
partie \reff{3.2}.
\med

Il est immédiat que les \dacs{}
sont compatibles avec l'addition et la multiplication scalaire.
Pour la compatibilité avec la multiplication entre développements, le seul point moins évident est de montrer qu'un produit $a(x)g\big(\ts\frac x\e\big)$, $a\in\H(r_0)$, $g\in\G(V)$ 
admet un \dac{}. Ceci est une conséquence de la formule
\rf{prodcomb} et du fait que $\S$ et $\T$ sont des endomorphismes.
La définition \rf{dprod} était faite pour que l'on ait
$a(x)\,g\big(\ts\frac x\e\big)\sim I(a)(x,\e)\,I(g)(x,\e)$.
\bigskip

{\noi\bf Composition} \sep
\med\\
Les \dacs{} sont aussi compatibles avec la composition à  gauche et à  droite par une fonction holomorphe, comme l'exprime la proposition suivante. L'énoncé (a) concerne la composition à  gauche par une fonction de trois variables, mais dans le cas d'un \dac{} sans terme en $\e^0$. L'énoncé (b) traite le cas de la composition à  gauche  sans cette restriction, mais par une fonction d'une variable uniquement. Ces deux énoncés sont donc complémentaires. Pour la composition à  droite, nous n'avons considéré que les fonctions de la seule variable $x$ pour une raison de simplicité, mais il est possible de généraliser le résultat au cas d'une fonction $\f$ des deux variables $x$ et $\e$, telle que $\f(0,0)=0$ et $\frac{\partial\f}{\partial x}(0,0)=1$. La version Gevrey de cette généralisation est faite dans la partie \reff{4.5}, théorème \reff{t4.7}.
Nous n'avons pas mis au point d'énoncé à  propos d'un changement de variable portant sur $\e$ car nous n'en avons pas l'utilité.
\propo{l2.4}{
\be[\rm(a)]\item
 Soit $P(x,z,\e)$ une fonction holomorphe définie quand
 $ |z| < r$, $\e\in S_2=S(\a_2,\b_2,\e_0)$ et $x\in V(\a_1,\b_1,r_0,\mu|\e|)$ 
telle que tous les coefficients $P_n$ du développement 
$P(x,z,\e)=\sum_{n\geq0}P_n(x,\e)z^n$
admettent un \dac{} $P_n(x,\e)\sim\widehat P_n(x,\e)$ quand $S_2\ni\e\to0$,
$x\in V(\a_1,\b_1,r_0,\mu|\e|)$. 
Soit $y(x,\e)={\cal O}(\e)$ une fonction admettant un 
\dac{} $\yc(x,\e)$ quand $S_2\ni\e\to0$ et 
$x\in V(\a_1,\b_1,r_0,\mu|\e|)$ sans termes en $\e^0$. On suppose que 
$\sup_{x,\e}|y(x,\e)|<r$.
Alors la fonction $u:(x,\e)\mapsto P(x,y(x,\e),\e)$ admet le \dac{}
$$
\P(\yc)(x,\e)=\sum_{n\geq0}\widehat P_n(x,\e)\yc(x,\e)^n.
$$
\item
 On considère une fonction holomorphe $y$ définie pour 
$\e\in S_2$ et $x\in V=V(\a_1,\b_1,r_0,\mu\norm\e)$, à  valeurs dans un ensemble 
$W\subset\C$
 et  admettant un \dac{} quand $\e\to0$. Soit $f$ une fonction holomorphe sur un 
voisinage de l'adhérence de $W$. 
Alors la fonction $z=f\circ y$ admet un \dac{} quand $S_2\ni\e\to0$, $x\in V$.
\item
Soit $\f$ une fonction holomorphe définie quand $\norm x<x_1$ 
telle que $\f(0)=0$ et $\f'(0)=1$ et soit $z=z(u,\e)$ une fonction ayant un \dac{}
$\sum_{n\geq0}\Big(a_n(u)+g_n\big(\ts\frac u\e\big)\Big)\e^n$ quand 
$S_2\ni\e\to0$ et $u\in V(\a_1,\b_1,r_0,\mu|\e|)$, avec $a_n\in\cch$ et $g_n\in\G$. 
Alors pour tout $\~\a_1,\~\b_1$ avec $\a_1<\~\a_1<\~\b_1<\b_1$ et tout $\tilde\mu<\mu$
il existe $\~r,\~\e_0>0$ tels que la fonction 
$y:(x,\e)\mapsto z(\f(x),\e)$ admet un \dac{} quand $S(\a_2,\b_2,\~\e_0)\ni\e\to0$ et 
$x\in V(\~\a_1,\~\b_1,\~r,\~\mu|\e|)$.
\ee
}
\rq
 Dans le (a), l'hypothèse ``$y$ bornée par $r$'' n'est pas essentielle~: 
il suffit de réduire le domaine de $u$ si elle n'est pas satisfaite.
\med\\
\pr{Preuve}
(a) Pour tout $N\in\N^*$, la somme finie $\ds\sum_{0\leq n\leq N-1}P_n(x,\e)y(x,\e)^n$ admet
un \dac{} (compatibilité avec produit et somme).
Il reste à  vérifier qu'il existe  une constante $L$ telle que le reste 
est borné par $L \norm\e^N$. Ceci est évident d'après les hypothèses.
\med\\
(b) Quitte à  modifier $f$ et $y$, 
on peut supposer que $a_0(0)=0$. Moyennant un développement de Taylor,
il suffit de montrer que $f\big(a_0(x)+g_0(\tfrac x\e)\big)$ admet un \dac{}.

Notons $h(u,v)=f(u+v)$. Il suffit donc de montrer que $h\big(a_0(x),g_0(\tfrac x\e)\big)$ admet
un DAC quand $\e$ tend vers $0$. Pour montrer ceci, on remarque qu'on peut écrire
$$h(x,y)=h(x,0)+h(0,y)-h(0,0)+xy\,k(x,y)$$
avec une certaine fonction holomorphe $k$ des deux variables $x,y$.
On a donc
$$
\begin{array}{rl}
h\big(a_0(x),g_0(\tfrac x\e)\big)=\!\!&
h(a_0(x),0)+h\big(0,g_0(\tfrac x\e)\big)-h(0,0)\med\\
&+ 
a_0(x)g_0(\tfrac x\e)\,k\big(a_0(x), g_0(\tfrac x\e)\big).
\end{array}
$$
Puisque $a_0(0)=0$, le produit $ a_0(x)g_0\big(\tfrac x\e\big)$ est de la forme $\O(\e)$~;
 on  obtient donc un \dac{} pour $h\big(a_0(x),g_0(\tfrac x\e)\big)$ en itérant cette procédure.

En passant, remarquons que le terme principal du \dac{} de $f(y(x,\e))$ (sans la réduction à  $a_0(0)=0$) a pour partie lente 
$f(a_0(x))$ et pour partie rapide $f\big(a_0(0)+g_0(\tfrac x\e)\big)-f(a_0(0))$.
\med\\
(c) Si $\~r,\~\e_0$ sont assez petits, alors $\f\big(V(\~\a_1,\~\b_1,\~r,\~\mu|\e|)\big)
\subset V(\a_1,\b_1,r_0,\mu|\e|)$ si $\e\in S_2,\norm\e<\~\e_0$. Il suffit alors
de montrer que $b\big(\frac{\f(x)}\e\big)$ admet un \dac, si $b$ est dans $\cal G$.
Il convient d'introduire les fonctions $\psi$ et $h$ définies par  
$\frac1{\f(x)}-\frac1x=h(x)$  et $\psi(x,t)=x/(1+txh(x))$. La fonction $h$ se prolonge en une fonction holomorphe définie quand $\norm x < x_1$, notée encore 
$h$ par abus de notation, et $\psi(x,0)=x,\psi(x,1)=\f(x)$. 
Le développement de Taylor de $b\big(\frac{\f(x)}\e\big)=b\big(\frac{\psi(x,1)}\e\big)$
par rapport à  $t$ donne pour tout $N\in\N$
$$b\big(\tfrac{\f(x)}\e\big)= \sum_{n=0}^{N-1} \tfrac{1}{n!}\,\tfrac{\partial^n}
   {\partial t^n}b\big(\tfrac{\psi(x,t)}\e\big)\mid_{t=0} + 
                   \tfrac{1}{(N-1)!}\,\int_0^1\tfrac{\partial^N}
   {\partial t^N}b\big(\tfrac{\psi(x,t)}\e\big)\mid_{t=\tau} (1-\tau)^{n-1}
 \,d\tau.  $$
En utilisant le fait que $\frac\partial{\partial t}\big[f\big(\frac{\psi(x,t)}\e\big)\big]=
\e h(x)(\De f)\big(\frac{\psi(x,t)}\e\big)$ avec l'opérateur 
$\De$ défini par $(\De f)(X)=-X^2f'(X)$, on obtient
\eq{compdroite}{
\begin{array}{rl}
b\big(\tfrac{\f(x)}\e\big)=\!\!&\ds \sum_{n=0}^{N-1} \tfrac{\e^n}{n!}\,h(x)^n(\De^nb)\big(\tfrac x\e\big)\med\\
&\ds  + 
\tfrac{\e^N}{(N-1)!}h(x)^N\,\int_0^1 (\De^Nb)\big(\tfrac{\psi(x,\tau)}\e\big) (1-\tau)^{n-1}
\,d\tau
\end{array}
}
et on peut vérifier que le dernier terme est ${\cal O}(\e^N)$. 
La compatibilité des \dac{} avec l'addition et la multiplication entraîne alors l'existence d'un \dac{} pour $b\big(\frac{\f(x)}\e\big)$.
\ep
\bigskip

{\noi\bf Dérivation} \sep
\lem{derivasympt}{Soit $y(x,\e)$
une fonction définie quand $\e\in S_2$ et $x\in V(\a_1,\b_1,$ $r_0,\mu\norm\e)$ telle que
$y(x,\e)\sim\sum_{n\geq0}\gk{a_n(x)+g_n\big(\ts\frac x\e\big)}\e^n=:\yc(x,\e)\in\cch(r_0,V)$
quand $S_2\ni\e\to0$.
On suppose que $g_0(X)\equiv0$.

Alors, pour tout $\~r_0$ entre $0$ et $r_0$ et tout $\~\mu<\mu$, on a
$$\frac{dy}{dx}(x,\e)\sim \frac{d\yc}{dx}(x,\e)=\ds\sum_{n\geq0}
      \gk{a_n'(x)+g_{n+1}'\big(\ts\frac x\e\big)}\e^n $$
quand $S_2\ni\e\to0$ et $x\in  V(\a_1,\b_1,\~r_0,\~\mu\norm\e)$, avec $\ds\frac{d\yc}{dx}(x,\e)\in\cch(r_0,\~V)$, où $\~V:=V(\a,\b,\infty,\~\mu)$.} 
\pr{Preuve}
Notons $\d=\min(|\e|(\mu-\~\mu),r_0-\~r_0)$ et, pour $N\in\N$ arbitraire,
\eq{rest}{
R_N(x,\e)=y(x,\e)-\ds\sum_{n<N}\gk{a_n(x)+g_n\big(\ts\frac x\e\big)}\e^n.
}
La formule de Cauchy pour la dérivée donne 
$$
\bigg|\frac{dR_{N+1}}{dx}(x,\e)\bigg|=\bigg|\frac1{2\pi i}\int_{|u-x|=\d}\frac{R_{N+1}(u,\e)}{(u-x)^2}du\bigg|\leq\frac1\d\max_{|u-x|=\d}|R_{N+1}(u,\e)|,
$$
ce qui donne $\ds\frac{dR_{N+1}}{dx}(x,\e)=O(|\e|^{N})$. Puisque les fonctions $a'_N$ et $g'_{N+1}$ sont elles-mêmes bornées sur $D(0,\~r_0)$, resp. sur $\~V$,  on en déduit $\ds\frac{dR_{N}}{dx}(x,\e)=O\big(|\e|^{N}\big)$.
\ep

\rq Le lemme \reff{derivasympt} n'est pas valable dans le cadre réel. 
Par exemple, il est bien connu qu'il existe de petites
fonctions avec des dérivées non bornées. 
Cependant on a le résultat suivant en réel~: 
si la dérivée d'une fonction
ayant un \dac{} admet elle aussi un \dac{}, alors la formule \rf{int} ci-après implique que
le \dac{} de la dérivée peut être obtenu en dérivant terme à  terme le \dac{}
de la fonction.
\bigskip

{\noi\bf Intégration} \sep
\med\\
Une difficulté pour la compatibilité avec l'intégration provient du fait que
les fonctions de $\cal G$ ayant un terme avec $1/X$ dans leur développement asymptotique à  l'infini ne possèdent pas de primitive dans $\cal G$.
Lorsque ces termes sont tous nuls, l'intégration ne pose pas de problème, comme l'indique l'énoncé qui suit. 
\propo{p40}{Considérons un \dac{}
$y(x,\e)\sim\sum_{n\geq0} \Big(a_n(x)+g_n\big(\ts\frac x\e\big)\Big)\e^n$ défini pour $\e\in S_2$ et $x\in V(\a_1,\b_1,r_0,\mu\norm\e)$, tel que toutes les fonctions $g_n$ satisfont $g_n(X)={\cal O}(X^{-2})$ quand $X\to\infty$.

Alors la fonction $(x,\e)\mapsto\ds\int_r^x y(t,\e)\,dt$ admet un \dac. Précisément, on a 
\eq{int}{\begin{array}c \ds\int_r^x y(t,\e)\,dt\sim \hat Y(x,\e)-\hat Y(r,\e),\med\\ \mbox{ où }
   \hat Y(x,\e)=A_0(x)+\ds\sum_{n=1}^{\infty} \Big(A_n(x)+
  G_{n-1}\big(\tfrac x\e\big)\Big)\e^n
\end{array}  
}
avec $A_n(x)=\ds\int_r^x a_n(t)\,dt$ et $G_n(X)=-\ds\int_X^\infty g_n(T)\,dT$. }
 Ici on a identifié
$\hat Y(r,\e)$  avec la série formelle obtenue en remplaçant $G\big(\frac r\e\big)$ par son développement asymptotique quand $\e\to0$. La preuve est immédiate~: on a bien
$A_n\in{\cal A}$ et, d'après l'hypothèse, $G_n\in{\cal G}$. 
Avant d'énoncer le résultat dans le cas général, nous devons encore introduire une notation.
Soit $\ell$ une fonction analytique dans le quasi-secteur $V=V(\a,\b,\infty,\mu)$ telle que sa dérivée $\ell'$ 
a un développement asymptotique à  l'infini commençant par $\frac1X$~: 
$$
\ell'(X)\sim\sum_{n\geq1}c_nX^{-n}\mbox{ avec }c_1=1.
$$
 On peut choisir par exemple $\ell(X)=\log(X-\g)$ avec $\g\notin V$.
Dans une situation réelle au départ, on pourra prendre $\ell(X)=\frac12\log(X^2+L^2)$ avec $L$ assez grand. Nous serons aussi amenés à  considérer une fonction $\ell$ dépendant des deux variables $X$ et $\e$. Pour simplifier, et parce que cela nous suffira, nous ne considérerons $\ell$ que sous la forme de la somme d'une fonction ne dépendant que de $X$ et d'une fonction ne dépendant que de $\e$. Dans l'application à  la résonance, nous choisirons   $\ell(X,\e)=\frac1p\log(X^p+1)+\log(\e)=\frac1p\log(x^p+\eps)$.
L'énoncé dans le cas général est le suivant.
\propo{p4}{
Étant donné un \dac{} $y(x,\e)\sim\sum_{n\geq0}\Big(a_n(x)+g_n\big(\ts\frac x\e\big)\Big)\e^n$,
notons $\hat R$ la série des résidus
des $g_n(X)$~: $\hat R(\e)=\sum_{n\geq0} g_{n1}\e^n$.
Alors on a
$\ds\int_r^x y(t,\e)\,dt \sim \hat Y(x,\e)-\hat Y(r,\e)$,
avec
\eq{int2}{
\begin{array}{rcl}
\hat Y(x,\e)\!\!&=\!\!&\e \hat R(\e) \Big(\ell\big(\tfrac x\e\big)-\ell\big(\tfrac r\e\big)\Big)+ 
  A_0(x)+
\med\\
&&
 \ds\sum_{n=1}^{\infty} \Big(A_n(x)+H_{n-1}\big(\tfrac x\e\big)\Big)\e^n,
\end{array}
}
où $A_n(x)=\ds\int_r^x a_n(t)\,dt$
et $H_{n}(X)=-\ds\int_X^\infty  \gk{g_n(T)-g_{n1}\ell'(T)}\,dT$. 
}
Ici aussi, nous avons identifié $\hat Y(r,\e)$ avec la série formelle
obtenue en  remplaçant $H_n\big(\frac r\e\big)$ par
son développement quand $\e$ tend vers $0$.
\med\\
\pr{Preuve}Le théorème de Borel-Ritt classique (voir ci-dessous) nous fournit une fonction  $R(\e)$ ayant $\hat R(\e)$ pour développement asymptotique. La différence
$y(x,\e)-  R(\e) \ell'\big(\frac x\e\big)$ satisfait la condition de la proposition \reff{p40}, donc son intégrale admet un \dac. Le ``\dac{} généralisé'' pour $y$ s'en déduit.\ep

On a aussi un énoncé analogue au théorème classique de Borel-Ritt. Ce théorème affirme, pour toute suite $(a_n)_{n\in\N}$ de nombres complexes et tout secteur $S(\a,\b,\e_0)$, l'existence d'une fonction $a=a(\e)$ définie sur $S$ et admettant la série formelle $\sum_{n=0}^\infty a_n\e^n$ pour développement asymptotique. Le résultat est aussi vrai lorsque  $(a_n)_{n\in\N}$ est une suite de fonctions analytiques d'une variable complexe $x$. Dans le cas de nos \dacs, l'énoncé est le suivant.
\lem{borel-ritt}{
(Borel-Ritt)
Soit $V=V(\a,\b,\infty,\mu)$ un quasi-secteur infini ($\mu>0$ ou $<0$),
$S_2=S(\a_2,\b_2,\e_0)$ un secteur fini, $r_0>0$ et soit $\a_1<\b_1$ tels que
$\a\leq\a_1-\b_2<\b_1-\a_2\leq\b$. Étant donnée une série formelle combinée
$\yc(x,\e)=\sum_{n\geq0}\gk{a_n(x)+g_n\big(\ts\frac x\e\big)}\e^n\in\cch(r_0,V)$, il existe
une fonction $y(x,\e)$ holomorphe définie pour $\e\in S_2$ et $x\in V(\a_1,\b_1,r_0,\mu\norm\e)$
telle que $y(x,\e)\sim\yc(x,\e)$ quand $\e\to0$.
}
\pr{Preuve}
Il suffit d'utiliser le théorème de Borel-Ritt pour des développements asymptotiques uniformes classiques deux fois : une fois pour $\sum a_n(x)\e^n$, une fois pour $\sum g_n(X)\e^n$.
\ep
\sub{2.3}{Développements combinés et ``matching''}
Notre notion de développement combiné mélange la notion classique de développement asymptotique au sens de Poincaré que nous appellerons ``extérieur'', de la forme $y(x,\e)\sim\sum_{n\geq0}c_n(x)\e^n$, et celle de développement dit ``intérieur'' de la forme $y(\e X,\e)\sim\sum_{n\geq0}h_n(X)\e^n$. Ces développements intérieurs et extérieurs occupent une place centrale dans la méthode de recollement des développements asymptotiques (méthode des ``matched asymptotic expansions'' en anglais).
Bien que les \dacs{} soient de nature différente, les liens avec les développements intérieurs et extérieurs sont étroits.

Nous montrons d'une part qu'une fonction ayant un \dac{} admet aussi un développement intérieur et un développement extérieur, et que ces deux développements ont une région commune de validité. En d'autres termes, 
une preuve d'existence d'un \dac{}  permet de donner un fondement solide à la méthode de recollement.

D'autre part la réciproque est vraie~: si la méthode de recollement est valide, \ie si une fonction  a des développements intérieur et extérieur avec une région commune de validité, et si de plus ces développements  vérifient une propriété supplémentaire, alors la fonction a aussi un \dac.
Le premier résultat est le suivant.
\propo{combextint}{
Soient $(a_n)_{n\in\N}$ une famille de fonctions de $\H(r_0)$ et $(g_n)_{n\in\N}$ une famille de fonctions de $\G(V)$ avec $V=V(\a,\b,\infty,\mu)$.
On note leurs développements $a_n(x)=\ds\sum_{m=0}^\infty a_{nm}x^m$ et $g_n(X)\sim\ds\sum_{m>0}g_{nm}X^{-m}$.
Supposons que $$y(x,\e)\sim \ds\sum_{n\geq0}\gk{a_n(x)+g_n\big(\ts\frac x\e\big)}\e^n$$ 
quand $S_2\ni\e\to0$ et $x\in V(\a_1,\b_1,r_0,\mu\norm\e)$ au sens de la définition 
\reff{d2.2}. 

Alors, pour $x\in S(\a_1,\b_1,r_0)$ fixé, on a
\eq {cc}{
y(x,\e)\sim\sum_{n\geq0}c_n(x)\e^n \mbox{ \ quand }S_2\ni\e\to0,
}
où $c_n(x)=a_n(x)+\ds\sum_{0\leq l\leq n-1} g_{l,n-l}x^{l-n}$. De plus, pour tout $r>0$, ce développement est uniforme par rapport à  $x$ sur l'ensemble des $x\in S(\a_1,\b_1,r_0)$ tels que $\norm x > r$.

De même, si $X\in V$ et $\a_3,\b_3,\e_3$ sont tels que $\e\in S(\a_3,\b_3,\e_3)$ implique
$\e\in S_2$ et $\e X \in V(\a_1,\b_1,r_0,\mu\norm\e)$, alors on a
\eq{hhh}{
y(\e X,\e)\sim\sum_{n\geq0}h_n(X)\e^n\mbox{ quand }S(\a_3,\b_3,\e_3)\ni\e\to0,
}
où $h_n(X)=g_n(X)+\ds\sum_{0\leq l\leq n} a_{n-l,l}X^l$. Ce développement est uniforme par rapport à  $X$ sur des parties compactes de $V$ vérifiant la condition précédente.
}
\rqs
1.\ 
Conformément à la littérature, 
nous appellerons le premier développement \rf{cc} le {\sl développement extérieur}
et le second  \rf{hhh} le {\sl développement intérieur}. Chaque fonction $c_n$ du développement extérieur peut avoir une singularité en $x=0$ mais seulement un pôle d'ordre au plus $n$~; de même chaque fonction $h_n$ du développement intérieur a une croissance polynomiale d'ordre au plus $n$ lorsque $X\to\infty$.
Le {\sl restreint index} au sens de Wasow \cit{wa1}, chapitre VIII est donc égal à 1.
\med\\
2.\ 
On peut montrer que, pour tout $\kappa\in\,]0,1[$, le  développement extérieur \rf{cc} est uniforme sur $\norm x>\norm\e^\kappa$,
et que le développement intérieur \rf{hhh} est uniforme sur $\norm X < \norm\e^{-\kappa}$, ce qui justifie
la méthode de ``matched asymptotic expansions'' lorsqu'un \dac{} existe. 
Dans chacun de ces cas, pour obtenir une approximation à  l'ordre $N$, il convient alors de sommer jusqu'à  l'ordre $\frac N{1-\kappa}$.
Cependant, il est souvent préférable d'avoir des approximations uniformes à sa disposition sur tout le domaine d'étude. Ceci semble même indispensable si on veut  obtenir des estimations de type Gevrey.
\med\\
3.\ 
Dans les cas où on peut démontrer indirectement 
l'existence d'un développement combiné pour une fonction
$y(x,\e)$, mais qu'on ne connaît pas encore les fonctions $a_n$ et $g_n$, une 
méthode pour les déterminer est d'appliquer la proposition précédente. Pour $x$ 
fixé différent  de $0$, on calcule le développement extérieur $y(x,\e)\sim\sum_{n\geq0}c_n(x)\e^n$, puis on rejette les termes avec des puissances négatives. On obtient ainsi les $a_n(x)$. Ensuite, on calcule
le  développement intérieur $y(\e X,\e)\sim\sum_{n\geq0}h_n(X)\e^n$ et on rejette les termes de puissances positives de $X$, ce qui donne les $g_n(X)$.

Dans des situations concrètes, le calcul des développements extérieur et intérieur mène souvent à  des équations de récurrence pour leurs coefficients. Ceci permet donc le calcul des $a_n,g_n$ sans avoir à utiliser les formules techniques pour la multiplication de séries formelles combinées.

Dans le cas des équations différentielles singulièrement perturbées, comme le remarquent Emmanuel Isambert et Véronique Gautheron dans \cit{gi}, le calcul du  développement extérieur ne fait intervenir que des opérations algébriques (une fois donnés les développement de Taylor des coefficients de l'équation) alors que celui du développement intérieur nécessite d'intégrer des équations différentielles linéaires, puis de choisir la constante d'intégration pour que la  solution ait un comportement asymptotique bien déterminé, ce qui introduit de la transcendance. Dans \cit{i},  Emmanuel Isambert  appelle ainsi ces développements extérieur et intérieur respectivement {\sl algebraic} et {\sl transcendental expansions}.
\med\\
\pr{Preuve}
Fixons $N\in\N^*$ et reprenons
 la notation \rf{rest}.
Notons de plus 
$$
r_{lk}(X)=g_l(X)-\ds\sum_{0<m<k}g_{lm}X^{-m}.
$$
Par hypothèse, il existe des constantes positives $C_N,\,A_{kn}$ et $C_{lk}$ telles que
$$
\forall\e\in S_2\;\forall x\in V_{1,\e}:=V(\a_1,\b_1,r_0,\mu|\e|)\qquad |R_N(x,\e)|\leq C_N|\e|^N,
$$
\eq{AA}{
\forall x\in V_{1,\e}\qquad \Big|a_k(x)-\sum_{l<n}a_{kl}x^l\Big|\leq A_{kn}|x|^n
}
et
\eq{34b}{
\forall X\in V\qquad|r_{lk}(X)|\leq C_{lk}|X|^{-k}.
}
Un calcul élémentaire donne
$$
y(x,\e)-\sum_{n<N}c_n(x)\e^n=R_N(x,\e)+\sum_{l<N}g_l\big(\ts\frac
x\e\big)\e^l-\ds\sum_{0<n<N}\bigg(\sum_{l<n}g_{l\,n\!-\!l}x^{l-n}\bigg)\e^l
$$
\eq{34c}{
=R_N(x,\e)+\sum_{l<N}r_{l\,N\!-\!l}(\ts\frac x\e\big)\e^l,
}
donc, quand $\norm x > r$,
\eq{cex}{
\Big|y(x,\e)-\sum_{n<N}c_n(x)\e^n\Big|\leq\bigg( C_N+\sum_{l<N}C_{l\,N\!-\!l}\,r^{l-N}\bigg)|\e|^N.
}
De même, on a 
$$
\begin{array}{rl}
y(\e X,\e)-\sum_{n<N}h_n(X)\e^n\!\!&\ds=R_N(\e X,\e)+\sum_{n<N}\bigg(a_n(\e X)-\sum_{l\leq n}a_{n\!-\!l\,l}X^l\bigg)\e^n\med\\&\ds
=R_N(\e X,\e)+\sum_{k<N}\bigg(a_k(\e X)-\sum_{l< N-k}a_{kl}\,\e^lX^l\bigg)\e^k
\end{array}
$$
d'où, pour tout $R>0$ et pour $\norm X\leq R$,
\eq{cin}{
\Big|y(\e X,\e)-\sum_{n<N}h_n(X)\e^n\Big|\leq\bigg(C_N+\sum_{k<N}A_{k\,N\!-\!k}R^{N-k}\bigg)|\e|^N.
}
\ep

Concernant la réciproque, l'énoncé est le suivant.
\propo{p3.10bis}{
Soit $y$ une fonction définie pour $\e\in S_2=S(\a_2,\b_2,\e_0)$ et 
$x\in V(\e)= V(\a_1,\b_1,r_0,\mu\norm\e)$. 
On suppose qu'il existe des nombres réels $a,b,\kappa$ avec $0<a<b$ et $0<\kappa<1$, et 
pour chaque $n\in\N$ une fonction $c_n$, $c_n(x)=P_n\big(\frac1x\big)+a_n(x)$, 
$P_n$ polynomial sans terme constant, $a_n\in\H(r_0)$ et une fonction 
$h_n=Q_n+g_n$, $Q_n$ polynomial et $g_n\in\G(V)$, $V= V(\a,\b,\infty,\mu)$,
$\a\leq\a_1-\b_2<\b_1-\a_2\leq\b$, avec les 
propriétés suivantes.
\be\item
Pour tout $N\in\N$, il existe une constante $C>0$ telle que
\eq{ccuni}{\norm{y(x,\e)-\sum_{n=0}^{N-1}c_n(x)\e^n}\leq C \norm\e^{N(1-\kappa)}}
pour tout $\e\in S_2$ et $x\in V(\e)$ avec $|x|> a|\e|^\kappa$ et
\eq{hhuni}{\norm{y(\e X,\e)-\sum_{n=0}^{N-1}h_n(X)\e^n}\leq C \norm\e^{N\kappa}}
pour tout $\e\in S_2$ et $X\in V$ tels que $\e X\in V(\e)$ avec $|X|<b|\e|^{\kappa-1}$. 

\item Pour tout $n$, les polynômes $P_n$ et $Q_n$ sont de degré inférieur à $n$.\ee

Alors $y$ admet un \dac{}  pour $\e\in S_2$ et $x\in V(\e)$ ; précisément
$$y(x,\e)\sim \sum_{n=0}^{\infty}\Big(a_n(x)+g_n\big(\tfrac x\e\big)\Big)\e^n .$$}

\rqs 1.\ Comme il y a une région commune aux développements \rf{ccuni} et \rf{hhuni},
les développements doivent être compatibles, ce que montre la preuve, \cf\rf{compa}.
\med\\
2.\ On ne peut pas avoir mieux que $\norm\e^{n(1-\kappa)}$ dans le reste de \rf{ccuni}
et $\norm\e^{n\kappa}$ dans celui de \rf{hhuni} en général, car les premiers termes négligés
ont cette taille lorsque $P_N$ et $Q_N$ sont de degré $N$.
\med\\
3.\ Cet énoncé est un cas particulier d'un théorème général du livre \cit e de 
W.\ Eckhaus.
Dans le procédé classique, on établit d'abord des développements intérieurs et 
extérieurs sur des domaines qui croissent quand $\e\to0$ et qui ont une intersection
non vide. Ensuite, on construit des développement dites \og composites \fg \ dont
nos \dacs{} sont donc un exemple, \cf aussi la partie \reff{8.}.
\med\\
\pr{Preuve}Notons $c_n(x)=\sum_{m=-n}^{+\infty}c_{nm}x^m$ et $M_{nN}$ le plus grand entier $M$ tel que
$M\kappa+n\leq N(1-\kappa)$. Alors \rf{ccuni} entraîne que, pour tout 
$N\in\N$, il existe $C_2>0$ tel que 
$${\norm{y(x,\e)-\sum_{n=0}^{N-1}\sum_{m=-n}^{M_{nN}}c_{nm}x^m\e^n}\leq 
       C_2 \norm\e^{N(1-\kappa)}}$$
quand $\e\in S_2$ et $x\in V(\e)$, $a\norm\e^\kappa<\norm x<b\norm\e^\kappa$.
Pour tout entier $S$, on peut donc trouver une constante $C_3$ telle que
\eq{ccmod}{\norm{y(x,\e)-\sum_{n\geq0,m\geq-n,m\kappa+n< S}c_{nm}x^m\e^n}\leq 
       C_3 \norm\e^S}
quand $\e\in S_2$ et $x\in V(\e)$, $a\norm\e^\kappa<\norm x<b\norm\e^\kappa$.

De manière analogue, en notant $h_n(X)\sim \sum_{m=-n}^{+\infty} z_{nm}X^{-m}$,
on trouve en remplaçant $X$ par $\frac x\e$ que, pour tout entier $S$, il existe une constante $C_4$
telle que 
\eq{hhmod}{\norm{y(x,\e)-\sum_{p\geq0,q\geq-p,-q(\kappa-1)+p<S}z_{pq}x^{-q}\e^{p+q}}\leq 
       C_4 \norm\e^S}
quand $\e\in S_2$ et $x\in V(\e)$, $a\norm\e^\kappa<\norm x<b\norm\e^\kappa$.

Comme \rf{ccmod} et \rf{hhmod} déterminent de manière unique les coefficients $c_{nm}$
et $z_{pq}$, ceux-ci doivent coïncider, \ie $c_{nm}=z_{n+m,-m}$ pour tout
$n\in\N$ et $m\in\Z$, $m\geq -n$. On a donc  l'égalité formelle
\eq{compa}{\sum_{n=0}^\infty \hat h_n\big(\tfrac x\e\big)\e^n = \sum_{n=0}^\infty 
   \hat c_n(x)\e^n ,}
dans laquelle $\hat h_n$ et $\hat c_n$ désignent les séries associées à $h_n$ 
et $c_n$.

Considérons maintenant la somme $Y_N(x,\e)=\ds\sum_{n=0}^{N}\Big(a_n(x)+g_n\big(\tfrac x\e\big)\Big)\e^n$.
Quand $a\norm\e^\kappa<\norm x$, on trouve avec 
$\norm{g_n\big(\tfrac x\e\big)-\ds\sum_{q=1}^{N-n-1}z_{nq}x^{-q}\e^{q}}
       \leq C_5\norm\e^{(N-n)(1-\kappa)}$ et donc avec $z_{nq}=c_{n+q,-q}$ 
$$
\begin{array}r
\big|y(x,\e)-Y_N(x,\e)\big|\leq\norm{y(x,\e)-\ds\sum_{n=0}^{N-1}
    \gk{a_n(x)+\sum_{m=1}^{n}c_{n,-m}x^{-m}}\e^n}+\med\\
    C_6\norm\e^{N(1-\kappa)}.
\end{array}
$$
Ceci implique
\eq{ab-a}{
\begin{array}{rl}
\big|y(x,\e)-Y_N(x,\e)\big|\!\!&\leq\norm{y(x,\e)-\sum_{n=0}^{N-1}c_n(x)\e^n}+
C_6\norm\e^{N(1-\kappa)}
\med\\
&\leq C_7  \norm\e^{N(1-\kappa)} 
\end{array}   
}

quand $\e\in S_2$, $x\in V(\e)$,  $a\norm\e^\kappa<\norm x$.

En utilisant les développements des $a_n$, on trouve de manière analogue, que
$$\norm{y(x,\e)-Y_N(x,\e)}\leq C_8\norm\e^{N\kappa}$$
aussi quand  $\e\in S_2$, $x\in V(\e)$,  $\norm x < b \norm\e^\kappa$.
Ensemble avec \rf{ab-a}, ceci démontre que, pour tout $N$, il existe
une constante $C_9$ telle que pour tout $\e\in S_2$ et $x\in V(\e)$ on a
$\norm{y(x,\e)-Y_N(x,\e)}\leq C_9\norm\e^{N \lambda}$ avec $\lambda=\min(\kappa,1-\kappa)$.

L'énoncé à démontrer correspond à $\norm\e^{N}$ au lieu de $\norm\e^{N \lambda}$ dans cette dernière inégalité. On l'obtient en deux temps. 
D'une part cette dernière assertion peut aussi s'écrire : 
il existe $C_{10}$ avec
$\norm{y(x,\e)-Y_S(x,\e)}\leq C_{10}\norm\e^{N}$, si $S\lambda>N$. D'autre part le fait que toutes les fonctions $a_n(x)$ et $g_n\big(\tfrac x\e\big)$ 
soient bornées sur l'ensemble des $x,\e$ en question
entraîne qu'il existe une constante $C_{11}$ telle que $\norm{Y_S(x,\e)-Y_N(x,\e)}\leq C_{11}\norm\e^{N}$.
\ep
\sub{prolodac}{Prolongements de développements combinés}
En relation avec les développements intérieurs de la méthode de recollement, nous avons aussi des résultats de
prolongement de \dac, qui nous seront bien utiles pour les solutions d'équations différentielles. 
Le premier résultat dit
essentiellement qu'une fonction ayant un \dac{} pour $x$ dans un quasi-secteur, 
dont le développement intérieur
existe sur un quasi-secteur plus grand, admet le \dac{} aussi sur le grand quasi-secteur. 
Le résultat précis est le suivant.
\propo{p3.12}{
Soit $y$ une fonction définie pour $\e\in S_2=S(\a_2,\b_2,\e_0)$ et 
$x\in V_1(\e)=V(\a_1,\b_1,r_0,\mu\norm\e)$
et ayant un \dac{} $\ts\sum_{n\geq0}\Big(a_n(x)+g_n\big(\ts\frac x\e\big)\Big)\e^n$, 
quand $S_2\ni\e\to0$ et $x\in V_1(\e)$, avec  
 $a_n\in\H(r_0)$ et $g_n\in\G(V)$,  où $V=V(\a,\b,\infty,\mu)$, 
$\a=\a_1-\b_2$ et $\b=\b_1-\a_2$. 
Soit $\nu>\mu$. Dans le cas où $\nu>|\mu|$, on pose $\W=D(0,\nu)$, sinon on pose $\W=V(\a,\b,-\mu+\g,\nu)$ avec $\g>0$ arbitrairement petit.

On suppose que la fonction $Y:(X,\e)\mapsto y(\e X,\e)$ peut être 
prolongée analytiquement sur $\W\times S_2$  et qu'elle admet un
développement asymptotique  $Y(X,\e)\sim\sum_{n=0}^\infty h_n(X)\e^n$ quand $\e$ 
tend vers $0$, uniformément sur $\W$. 

Alors $y$ peut être prolongée analytiquement sur l'ensemble des $(x,\e)$ avec 
$\e\in S_2$ et avec $x\in V(\a_1,\b_1,r_0,\nu\norm\e)$ et y admet un \dac{} quand $\e\to0$.
}
\figu{f3.2}{
\vspace{-7mm}
\epsfxsize10cm\epsfbox{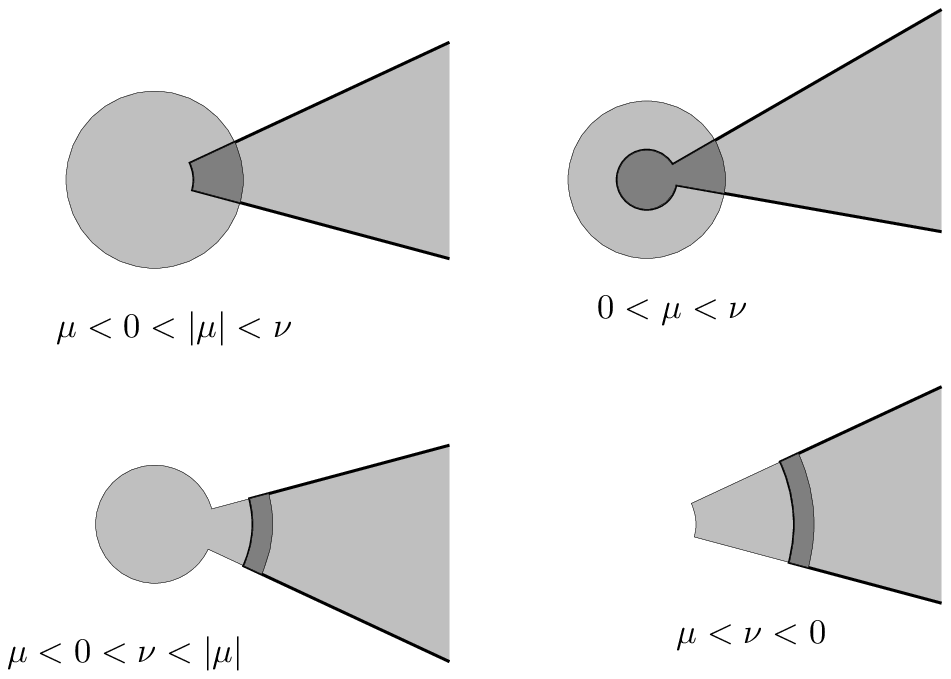}\vspace{-3mm}
}{Quelques domaines $V$ et $\W$ suivant les signes de $\mu$, $|\mu|-\nu$ et $\nu$. En trait gras le bord de $V$, en trait fin celui de $\W$, en gris foncé leur intersection}

\rqs
1.\ 
Le domaine $\W$ a été choisi de façon à avoir $\W$ borné, $\W\cap V\neq\emptyset$ et 
$V\cup\W = V(\a,\b,\infty,\nu)$.
Les signes de $\mu$ et $\nu$ sont arbitraires. Nous utiliserons ce résultat en particulier dans le cas $\mu<0<\nu$.
\med\\
2.\ L'hypothèse sur le domaine et le développement asymptotique de $Y$ pourrait être
un peu affaiblie (le domaine par rapport à  $X$ pourrait dépendre de l'argument de $\e$), 
mais la version présentée suffit pour nos applications aux équations différentielles.
\med\\
3.\ Il est possible de montrer ce résultat en utilisant les propositions \reff{combextint} et \reff{p3.10bis}, mais nous préférons présenter une preuve indépendante. Une raison pour ce choix est que cette preuve servira pour l'analogue Gevrey, la proposition \reff{matching-gevrey}. Par contre nous n'avons pas d'analogue Gevrey de la proposition \reff{p3.10bis}.
\med

\pr{Preuve}
Le développement de $Y$ de l'hypothèse et le développement intérieur 
correpondant au \dac{} de $y$ donné par la proposition \reff{combextint} coexistent sur
un certain ouvert, donc coïncident par l'unicité d'un développement asymptotique. Ainsi, 
les $h_n(X)$ de l'énoncé sont nécessairement les prolongements analytiques
des coefficients de ce développement intérieur.

L'hypothèse entraîne que $y(x,\e)$ peut être prolongée analytiquement sur l'ensemble des 
$(x,\e)$ tels que $\e\in S_2$ et $x\in \~V_1(\e)=V(\a_1,\b_1,r_0,\nu\norm\e)$, précisément en posant 
$y(x,\e)=Y\big(\tfrac x\e,\e\big)$.
Nous reprenons maintenant les notations de la preuve précédente. 
Il s'agit de montrer que le reste  $R_N(x,\e)$ donné par \rf{rest} est majoré par une 
constante fois $|\e|^N$, aussi pour $x\in \~V_1(\e)\setminus V_1(\e)$.
L'hypothèse sur $\W$ assure que pour tout $\e\in S_2$ et tout $x\in \~V_1(\e)\setminus V_1(\e)$, on a $x/\e\in\W$.
 Par hypothèse, il existe $D_N$ tel que 
$$
\Big|y(x,\e)-\sum_{n<N}h_n\big(\tfrac x\e\big)\e^n\Big|\leq D_N|\e|^N.
$$
Par ailleurs, l'égalité au-dessus de \rf{cin} peut s'écrire
$$
R_N(x,\e)=y(x,\e)-\sum_{n<N}h_n\big(\tfrac x\e\big)\e^n-\sum_{k<N}\bigg(a_k(x)-
   \sum_{l< N-k}a_{kl}\,x^l\bigg)\e^k.
$$
De plus, quitte à  modifier les constantes $A_{kn}$, l'inégalité \rf{AA} est valable pour tout $x\in D(0,r_0)$, donc en particulier pour $x\in \~V_1(\e)\setminus V_1(\e)$. 
Ceci montre que pour tout $\e\in S_2$ et pour tout 
$x\in \~V_1(\e)\setminus V_1(\e)$
\eq{rn}{
|R_N(x,\e)|\leq\bigg(D_N+\sum_{k<N}A_{k\,N\!-\!k}M^{N-k}\bigg)|\e|^N
}
avec $M=\sup_{X\in\W}\norm X=|\mu|$ ou $\nu$.
\ep
\\
Le deuxième résultat de prolongement concerne le prolongement vers l'extérieur.
\propo{p3.11bis}{
Soit  $0<r_0<\~r_0$  et soit $y$ une fonction définie pour $\e\in S_2=S(\a_2,\b_2,\e_0)$ et
$x\in\~V_1(\e)=V(\a_1,\b_1,\~r_0,\mu\norm\e)$. On suppose que $y$ a un \dac{}
$\ts\sum_{n\geq0}\gk{a_n(x)+g_n\big(\ts\frac x\e\big)}\e^n$,
quand $S_2\ni\e\to0$ et $x\in V_1(\e)=V(\a_1,\b_1,r_0,\mu\norm\e)$, avec  $a_n\in\H(r_0)$ et $g_n\in\G(V)$,
$V=V(\a,\b,\infty,\mu)$ tels que
$\a\leq\a_1-\b_2<\b_1-\a_2\leq\b$.

On suppose de plus que $y$ admet un
développement asymptotique  $y(x,\e)\sim\sum_{n=0}^\infty c_n(x)\e^n$ quand $\e$
tend vers $0$, uniformément pour $x\in V(\a_1,\b_1,r_0-\g,\~r_0)$
avec $\g>0$ arbitrairement petit.

Alors \rf{defcomb} est satisfait pour tout $\e\in S_2$ et tout $x\in\~V_1(\e)$.
}
\rq Par abus de langage, nous dirons que $y$ a un \dac{} pour $\e\in S_2$ et $x\in\~V_1(\e)$, bien que les  fonctions
$a_n$ ne soient \apriori{} pas définies sur tout le disque $D(0,\~r_0)$.
\med

\pr{Preuve}
D'abord, on peut utiliser la proposition \reff{combextint} sur le quasi-secteur
$V(\a_1,\b_1,$ $r_0-\g,r_0)$ et, en comparant \rf{cc} avec la deuxième hypothèse, on 
obtient que les fonctions $c_n$ de l'hypothèse sont des prolongements analytiques de 
celles
de la proposition. On peut donc aussi prolonger les fonctions $a_n$ analytiquement sur 
$D(0,r_0)\cup V(\a_1,\b_1,r_0-\g,\~r_0)$.
 
Il faut encore majorer $R_N$ pour $x\in\~V_1(\e)\setminus V_1(\e)$. Par \rf{34c}, on a
\eq{322}{
R_N(x,\e)=y(x,\e)-\sum_{n<N}c_n(x)\e^n
-\sum_{l<N}r_{l\,N\!-\!l}\big(\ts\frac x\e\big)\e^l.
}
et par \rf{34b} $|r_{l\,N\!-\!l}(X)|\leq C_{l\,N\!-\!l}|X|^{l-N}$. Par hypothèse, il existe $A_N>0$ tel que
$$
|y(x,\e)-\sum_{n<N}c_n(x)\e^n|\leq A_N|\e|^N
$$ 
pour tout $\e\in S_2$
 et tout $x\in V_1(\e)\setminus\~V_1(\e)$.
On obtient ainsi $|R_N(x,\e)|\leq C|\e|^N$ avec $C=A_N+
  \sum_{l<N}C_{l\,N\!-\!l}r_0^{l-N}$.
\ep
\sub{divers}{Quotients de \dacs{} et une extension}
{\noi\bf Quotients de \dac} \sep
\med\\
Nous étudions ici sous quelle condition l'inverse d'une fonction ayant un \dac{} admet un \dac{}.
Soit $y=y(x,\e)$ une fonction définie et analytique pour $\e\in S=S(-\d,\d,\e_0)$ et $x\in V(\a,\b,r_0,\mu|\e|)$ et ayant un \dac{} 
$y(x,\e)\sim\sum_{n\geq0}\gk{a_n(x)+g_n\big(\ts\frac x\e\big)}\e^n$
quand $\e\to0$. Nous proposons un énoncé un peu plus général, qui s'avère plus utile en pratique~: 
 sous quelle condition existe-t-il $k\in\N$ tel que la fonction $(x,\e)\mapsto \e^k/y(x,\e)$ admette un \dac.

D'après la proposition \reff{combextint}, $y$ admet un développement intérieur
$${
y(\e X,\e)\sim\sum_{n\geq0}h_n(X)\e^n \mbox{ \ quand }S_2\ni\e\to0,
}$$
uniformement  par rapport à  $X$ sur des compacts de $S(\a_1,\b_1,\infty)$ 
où $\a_1=\a-\d,$ $\b_1=\b+\d$ et
\eq{hhhh}{h_n(X)=g_n(X)+\ds\sum_{0\leq l\leq n} a_{l,n-l}X^{n-l}.} 
Pour tout $n$, notons $\ds\sum_{m=-n}^{+\infty}h_{nm}X^{-m}$ le développement asymptotique à l'infini de $h_n$ et $\vu_n=\val(h_n)$ le minimum des $m\geq-n$ tels que $h_{nm}\neq0$.
Si $h_n$ est plate, on note $\vu_n=\val(h_n)=+\infty$.
Nous disons que $y$ est {\em dégénérée}, si elle est plate ou s'il existe $N\in\N$
tel que $h_0=...=h_{N-1}=0$ et $h_N\neq0$ est plate.
Si $y$ n'est pas dégénérée, notons $C(y)$ le couple $(N,M)$ de 
$N\in\N$ tel que $h_0=...=h_{N-1}=0$, $h_N\neq0$ et $M=\val(h_N)\geq-N$. 
\coro{quot}{
Avec les notations précédentes, les trois conditions suivantes sont équivalentes.
\be[\rm(a)]
\item 
Il existe $k\in\N$, $\~\e_0<\e_0$, $\~r_0<r_0$ et $\~\mu\leq\mu$ tels que la fonction $(x,\e)\mapsto \e^k/y(x,\e)$ admette un \dac{} quand $\e\to0$ dans $\~S=S(-\d,\d,\~\e_0)$ et $x\in V(\a,\b,\~r_0,\~\mu)$.
\item
$y$ est non dégénérée et, si $C(y)=(N,M)$, on a $\val(h_n)\geq M - n + N$
pour tout $n\geq N$.
\item
Il existe $k\in\N$, $\ell\in\Z$ tels que la fonction 
$(x,\e)\mapsto \e^{-k}x^{\ell}y(x,\e)$ admet un \dac{} dont le premier coefficient lent
$\~a_0$ satisfait $\~a_0(0)\neq 0$. 
\ee
}
\rqs 
1.\ 
La deuxième condition du (b) peut aussi
être écrite en utilisant le développement extérieur et ceux de ses coefficients. On a la relation 
\rf{compa}.
\med\\
2.\ Graphiquement, la deuxième condition dans (b) signifie que les points du plan 
 de coordonnées $(n,m)$ tels que $h_{nm}\neq0$ (le \og support \fg{} du développement intérieur) sont tous dans la partie du plan à droite et au-dessus des droites verticale et de pente $-1$ passant pas $C(y)$. Puisque $h_n$ a une partie polynomiale de degré au plus $n$, on sait déjà que ce support est dans le quart de plan à droite et au-dessus de l'axe des ordonnées et de la deuxième bissectrice. Le changement de variable $y\to z:(x,\e)\mapsto \e^{-k}x^{\ell}y(x,\e)$ induit justement une translation de $-C(y)=(-N,-M)$ sur les supports, avec $N=k-\ell$ et $M=\ell$.
 
\figu{f3.3}{
\vspace{-3mm}
\epsfxsize10cm\epsfbox{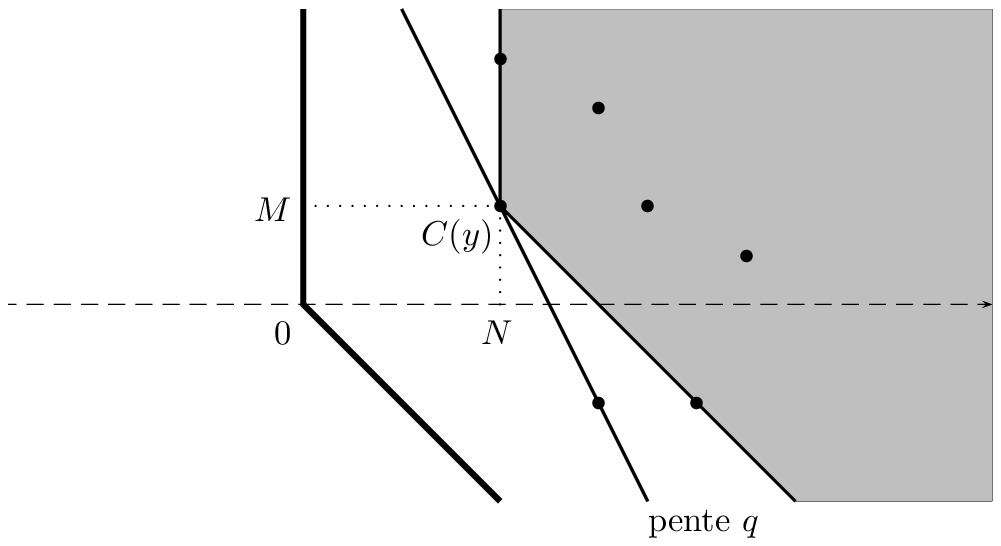}\vspace{-5mm}
}{En grisé, la partie du plan contenant le support du développement intérieur de $y$. En gras, le bord de la partie analogue pour $z$}
\noi3.\ 
La démonstration donne aussi une procédure pour obtenir  ce \dac~: par la translation précédente, on se ramène au cas où $y$ a un premier terme $a_0$ lent non nul en $x=0$ et on compose le \dac{} par l'inversion $u\mapsto1/u$.
\med\\
\pr{Preuve}
Nous montrons les implications (b)$\Rightarrow$(c)$\Rightarrow$(a)$\Rightarrow$(b).
Supposons que la  condition (b) est satisfaite. 

Si $M<0$,
considérons d'abord $z(x,\e)=\gk{\frac\e x}^{-M}y(x,\e)$. Comme produit de 
deux fonctions admettant des \dacs{}, $z$ a un \dac{} sur le même domaine que $y$.
Le développement intérieur correspondant est celui de $X^{M}y(\e X,\e)$ et satisfait 
donc une condition analogue à (b) avec $(N,0)$ à la place de $(N,M)$.

Si $M>0$, considérons $z(x,\e)=x^{M}y(x,\e)$. Comme avant, $z$ admet un \dac{} et 
le développement intérieur  correspondant est celui de $\e^{M}X^{M}y(\e X,\e)$ 
et satisfait 
donc une condition analogue à (b) avec $(N+M,0)$ à la place de $(N,M)$.
Les deux cas $M>0$ et $M<0$ peuvent donc être réduits au cas $M=0$.

Si $M=0$, alors la formule \rf{hhhh} et la condition sur les $h_n$
montrent que $a_{s m}=0$ pour $0\leq s<N$ et $m\geq0$, et que $a_{N0}\neq0$. 
Comme les fonctions $a_s$ sont analytiques,
ceci implique  $a_s=0$ pour $s=0,...,N-1$. 
Alors la fonction $\~y:(x,\e)\mapsto\e^{-N}y(x,\e)$ a aussi un \dac{}
sur le même domaine que $y$ et satisfait la condition (c).
En résumé, dans les trois cas $y$ satisfait la  condition (c) avec $k=N+M$ et $\ell=M$.
\med

Supposons maintenant la  condition (c) satisfaite et posons $z(x,\e)=\e^{-k}x^\ell y(x,\e)$. Pour $\~r_0$ et $\~\e_0$ assez petits et $\~\mu\leq\mu$ convenable, la fonction $z$ ne s'annule pas si $\e\in\~S=S(-\d,\d,\~\e_0)$ et $x\in V(\a,\b,\~r_0,\~\mu)$. La proposition \reff{l2.4} (b) s'applique avec la fonction $f:u\mapsto1/u$ et on en déduit que $1/z$ a un \dac.
Dans le cas où $\ell\geq0$, la fonction $(x,\e)\mapsto x^\ell/z(x,\e)$ a donc, elle aussi, un \dac, ce qui prouve (a). Dans le cas  $\ell<0$, c'est la fonction $(x,\e)\mapsto\big(\frac x\e\big)^\ell/z(x,\e)$ qui a un \dac, ce qui donne (a) avec $k-\ell$ au lieu de $k$.
\med

Supposons enfin que la  condition (a) est satisfaite, posons $\~y(x,\e)=\e^k/y(x,\e)$ et notons  $\~h_n$ les coefficients du développement intérieur de $\~y$. Puisque $(y\~y)(x,\e)=\e^k$, c'est qu'il existe un premier terme non identiquement nul $h_r$ dans le développement intérieur de $y$ et un premier terme $\~h_s$ pour $\~y$, avec $r+s=k$ et $h_r\~h_s=1$. Chacune de ces fonctions est à croissance polynomiale lorsque $X\to\infty$ donc aucune ne peut être plate. Ainsi les deux fonctions $y$ et $\~y$ sont non dégénérées. On note $(N,M)=C(y)$ et $(\~N,\~M)=C(\~y)$.

Par l'absurde, supposons qu'il existe $n>N$
tel que $\val(h_n)<M-n+N$. 
Notons $q=\min\big\{\frac{\val(h_s)-M}{s-N}\tq s>N\big\}<-1$ et
${\cal M}=\{s\geq N\tq \val(h_s)-s q= M-Nq\}$~; c'est un ensemble de cardinal au moins 2 (contenant au moins $N$ et un $s$ réalisant le minimum $q$) et fini (puisque $\val(h_s)\geq-s$).

Si $\~y$ satisfait la  condition (b), on pose  $\~q=-1$, sinon $\~q$ est l'analogue de $q$ pour $y$. Quitte à échanger les rôles de $\~y$ et $y$, on peut supposer sans perte que $q\leq\~q$.
Notons alors $K=\min\{\val(\~h_s)-s q \tq s>\~N \}$ et $\cal N$ l'ensemble fini non vide
des $s\in\N$ tels que $\val(\~h_s)-s q=K$. Remarquons que ${\cal N}=\{\~N\}$ si $\~q>q$ et
que le cardinal de $\cal N$ est au moins 2 si $\~q=q$.

Rappelons que le minimum de  $\cal M$ est $n_1=N$, et notons $n_2=\max\cal M$~; notons $\~n_1$ et $\~n_2$ le minimum et le maximum de $\cal N$. 
Considérons le développement intérieur
de $p=y\~y$~: ${p(\e X,\e)\sim\sum_{n\geq0}p_n(X)\e^n}$.
On a donc $p_n=\sum_{r+s=n}h_r\~h_s$ pour tout $n\geq0$. 
Si $n=n_1+\~n_1$ alors $h_{n_1}\~h_{\~n_1}$ a la valuation 
$M+\~n_1+K= M-Nq+n q+K$~;  si
$r+s=n=n_1+\~n_1$ avec $r\neq n_1$, alors la valuation de $h_r\~h_s$ est strictement supérieure
à ce nombre à cause du choix de $n_1$ et $\~n_1$. On obtient donc $p_{n_1+\~n_1}\neq 0$.
De manière analogue, on obtient aussi $p_{n_2+\~n_2}\neq 0$.
Puisque $n_1<n_2$ et $\~n_1\leq \~n_2$, ceci contredit l'hypothèse que le produit
$y\~y$ est réduit au monôme $\e^k$.
\ep
\med

{\noi\bf Une extension} \sep
\med\\
Appelons provisoirement {\sl point tournant} d'une fonction $y$ un point au voisinage duquel $y$ admet un \dac.
Pour terminer cette partie, nous voudrions discuter le cas d'un domaine
 ayant {\sl deux} points tournants sur sa frontière. 
L'énoncé qui suit montre qu'il est inutile de généraliser la notion de \dacs{}  pour des développements uniformes
sur des domaines ayant plusieurs points tournants sur leur frontière,
 car on peut se ramener au cas d'un seul point 
tournant. Pour simplifier, nous nous sommes placés dans le cas d'un intervalle réel.
\propo{deuxpt}{
Soit $a<b<c<d$ quatre nombres réels et soit 
$y:\,]a,d[\times]0,\e_1]\to\R$ une fonction admettant un \dac{} 
$$
y(x,\e)\sim\sum_{n=0}^\infty\Big(a_n(x)+g_n\big(\tfrac{x-a}\e\big)\Big)\e^n
$$
quand $\e\to0$, uniformément sur  $]a,c[$, avec $a_n$ holomorphe dans
un voisinage de $[a,c]$ et $g_n\in\G(S)$, $S=S(-\d,\d,\infty)$
avec un certain $\d>0$.

On suppose que $y$ admet aussi un \dac{}
$$
y(x,\e)\sim\sum_{n=0}^\infty\Big(b_n(x)+h_n\big(\tfrac{d-x}\e\big)\Big)\e^n
$$
quand $\e\to0$, uniformément sur $]b,d[$, avec $b_n$ holomorphe dans un
voisinage de $[b,d]$ et $h_n\in\G(S)$.

Alors $y$ admet un développement asymptotique
\eq{2pt}{y(x,\e)\sim \sum_{n=0}^\infty (c_n(x)+g_n(\tfrac{x-a}\e)
   +h_n(\tfrac {d-x}\e))\e^n}
quand $\e\to0$, uniforme sur $]a,d[$ avec les $g_n,h_n$ des formules précédentes
et avec des fonctions $c_n$ holomorphes sur un voisinage de $[a,d]$.

Précisément, si $g_n(X)\sim\sum_{m=1}^\infty g_{nm}X^{-m}$ et $h_n(X)\sim\sum_{m=1}^\infty h_{nm}X^{-m}$
quand $X\to\infty$,
alors $c_n(x)=a_n(x)-\sum_{\ell=0}^{n-1} h_{\ell\,n-\ell} (d-x)^{\ell-n}$ quand $x$ est dans un voisinage de
$[a,c]$ et $c_n(x)=b_n(x)-\sum_{\ell=0}^{n-1} g_{\ell\,n-\ell} (x-a)^{\ell-n}$
quand $x$ est dans un voisinage de $[b,d]$.
}
\rq 
Les fonctions $c_n$ sont donc les parties non polaires des fonctions
$b_n$ au point  $x=a$ et les parties non polaires des fonctions $a_n$ au point $x=d$.
\med
\\
\pr{Preuve succincte}
Par le théorème de Borel-Ritt classique, on construit deux fonctions $g,h$ avec
$g(X,\e)\sim\sum_{n\geq0}g_n(X)\e^n$ et $h(X,\e)\sim\sum_{n\geq0}h_n(X)\e^n$ uniformément
sur $S$. 
On considère alors la différence
$z(x,\e)=y(x,\e)-g\big(\tfrac{x-a}\e,\e\big)-h\big(\tfrac{d-x}\e,\e\big)$.
La proposition \reff{combextint}, appliquée à $h$ respectivement $g$, 
montre qu'elle admet deux développements 
uniformes lents sur $[a,c]$ et sur
$[b,d]$. Par unicité des développements asymptotiques, ils doivent coïncider sur $[b,c]$. Ceci implique que leurs coefficients doivent être des prolongements les uns des autres et on obtient l'énoncé.
\ep
%
%
%
\sec{2bis} {Développements asymptotiques combinés :\\ étude Gevrey}
Dans la suite de l'article (sections \reff{5.} et \reff{6.}) nous appliquerons les \dacs{} à  des problèmes d'équations différentielles singulièrement perturbées. Nous verrons à  cette occasion que les notions de séries formelles combinées Gevrey et de \dacs{} de type Gevrey jouent un rôle clé. Cette notion Gevrey a déjà  joué un rôle essentiel dans la théorie classique des séries formelles et des développements asymptotiques dans les applications à  la perturbation singulière \cit{crss,bfsw}.
\`A partir de maintenant, 
 $p$ désigne un entier positif. En général (sauf dans la partie \reff{5.1}) cet entier sera supérieur ou égal à  $2$.
\sub{2bis.1} {Séries formelles combinées Gevrey}
\df{d2.3.1}{
Soit $r_0>0$, $\mu\in\R$ , $V=V(\a,\b,\infty,\mu)$ donné par \rf{V1} ou \rf{V2} et soit
$$
\yc(x,\e)=\sum_{n\geq0}\gk{a_n(x)+g_n\big(\ts\frac x\e\big)}\e^n\in\cch(r_0,V) ,
$$
où $g_n(X)\sim\ds\sum_{m>0}g_{nm}X^{-m}$ quand $V\ni X\to\infty$.
Nous disons que $\yc$ est {\sl Gevrey d'ordre $\usp$ et de type $(L_1,L_2)$,}
s'il existe une constante $C$ telle que pour tout $n\in\N$ on a
$\ds\sup_{\norm x< r_0}\norm{a_n(x)}\leq C L_1^n \Gamma
\big(\tfrac np+1\big)$ et pour tout $n,M\in\N$ et tout $X\in V$ on a
\eq{gnm}{  
\norm X^M\bigg|g_n(X)-\sum_{m=1}^{M-1}g_{nm}X^{-m}\bigg|
\leq C L_1^n L_2^M\Gamma\big(\ts\frac {M+n}p+1\big).  }
Nous disons que des fonctions $g_n(X),$ $n=0,1,...$, holomorphes sur $V$
{\sl admettent des développements asymptotiques Gevrey
d'ordre $\usp$ compatibles (de type $(L_1,L_2)$)}, s'il existe une constante $C$ telle
que (\reff{gnm}) est satisfaite pour tout $n,M\in\N$ et $X\in V$.
}
Nous avons inclus la notion de type dans la définition, mais nous utiliserons peu cette notion dans le présent mémoire.

Il serait possible de formuler (\reff{gnm}) avec des produits 
$\Gamma\big(\frac np+1\big)\Gamma\big(\frac Mp+1\big)$ en vertu des inégalités
\eq{gamma}{\Gamma(a+1)\Gamma(b+1)\leq\Gamma(a+b+1)\leq2^{a+b}\Gamma(a+1)\Gamma(b+1),
\ \mbox{quand}\ a,b>0}
mais nous préférons l'écriture (\reff{gnm}), qui simplifie certaines preuves.
On déduit facilement  de (\reff{gnm}) (avec successivement $M=m$ et $M=m+1$) que 
\eq{gnmc}{\norm{g_{nm}}\leq C L_1^nL_2^m\Gamma\big(\ts\frac {m+n}p+1\big)}
pour tout $n,m\in\N$. 
Pour $M\geq1$, les inégalités (\reff{gnm}) sont équivalentes à 
$\norm X\norm{(\T^{M-1}g_n)(X)}\leq C L_1^n L_2^M\Gamma\big(\ts\frac {M+n}p+1\big) $
avec l'opérateur $\T$ de (\reff{4}). 

Il nous sera utile de ne pas avoir
à  montrer les inégalités (\reff{gnm}) pour $X$ proche de 0, d'où la remarque qui suit.
\begin{rem}\lb{remgnm}\sep
{\sl 
Les fonctions $g_n(X),$ $n=0,1,...$,  admettent des développements asymptotiques Gevrey
d'ordre $\usp$ compatibles de type $(L_1,L_2)$ si et seulement s'il 
existe deux constantes $C,T>0$ telles que 
\eq{gnmrem1}{
\sup_{X\in V}\norm{g_n(X)}\leq C L_1^n \Gamma
\big(\tfrac np+1\big)
}
pour tout $n\in\N$ et telles que
\eq{gnmrem2}{
\bigg|g_n(X)-\sum_{m=1}^{M-1}g_{nm}X^{-m}\bigg|
\leq C L_1^n L_2^M\Gamma\big(\ts\frac {M+n}p+1\big)\norm X^{-M}
}
pour tout $n,M\in\N$ et tout $X\in V$ avec $\norm X > T$.
}
\end{rem}
\pr{Preuve}
L'implication directe est évidente (pour \rf{gnmrem1}, choisir $M=0$). 
Pour la réciproque, il s'agit de montrer, à  partir de (\reff{gnmrem1}) et (\reff{gnmrem2}),
 que (\reff{gnm}) est aussi vrai pour les $X\in V$  avec $\norm X\leq T$. 
Tout d'abord (\reff{gnmrem2}) implique
(\reff{gnmc}). En utilisant (\reff{gamma})
on obtient
$$\begin{array}{rcl}
\norm X^{M}\bigg|g_n(X)-\ds\sum_{m=1}^{M-1}g_{nm}X^{-m}\bigg|&\leq&
T^M  C L_1^n \Gamma\big(\tfrac np+1\big)+\\
&&\ds\sum_{m=1}^{M-1}C L_1^nL_2^m
\Gamma\big(\ts\frac {m+n}p+1\big)T^{M-m}\med\\
&\leq&
\~ C L_1^n L_2^M\Gamma\big(\ts\frac {M+n}p+1\big)
\end{array} $$
avec $\~ C=C\sum_{k\geq1}\gk{\frac T{L_2}}^k\frac1\Gamma\big(\frac kp+1\big)$.
\ep

Les inégalités \rf{gnm} sont également 
indispensables pour la compatibilité de la nouvelle notion 
avec les opérations élémentaires. Nous laissons les cas de l'addition et de la
dérivation comme exercices~; l'énoncé (et sa preuve) 
pour la dérivation est analogue au lemme \reff{derivasympt}.
Nous démontrons ci-dessous la compatibilité avec la multiplication.

\fp
Considérons d'abord le cas classique d'un produit 
$\widehat c(x,\e)=\widehat a_1(x,\e)\widehat a_2(x,\e)$ de 
deux séries formelles $\widehat a_j(x,\e)=\sum_{n\geq0} a_{jn}(x)\e^n$ Gevrey d'ordre $\usp$.
Si $\ds\sup_{|x|<r_0}\norm{a_{jn}(x)}$ $\leq C_j L_j^n \Gamma\big(\tfrac np+1\big)$ pour $j=1,2$ et tout $n\in\N$,
on obtient pour les coefficients de $\widehat c$
$$
\ds\sup_{|x|<r_0}\norm{c_n(x)}\leq C_1C_2\max(L_1,L_2)^n\sum_{l=0}^n
\Gamma\big({\tfrac lp+1\big)}\Gamma\big({\tfrac {n-l}p+1}\big)
$$
et on obtient le caractère Gevrey de $\widehat c$ 
de l'inégalité 
\eq{prodgamma}{\sum_{\nu=0}^n\Gamma\big(\tfrac{n-\nu}p+1\big)\Gamma\big(\tfrac{\nu}p+1\big)\leq K_p
\Gamma\big(\tfrac{n}p+1\big)}
avec la constante $K_p=2\Big(1+\Gamma\big(\frac{1}p+1\big)+....+\Gamma\big(\frac{p-1}p+1\big)\Big)+2p$.

Considérons maintenant un produit $\widehat h(x,\e)=\widehat g_1(x,\e)g_2(x,\e)$ 
de deux séries formelles
$\widehat g_j(x,\e)=\sum g_{jn}(\tfrac x\e)\e^n$ Gevrey d'ordre $\usp$ et supposons 
d'après la remarque \reff{remgnm} que 
$\sup\norm{g_{jn}(X)}\leq C_j L_j^n \Gamma\big(\tfrac np+1\big)$ pour tout $n$.
On montre comme avant que les coefficients $h_n$ de $\widehat h$ satisfont
$$
\sup\norm{h_n(X)}\leq K_p C_1C_2 \max(L_1,L_2)^n\Gamma\big(\tfrac np+1\big).
$$
Il faut encore montrer que les $h_n(X)$ admettent des développements asymptotiques
Gevrey compatibles, \ie (\reff{gnmrem2}).

On suppose donc (sans perte avec les mêmes constantes)
$$
\norm X \norm{\T^{M-1}g_{jn}(X)}\leq C L_1^n L_2^M\Gamma\big(\ts\frac {M+n}p+1\big)
$$
pour  tout $n,M$ avec $M\geq1$
ainsi que la même majoration pour les $g_{jnm}$. On a à  estimer
$\norm X \norm{\T^{M-1}h_{n}(X)}$ avec $h_n=\sum_{\nu=0}^ng_{1\nu}g_{2,n-\nu}$.
En utilisant plusieurs fois la formule $\T(f\~ f)=\T(f)\~ f + f_1\~ f$
pour chaque terme de cette somme, on obtient
$$
\norm X \norm{\T^{M-1}h_{n}(X)}\leq
C^2 L_1^n L_2^M \sum_{\nu=0}^n\sum_{m=1}^{M} 
\Gamma\big(\tfrac {\nu+m}p+1\big)\Gamma\big(\tfrac {n-\nu+M-m}p+1\big).
$$
De manière analogue à  (\reff{gamma}), on montre l'existence d'une constante $\~ K$
telle qu'on a pour tout $n,M$
$$
\sum_{\nu=0}^n\sum_{m=0}^{M} 
\Gamma\big(\tfrac {\nu+m}p+1\big)\Gamma\big(\tfrac {n-\nu+M-m}p+1\big)\leq \~ K
\Gamma\big(\tfrac {M+n}p+1\big).
$$
Ceci complète la majoration des termes $\norm X \norm{\T^{M-1}h_{n}(X)}$.

Il reste à  traiter le cas mixte.
Soit $\widehat y(x,\e)=\sum_{n\geq0} a_n(x)\e^n$ et $\widehat z(x,\e)=$ $ \sum_{n\geq0}$ $
g_n\big(\frac x\e\big)\e^n$ deux séries formelles combinées Gevrey d'ordre $\usp$ de type 
$(L_1,L_2)$ selon la définition \reff{d2.3.1}.
Notons $a_n(x)=\sum_{m\geq0}a_{nm}x^m$ quand $\norm x<r_0$ et $g_n(X)\sim\sum_{m\geq1}
g_{nm}X^{-m}$ quand $X\rightarrow\infty$.
D'après la définition du produit de séries formelles combinées,
on a
\eq{prodform}{
(\widehat y\cdot\widehat z)(x,\e)=\sum_{n\geq0}\e^n\sum_{\nu=0}^n 
   \big(I(a_\nu)I(g_{n-\nu})\big)(x,\e)=\sum_{n\geq0}\big(b_n(x)+h_n(\tfrac x\e)\big)\e^n,
   } 
où 
$$
b_n(x)=\ds\sum_{\nu+l+m=n}g_{l,m}(\S^ma_\nu)(x)~\mbox{ et }~h_n(X)=\ds\sum_{\nu+l+m=n}
a_{\nu m}(\T^mg_l)(X).
$$
D'après l'hypothèse, on a
$$
\norm{a_{\nu m}}\leq C L_1^\nu \gk{\frac1{r_0}}^m\Gamma\big(\frac\nu p+1\big),
\quad
\sup \norm{\S^ma_\nu}\leq C L_1^\nu \big({\frac2{r_0}}\big)^m\Gamma\big(\frac\nu p+1\big),
$$
 ainsi que 
$$
\norm{g_{lm}}\leq C L_1^l L_2^m\Gamma\big(\frac{m+l}p+1\big)~\mbox{ et  }~\sup\norm{\T^mg_l(X)}\leq C L_1^l L_2^m\Gamma\big(\frac{m+l}p+1\big).
$$
Ceci implique que 
$$
\sup\norm{b_n(x)}\leq C^2\ds\sum_{\nu+l+m=n}L_1^lL_2^m
\Gamma\big(\tfrac{m+l}p+1\big)L_1^\nu \big({\tfrac2{r_0}}\big)^m\Gamma\big(\tfrac\nu p+1\big).
$$ 
Quitte à  agrandir
$L_1$, on peut supposer que $q:=\frac{2L_2}{r_0L_1}<1$ et on obtient
\begin{eqnarray*}
\sup\norm{b_n(x)}&\leq& C^2 L_1^n\sum_{\nu+l+m=n}q^m \Gamma\big(\tfrac{m+l}p+1\big)
\Gamma\big(\tfrac\nu p+1\big)\\
&\leq&\frac{C^2 L_1^n}{1-q}\sum_{\nu=0}^n
\Gamma\big(\tfrac{n-\nu}p+1\big)\Gamma\big(\tfrac{\nu}p+1\big).
\end{eqnarray*}
La majoration désirée de $b_n$ est donc de nouveau une conséquence de l'inégalité
(\reff{prodgamma}).
La majoration de $h_n(X)$ étant la même, il reste à  vérifier que les $h_n(X)$
admettent des développement Gevrey d'ordres $\usp$ compatibles. 
Comme ceci équivaut à  majorer 
$$
X\T^{M-1}h_n(X)=
\ds\sum_{\nu+l+m=n} a_{\nu m}X\T^{m+M-1}g_l(X),
$$
 ceci est de nouveau analogue à  la
démarche précédente.

\ft\med
Contrairement au cas de l'addition, de la multiplication et de la dérivation, la compatibilité de la notion de série formelle combinée Gevrey avec la composition à  gauche ou à  droite par une fonction analytique serait très fastidieuse et délicate à  montrer en n'utilisant que la définition. 
Elle sera démontrée dans la partie \reff{4.5} en utilisant la proposition \reff{p3.11} et le théorème-clé \reff{t3.1}.
\sub{2bis.2} {Développements asymptotiques combinés Gevrey}
La définition d'un \dac{} Gevrey est assez proche de la
définition \reff{d2.2}.
\df{d2.3.2}
{Soit $V=V(\a,\b,\infty,\mu)$ un quasi-secteur infini ($\mu$ positif ou négatif),
$S_2=S(\a_2,\b_2,\e_0)$ un secteur fini et soit $\a_1<\b_1$ tels que
$\a\leq\a_1-\b_2<\b_1-\a_2\leq\b$. Soit
$y(x,\e)$
une fonction holomorphe définie quand $\e\in S_2$ et $x\in V(\a_1,\b_1,r_0,\mu\norm\e)$
et 
$\yc(x,\e)=\sum_{n\geq0}\gk{a_n(x)+g_n\big(\ts\frac x\e\big)}\e^n\in\cch(r_0,V)$.
Nous disons alors que {\sl $y$ admet $\yc$ comme \dac{} Gevrey 
d'ordre $\usp$ et de type $(L_1,L_2)$}
et nous écrivons {\sl $y(x,\e)\sim_{\frac1p}\yc(x,\e)$ quand $S_2\ni\e\to0$  
 et $x\in V(\a_1,\b_1,r_0,\mu\norm\e)$},
si $\yc(x,\e)$ est Gevrey d'ordre $\usp$ de type $(L_1,L_2)$ au sens de définition \reff{d2.3.1} et
s'il existe une constante $C$, telle que pour tout $N$, pour tout
$\e\in S_2$ et tout $x\in V(\a_1,\b_1,r_0,\mu\norm\e)$,
\eq{defcombgevrey}{
\norm{y(x,\e)-\sum_{n=0}^{N-1}\gk{a_n(x)+g_n\big(\ts\frac x\e\big)}\e^n}\leq 
C L_1^N \Gamma\big(\tfrac N p +1\big) \norm\e^N.
}
}
\rq 
Dans le cas d'une couronne, cette définition d'un \dac{} Gevrey est équivalente
à celle de développement asymptotique monomial Gevrey~; \cf \cit{cms},
définition 3.6. Pour voir ceci, on procède de la même manière que dans la remarque 1
après notre définition \reff{d2.2}. La condition \rf{gnmrem2} est automatiquement
satisfaite dans ce cas, car les $g_n(X)$ sont holomorphes. 
Par contre, la notion de {\sl sommabilité} monomiale de \cit{cms}
n'a pas été généralisée pour les \dacs{} dans notre mémoire.
Un premier pas dans cette direction est constitué par la proposition \reff{watson}.
\med

La compatibilité de la définition \reff{d2.3.2} avec l'addition et la dérivation est laissée en exercice. 
Concernant l'intégration, les énoncés précis sont exactement les mêmes que dans les propositions \reff{p40} et \reff{p4} en remplaçant le symbole $\sim$ par $\sim_\usp$ dans les hypothèses et dans les conclusions. En effet la majoration \rf{gnmc} montre que la suite $(g_{n 1})_{n\in\N}$ est Gevrey et le théorème de Borel-Ritt-Gevrey classique \cit{r1} (voir le début de la preuve du lemme \reff{brg1} dans la suite) fournit une fonction $R$ asymptotique Gevrey à  la série  $\hat R(\e)=\sum_{n=0}^\infty g_{n1}\e^n$. Le reste des preuves est identique.

Nous détaillons ci-dessous la compatibilité avec la multiplication. Nous ne la démontrons que pour le cas d'un produit ``mixte''
de $y(x,\e)\sim_{\frac1p}\widehat y(x,\e)=\sum_{n\geq0} a_n(x)\e^n$ et 
$z(x,\e)\sim_{\frac1p}\widehat z(x,\e)= \sum_{n\geq0}g_n\big(\frac x\e\big)\e^n$,
qui est le cas le plus intéressant. \fp Écrivons 
$$
y(x,\e)=\sum_{n=0}^{N-1} a_n(x)\e^n + P_N(x,\e)
$$ 
et
$$
z(x,\e)=\sum_{n=0}^{N-1} g_n\big(\tfrac x\e\big)\e^n + Q_N(x,\e)
$$
où les deux restes sont majorés par $CL_1^N\Gamma\big(\frac Np+1\big)\norm\e^N$.
La première majoration ci-dessous est classique~; c'est d'ailleurs la même pour les deux cas $a_1a_2$ et $g_1g_2$ non explicités.
$$
\begin{array}{rcl}\norm{R_N(x,\e)}&:=&
\norm{(y\cdot z)(x,\e)-\ds\sum_{n=0}^{N-1}\e^n\sum_{\nu=0}^n a_\nu(x)
g_{n-\nu}\big(\tfrac x\e\big)}\med\\
&=&\norm{\ds\sum_{n=0}^{N-1}P_n(x,\e)g_{N-n}\big(\tfrac x\e\big)\e^{N-n}
+y(x,\e)Q_N(x,\e)}\\
&\leq &C^2 L_1^N \norm\e^N\ds\sum_{n=0}^{N}
\Gamma\big(\ts\frac np+1\big)\Gamma\big(\frac {N-n}p+1\big)
\end{array}
$$
et (\reff{prodgamma}) implique $\norm{R_N(x,\e)}\leq K_p C^2 L_1^N\Gamma\big(\frac Np+1\big)$
avec la constante $K_p$ indiquée.

En utilisant la notation de (\reff{prodform}), il nous reste à  majorer
$${\sum_{n=0}^{N-1}\e^n\sum_{\nu=0}^n a_\nu(x)
g_{n-\nu}\big(\tfrac x\e\big)-\sum_{n=0}^{N-1} \Big(b_n(x)+h_n\big(\tfrac x\e\big)\Big)\e^n   }.$$
D'après (\reff{prodcomb}), utilisé sous la forme 
$$
a(x)g\big(\tfrac x\e\big)=
\sum_{l=0}^{m-1}\big(a_l(\T^lg)+g_l(\S^la)\big)(x,\e)\e^l+
 (\S^ma)(\T^mg)(x,\e)\e^m,
 $$
il reste à  majorer
$\ds\sum_{\nu+l+m=N}(\S^ma_\nu)(\T^mg_l)(x,\e)$
ce qui est analogue au cas du produit de séries formelles Gevrey traité dans la partie
\reff{2bis.1}. On obtient enfin comme désiré
$$
\norm{(y\cdot z)(x,\e)-\sum_{n=0}^{N-1} \Big(b_n(x)+h_n\big(\tfrac x\e\big)\Big)\e^n}\leq
 K_p C^2 L_1^N\Gamma\big(\tfrac Np+1\big)\norm\e^N
$$
avec la constante $K_p$ indiquée dans (\reff{prodgamma}).
\ep

Le résultat concernant les quotients de \dacs{} Gevrey et sa preuve
sont identiques au corollaire \reff{quot} et à sa preuve 
concernant les \dacs{} \og ordinaires \fg{} 
en ajoutant \og Gevrey\fg{} au mot \dac{}. Il suffit d'utiliser le 
théorème \reff{t4.7} (b) de la partie \reff{4.5}
au lieu de la proposition \reff{l2.4} (b) pour obtenir le \dac~de la composée d'un \dac~par l'inversion.
\med

L'analogue de la proposition \reff{combextint} pour l'asymptotique Gevrey est aussi
vrai.
\propo{matching-gevrey}{Soit $a_n(x)=\ds\sum_{m\geq0}a_{nm}x^m\in\H(r_0)$ et  $g_n(X)\sim\ds\sum_{m>0}g_{nm}X^{-m}\in\G(V)$.
Supposons que 
$$
y(x,\e)\sim_{\frac1p} \ds\sum_{n\geq0}\gk{a_n(x)+g_n\big(\ts\frac x\e\big)}\e^n
$$ 
quand $S_2\ni\e\to0$ et $x\in V(\a_1,\b_1,r_0,\mu\norm\e)$ au sens de la définition 
\reff{d2.2}. 
Alors, pour tout $r>0$, on a
$$
y(x,\e)\sim_{\frac1p}\sum_{n\geq0}c_n(x)\e^n \mbox{ \ quand }S_2\ni\e\to0
$$
uniformément par rapport à  $x$ sur l'ensemble des $x\in S(\a_1,\b_1,r_0)$ tels que $\norm x > r$, avec $c_n(x)=a_n(x)+\ds\sum_{0\leq l\leq n-1} g_{l,n-l}x^{l-n}$.

De même, étant donnés une partie compacte $K$ de $V$ et $\a_3,\b_3,\e_3$ tels que $X\in K$ et $\e\in S(\a_3,\b_3,\e_3)$ impliquent
$\e\in S_2$ et $\e X \in V(\a_1,\b_1,r_0,\mu\norm\e)$, on a
$$
y(\e X,\e)\sim_{\frac1p}\sum_{n=0}^\infty h_n(X)\e^n\mbox{ quand }S(\a_3,\b_3,\e_3)\ni\e\to0
$$
 uniformément pour $X\in K$, où $h_n(X)=\ds\sum_{0\leq l\leq n} a_{n-l,l}X^l  +g_n(X)$. 
}
\pr{Preuve}
 On adapte la démonstration de la proposition  \reff{combextint} en donnant aux constan\-tes
des formes spécifiques à  l'asymptotique Gevrey. D'après la définition, en particulier (\reff{gnm}), on a 
$C_N=C L_1^N \Gamma\big(\frac Np+1\big)$  et $C_{lk}=C L_1^l L_2^k \Gamma\big(\frac {k+l}p+1\big)$. On a aussi 
$$
||a_k||:=\ds\sup_{|x|<r_0}\ts\norm {a_k(x)}\leq C L_1^k\Gamma\big(\frac kp+1\big).
$$
Par ailleurs, en utilisant 
les inégalités de Cauchy, on a  $|a_{kl}|\leq||a_k||\;|r_0|^{-l}$. La fonction $b:x\mapsto a_k(x)-\ds\sum_{l<n}a_{kl}x^l$ est donc bornée par $(n+1)||a_k||$ sur le disque $|x|<r_0$. 
Ainsi, la fonction (analytique) $x\mapsto x^{-n}b(x)$ est bornée par $(n+1)||a_k||r^{-n}$ sur le cercle $|x|=r$ pour tout $r<r_0$. Par le principe du maximum, on en déduit que $|b(x)|\leq(n+1)||a_k||r_0^{-n}$. En résumé, on peut choisir
$A_{kn}=C L_1^k\Gamma\big(\frac kp+1\big)(n+1)r_0^{-n} $ dans \rf{AA}.

Pour le développement extérieur, il suffit donc de majorer dans \rf{cex}
$$
C_N+\sum_{l<N} C_{l,N-l}r^{l-N}\leq C \Gamma\big(\tfrac Np+1\big)\,
\sum_{l\leq N}L_1^l L_2^{N-l}r^{l-N}\leq C\big(L_1+\ts\frac{L_2}r\big)^N \Gamma\big(\tfrac Np+1\big).
$$
Pour le développement intérieur, on majore dans \rf{cin} en utilisant (\reff{gamma})
\begin{eqnarray}\lb{Cn}
C_N+\sum_{k<N} A_{k,N-k}R^{N-k}\!\!&\ds\leq& C 
\sum_{k\leq N} L_1^k\Gamma\big(\tfrac kp+1\big)(N-k+1)\gk{\tfrac R{r_0}}^{N-k}\nonumber\\&\ds
\leq&
\~ C L_1^N \Gamma\big(\tfrac Np+1\big)
\end{eqnarray}
avec 
\eq{ct}{
\~ C=C \ds\sum_{\nu\geq0}\ts(\nu+1)\gk{\frac R{r_0L_1}}^\nu
\frac1\Gamma\big(\tfrac \nu p+1\big).
}

\vspace{-5mm}
\ep

En relation avec  la remarque 2 qui suit la proposition  \reff{combextint}, nous n'avons pas trouvé d'estimation Gevrey utile pour $|x|>|\e|^\kappa$ ou $|X|<|\e|^{-\kappa}$, $\kappa\in\,]0,1[$. Nous n'avons pas non plus d'analogue Gevrey de la proposition \ref{p3.10bis}. En revanche, de tels analogues Gevrey existent pour
les proposition \reff{p3.12} et \reff{p3.11bis}~: 
sous des conditions similaires sur les développements intérieurs, 
resp.\ extérieurs, le domaine de validité d'un \dac{} Gevrey peut être 
étendu.
\propo{p4.5}{
Soit $y$ une fonction définie pour $\e\in S_2=S(\a_2,\b_2,\e_0)$ et $x\in V_1(\e)=V(\a_1,\b_1,r_0,\mu\norm\e)$
et ayant un \dac{}  Gevrey d'ordre $\usp$ au sens de définition \reff{d2.3.2}~:
$$
y(x,\e)\sim_\usp\sum_{n\geq0}\gk{a_n(x)+g_n\big(\ts\frac x\e\big)}\e^n
$$ 
quand $S_2\ni\e\to0$ et  $x\in V_1(\e)$, avec  
 $a_n\in\H(r_0)$ et $g_n\in\G(V)$,   $V=V(\a,\b,\infty,\mu)$ tels que
$\a\leq\a_1-\b_2<\b_1-\a_2\leq\b$.

Soit $\nu>\mu$. Dans le cas où $\nu>|\mu|$, on pose $\W=D(0,\nu)$, sinon on pose $\W=V(\a,\b,-\mu+\g,\nu)$ avec $\g>0$ arbitrairement petit.

On suppose que la fonction $Y:(X,\e)\mapsto y(\e X,\e)$ peut être 
prolongée analytiquement sur $\W\times S_2$  et qu'elle admet un
développement asymptotique  Gevrey 
$$
Y(X,\e)\sim_\usp\sum_{n=0}^\infty h_n(X)\e^n
$$
$\mbox{ quand }\e\mbox{ tend vers }0,$
uniformément sur $\W$. 

Alors $y$ peut être prolongée analytiquement sur l'ensemble des $(x,\e)$ avec 
$\e\in S_2$ et avec $x\in V(\a_1,\b_1,r_0,\nu\norm\e)$ et y admet un \dac{} 
Gevrey d'ordre $\usp$ quand $\e\to0$.}

\pr{Preuve}
Avec les notations des preuves précédentes, il s'agit d'une part de donner une majoration Gevrey à  $R_N(x,\e)$  donné par \rf{rest}, pour $x$ dans $\~V_{1}(\e)\setminus V_1(\e)$,
$\~V_1(\e)=V(\a_1,\b_1,r_0,\nu\norm\e)$
et d'autre part de montrer que les fonctions $g_n$ peuvent être prolongées analytiquement
et ont des développements Gevrey compatibles pour 
$X\in V(\a_1-\b_2,\b_1-\a_2,\infty,\nu)$. 

On majore $R_N$ en partant de \rf{rn} et en procédant comme pour \rf{Cn}. Par hypothèse, on peut choisir $D_N$ de la forme  
$D_N=C L_1^N \Gamma\big(\frac Np+1\big)$.
Avec la même constante $A_{kn}=C L_1^k\Gamma\big(\frac kp+1\big)(n+1)r_0^{-n} $ que précédemment, l'inégalité \rf{AA} est valide pour $|x|<r_0$, d'où
$$
|R_N(x,\e)|\leq\bigg(D_N+\sum_{k<N}A_{k\,N\!-\!k}M^{N-k}\bigg)|\e|^N\leq
\~ C L_1^N \Gamma\big(\tfrac Np+1\big) 
$$
avec $\~C$ donné par \rf{ct} et $M=\sup\{\norm x; x\in\W\}$.

Pour le caractère Gevrey compatible des fonctions $g_n$, on utilise la remarque 
\reff{remgnm}. Par hypothèse la série formelle $\sum_{n\geq0} h_n(X)\e^n$ est Gevrey d'ordre 
$\usp$ sur $\W$ , \ie il existe des constantes $H,L_1$ telles que pour tout $n\in\N$
\eq{gnmremh}{
\sup_{\norm X\in \W}\norm{h_n(X)}\leq H L_1^n \Gamma\big(\tfrac np+1\big).
}
Par définition d'un \dac{} Gevrey, la série formelle $\sum a_k(x)\e^k$ est aussi Gevrey
d'ordre $\usp$~; en diminuant au besoin  $L_1$, on en déduit qu'il existe une constante $A$ telle que 
$\ds\sup_{\norm x< r_0}\norm{a_k(x)}\leq A L_1^k \Gamma
\big(\tfrac kp+1\big)$ pour tout $k$.
Ceci implique par les inégalités de Cauchy que les coefficients de Taylor des $a_k$ 
satisfont $\norm{a_{kl}}\leq  A L_1^k \Gamma\big(\tfrac kp+1\big)r_0^{-l}$.
Utilisons maintenant la formule  $h_n(X)=\ds\sum_{0\leq l\leq n} a_{n-l,l}X^l  +g_n(X)$ de la proposition \reff{combextint}.
Elle implique en utilisant \rf{gamma} que 
$$
\norm{g_n(X)}\leq H L_1^n \Gamma\big(\tfrac np+1\big) + \sum_{0\leq l\leq n} 
  A L_1^{n-l} \Gamma\big(\tfrac {n-l}p+1\big)r_0^{-l}M^l
  \leq \tilde C L_1^n \Gamma\big(\tfrac np+1\big) 
$$ 
avec $\tilde C = H + \sum_{l\geq0} \gk{\tfrac M{L_1r_0}}^l\tfrac1\Gamma\big(\tfrac lp+1\big)$
pour $X\in\W$ et $n\in\N$. 
Les $g_n$ satisfont donc \rf{gnmrem1} pour $X\in\W$, et aussi 
pour $X\in V$ par hypothèse, donc pour  $X\in V(\a_1-\b_2,\b_1-\a_2,\infty,\nu)$. 
Enfin \rf{gnmrem2} est aussi satisfait par hypothèse.
Les hypothèses de la remarque \reff{remgnm} sont donc satisfaites et on peut conclure.
\ep
\propo{p4.5bis}{
Soit $0<r_0<\~r_0$ et soit $y$ une fonction définie pour $\e\in S_2=S(\a_2,\b_2,\e_0)$ et
$x\in V_1(\e)=V(\a_1,\b_1,\~r_0,\mu\norm\e)$. On suppose que 
$y(x,\e)\sim_\usp\ts\sum_{n\geq0}\gk{a_n(x)+g_n\big(\ts\frac x\e\big)}\e^n$,
quand $S_2\ni\e\to0$ et $x\in V_1(\e)=V(\a_1,\b_1,r_0,$ $\mu\norm\e)$, au sens de la définition \reff{d2.3.2}.

On suppose de plus qu'il existe $\g>0$ tel que $y$ admet un
développement asymptotique  Gevrey $y(x,\e)\sim_{\usp}\sum_{n=0}^\infty c_n(x)\e^n$ quand 
$\e$ tend vers $0$, uniformément pour $x\in V(\a_1,\b_1,r_0-\g,\~r_0)$.

Alors \rf{defcombgevrey} est satisfait pour tout $\e\in S_2$ et tout $x\in\~V_1(\e)$.
}
\rq 
De même que pour la proposition \reff{p3.11bis}, nous dirons par abus de langage que $y$ a
un \dac{} Gevrey pour $\e\in S_2$ et $x\in V_1(\e)$, même si les  fonctions
$a_n$ ne sont pas définies sur $D(0,r_0)$ en entier.
\med

\pr{Preuve}
D'après la proposition \reff{p3.11bis}, les fonctions $a_n$ peuvent être prolongées
analytiquement sur l'union de $D(0,r_0)\cup V(\a_1,\b_1,r_0-\g,\~r_0)$.
Il faut encore vérifier \rf{defcombgevrey} pour $x\in V_1(\e)$ avec 
$|x|\geq r_0$, \ie  majorer
$R_N$ donné par \rf{322}. D'après \rf{gnm}, il existe $C,L_1,L_2$ tels que pour tout $l<N\in\N$
$$
|r_{l\,N-l}(X)|\leq CL_1^lL_2^{N-l}\Gamma\big(\ts\frac Np+1\big)|X|^{-N+l}.
$$
Quitte à augmenter  les constantes $C$ et $L_1$, on a par hypothèse
$$
\Big|y(x,\e)-\sum_{n<N}c_n(x)\e^n\Big|\leq CL_1^N\Gamma\big(\ts\frac Np+1\big)|\e|^N.
$$
On en déduit
\begin{eqnarray*}
|R_n(x,\e)|
&\leq& C\lp L_1^N+\sum_{0\leq l<N}L_1^lL_2^{N-l} r_0^{l-N}\rp\Gamma\big(\ts\frac
Np+1\big)|\e|^N\\
&\leq& C\Big(L_1+\frac{L_2}{r_0}\Big)^N\Gamma\big(\ts\frac Np+1\big)|\e|^N.
\end{eqnarray*}
\vspace{-5mm}
\ep
\sub{2bis.3}{Fonctions plates Gevrey}
Comme pour les développements asymptotiques classiques, 
les fonctions plates sont particulièrement intéressantes.
Ici, on peut définir deux notions  de platitude, suivant que l'on demande aux fonctions $a_n$ et $g_n$ d'être identiquement nulles, ou seulement à  leurs coefficients $a_{nm}$ et $g_{nm}$.
Une fonction analytique étant déterminée par les coefficients de sa série de Taylor,
la situation pour les $g_n$ et pour les $a_n$ n'est pas symétrique.
\df{d14}{
Avec les notations de la définition \reff{d2.3.2}, on dit qu'une fonction
$y(x,\e)$ est {\sl plate au sens fort}
si elle admet un \dac{} et 
si toutes les fonctions $a_n$ et
$g_n$ de la série formelle correspondante sont  identiquement nulles. On dit qu'elle est
{\sl plate au sens faible} si elle a un \dac{} et si toutes les fonctions $a_n$ sont nulles et tous 
les coefficients $g_{nm}$ des développements asymptotiques (\reff{gnm}) des fonctions
$g_n$ s'annulent.
}
Un des points clés de la théorie classique des développements asymptotiques Gevrey est la relation
entre les fonctions plates Gevrey et les fonctions exponentiellement petites.
\propo{p2.3.1}{
Soit $y$ une fonction définie et analytique lorsque $\e\in S_2$ et $x\in  V(\a_1,\b_1,r_0,\mu\norm\e)$.
\be[\rm(a)]\item
Si $y$  est plate au sens fort, alors il existe $A,C>0$ tels que
\eq a{
\norm{y(x,\e)}\leq C\exp(- A/\norm\e^p)
}
pour $\e\in S_2$ et $x\in  V(\a_1,\b_1,r_0,\mu\norm\e)$.
\smallskip\\
Réciproquement, si $y$ satisfait \rf a, alors $y$ est plate au sens fort.
\item
Si  $y$ est plate au sens faible, alors il existe $B,C>0$ tels que 
\eq b{
\norm{y(x,\e)}\leq C\exp(-B\norm x^p/\norm\e^p)
}
pour $\e\in S_2$ et $x\in  V(\a_1,\b_1,r_0,\mu\norm\e)$.
\smallskip\\
Réciproquement,  si une fonction $y$ satisfait \rf b, et si de plus $y$ admet un \dac{}, alors $y$ est plate au sens faible.
\ee
}
\rq 
L'hypothèse  que $y$ admette un \dac{} est indispensable pour la réciproque du sens faible, comme le montre par exemple la fonction $y(x,\e)=\e^{1/2}\exp\big\{-\big(\frac x\e\big)^p\big\}$.
\med\\
\pr{Preuve}
Le (a) est classique. On a pour tout $N$ la majoration
\eq l{
\norm{y(x,\e)}\leq C L_1^N \Gamma\big(\tfrac N p +1\big) \norm\e^N. 
}
Il suffit de choisir $N$ proche de sa valeur optimale 
$N\sim p \gk{\tfrac{1}{\norm{\e}L_1}}^p$ pour obtenir la majoration avec $A=\frac1{L_1}$.

Réciproquement, une fonction satisfaisant \rf a satisfait \rf l pour tout $L_1>\frac1{A}$.
\med\\
(b) 
Comme $y$ est bornée d'après la définition \reff{d2.3.2}, il suffit de montrer l'inégalité
quand $x\in V_1(\e)=V(\a_1,\b_1,r_0,\mu\norm\e)$ satisfait $\norm{x}\geq K\norm{\e}$, où
$K>0$ est assez grand. Quitte à  agrandir $L_2$, on peut supposer que $r_0+\mu\e_0
\leq L_2/L_1$ et donc $\norm x\leq L_2/L_1$ pour tout $x\in V_1(\e)$.

D'après l'hypothèse, on a pour tout entier positif $N$
$$
\norm{y(x,\e)}\leq\sum_{n=0}^{N-1}\norm{g_n\big(\ts\frac x\e\big)}|\e^n|+ 
C (L_1\norm\e)^N \Gamma\big(\tfrac N p +1\big)
$$
et, pour tout $n$ et pour tout $M$ entier positif (cf.\ (\reff{gnm}))
$$
\norm{g_n(\tfrac x\e)}\norm\e^n\leq
C L_1^n (L_2/\norm x)^M \norm\e^{M+n}\Gamma\big(\tfrac{M+n}p+1\big).
$$
L'inégalité $L_1\leq L_2/\norm x$ et le choix de $M+n=N$ mènent à 
$$
\norm{y(x,\e)}\leq C (N+1)(L_2\norm{\e/x})^N\Gamma\big(\tfrac N p +1\big)
$$
pour tout entier positif $N$. En choisissant $N$ proche de sa valeur optimale
$N\sim p \gk{\tfrac{\norm{x}}{\norm{\e}L_2}}^p$  (ici on utilise que  
$\norm{x}\geq K\norm{\e}$ avec $K$ assez grand) on obtient
l'existence d'une constante $\~C$ telle que
$\norm{y(x,\e)}\leq \~C \norm{\tfrac x \e}^p \exp\gk{-\tfrac{\norm x^p}
{\norm\e^p L_2^p}}$. Ceci implique l'énoncé pour tout $0<B<L_2^{-p}$.

Réciproquement, si $y$ admet un \dac{} et satisfait \rf b, alors non seulement la partie lente de son \dac{}  est nulle, mais aussi son développement extérieur. La première partie de la proposition \reff{combextint} montre que les fonctions de sa partie rapide ont aussi un développement asymptotique nul, \cf \rf{cc}.
\ep

Dans le cadre des \dacs, nous avons aussi un équivalent du lemme de Watson classique :
si une fonction est plate sur des secteurs en $\e$ et en $x$ 
suffisamment grands, alors c'est la fonction nulle.
\propo{watson}{
Soit $\e_0,r_0,\mu>0$, $\w>\psi+\pi/p$ et $\d>\g+\pi/p$. 
On note
$$
{\cal D}=\{(x,\e)\tq \e\in S(\psi,\w,\e_0),\,x\in V(\g+\arg\e,\d+\arg\e,r_0,\mu|\e|)\}.
$$ 
Soit $y(x,\e)$ une fonction holomorphe
définie sur $\cal D$ telle que, pour tous 
$\psi\leq \a_2<\b_2\leq \w$ et $\a_1<\b_1$ satisfaisant
$\g+\b_2\leq\a_1<\b_1\leq\d+\a_2$, ses restrictions à  $\e\in S(\a_2,\b_2,\e_0)$ et à 
$x\in V(\a_1,\b_1,r_0,\mu\norm\e)$ ont un \dac{}
Gevrey d'ordre $\usp$ plat au sens faible. 
Alors $y$ est la fonction $0$. 
}
\pr{Preuve}D'après la proposition précédente, il existe des constantes $C,B>0$ telles que
la fonction $y$ satisfait $\norm{y(x,\e)}\leq C e^{-B\norm x^p/\norm\e^p}$ 
pour tout $(x,\e)\in{\cal D}$. \Apriori{}, les constantes dépendent des domaines choisis
pour les restrictions, mais comme on peut couvrir $\cal D$ par un nombre fini de tels
domaines, il suffit de prendre le maximum des majorations.

Dans la suite, on peut supposer sans perte que $\d-\g>\w-\psi$.
On prend $\a_2=\psi,\ \b_2=\w$, $\a_1=\g+\w$ et $\b_1=\d+\psi$ ;
la restriction $\~ y$ de $y$ à  $\e\in S(\a_2,\b_2,\e_0)$,
$x\in S(\a_1,\b_1,r_0)$ 
satisfait donc 
$\norm{\~ y(x,\e)}\leq C e^{-B\norm x^p/\norm\e^p}$.
Pour tout $x\in S(\a_1,\b_1,r_0)$ fixé, cette fonction est holomorphe et 
exponentiellement petite sur le secteur
$S(\a_2,\b_2,\e_0)$ d'angle d'ouverture $>\pi/p$. D'après le lemme de Watson classique,
on a donc $\~ y=0$. Le théorème d'identité des fonctions holomorphes implique
donc l'énoncé de la proposition.\ep

Étant donnés une série formelle $\sum a_n(x)\e^n$ Gevrey d'ordre $\usp$ et des coefficients
$g_{nm}$ satisfaisant la condition nécessaire (\reff{gnmc}), 
il existe donc au plus une fonction sur un domaine
$\cal D$ tel que décrit dans la proposition \reff{watson} ayant un \dac{} Gevrey d'ordre $\usp$
correspondant à  ces données. Dans une théorie de la resommation des \dacs{} 
qui reste à  faire, une telle fonction, si elle existe, pourrait donc être appelée la {\sl somme} de la série
doublement formelle $\sum_{n\geq0}\Big(a_n(x)+\widehat g_n\big(\tfrac x\e\big)\Big)\e^n$, 
$\widehat g_n(X)=\sum g_{nm}X^{-m}$ sur $\cal D$.
\sub{2bis.4}{Théorèmes de type Borel-Ritt}
De même que dans la théorie classique de l'asymptotique Gevrey, 
il nous sera utile de construire des fonctions ayant un \dac{} Gevrey prescrit. 
Comme dans la théorie classique, ce ne sera possible que
si la taille des secteurs est inférieure à  $\pi/p$. 
Nous présentons deux résultats de 
ce type. Dans le premier, le lemme \reff{brg1},  les fonctions $a_n,g_n$ sont données.
Dans le deuxième, le corollaire \reff{brg2}, ce sont les fonctions $a_n$ et les coefficients $g_{nm}$ des fonctions $g_n$ qui sont donnés.
\lem{brg1}{
Soit $r_0>0$, soit $\a,\b,\a_1,\b_1,\a_2,\b_2,\mu\in\R$ vérifiant $\a_1<\b_1$, $\a_2<\b_2$ et  et $\a\leq\a_1-\b_2<\b_1-\a_2\leq\b$. Soit $V=V(\a,\b,\infty,\mu)$ un quasi-secteur infini et $S_2=S(\a_2,\b_2,\e_0)$ un secteur fini. Enfin soit
$$
\yc(x,\e)=\sum_{n\geq0}\gk{a_n(x)+g_n\big(\ts\frac x\e\big)}\e^n\in\cch(r_0,V)
$$
une série formelle combinée Gevrey d'ordre $\usp$ selon la définition \reff{d2.3.1}.

\noindent On suppose que $\b_2-\a_2<\pi/p$. 

\noindent Alors il existe une fonction holomorphe $y$ définie pour $\e\in S_2$ et $x\in V(\a_1,\b_1,r_0,\mu\norm\e)$ qui admet $\yc$ comme \dac{} Gevrey d'ordre $\usp$.
}
\rq
La preuve montre que, si $\yc$ a pour type $(L_1,L_2)$, alors il existe $y$ ayant $\yc$
pour \dac{} de type $\~L_1$ pour tout $\~L_1>\frac{L_1}{\cos(p\psi)}$, avec $\psi=\frac12(\b_2-\a_2)$.

\pr{Preuve}
D'après l'hypothèse, on a pour tout 
$n\in\N$, $\ds\sup_{\norm x< r_0}\norm{a_n(x)}\leq C L_1^n 
\Gamma\big(\tfrac np+1\big)$ et aussi 
$\ds\sup_{\norm x< r_0}\norm{g_n(x)}\leq C L_1^n 
\Gamma\big(\tfrac np+1\big)$. 

Nous rappelons d'abord 
le théorème de Borel-Ritt-Gevrey, \ie la construction classique
d'une fonction ayant un développement asymptotique Gevrey d'ordre $\usp$ prescrit. Étant donnée une série $\widehat a(x,\e)=\ds\sum_{n=0}^\infty a_n(x)\e^n$ 
Gevrey d'ordre $\usp$ et de type $L_1$, on pose 
$\ds\breve a(x,t)=\sum_{n=0}^\infty a_{n}(x)\tfrac1\Gamma\big(\tfrac np+1\big)t^{n}$ ;
elle est appelée {\sl la transformée de Borel formelle } de $\widehat a(x,\eps)$. 
Nous l'avons un peu modifiée par rapport aux présentations les plus répandues : nous n'avons pas introduit de décalage, car cela simplifie notre présentation.
Cette série a un rayon de convergence au moins $1/L_1$.  
Soit $\f=\frac12(\a_2+\b_2)$, $\psi=\frac12(\b_2-\a_2)$,
$0<\rho<1/L_1$ et $T=\rho e^{i\f}$. 
À présent on pose
\eq{laplacetronq}
{\~ a(x,\e)={\cal L}_{T,p}(\breve a)(x,\e):=
  \e^{-p}\int_0^{T} e^{-t^p/\e^p}\,\breve a(x,t)\,d\,(t^p)\ \ ;}
elle est appelée la {\sl transformée de Laplace tronquée} de $\breve a$
(elle aussi un peu modifiée). Remarquons, qu'elle est définie pour tout $\e\in\C^*$ et $x\in D(0,r_0)$.
Les transformations de Borel formelle et de Laplace tronquée sont importantes en théorie Gevrey 
car elles transforment une série formelle Gevrey en une fontion l'admettant comme développement asymptotique
Gevrey. Ceci est basé sur la formule suivante.
$$
{\cal L}_{\infty e^{i\f},p}(t^n)=\Gamma\big(\tfrac n p +1\big) \e^{n}\ \  
\mbox{quand}\ \  \norm{\arg\e-\f}<\tfrac\pi p.
$$
Pour montrer que $\~ a(x,\e)$ admet $\widehat a(x,\e)$ comme développement 
asymptotique Gevrey d'ordre $\usp$, on écrit pour un $N$ donné
$\ds\breve a(x,t)=P_N(x,t)+ r_N(x,t),\ \mbox{où}$
$$P_N(x,t)=\sum_{n=0}^{N-1} a_{n}(x)\tfrac1\Gamma\big(\tfrac np+1\big)t^{n}\ \mbox{et}\ 
r_N(x,t)=\sum_{n=N}^\infty a_{n}(x)\tfrac1\Gamma\big(\tfrac np+1\big)t^{n}.$$
Alors pour $\e\in S_2 $ et pour $\norm{x}\leq r_0$, on a
$$
\~ a(x,\e)-\sum_{n=0}^{N-1} a_n(x)\e^n=-I+J\ \ \mbox{}
$$
avec $I=\e^{-p}\ds\int_T^{\infty e^{i\f}}e^{-t^p/\e^p}P_N(x,t)\,d(t^p)$ et 
$J=\e^{-p}\ds\int_0^{T}e^{-t^p/\e^p}r_N(x,t)\,d(t^p)$. 
Puisque 
$$
\norm{P_N(x,t)}\leq\rho^{-N}|t|^NC\ds\sum_{n=0}^{N-1}(L_1\rho)^n\,
\ \ \mbox{quand}\ \norm t\geq \rho\,,
$$ 
on a d'abord, avec la constante $K=\ds\frac C{1-L_1\rho}$ ,
\begin{eqnarray*}
\norm I
&\leq& K\norm\e^{-p}\rho^{-N}\int_0^{+\infty}\exp\gk{-\frac{s^p\cos(p\psi)}{\norm\e^p}}
s^{N}\,d(s^p)\\
&=&\frac K{\cos(p\psi)}\,(\rho\cos(p\psi)^{1/p})^{-N}
\Gamma\big(\tfrac{N}p+1\big)\norm{\e}^{N}.
\end{eqnarray*}
Avec la même constante $K$, on a l'inégalité
$\norm{r_N(x,t)}\leq K \rho^{-N}\norm t^{N}$
quand $t\in[0,T]$. On obtient donc de manière analogue
\begin{eqnarray*}
\norm J
&\leq& K\norm\e^{-p}\int_0^\rho \exp\gk{-\frac{s^p\cos(p\psi)}{\norm\e^p}}
\rho^{-N}s^{N}\,d(s^p)\\
&\leq& \frac K{\cos(p\psi)}\,\Big(\rho\big(\cos(p\psi)\big)^{1/p}\Big)^{-N}
\Gamma\big(\tfrac{N}p+1\big)\norm{\e}^{N}.
\end{eqnarray*}
Ceci démontre 
$\norm{\~ a(x,\e)-\ds\sum_{n=0}^{N-1} a_n(x)\e^n}\leq 
  \~ C {\~ L}^N\Gamma\big(\tfrac Np+1\big)\norm{\e}^N$
avec les constantes $\~ C=\ds\frac{2K}{\cos(p\psi)}=\frac{2C}{(1-L_1\rho)\cos(p\psi)}$ et  $\~ L=\Big(\rho\big(\cos(p\psi)\big)^{1/p}\Big)^{-1}$, autrement dit l'asymptotique Gevrey d'ordre $\usp$ de $\~ a$, 
uniformément pour $x\in D(0,r_0)$.

De manière analogue, 
 on construit $\~ g(X,\e)$ ayant $\widehat g=\ds\sum_{n=1}^\infty
g_n(X)\e^n$ comme développement asymptotique Gevrey d'ordre $\usp$
uniformément pour $X\in V(\a,\b,\infty,\mu)$.
Comme fonction $y$ ayant le \dac{} cherché, on peut alors choisir la fonction donnée par
$y(x,\e)=\~ a(x,\e)+\~ g \big(\tfrac x\e,\e\big)$.
\ep

Concernant le deuxième résultat, où seuls les coefficients $g_{nm}$ sont donnés avec la condition
\rf{gnmc},  il nous faut d'abord construire des fonctions $g_n(X)$
satisfaisant \rf{gnm}, ce que nous faisons dans l'énoncé ci-dessous.
\lem{brg2-lem}{
Soit $g_{nm}$, $n,m\in\N$ des nombres complexes vérifiant \rf{gnmc}. Soit $\a<\b<\a+\pi/p$ et $\mu\in\R$. Alors il existe $C',L'_1,
L'_2>0$ et
une suite de fonctions $g_n$ définies sur $V(\a,\b,\infty,\mu)$ telles que
 pour tout $X\in V(\a,\b,\infty,\mu)$
$$
\norm X^{M}\bigg|g_n(X)-\sum_{m=1}^{M-1}g_{nm}X^{-m}\bigg|
\leq C' L_1^{\prime\ n}  L_2^{\prime\ M}\Gamma\big(\ts\frac {M+n}p+1\big), 
$$
\ie $\sum_{n\geq0} g_n\big(\frac x\e\big)\e^n$ est Gevrey d'ordre $\usp$
au sens de la définition \reff{d2.3.1}.
}
\coro{brg2}{
Soit $g_{nm}$, $n,m\in\N$ des nombres complexes satisfaisant 
(\reff{gnmc}) et $a_n\in\H(r_0)$
des fonctions holomorphes telles qu'il existe
des constantes $C,L_1$ avec $\sup_{|x|<r_0}\norm{a_n(x)}\leq CL_1^n\Gamma\big(\tfrac{n}p+1\big)$
pour tout $n\in\N$. 
Soit $\a<\b<\a+\pi/p$,  $\a_2<\b_2<\a_2+\pi/p$,
$S_2=S(\a_2,\b_2,\e_0)$ un secteur fini, soit $\a_1<\b_1$ tels que
$\a\leq\a_1-\b_2<\b_1-\a_2\leq\b$ et soit $\mu\in\R$.

Alors il existe une suite de fonctions $(g_n^{\a,\b})_{n\in\N}$ définies sur 
$V(\a,\b,\infty,\mu)$ sa\-tisfaisant (\reff{gnm}) avec les coefficients
$g_{nm}$ donnés et
une fonction holomorphe $y(x,\e)$ définie quand $\e\in S_2$ et
$x\in V(\a_1,\b_1,r_0,\mu\norm\e)$ qui admet 
$\yc(x,\e):=\sum_{n\geq0}\gk{a_n(x)+g_n^{\a,\b}\big(\ts\frac x\e\big)}\e^n$ 
comme \dac{} Gevrey d'ordre $\usp$.
}
Le corollaire est une conséquence immédiate des deux lemmes précédents. 
\med\\
\pr{Preuve de lemme \reff{brg2-lem}} 
 Les séries $\widehat g_n(X)=\ds\sum_{m=n+1}^\infty g_{n,m-n}X^{-m}$ sont Gevrey d'ordre
$\usp$~; leurs coefficients satisfont précisément $\norm{g_{n,m-n}}\leq C 
\gk{\frac{L_1}{L_2}}^n\,L_2^m\Gamma\big(\frac mp+1\big)$ pour tout $n,m$.
La méthode présentée dans la démonstration du lemme précédent s'applique
à  chaque $\widehat g_n$, en remplaçant la variable $\e$ par $X^{-1}$ et
la constante $C$ par $C\gk{\frac{L_1}{L_2}}^n$.
Précisément, soit $T>\norm\mu$ et soit $\rho\in\,]0,1/L_2[$.
Posons $c=\cos\big(\frac{p(\b-\a)}2\big)$.

Alors on définit les fonctions $\~ g_n$ par
$$\~g_n(X)= X^{p}\int_0^{T} e^{-t^pX^p}\,\sum_{m=n+1}^\infty g_{n,m-n}\tfrac1\Gamma(\tfrac mp+1)t^m\,d\,(t^p)$$
sur  $V=V(\a,\b,\infty,\mu)$. Pour tout $M\geq n$ et tout $X\in V$, $\norm X\geq T$ on a
d'après les majorations dans la preuve du lemme précedent
$$
\norm{\~ g_n(X)-\sum_{m=n+1}^{M-1}g_{n,m-n}X^{-m}}
\leq \frac {2C}{(1-L_2\rho)c}\gk{\frac{L_1}{L_2}}^n
(\rho c^{1/p})^{-M}\Gamma\big(\tfrac Mp+1\big)\norm X^{-M}.
$$
En posant
\eq C{
\~ C=\frac{2C}{(1-L_2\rho)c},\qquad\~ \rho=\rho c^{1/p},
} 
en remplaçant $M$ par $M+n$ et 
en multipliant par $\norm X^n$, ceci implique que les fonctions
$g_n:X\mapsto\~ g_n(X)X^n$ satisfont
$$
\norm{g_n(X)-\sum_{m=1}^{M-1}g_{nm}X^{-m}}
\leq \~ C \gk{\frac{L_1}{L_2\~\rho}}^n\,\~\rho^{-M}
\Gamma\big(\tfrac{M+n}p+1\big)\norm X^{-M},
$$
\ie\ \rf{gnmrem2} est satisfaite pour $\norm X>T$ (avec les constantes $\~C$ de \rf C, $\~L_1=\frac{L_1}{L_2\~\rho}$ et $\~L_2=\frac1{\~\rho}$ au lieu de $C$, $L_1$ et $L_2$).

D'après la remarque \reff{remgnm}, il reste à  majorer $g_n(X)$ pour 
$X\in V(\a,\b,\infty,\mu)$, $\norm X\leq T$. Une majoration analogue à  celle
de $J$ dans la démonstration précédente montre que
\begin{eqnarray*}
\norm{\~ g_n(X)}
&\leq& \norm X^{p}\int_0^\rho \exp( s^p\norm X^p)\,
\frac{2C}{1-L_2\rho}\gk{\frac{L_1}{L_2}}^n\gk{\frac s\rho}^{n}\,d(s^p)\\
&\leq&
 \~ C\exp(\rho^pT^p)\gk{\frac{L_1}{L_2}}^n,
\end{eqnarray*}
où $\~ C$ est la constante définie dans \rf C. On a donc $\norm{g_n(X)}\leq \~ C\exp(\rho^pT^p)\big(\frac{L_1T}{L_2}\big)^n$ pour $|X|\leq T$. Ceci implique enfin (\reff{gnmrem1}) et complète ainsi la démonstration du lemme.
\ep

\sub{2bis.5}{Bons recouvrements cohérents}
À présent, nous allons utiliser les résultats de la partie \reff{2bis.4} pour associer, à  une fonction $y(x,\e)$ ayant un 
\dac{} Gevrey d'ordre $\usp$, des fonctions définies sur d'autres
secteurs dans la variable $\e$ et d'autres quasi-secteurs en $x$ ayant des 
\dacs{} Gevrey d'ordre $\frac1p$
de telle manière que, d'une part  les séries formelles 
$\sum_n a_n(x)\e^n$ et $\sum_n \widehat g_n(X)\e^n$
soient les mêmes pour toutes ces fonctions et
d'autre part l'union des domaines de définition contienne 
l'ensemble des $(x,\e)$ tels que $\norm\e<\e_0$ et $-\mu\norm\e<\norm x<r_0$.
Si $\mu>0$, cette union est donc simplement $D(0,\e_0)^*\times D(0,r_0)$.
D'après la proposition \reff{p2.3.1}, les différences de ces fonctions sont donc exponentiellement petites.

Soit $V=V(\a,\b,\infty,\mu)$ un quasi-secteur infini, soit $S=S(\psi,\w,\e_0)$
un secteur fini avec $0<\w-\psi<\b-\a$, soit $\kappa<\l$ tels que $\a<\kappa-\w<\l-\psi<\b$ et 
soit $y(x,\e)$ une fonction holomorphe définie quand $\e\in S$ et
$x\in V(\kappa,\l,r_0,\mu\norm\e)$ telle que 
\eq{ydonne}{y(x,\e)\sim_{\frac1p}\sum_{n=0}^\infty \Big(a_n(x)+
g_n\big(\tfrac x\e\big)\Big) \e^n \mbox{ quand } S\ni\e\rightarrow0} 
où $g_0,g_1,...$ sont des fonctions de $\G(V)$ satisfaisant (\reff{gnm}) avec des coefficients
$g_{nm}$. 

Pour associer une famille de fonctions à  $y$, on construit d'abord leurs domaines de 
définition. D'abord, on choisit des secteurs infinis $V^j=V(\a^j,\b^j,$ $\infty,\mu)$,
$j=2,...,J$, où $J$ est un entier supérieur ou égal à  $2$,  
tels que leurs angles d'ouverture soient inférieurs à  $\pi/p$ et tels
qu'en les complétant par $V^1:=V$, ils forment  un {\sl bon recouvrement de $\C$} si 
$\mu>0$, ou de la couronne infini $C(\norm\mu,\infty)=\{X\in\C\mid\norm\mu<
   \norm X<\infty\}$ si $\mu\leq0$. 
Ceci veut dire qu'on a $\bigcup_{j=1}^J V^j = \C$, resp.\ $C(\norm\mu,\infty)$,  
ainsi que $V^j\cap V^m$ borné si $j\not\in\{m-1,m,m+1\}$, \ie $\a^j<\b^{j-1}<\a^{j+1}$
 (par convention on pose $V^{J+1}=V^1$ et $V^{0}=V^J$, autrement dit les 
indices supérieurs sont pris mod $J$). 
\figu{f4.1}{
\vspace{-1.1cm}
\epsfxsize10.2cm\epsfbox{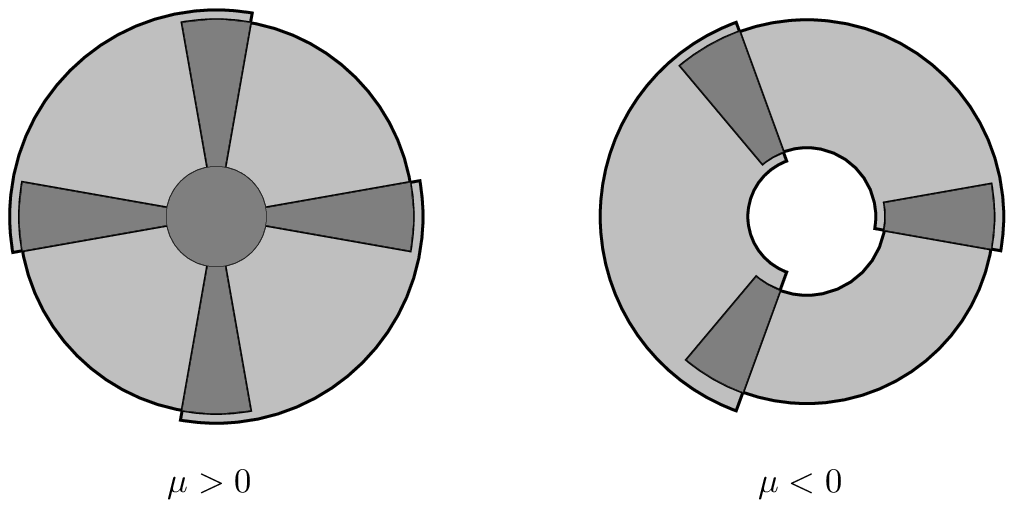}
\vspace{-5mm}
}{Exemples de recouvrements par des quasi-secteurs $(V_l^j(\e))_{0<j\leq J}$ pour $l$ fixé. Les rayons ont été choisis un peu différents pour que les figures soient plus lisibles}
Pour avoir de grands angles d'ouverture pour les quasi-secteurs en $x$,
il convient de choisir de petits secteurs en $\e$, pour que $X=x/\e$ soit encore
dans les $V^j$; \apriori{} le secteur donné $S$ ne sera donc pas inclus dans cette famille. 
Soit $\d>0$ strictement inférieur au minimum des demi-angles
d'ouverture des intersections $V^j\cap V^{j+1}$, \ie 
$2\d<\min\{\b^j-\a^{j+1}, j=1,...,J\}$ (avec la convention 
$\b^J-\a^{J+1}=\b^J-\a^1-2\pi$). On choisit ensuite une famille de secteurs 
finis $S_l=S(\a_l,\b_l,\e_0)$, $l=1,...,L$ 
tels qu'ils forment un {\sl bon recouvrement du disque
épointé $D(0,\e_0)^*$}, \ie leur union est $D(0,\e_0)^*$ et 
$S_l\cap S_m=\emptyset$
si $l\not\in\{m-1,m,m+1\}$, et tels que les angles d'ouverture satisfassent $\b_l-\a_l\leq\d$.
On prend ici et dans la suite les indices du bas mod $L$.
On pose $\f_l=(\a_l+\b_l)/2$ les directions bissectrices des secteurs $S_l$.
Enfin, on choisit les quasi-secteurs
$V_l^j(\e)=V(\a_l^j,\b_l^j,r_0,\mu\norm\e)$
avec $\a_l^j:=\a^j+\f_l+\d$ et $\b_l^j:=\b^j+\f_l-\d$. On  vérifie facilement
que, pour tout $\e\in S_l$, 
les familles des $V_l^j(\e)$, $j=1,...,J$ forment un bon recouvrement de
$D(0,r_0)$ si $\mu>0$,  de la couronne $C(-\mu\norm\e,r_0)$ si $\mu\leq 0$,
et que $\e\in S_l$ et $x\in V_l^j(\e)$ impliquent $x/\e\in V^j$.
\med

En résumé, la collection des  $V^j,S_l,V_l^j(\e)$ construits est un {\sl bon recouvrement 
cohérent}
au sens de la définition suivante. Pour pouvoir combiner les cas $\mu$ positif et $\mu$ négatif,
on utilise $C(\sigma,\rho)=\{x\in\C\mid \sigma<\norm x<\rho\}$ pour $\rho>0$, $\sigma\in\R$, $\sigma<\rho$.
\df{bonreccons}{
On appelle la collection des $V^j,S_l,V_l^j(\e)$, $j=1,...,J$, $l=1,...,L$
un {\sl bon recouvrement cohérent}
si $V^j=V(\a^j,\b^j,\infty,\mu)$, $j=1,...,J$, forment un bon recouvrement de $C(-\mu,\infty)$, si 
$S_l=S(\a_l,\b_l,\e_0)$, $l=1,...,L$, forment un bon recouvrement de 
$D(0,\e_0)^*$, s'il existe
$\d$ tel que
$$
\max\{\b_l-\a_l,l=1,...,L\}\leq\d<\ts\frac12 \min\{\b^1-\a^2,...,\b^{J-1}-\a^J,
\b^J-\a^1-2\pi\}
$$
et si $V_l^j(\e)=V(\a_l^j,\b_l^j,r_0,\mu\norm\e)$ avec $\a_l^j= \a^j+\f_l+\d$ et 
$\b_l^j= \b^j+\f_l-\d$, où $\f_l=(\a_l+\b_l)/2$.
On appelle le nombre $\max(\b_l-\a_l,l=1,...,L)$ la {\sl finesse} du  recouvrement.
}
Le lemme \reff{brg1} établit l'existence de fonctions $y_l^1$ holomorphes
quand $\e\in S_l$ et $x\in V_l^1(\e)$ telles que
$y_l^1(x,\e)\sim_{\frac1p} \sum_{n=0}^\infty \big(a_n(x)+g_n(\tfrac x\e)\big)\e^n$
avec les fonctions $a_n(x)$ et $g_n(X)$ de (\reff{ydonne}).
De plus, quitte à  diminuer au besoin $\mu$ et $\e_0$, le corollaire \reff{brg2} implique l'existence de suites de fonctions
$\big(g_n^j(X)\big)_{n\in\N}$, $j=2,...,J$, holomorphes sur $V^j$ satisfaisant (\reff{gnm}) avec les $g_{nm}$
de (\reff{ydonne}) et certaines constantes $C,L_1,L_2$ et des fonctions $y_l^j(x,\e)$ holomorphes sur l'ensemble
des $\e\in S_l$, $x\in V_l^j(\e)$ telles que 
\eq{yconstr}{
y_l^j(x,\e)\sim_{\frac1p}\sum_{n=0}^\infty \Big(a_n(x)+g_n^j\big(\tfrac x\e\big)\Big)\e^n.
} 
Pour compléter les $g_n^j$, on pose $g_n^1=g_n$ pour tout $n$.

On a ainsi une famille $y_l^j$ de fonctions définies sur un bon recouvrement
cohérent ayant des \dacs{} contenant les
mêmes séries formelles que la fonction $y$ donnée.
La proposition \reff{p2.3.1} implique  qu'il existe des constantes $A,B,C>0$ telles
que
\eq{expo-l}{
\norm{y_{l+1}^j(x,\e)-y_l^j(x,\e)}\leq Ce^{-A/\norm\e^p},
}
quand $j=1,...,J,\,l=1,...,L,\,\e\in S_{l+1}\cap S_l,\,
      x\in V_{l+1}^j(\e)\cap V_l^j(\e),$
\eq{expo-j}{
\norm{y_{l}^{j+1}(x,\e)-y_l^j(x,\e)}\leq Ce^{-B\norm x^p/\norm\e^p},
}
quand $j=1,...,J,\,l=1,...,L,\,\e\in S_l,\,
      x\in V_{l}^{j+1}(\e)\cap V_l^j(\e)$
ainsi que
\eq{expo-y}{
\norm{y_{l}^{1}(x,\e)-y(x,\e)}\leq Ce^{-A/\norm\e^p},
}
quand $l=1,...,L,\,\e\in S_l\cap S,$ si cette intersection est non vide,
et $x\in V_{l}^{1}(\e)\cap V(\kappa,\l,r_0,\mu\norm\e)$. 
Nous avons ainsi démontré ce qui suit.
\propo{p3.11}{
Soit $V=V(\a,\b,\infty,\mu)$ un quasi-secteur infini, soit $S=S(\psi,\w,\e_0)$
un secteur fini avec $\w-\psi<\b-\a$, soit $\kappa<\l$ tels que $\a<\kappa-\w<\l-\psi<\b$ et 
soit $y(x,\e)$ une fonction holomorphe définie quand $\e\in S$ et
$x\in V(\kappa,\l,r_0,\mu\norm\e)$ ayant un \dac{} Gevrey d'ordre $\usp$.

Alors 
il existe un bon recouvrement cohérent $V^j,S_l,V_l^j(\e)$, $j=1,...,J$, 
$l=1,...,L$ 
et une famille de fonctions holomorphes bornées $y_l^j(x,\e)$ définies
quand $\e\in S_l$ et $x\in V_l^j(\e)$ telle que les différences des fonctions 
sont exponentiellement petites, i.e.\ il existe des constantes $A,B,C>0$ telles que
\rft{expo-l}{expo-y} sont satisfaites.
}
Le partie suivante est consacrée à la démonstration d'une réciproque de cette
proposition.
%
%
%
%
\sec{3.}{Un théorème de type Ramis-Sibuya}
Nous présentons ici un résultat clé de ce mémoire~: des fonctions $y_l^j$ définies et analytiques sur un bon recouvrement cohérent au sens de la définition \reff{bonreccons} et  satisfaisant
\rff{expo-l}{expo-j} ont des \dacs{} Gevrey d'ordre $\usp$ avec les mêmes séries formelles
$\sum_n a_n(x)\e^n$ et $\sum_n \widehat g_n(X)\e^n$.
Le résultat est présenté pour  $\mu$ positif ou négatif, mais pour simplifier l'exposition, les sections \reff{3.3}, \reff{3.4} et \reff{secgn}, qui détaillent
les parties lentes et rapides de ces \dacs{} ne concernent que  $\mu$ positif. Les
 modifications à  apporter pour $\mu$ négatif sont présentées dans la section \reff{5.5}.
Dans la dernière section \reff{4.5}, nous utilisons ce résultat pour montrer la stabilité des \dacs{} Gevrey par composition avec une fonction analytique.
\sub{3.2}{Énoncé du théorème et début de la preuve}
\theo{t3.1}{Soit $J,L$ deux entiers naturels non nuls, $\mu$ un nombre réel positif ou négatif, et soit $V^j=V(\a^j,\b^j,\infty,\mu),\ S_l=S(\a_l,\b_l,\eta_0),\ 
V_l^j(\e)=V(\a_l^j,\b_l^j,r_0,\mu\norm\e)$, $j=1,...,J$, $l=1,...,L$, un bon recouvrement cohérent de finesse $\d$ au sens de la définition \reff{bonreccons}. 
Soit  $\~r_0>r_0$ et $\~\mu>\mu$. On pose $\~V_l^j(\e)=V(\a_l^j-2\d,\b_l^j+2\d,\~r_0,\~\mu\norm\e)$.
On suppose qu'il existe des fonctions holomorphes bornées $y_l^j(x,\e)$ définies
quand $\e\in S_l$ et $x\in\~V_l^j(\e)$ et des constantes $A,B,C$ positives 
telles que les différences satisfont les majorations suivantes~:
\eq6{
\left|y_{l+1}^j(x,\e)-y_{l}^j(x,\e)\right|\leq C\exp\big(-\ts\frac A{|\e|^p}\big)
}
si $\e\in S_l\cap S_{l+1}$ et $x\in \~V_l^j(\e)\cap \~V_{l+1}^j(\e)$ et
\eq7{
\big|y_{l}^{j+1}(x,\e)-y_{l}^j(x,\e)\big|\leq C\exp\big(-B\big|\ts\frac x\e\big|^p\big)
}
si $\e\in S_l$ et $x\in\~V_l^j(\e)\cap\~V_l^{j+1}(\e)$.
Alors les restrictions des fonctions $y_l^j$ sur l'ensemble des $(x,\e)$
tels que $\e\in S_l$, $x\in V_l^j(\e)$
admettent un \dac{} Gevrey d'ordre $\usp$ au sens de la définition \reff{d2.3.2}.
\medskip

Précisément, il existe une suite $(a_n)_{n\in\N}, a_n\in\H(r_0)$, pour chaque $j=1,...,J$ une suite $(g_n^j)_{n\in\N}, g_n^j\in\G(V^j)$ et une constante $\~ C$ 
telles que  pour tout couple d'indices $(j,l)$, tout $\e\in S_l$, tout 
$x\in V_l^j(\e)$,  et tout entier $N>0$
\eq8{
\Big|y_{l}^j(x,\e)-\sum_{n=0}^{N-1}\lp a_n(x)+g_n^j\big(\ts\frac x\e\big)\rp\e^n\Big|\leq \~ C\,\~ A^{-N/p}\Gamma\big(\ts\frac Np+1\big)\,|\e|^N\ ;
}
ici et dans la suite, on a posé $\~ A=\min(A,Br_0^p)$ avec les constantes $A,B$ de \rf{6},\rf{7}.
De plus, les suites de fonctions $(g_n^j)_{n\in\N}$, ont des développements asymptotiques
Gevrey d'ordre $\usp$ compatibles à  l'infini de type $(\~ A^{-1/p},B^{-1/p})$ 
avec des coefficients indépendants de $j$, 
i.e.\ 
il existe une suite double $(g_{nm})_{n,m\in\N}$ de nombres complexes et une constante
$\~ C$ telles que pour tout $M,n\in\N,\,j\in\{1,...,J\}$ et tout $X\in V^j$ 
\eq9{\norm X^{M}
\Big|g_n^j(X)-\sum_{m=1}^{M-1}g_{nm}X^{-m}\Big|\leq 
     \~ C\, \~ A^{-n/p}\,B^{-M/p}\Gamma\big(\ts\frac {M+n}p+1\big).}
}
\newcommand{\ext}{{\mbox{\scriptsize\rm ext}}}
\newcommand{\inte}{{\mbox{\scriptsize\rm int}}}
\newcommand{\ye}{y_l^\ext}
\newcommand{\yep}{y_{l+1}^\ext}
\newcommand{\yek}{y_k^\ext}
\newcommand{\yepk}{y_{k+1}^\ext}
\newcommand{\yi}{y_l^{j\,\inte}}
\newcommand{\yiun}{y_l^{1\,\inte}}
\newcommand{\yiJ}{y_l^{J\,\inte}}
\newcommand{\yip}{y_{l+1}^{j\,\inte}}
\newcommand{\yipj}{y_l^{j+1\,\inte}}
\rq L'asymptotique des fonctions $y^j_l$ ne peut être valide que dans des quasi-secteurs $V_l^j(\e)$ d'angle d'ouverture plus petit que leurs
quasi-secteurs de définition $\~V_l^j(\e)$. En effet les fonctions $g_n^j$ sont construites à  partir des $y_l^j$, 
 donc les quasi-secteurs  $V^j$ sur lesquels elles sont  définies  doivent déjà  avoir un angle d'ouverture plus petit que les $\~V_l^j(\e)$ ; ensuite la condition que les quotients $x/\e$
doivent être dans le domaine $V^j$ de $g_n^j$ entraîne que l'asymptotique ne peut
être valable que dans des quasi-secteurs $V_l^j(\e)$ d'angle d'ouverture encore 
plus petite que $V^j$. Précisément, l'ouverture des $\~V_l^j(\e)$ est $\b^j-\a^j+2\d$, celle des $V^j$ est $\b^j-\a^j$, celle des $V_l^j(\e)$ est $\b^j-\a^j-2\d$ et on a les
implications suivantes
\eq x{
\mbox{ si }\e\in S_l\mbox{ et }x\in V_l^j(\e),\mbox{ alors }x/\e\in V^j
}
et réciproquement
\eq X{
\mbox{ si }X\in V^j\mbox{ et }\e\in S_l\mbox{ avec }|\e X|<\~r_0,
  \mbox{ alors }\e X\in \~V_l^j(\e).
}
En ce qui concerne les constantes  $\~r_0$ et
$\~\mu$, elles peuvent être arbitrairement proches de  $r_0$ et
$\mu$ respectivement mais pas égales. 
\med\\
{\sl Principe de la preuve} \sep
L'idée est d'écrire chaque fonction $y_l^j$ comme la somme de deux fonctions $\ye$ et $\yi$ qui ont chacune un développement asymptotique en puissances de $\e$ Gevrey d'ordre $\usp$, l'une $\ye$ indépendante de $j$ et ayant 
pour développement la partie lente $\sum a_n(x)\e^n$, et l'autre $\yi$ ayant pour développement la partie rapide $\sum g_n^j\big(\tfrac x\e\big)\e^n$. 
Pour construire ces deux fonctions, nous utilisons la {\sl formule de Cauchy-Heine} par rapport à  la variable $x$ en  choisissant des points au bord des intersections 
$\~V_{l}^{j,j+1}(\e):=\~V_l^j(\e)\cap\~V_l^{j+1}(\e)$.
Dans le cas $\mu\geq0$, il suffit de choisir des points $x_l^j$ avec $\norm{x_l^j}=r_0$,
voir ci-dessous les formules \rf e et \rf i~; dans le cas $\mu<0$, on suppose sans perte que $\mu<\~\mu<0$ et il faut choisir, de 
plus, des points $\~x_l^j(\e)$ avec $\norm{\~x_l^j(\e)}=\norm{\~\mu\e}$, voir  
les formules \rf{emod} et \rf{imod} dans la partie \reff{5.5}.
 Ces fonctions dépendent des points $x_l^j$, mais pas des $\~x_l^j(\eta)$. 
Pour des raisons techniques, les fonctions $\ye$ et $\yi$ n'ont le comportement voulu que sur des secteurs excluant les points $x_l^{j-1}$ et $x_l^j$. Cependant, ni les fonctions $y_l^j$, ni les fonctions $a_n$ et $g_n^j$ ne dépendent des points $x_l^j$, si bien qu'en modifiant au besoin les points $x_l^k$ on obtiendra le comportement voulu pour $y_l^j$ sur tout le quasi-secteur $V_l^j(\e)$.
\med\\
Des notations similaires à  $\~V_{l}^{j,j+1}(\e)$ seront utilisées partout dans la preuve.
\med\\
Quitte à  choisir pour $A$ la constante $\~ A=\min(A,Br_0^p)$ de l'énoncé, on peut supposer dans toute la partie \reff{3.}
\eq z{
A\leq B r_0^p ;
}
\figu{f5.1}{
\vspace{-2mm}
\epsfxsize11.6cm\epsfbox{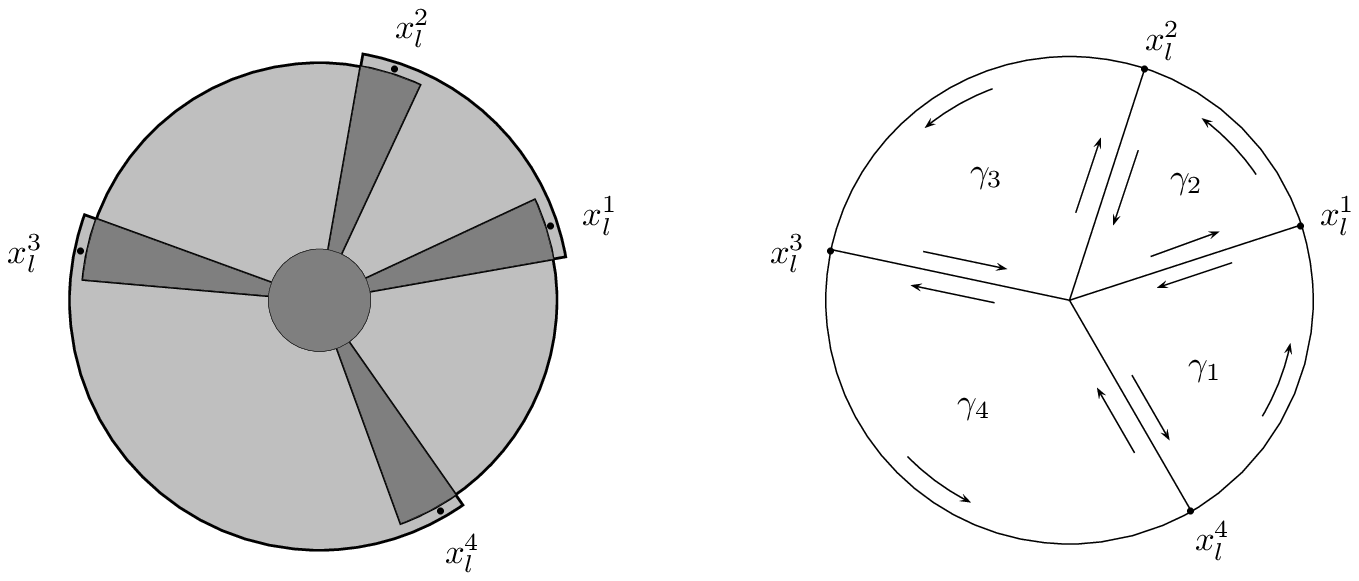}
\vspace{-8mm}
}{A gauche, le recouvrement $\big(\~V_{l}^{j}(\e)\big)_{0<j\leq J}$, à droite les chemins $\g^j$}
Dans la suite, nous ne traitons que le cas $\mu\geq0$~; les modifications nécessaires 
pour le cas $\mu<0$ sont indiquées dans la partie \reff{5.5}.
Pour chaque $j\in\{1,...,J\}$, soit $\psi^j\in\,]\a^{j+1}-\d,\b^j+\d[$.
Pour chaque $(j,l)\in\{1,...,J\}\times\{1,...,L\}$ on pose
$x_l^j=\~r_0 e^{i(\psi^j+\f_l)}$. Ainsi on a  
$x_l^j\in\cl\lp \~V_{l}^{j,j+1}(\e)\rp$ pour tout $|\e|<\e_0$.
Notons $\g^j$ le lacet allant de $0$ à  $x_l^{j-1}$ le long du segment, puis de $x_l^{j-1}$ à  $x_l^j$ le long d'un chemin arbitrairement proche de l'arc de cercle de rayon $r_0$ et enfin de $x_l^j$ à  $0$ le long du segment. 
Pour $\e\in S_l$ et $x\in\~V_l^j(\e)$ avec $\arg x_l^{j-1}<\arg x<\arg x_l^j$, puisque seul $\g^j$ entoure $x$ (une fois et dans le sens positif), la formule des résidus donne
\eq r{
\int_{\g^j}\frac{y_l^j(u,\e)}{u-x}\,du=2\pi i\,y_l^j(x,\e)~\mbox{ et }~
\int_{\g^k}\frac{y_l^k(u,\e)}{u-x}\,du=0~\mbox{ si }~k\neq j.
}
En réordonnant la somme de ces intégrales, on obtient la {\sl formule de Cauchy-Heine}
$$
y_l^j(x,\e)=\ye(x,\e)+\yi(x,\e)
$$
avec
\eq e{
\ye(x,\e)=\frac1{2\pi i}\ds\sum_{k=1}^J\int_{x_l^{k-1}}^{x_l^k}\frac{y_l^k(u,\e)}{u-x}\,du
}
où l'intégration se fait le long d'un chemin arbitrairement proche de l'arc de cercle de rayon $\~r_0$, et
\eq i{
\yi(x,\e)=\ds\frac1{2\pi i}\ds\sum_{k=1}^J\int_0^{x_l^k}\frac{y_l^{k+1}(u,\e)-y_l^k(u,\e)}{u-x}\,du.
}

Par sa définition (\reff{e}), la fonction $\ye(x,\e)$ est holomorphe bornée sur 
$D(0,r_0)\times S_l$. Par la relation $\yi(x,\e)=y_l^j(x,\e)-\ye(x,\e)$, la fonction
$\yi$ est holomorphe bornée quand $\e\in S_l$, $x\in \~V_l^j(\e)$, $\norm x\leq r_0$ ;
la formule (\reff{i}) prolonge $\yi$ dans le secteur {\sl infini}
$S(\arg x_l^{j-1},\arg x_l^{j},
\infty)$ en $x$. 
Concernant les fonctions $\ye$, nous démontrerons dans la section \reff{3.3} le résultat suivant.
\lem{lm3.2}{
Il existe une suite $(a_n)_{n\in\N}, a_n\in\H(r_0)$ et une constante $C>0$ telles que,  pour  tout entier $N>0$, tout indice $l\in\{1,...,L\}$ et tout couple $(x,\e)\in D(0,r_0)\times S_{l}$ 
$$
\Big|\ye(x,\e)-\sum_{n=0}^{N-1} a_n(x)\e^n\Big|\leq C A^{-N/p} \Gamma\big(\ts\frac Np+1\big)\,|\e|^N.
$$
}
Au vu de la dépendance des $x_l^j$ par rapport à  $l$, il convient d'introduire
les fonctions
$$Y_l^j(X,\e):=\yi(\e X,\e)$$ 
sur l'ensemble des $(X,\e)$ tels que $\e\in S_l$ et, ou bien $\arg(X)\in\,]\psi^{j-1}+\d,\psi^{j}-\d[$,
ou bien ($\norm{\e X}\leq r_0$ et $X\in V^j$). On note cet ensemble $E_l^j$ ; il dépend
du choix des $\psi^j$. Pour ces fonctions, on a le résultat suivant.
\lem{lm3.5}{
Il existe une constante $C>0$ et,
pour chaque $j=1,...,J$, une suite $(g_n^j)_{n\in\N}$ de fonctions holomorphes dans $V^j$,
tendant vers 0 quand 
$\arg X\in ]\psi^{j-1}+\d,\psi^{j}-\d[,$ $\norm X\rightarrow\infty$
telles que pour tout $l\in\{1,...,L\}$, tout $N\in\N^*$ et tout $(X,\e)$
avec $\e\in S_l$ et $\arg X\in\,]\psi^{j-1}+\d,\psi^{j}-\d[ $ on a 
\eq g{
\Big|Y_{l}^j(X,\e)-\sum_{n=0}^{N-1}g_n^j(X)\e^n\Big|\leq C A^{-N/p}
   \Gamma\big(\ts\frac Np+1\big)\,|\e|^N.  } 
De plus, pour tout compact $K$ de $V^j$, il existe $\e_0>0$ tel que pour tout
$l\in\{1,...,L\}$, tout $N\in\N^*$, tout $\e\in S_l$ avec $\norm\e<\e_0$ 
et tout $X\in K$, on ait l'inégalité \rf g.
}
Les fonctions ainsi que la constante dans ce lemme dépendent \apriori{} des 
$\psi^j,\ j=1,...,J$. Ce lemme est démontré dans la partie \reff{3.4}.

Comme $y^j_l(x,\e)=\ye(x,\e)+Y_l^j\big(\tfrac x\e,\e\big)$,
on obtient ainsi presque un \dac. Précisément il existe une constante
$\~ C$ telle que pour tout $j,l,N$ et tout $\e\in S_l$, $x\in S(\psi^{j-1}+\f_l+2\d,
  \psi^j+\f_l-2\d,r_0)$
\eq{temp}{
\Big|y_{l}^j(x,\e)-\sum_{n=0}^{N-1}\lp a_n(x)+g_n^j\big(\ts\frac x\e\big)\rp\e^n\Big|
     \leq \~ C\,A^{-N/p}\Gamma\big(\ts\frac Np+1\big)\,|\e|^N.
}
Puisque les fonctions $g_n^j$ tendent vers 0 quand $\norm X\rightarrow\infty$,
celles-ci et les $a_n$ sont déterminées uniquement par les inégalités (\reff{temp}) d'après
la remarque 3 après la définition \reff{d2.2}.
Par conséquent, les $a_n$ et les $g_n^j$ sont indépendantes du choix des $\psi^j$, 
\ie des $x_l^j$, et les inégalités (\reff{temp}) sont donc valables dans toute union
finie de secteurs en $x$ qu'on peut obtenir par des choix différents de ces points.
En utilisant la dernière phrase du lemme, on peut compléter en des quasi-secteurs
dans le cas $\mu>0$.

Ceci démontre \rf8 ; on ne peut pas encore parler d'un \dac, car on n'a pas encore la
propriété requise que les $g_n^j(X)$ admettent des développements asymptotiques quand
$X\rightarrow\infty$.
Les inégalités (\ref9), \ie le fait que les $g_n^j(X),n\in\N$ admettent 
des développements asymptotiques Gevrey d'ordre $\usp$ compatibles, 
seront démontrées dans section \reff{secgn}.
Ceci complétera la démonstration du théorème.
\sub{3.3}{La partie lente}
La formule \rf e définit pour chaque indice $l$ une fonction analytique bornée sur 
$D(0,r_0)\times S_l$,  et les différences satisfont, pour tout 
$x\in D(0,r_0)$ et tout $\e\in S_{l,l+1}:=S_l\cap S_{l+1}$
\begin{eqnarray*}
2\pi i\lp\yep(x,\e)-\ye(x,\e)\rp&=&\ds\sum_{k=1}^J\int_{x_{l+1}^{k-1}}^{x_l^k}\frac{y_{l+1}^k(u,\e)-y_l^k(u,\e)}{u-x}\,du+\\
&&\sum_{k=1}^J\int_{x_l^k}^{x_{l+1}^k}\frac{y_{l+1}^k(u,\e)-y_l^{k+1}(u,\e)}{u-x}\,du.
\end{eqnarray*}
Par hypothèse,  sur l'arc de $x_{l+1}^{k-1}$ à  $x_l^k$ on a 
$$
\norm{y_{l+1}^k(u,\e)-y_l^k(u,\e)}\leq C\exp\big(-\ts\frac A{|\e|^p}\big).
$$
Pour la deuxième intégrale, on utilise le fait que, suivant les valeurs de $u$ sur l'arc de $x_l^k$ à  $x_{l+1}^k$, une (au moins) des deux identités suivantes est valide~:
$y_{l+1}^k(u)-y_l^{k+1}(u)=\lp y_{l+1}^k(u)-y_l^{k}(u)\rp+\lp y_{l}^k(u)-y_l^{k+1}(u)\rp$
ou bien
$y_{l+1}^k(u)-y_l^{k+1}(u)=\lp y_{l+1}^k(u)-y_{l+1}^{k+1}(u)\rp+\lp y_{l+1}^{k+1}(u)-
y_l^{k+1}(u)\rp$.
En utilisant \rf{6} et \rf{7}, 
et avec des chemins d'intégration $\eps$-proches de l'arc de cercle, $\eps$ arbitrairement petit, 
on obtient
\begin{eqnarray*}
\norm{y_{l+1}^k(u,\e)-y_l^{k+1}(u,\e)}&\leq& C\lp\exp\big(-\ts\frac A{|\e|^p}\big)+
  \exp\big(-\ts\frac{B \~ r_0^p}{|\e|^p}\big)\rp\\
   &\leq& 2C\exp\big(-\ts\frac A{|\e|^p}\big),
\end{eqnarray*}
puisque $A\leq Br_0^p$ et $r_0<\~ r_0$.
Sur les chemins d'intégration on a $|u-x|\geq \~r_0-r_0 -\eps$.
Ainsi, puisque $\eps$ est arbitraire, 
on a pour tout $(x,\e)\in D(0,r_0)\times S_{l,l+1}$ 
\eq{ep}{
\left|\yep(x,\e)-\ye(x,\e)\right|\leq \tfrac{2C}{\~ r_0-r_0}
   \exp\big(-\ts\frac A{|\e|^p}\big).
}
En utilisant le lemme de Ramis-Sibuya classique (voir lemme \reff{RS} ci-dessous), 
on obtient que les fonctions $\ye$ ont un développement asymptotique commun Gevrey 
d'ordre $\usp$ uniformément dans $D(0,r_0)$. 
Pour la complétude du mémoire, nous donnons une preuve de ce résultat.
\lem{RS}{
On suppose que $S_l=S(\a_l,\b_l,\e_0),\,l=1,...,L$ 
forment un bon recouvrement du disque 
épointé $D(0,\e_0)^*$. Soit $f_l:S_l\rightarrow \C$, l=1,...,L, 
des fonctions holomorphes bornées par
une constante $C>0$. On suppose qu'il existe des constantes $D,A>0$ telles que
$\norm{f_{l+1}(\e)-f_l(\e)}\leq D e^{-A/\norm\e^p}$ pour tout $l=1,...,L$ et tout
$\e\in S_{l,l+1}$. 
Alors les fonctions $f_1,...f_L$ admettent un développement asymptotique Gevrey $\usp$ 
commun.
Précisément, il existe une constante $\~C$ et une suite $(c_n)_{n\in\N}$ telles que
$$
\norm{f_l(\e)-\sum_{n=0}^{N-1}c_n\e^n}\leq \~C A^{-N/p} 
    \Gamma\big(\tfrac Np+1\big)\norm\e^N
$$
pour tout $N\in\N^*$, tout $l$ et tout $\e\in S_l$.
}
\rq 
La constante $\~C$ ne dépend que des données $A,C,D,\e_0,\a_l,\b_l$. 
Les coefficients $c_n$ peuvent être exprimées par des intégrales ; voir la formule
\rf{cn} ci-dessous. Par conséquent, si les fonctions dépendent de paramètres de façon analytique,
resp.\ continue, si les fonctions sont uniformément bornées et si les constantes $D$ et $A$
sont indépendantes des paramètres, alors 
les $c_n$ sont aussi analytiques, resp.\ continus, par rapport à  ces 
paramètres et les majorations des restes dans le lemme sont uniformes.
\med\\
\pr{Preuve}
Fixons $\g>0$ tel que $\b_l-\a_{l+1}\geq2\g$ pour tout $l\in\{0,...,L\}$.
On pose $T_l=\e_0e^{i\phi_l}$ avec des $\phi_l\in[\a_{l+1},\b_l]$ à  choisir.
Par un réarrangement de formules de résidus analogues à  \rf r, on obtient
la formule de Cauchy-Heine classique
\eq y{
f_l(\e)=\frac1{2\pi i}\ds\sum_{k=1}^L\lp
\int_0^{T_k}\frac{f_{k+1}(v)-f_k(v)}{v-\e}\,dv
+
\int_{T_{k-1}}^{T_k}\frac{f_k(v)}{v-\e}\,dv
     \rp,
}
où, de même que précédemment, l'intégration de $T_{k-1}$ à  $T_k$ 
se fait le long d'un chemin arbitrairement proche de l'arc de cercle de rayon $\e_0$.

Pour $n\in\N$ posons
\eq{cn}{
c_n=\frac1{2\pi i}\ds\sum_{k=1}^L\lp
   \int_0^{T_k}\Big(
f_{k+1}(v)-f_k(v)\Big)\,\frac{dv}{v^{n+1}}
+
\int_{T_{k-1}}^{T_k}f_k(v)\,\frac{dv}{v^{n+1}}
\rp.
}
Fixons à  présent $N\in\N^*$ 
et posons
$$
r(\e)=f_l(\e)-\sum_{n=0}^{N-1}c_n\e^n.
$$
En remplaçant $\ds\frac1{v-\e}$ par
$\ds\sum_{n=0}^{N-1}\frac{\e^n}{v^{n+1}}+
 \frac{\e^N}{v^N(v-\e)}$ dans \rf y, on obtient 
\eq{r1r2}{
r(\e)=r_1(\e)+r_2(\e)
}
avec
$$
r_1(\e)=\frac{\e^N}{2\pi i}\,\sum_{k=1}^L\int_0^{T_k}
\frac{f_{k+1}(v)-f_k(v)}{v^N(v-\e)}\,dv
$$
et 
$$
r_2(\e)=\frac{\e^N}{2\pi i}\,\sum_{k=1}^L\int_{T_{k-1}}^{T_k}
\frac{f_k(v)}{v^N(v-\e)}\,dv.
$$
Supposons d'abord $\a_l+\g\leq\arg\e\leq\b_l-\g$ et $\norm\e\leq(1-\g)\e_0$.
Alors, en utilisant $\phi_{l-1}=\a_l$ et $\phi_l=\b_l$, on obtient 
$\norm{v-\e}\geq\norm v\sin\g$ sur tous les segments $[0,T_k]$. En tenant compte de l'hypothèse, ceci implique
$$
|r_1(\e)|
  \leq\frac 1{2\pi}\int_0^{\e_0}\frac{De^{-A/t^p}}{\sin\g}\frac{dt}{t^{N+1}}|\e|^N
  \leq\frac{D\Gamma\big(\frac Np\big)}{2\pi pA^{N/p}\sin\g}\,|\e|^N
$$
Pour l'autre partie du reste, on a 
$$
|r_2(\e)|\leq\frac
{C}{2\pi}\int_0^{2\pi}\frac{\e_0d\theta}{\e_0^{N+1}\g}|\e|^N=\frac
{C}{\e_0^N\g}|\e|^N.
$$
La combinaison des deux majorations donne l'énoncé avec une constante
$\~ C$, qui dépend de $C,D,\e_0,\gamma$ et de $\sup_N \eta_0^{-N} A^{N/p}
       \frac1\Gamma\big(\frac Np+1\big)$.

Dans le cas où $\arg\e-\a_l<\g$, on a $\b_{l-1}-\arg\e\geq\g$ et en utilisant un 
autre choix d'angle $\phi_k$, on obtient une majoration analogue du reste pour $f_{l-1}.$
Comme $\norm{f_{l-1}(\e)-f_l(\e)}\leq D e^{-A/\norm\e^p}$ et 
$\norm\e^{-N}e^{-A/\norm\e^p}\leq A^{-N/p}\Gamma\big(  \tfrac Np+1\big)$, on conclut
que la majoration du reste est aussi valable pour $f_l$.

Dans le cas $\b_l-\arg\e<\g$, on procède de manière analogue avec $f_{l+1}$ à  la place de
$f_l$ ; dans le cas $\e_0(1-\g)\leq \norm\e\leq \e_0$, la majoration du reste
est faite simplement en utilisant les majorations des $c_n$ et de $f_l$,
si on choisit $\~C$ assez grand.\ep
\sub{3.4}{La partie rapide}
Considérons à  présent, pour $l=1,...,L$, $j=1,...,J$,
les fonctions 

\begin{eqnarray}\lb{Ylj}
Y_l^j(X,\e)=\yi(\e X,\e)&=&
\ds\frac1{2\pi i}\ds\sum_{k=1}^J\int_0^{x_l^k}\frac{y_l^{k+1}(u,\e)-y_l^k(u,\e)}
       {u-\e X}\,du\nonumber\\
&=&y_l^j(\e X,\e)-\ye(\e X,\e)
\end{eqnarray}
définies sur
l'ensemble $E_l^j$ des $(X,\e)$ tels que $\e\in S_l$ et, ou bien
$\arg(X)\in ]\psi^{j-1}+\d,\psi^{j}-\d[$, ou bien
 $\norm{\e X}<r_0$, $X\in V^j$. Rappelons que l'ensemble $E_l^j$ dépend
du choix des $\psi^j$ et que les notations de la section \reff{3.2} sont utilisées.

Pour démontrer le lemme \reff{lm3.5}, nous montrons d'abord
\lem{lm3.4}{
Il existe une constante $\~ C>0$ 
telle que pour tout $(j,l)\in\{1,...,J\}\times\{1,...,L\}$ 
et tout $(X,\e)\in E_{l,l+1}^j$ 
on ait 
\eq{Ydif}{
\left|Y_{l+1}^j(X,\e)-Y_{l}^j(X,\e)\right|\leq \~ C\exp\big(-\ts\frac A{|\e|^p}\big).
}
}
\pr{Preuve}
Si $X\in V^j$ et $|\e X|<r_0$, alors $\e X\in \~ V_{l,l+1}^j(\e)$.
On utilise dans ce cas
$$
Y_{l+1}^j(X,\e)-Y_{l}^j(X,\e)=(y_{l+1}^j-y_{l}^j+\ye-\yep)(\e X,\e)
$$
et les majorations \rf6 et \rf{ep}.
Si  $\arg(X)\in ]\psi^{j-1}+\d,\psi^{j}-\d[$ et $\norm{\e X}\geq r_0$, 
on reprend la formule intégrale définissant $\yi$, qui donne
$$\begin{array}{l}
Y_{l+1}^j(X,\e)-Y_{l}^j(X,\e)\med\\=\ds\frac1{2\pi i}\ds\sum_{k=1}^J\lp\int_0^{x_{l+1}^k}\, \frac{(y_{l+1}^{k+1}-y_{l+1}^k)(u,\e)}{u-\e X}\,du
  -\int_0^{x_l^k}\frac{(y_{l}^{k+1}-y_{l}^k)(u,\e)}{u-\e X}\,du\rp\med\\\med
=\ds\frac1{2\pi i}\ds\sum_{k=1}^J\lp\int_0^{\xi_{l}^k}\frac{(y_{l+1}^{k+1}-y_{l}^{k+1}+y_{l}^k-y_{l+1}^k)(u,\e)}{u-\e X}\,du\
+\right.\\\ds
\left.\qquad\qquad\qquad\qquad
\int_{\xi_l^k}^{x_{l+1}^k}\frac{(y_{l+1}^{k+1}-y_{l+1}^k)(u,\e)}{u-\e X}\,du
+\int_{x_l^k}^{\xi_{l}^k}\frac{(y_{l}^{k+1}-y_{l}^k)(u,\e)}{u-\e X}\,du\rp
\end{array}$$
où $\xi_l^k=\~ r_0 e^{i(\psi^j+\gamma_l)}$, $\f_l<\gamma_l<\f_{l+1}$ 
désigne un point dans l'intersection $\~ V_l^{k,k+1}(\e) 
\cap \~ V_{l+1}^{k,k+1}(\e)$.
Sur chacun des chemins, le dénominateur $u-\e X$ est minoré en module par 
$r_0\sin(\d/2)$, 
le numérateur $y_{l+1}^{k+1}-y_{l}^{k+1}+y_{l}^k-y_{l+1}^k$ dans la première 
intégrale est majoré par $2C\exp\big(-\ts\frac A{|\e|^p}\big)$ d'après \rf6, 
et les numérateurs des deuxième et troisième intégrales sont majorés par
$C \exp\big(-\ts\frac{B r_0^p}{|\e|^p}\big)
   \leq C\exp\big(-\ts\frac A{|\e|^p}\big)$ d'après \rf7 et \rf z.
\ep

Le lemme \reff{lm3.5} découle directement du lemme précédent par le lemme de Ramis-Sibuya 
classique \reff{RS} ; la dépendance analytique des fonctions $g_n^j(X)$ 
obtenues  découle de la formule (\reff{cn}) définissant ces coefficients, comme indiqué
dans la remarque qui suit le lemme \reff{RS}. De la même façon, on déduit d'abord
de leurs définitions par des intégrales que toutes les $Y_l^j(X,\e)$ tendent vers
$0$ quand $X\rightarrow\infty$, $\arg(X)\in  ]\psi^{j-1}+\d,\psi^{j}-\d[$,
si $\e$ est dans un compact ne pas contenant 0 et ensuite que les $g_n^j(X)$
tendent vers 0 quand $X\rightarrow\infty$.

\sub{secgn}{Étude des fonctions $g_n^j$}
\fp
Rappelons que les fonctions $g_n^j$ du développement asymptotique de la fonction $Y_l^j$ du
lemme \reff{lm3.5} sont obtenues par l'application du lemme
\reff{RS} et donc définies par
\begin{eqnarray}\lb{defgnj}
g_n^j(X)&=&\frac1{2\pi i}\sum_{k=1}^L\Bigg(\int_{\e_{k-1}}^{\e_k}Y_k^j(X,v)\frac{dv}{v^{n+1}}+\med\\\nonumber
  &&\hspace{1.8cm}\int_0^{\e_k}\lp Y_{k+1}^j(X,v)-Y_k^j(X,v)
  \rp\frac{dv}{v^{n+1}}\Bigg)
\end{eqnarray}
où $\e_k\in\cl(S_{k,k+1})$, $k=1,...,L$ satisfont 
$\norm{\e_1}=...=\norm{\e_L}=:r\leq \e_0$,
et que les fonctions $Y_l^j$ sont définies par (\reff{Ylj}) ; la formule est donc valide
pour tout $X\in V^j$ tel que, ou bien $r\norm X<r_0$, ou bien $\psi^{j-1}<\arg X<\psi^j$. 
Rappelons également que
l'asymptotique (\reff{temp}) détermine les fonctions $g_n^j$ de manière unique et que celles-ci 
sont donc indépendantes du choix des $x_l^j$ dans (\reff{Ylj}) et du choix des $\eta_k$ 
dans (\reff{defgnj}). Notons enfin une majoration de ces fonctions qui est une 
conséquence de (\reff{g}) et de la liberté de choix des $x_l^j$ : 
il existe une constante $C$ telle que 
\eq{majgnj}{
\norm{g_n^j(X)}\leq C A^{-n/p} \Gamma\big(\tfrac np+1\big)
}
pour tout $n$ et $X\in V^j$.

Il reste à  montrer que les fonctions $g_n^j,\ n\in\N$ admettent des développements 
asymptotiques Gevrey $\frac1p$
compatibles au sens de la définition \reff{d2.3.1}. Ceci sera fait en 
adaptant la preuve du lemme \reff{RS}~: on majore d'abord
les différences $Y_l^{j+1}-Y^j_l$, ensuite $g_n^{j+1}-g_n^j$ et enfin les restes
des développements asymptotiques.

D'abord d'après \rf X, si $\eta\in S_l$ et $X\in V^{j,j+1}$ sont 
tels que $|\e X|<r_0$, alors $\e X\in \~ V_l^{j,j+1}(\e)$. 
De plus, puisque $\ye$ ne dépend pas de $j$, on a dans ce cas
$$
Y_l^{j+1}(X,\e)-Y_{l}^j(X,\e)=y_l^{j+1}(\e X,\e)-y_l^{j}(\e X,\e)
$$
et \rf7 s'applique.
On obtient, pour $X\in V^{j,j+1}$ et 
$\e\in S_{l}$ tels que $|\e X|< r_0$, la majoration
\eq{Yj}{
\big|Y_l^{j+1}(X,\e)-Y_{l}^j(X,\e)\big|\leq C\exp\lp- B|X|^{p}\rp.
}
On utilise maintenant (\reff{defgnj}) pour majorer $g_n^{j+1}-g_n^j$.
Contrairement au lemme \reff{lm3.5}, nous voulons ici un résultat asymptotique lorsque 
$X$ tend vers l'infini dans le secteur $V^{j,j+1}$. Pour cela, nous allons utiliser 
l'estimation exponentielle des différences $Y_l^{j+1}-Y_l^j$ donnée par (\reff{Yj}).
Par conséquent le module $r$ des points $\e_l$ ne peut plus 
être fixé, mais doit nécessairement dépendre de $X$.
On a
\begin{eqnarray*}
g_n^{j+1}(X)-g_n^j(X)&=&\frac1{2\pi i}\sum_{l=1}^L\Bigg(
\int_{\e_{l-1}}^{\e_l}\lp Y_{l}^{j+1}-Y_l^j\rp(X,v)\frac{dv}{v^{n+1}}+\\
&&\qquad
\int_{0}^{\e_l}\lp Y_{l+1}^{j+1}-Y_{l}^{j+1}-Y_{l+1}^j+Y_l^j\rp(X,v)\frac{dv}{v^{n+1}}\Bigg)
\end{eqnarray*}
où les intégrales de $\e_{l-1}$ à  $\e_l$ sont arbitrairement proches des arcs du cercle de rayon $r$. 
Sur ces arcs, on a bien $\e X\in \~V_l^{j,j+1}(\e)$. On choisit $r=r_0/\norm X$, quand
$\norm X\geq L:=\e_0/r_0+1$. Alors en utilisant (\reff{Yj}), 
on obtient qu'il existe une constante $\~ C$
telle que, pour tout $\norm X\geq L$ et tout $(n,j)$, la somme des intégrales de $\e_{l-1}$
à  $\e_l$ est inférieure à  $\~C\,r_0^{-n}\norm X ^n\exp(-B\,\norm X^p)$.

Aux intégrales de 0 à  $\e_l$, on peut appliquer (\reff{Ydif}) si on écrit le numérateur
$(Y_{l+1}^{j+1}-Y_{l}^{j+1})-(Y_{l+1}^j-Y_l^j)$ et on peut appliquer \rf{Yj} si on 
l'écrit $(Y_{l+1}^{j+1}-Y_{l+1}^{j})-(Y_{l}^{j+1}-Y_l^j)$. Ceci entraîne 
qu'il existe une constante $\bar C$ telle que la 
somme des intégrales de $0$ à  $\e_l$ est majorée par 
$$\bar C\int_0^{r}\min\gk{e^{-A/v^p},e^{-B\norm X^p}}v^{-n-1}\,dv\ ;$$
rappelons que $\norm X\geq L$ et $r=r_0/\norm X$.
En déterminant $q>0$ par l'équation $A/q^p=B$, ceci devient
$$
\bar C\gk{\int_0^{q/\norm X} e^{-A/v^p}v^{-n-1}\,dv+
      e^{-B\norm X^p}\int_{q/\norm X}^{r_0/\norm X}v^{-n-1}\,dv},
$$
où $\norm X\geq L\geq1.$ 
En faisant un changement de variable dans la première intégrale, on en déduit 
que cette somme d'intégrales de $0$ à  $\e_l$ est majorée par 
$$
\bar C\gk{\tfrac1pA^{-n/p}\,\Gamma\big(\tfrac np,B\norm X^p\big)+ \tfrac1nq^{-n}\,\norm X^n\,e^{-B\norm X^p}}\ ;
$$
ici $\Gamma(a,u)=\int_u^\infty e^{-t}t^{a-1}\,dt$, $\re a>0$, est la fonction Gamma
incomplète. Il est préférable de laisser ici la fonction Gamma incomplète car des
majorations plus simples sont différentes si $n$ est petit, respectivement grand.

En combinant les deux majorations, on en déduit l'existence d'une cons\-tante $C$
telle que 
\eq{majdgnj}{
\norm{g_n^{j+1}(X)-g_n^j(X)}\leq C\gk{A^{-n/p}\,\Gamma\big(\tfrac np,B\norm X^p\big)
        + q^{-n}\,\norm X^n\,e^{-B\norm X^p}
        }      }
pour tout $n,j$ et $X\in V^{j,j+1}$, $\norm X\geq L$.

Utilisons maintenant la formule de Cauchy-Heine dans la variable $Z=\frac1X$. 
Considérons le secteur $V^j$ et d'abord seulement les $X\in V^j$ avec $\norm X\geq L+1$
et $\a^j+\delta<\arg X<\b^j-\d$.
Pour chaque $k=1,...,J$ soit 
$X_k\in V^{k,k+1}$ tels que $|X_k|=L$ ; en particulier
$\arg X_{j-1}=\a^j$ et $\arg X_{j}=\b^j$.
Soit $Z_k=\frac1{X_k}$. En posant 
$$
g_{mn}=\frac1{2\pi i}\sum_{k=1}^J\lp\int_0^{Z_k}\lp g_n^{k+1}\big(\ts\frac1\zeta\big)-g_n^k\big(\ts\frac1\zeta\big)\rp\frac{d\zeta}{\zeta^{m+1}}+\int_{Z_k}^{Z_{k+1}}g_n^k\big(\ts\frac1\zeta\big)\frac{d\zeta}{\zeta^{m+1}}\rp,
$$
on obtient ainsi, par un calcul analogue à  celui de $r_1$ et $r_2$ 
dans la preuve du lemme \reff{RS}
\begin{eqnarray*}
X^M\gk{ g_n^j(X)-\sum_{m=0}^{M-1}g_{mn}X^{-m}}\!\!&=\!\!&\frac1{2\pi i}\sum_{k=1}^J
\lp\int_0^{Z_k}\frac{g_n^{k+1}\big(\ts\frac1\zeta\big)-g_n^k\big(\ts\frac1\zeta\big)}
{\zeta^M\lp\zeta-\ts\frac1X\rp}\,d\zeta\right.+\\
 &&\qquad\qquad\qquad\qquad
\left.\int_{Z_{k-1}}^{Z_{k}}\frac{g_n^k\big(\ts\frac1\zeta\big)}{\zeta^M\lp\zeta-\frac1X\rp}\,d\zeta\rp
\end{eqnarray*}
pour $X\in V^j$, $\norm X\geq L+1$, $\a^j+\delta<\arg X<\b^j-\d$. 
En utilisant la majoration \rf{majdgnj} et la minoration 
$\big|\zeta-\frac1X\big|\geq|\zeta|\sin\d$
sur les chemins d'intégration, on majore en module la somme des intégrales de $0$ à  $Z_k$ par
$$
\frac{J\,C}{2\pi\sin\d}\lp A^{-n/p}\int_0^{+\infty}\,
   \Gamma\big(\tfrac np,B\,t^{-p}\big)\,\frac{dt}{t^{M+1}}+q^{-n}\int_0^{+\infty} 
      e^{-B\,t^{-p}}\,\frac{dt}{t^{M+n+1}}\rp
$$
et, après le changement de variable $s=Bt^{-p}$ et en utilisant la formule
$$
\int_0^\infty s^{a-1}\Gamma(b,s)ds=\frac1a\,\Gamma(a+b),
$$
on majore cette somme d'intégrales
par 
$$
\frac{J\,C}{2\pi p\sin\d}\big(\tfrac pM+1\big)A^{-n/p}B^{-M/p}\,\Gamma\big(\tfrac{n+M}p\big).
$$
La majoration des intégrales de $Z_{k-1}$ à  $Z_{k}$ suit de manière plus simple 
de \rf{majgnj}. Leur somme est majorée par 
$$
J\,C\,A^{-n/p}\big(\tfrac1L-\tfrac1{L+1}\big)^{-1}\Gamma\big(\ts\frac np+1\big)\,L^{-M}.
$$
Comme il existe une constante $D$ telle que $\Gamma\big(\ts\frac np+1\big)\,L^{-M}\leq
D B^{-M/p}\,\Gamma\big(\tfrac{n+M}p+1\big)$ pour tout $n,M\in\N$, 
on en déduit l'existence d'une constante $C$ telle que
\eq{majcomp}{
\norm{X^M\bigg(g_n^j(X)-\sum_{m=0}^{M-1}g_{mn}X^{-m}\bigg)}\leq
    C\,A^{-n/p}B^{-M/p}\,\Gamma\big(\tfrac{n+M}p+1\big)
    }
pour tout $n,M$ et $X\in V^j$ avec $\norm X\geq L+1$ et $\a^j+\delta<\arg X<\b^j-\d$. 

Pour $X\in V^j$ avec $\arg X$ proche de $\a^j$ ou $\arg X$ proche de $\b^j$, 
on procède de manière analogue
à  la fin de la preuve du lemme \reff{RS}. Il suffit de traiter le premier cas : 
on considère d'abord les majorations
pour $g_n^{j-1}$ et on ajoute des termes dus aux majorations des différences
$g_n^j-g_n^{j-1}$. 
On utilise pour cela le fait 
qu'il existe une constante $D$
telle que $T^M\Gamma\big(\tfrac np,BT^p\big)\leq D B^{-M/p}\,\Gamma\big(\tfrac{n+M}p+1\big)$ et
$T^{n+M}e^{-BT^p}\leq D B^{-M/p}\,\Gamma\big(\tfrac{n+M}p+1\big)$ pour tout $n,M,T,B>0$ 
(utiliser le fait que pour $\a,\b>0$, 
la fonction $f(t)=t^\a\Gamma(\b,t)$ admet son maximum en un point $t$ tel que
$f(t)=\frac1\a t^{\a+\b}e^{-t}$).

Ainsi, on a enfin démontré qu'il existe une constante $C$ telle que \rf{majcomp}
est valable pour { tout} $n,j,M$ et {\sl tout} $X\in V^j$ avec $\norm X\geq L+1$. 
En vertu de la remarque \reff{remgnm} et de \rf{majgnj}, ceci suffit pour montrer
que les suites des $g_n^j$, $n=0,1,...$ admettent des développements Gevrey d'ordre
$\frac1p$ compatibles au sens de la définition \reff{d2.3.1}.
La preuve du théorème \reff{t3.1} dans le cas $\mu\geq0$ est à  présent complète.
\ep
\sub{5.5}{Le cas $\mu$ négatif} 
Nous indiquons dans la suite les modifications de la démonstration précédente 
nécessaires dans le cas $\mu<\~\mu<0$. 

Comme avant, pour chaque $j\in\{1,...,J\}$, soit $\psi^j\in\,]\a^{j+1}-\d,\b^j+\d[$.
Pour chaque $(j,l)\in\{1,...,J\}\times\{1,...,L\}$ on pose
$x_l^j=\~r_0 e^{i(\psi^j+\f_l)}$ et pour tout $\eta\neq0$ on pose
$\~x_l^j=\~x_l^j(\e)=\norm{\~\mu\e} e^{i(\psi^j+\f_l)}$ . Ainsi on a  
$x_l^j,\~x_l^j\in\cl\lp \~V_{l}^{j,j+1}(\e)\rp$ pour tout $|\e|<\e_0$.
Notons $\g^j$ un lacet inclus dans $\~V_l^{j,j+1}(\e)$ allant de $\~x_l^{j-1}$ à  $x_l^{j-1}$ le long du segment, 
puis de $x_l^{j-1}$ à  $x_l^j$ le long d'un chemin arbitrairement proche de l'arc de 
cercle de rayon $r_0$, ensuite de $x_l^j$ à  $\~x_l^j$ le long du segment et enfin de
$\~x_l^j$ à $\~x_l^{j-1}$  arbitrairement proche de l'arc de cercle de rayon 
$\norm{\~\mu\eta}$. 
Pour $\e\in S_l$ et $x\in\~V_l^j(\e)$ avec $\arg x_l^{j-1}<\arg x<\arg x_l^j$, 
on obtient comme avant
$${
\int_{\g^j}\frac{y_l^j(u,\e)}{u-x}\,du=2\pi i\,y_l^j(x,\e)~\mbox{ et }~
\int_{\g^k}\frac{y_l^k(u,\e)}{u-x}\,du=0~\mbox{ pour tout }~k\neq j
}$$
et en réordonnant la somme de ces intégrales, on obtient la {formule de Cauchy-Heine}
$$
y_l^j(x,\e)=\ye(x,\e)+\yi(x,\e)
$$
avec
\eq{emod}{
\ye(x,\e)=\frac1{2\pi i}\ds\sum_{k=1}^J\int_{x_l^{k-1}}^{x_l^k}\frac{y_l^k(u,\e)}{u-x}\,du
}
où l'intégration se fait le long d'un chemin arbitrairement proche de l'arc de cercle de rayon $\~r_0$, et
\begin{eqnarray}
\lb{imod}
\yi(x,\e)&=&\ds\frac1{2\pi i}\ds\sum_{k=1}^J\int_{\~x_l^k}^{x_l^k}
\frac{y_l^{k+1}(u,\e)-y_l^k(u,\e)}{u-x}\,du-\\
&&
\frac1{2\pi i}\ds\sum_{k=1}^J\int_{\~x_l^{k-1}}^{\~x_l^k}\frac{y_l^k(u,\e)}{u-x}\,du
\nonumber
\end{eqnarray}
où l'intégration se fait le long d'un chemin arbitrairement proche de l'arc de cercle de 
rayon $\norm{\~\mu\e}$, \cf figure \reff{f5.3}.
Puisque ni $y_l^j$ ni $\ye$ ne dépendent des  points $\~x_l^k$, les fonctions $\yi$ n'en dépendent pas non plus.
\figu{f5.3}{\vspace{-.4cm}
\epsfxsize5cm\epsfbox{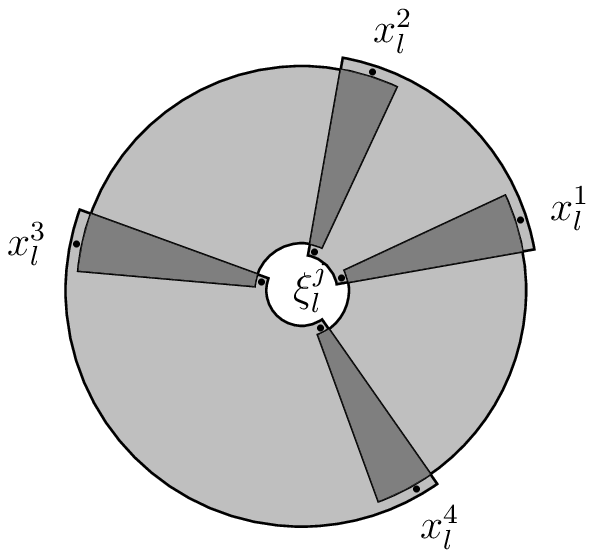}\vspace{-.5cm}\hspace{1.5cm}\epsfxsize5cm\epsfbox{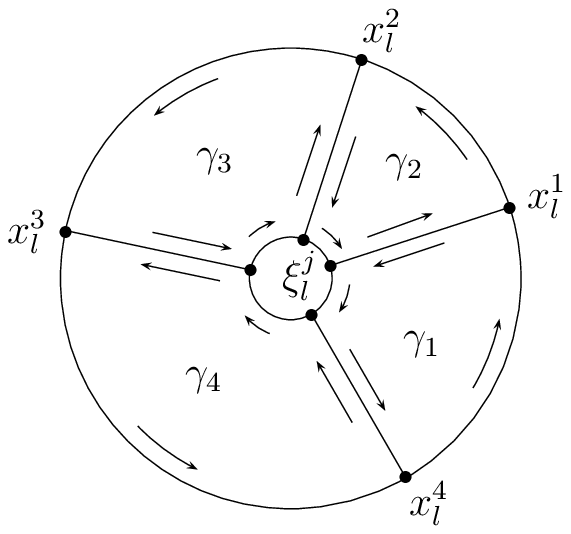}
}{Les chemins d'intégration pour $\mu<0$}
Le traitement des fonctions $\ye$ est identique au cas $\~\mu>0$ ; pour les fonction
$\yi$, il faut vérifier que les intégrales de $\~x_l^{k-1}$ à $\~x_l^{k}$, qui ne 
figuraient pas dans la formule \rf i du cas $\~\mu>0$, n'introduisent aucun problème.

D'abord, les fonctions $Y_l^j(X,\e)=\yi(\e X,\e)$ sont définies  
sur l'ensemble des $(X,\e)$ tels que $\e\in S_l$ et, ou bien 
($\arg(X)\in\,]\psi^{j-1}+\d,\psi^{j}-\d[$ et $\norm X>\~\mu$)
ou bien ($\norm X>\~\mu$, $\norm{\e X}\leq r_0$ et $X\in V^j$). 
On note cet ensemble $E_l^j$ ; il dépend
du choix des $\psi^j$ et de $\norm\e$.

Les lemmes \reff{lm3.4} et \reff{lm3.5} 
 restent valables. En effet, la formule exprimant la différence des fonctions 
dans le cas $\norm{\e X}\geq r_0$ devient
$$
\begin{array}l
Y_{l+1}^j(X,\e)-Y_{l}^j(X,\e)=
\med\\
\qquad\ds\frac1{2\pi i}\ds\sum_{k=1}^J\Bigg(\int_{\~\xi_l^k}^{\xi_{l}^k}\frac{(y_{l+1}^{k+1}-y_{l}^{k+1}+y_{l}^k-y_{l+1}^k)(u,\e)}{u-\e X}\,du+
\med\\
\qquad\qquad\qquad\ds
\int_{\xi_l^k}^{x_{l+1}^k}\frac{(y_{l+1}^{k+1}-y_{l+1}^k)(u,\e)}{u-\e X}\,du+
\int_{x_l^k}^{\xi_{l}^k}\frac{(y_{l}^{k+1}-y_{l}^k)(u,\e)}{u-\e X}\,du-
\med\\
\qquad\qquad\qquad\ds
 \int_{\~\xi_l^{k-1}}^{\~\xi_{l}^k}\frac{(y_{l+1}^{k}-y_{l}^k)(u,\e)}{u-\e X}\,du
\Bigg)
\end{array}
$$
où $\xi_l^k=\~ r_0 e^{i(\psi^j+\gamma_l)}$
et $\~\xi_l^k=\norm{\e\~\mu} e^{i(\psi^j+\gamma_l)}$, $\f_l<\gamma_l<\f_{l+1}$ 
désignent des points dans l'intersection $\~ V_l^{k,k+1}(\e) 
\cap \~ V_{l+1}^{k,k+1}(\e)$. 
Ici, les points $\~x_l^j$ respectivement $\~x_{l+1}^j$
de la formule \rf{imod} pour $Y_l^j$ respectivement 
$Y_{l+1}^j$ ont été modifiés en $\~\xi_l^j$ avant la soustraction.
Seule la dernière intégrale ne figurait pas avant~; elle est exponentiellement petite
par l'hypothèse \rf{6}.
Dans l'étude des fonctions $g_n^j$, rien ne change.

\sub{4.5}{Composition d'un \dac{} Gevrey par une fonction analytique}
Le théorème \reff{t3.1} et sa réciproque la proposition \reff{p3.11} 
permettent de démontrer rapidement une généralisation
de la proposition \reff{l2.4} aux \dacs{} Gevrey~: la composée --- à gauche ou à droite --- d'une fonction ayant un  \dac{} Gevrey par une fonction analytique a un \dac{} Gevrey.
Auparavant, nous rappelons le résultat analogue concernant les développements Gevrey classiques. D'une part cela permet d'indiquer la preuve utilisant le théorème de Ramis-Sibuya classique, d'autre part ce résultat sera utilisé dans la suite.
\propo{p5.6}
{ On considère $y$ une fonction holomorphe définie sur $D\times S$, où $D$ est un domaine et où $S$ est le secteur $S(\a,\b,\e_0)$, à  valeurs dans un ensemble $W\subset\C$ et  
admettant un développement Gevrey d'ordre $\usp$ en $\e$ uniforme par rapport à $x\in D$.
\be[\rm(a)]\item
 Soit $f$ une fonction holomorphe sur un voisinage de l'adhérence de $W$. 
Alors la fonction $z=f\circ y$ admet un  développement  Gevrey d'ordre $\usp$  en $\e\in S$ uniforme par rapport à $x\in D$.
\item
Soit $\f:D(0,x_1)\times D(0,\e_0)\to\C$ une fonction holomorphe bornée telle que $\f(0,0)=0$ et  $\frac{\partial\f}{\partial x}(0,0)=1$. 
Soit $K$ un compact de $D$ et soit $\~D\subset\C$ tel que $\f(x,\e)\in K$
quand $\e\in S$ et $x\in \~D$.
Alors la fonction $z:(x,\e)\mapsto y(\f(x,\e),\e)$ admet un développement Gevrey d'ordre $\usp$ en $\e$ uniforme par rapport à $x\in\~D$. 
\ee
}
\pr{Preuve}
Soit $y(x,\e)\sim_\usp \sum_{n=0}^\infty y_n(x)\e^n$ le 
développement asymptotique de $y$. On choisit des secteurs $S_l$, $l=2,...,L$
d'angle d'ouverture strictement inférieur à $\frac\pi p$, tels qu'en ajoutant $S_1:=S$ 
on obtienne un bon recouvrement de $D(0,\e_0)^*$.
Comme indiqué dans la preuve de lemme \reff{brg1}, il existe des fonctions
$y_l$ sur $D\times S_l$, $l=2,...,L$, admettant 
$\sum_{n=0}^\infty y_n(x)\e^n$ comme développement asymptotique Gevrey d'ordre $\usp$
uniformément sur $D$. En posant $y_1=y$ on obtient que les différences admettent
la série $0$ comme développement Gevrey d'ordre $\usp$ ; elles sont donc
exponentiellement petites : il existe des constantes $C,A>0$ telles que
$\norm{y_l(x,\e)-y_{l-1}(x,\e)}\leq C \exp(-A/\norm\e^p)$ quand $\e\in S_{l-1}\cap S_l$
et $x\in D$. 
Sous l'hypothèse de l'item (a), quand $\norm\e $ est assez petite, on peut 
définir $z_l=f\circ y_l$ et l'on a 
$$\norm{z_l(x,\e)-z_{l-1}(x,\e)}\leq C M \exp(-A/\norm\e^p)$$
où $M$ est le maximum de $|f'|$ sur l'adhérence de $W$. 
Par le lemme \reff{RS}, ceci implique (a).
On trouve une preuve plus détaillée dans \cit{crss}.

Sous l'hypothèse de (b), on pose de manière analogue $z_l(x,\e)=y_l(\f(x,\e),\e)$
et  on vérifie qu'elles sont définies quand $x\in\~D$ et $\e\in S_l$ est 
assez petit. Alors on a évidemment
$\norm{z_l(x,\e)-z_{l-1}(x,\e)}\leq C \exp(-A/\norm\e^p)$ quand $\e\in S_{l-1}\cap S_l$
et $x\in \~D$ et on conclut comme avant.
\ep
\theo{t4.7}{
On considère $y$ une fonction holomorphe définie pour $\e\in S_2=S(\a_2,\b_2,\e_0)$ et $x\in V_1(\e)=V(\a_1,\b_1,r_0,\mu\norm\e)$, à  valeurs dans un ensemble $W\subset\C$ et  
admettant un \dac{} Gevrey d'ordre $\usp$.
\be[\rm(a)]\item
Soit $P(x,z,\e)$ une fonction holomorphe définie quand
 $ |z| < r$, $\e\in S_2=S(\a_2,\b_2,\e_0)$ et $x\in V_1(\e)$ 
admettant un \dac{} Gevrey d'ordre $\usp$~: 
$P(x,z,\e)\sim_{\usp}\hat P(x,z,\e)=\sum_{n\geq0}\gk{A_n(x,z)+G_n(\tfrac x\e,z)}\e^n$
quand $S_2\ni\e\to0$, $x\in V_1(\e)$, uniformément pour $\norm z <r$,
\ie les constantes des définitions \reff{d2.3.2} et \reff{d2.3.1} sont indépendantes
de $z$ avec $\norm z<r$. Supposons de plus que $y(x,\e)=\O(\e)$.

Alors la fonction $u:(x,\e)\mapsto P(x,y(x,\e),\e)$ admet un \dac{} Gevrey d'ordre $\usp$.
\item
  Soit $f$ une fonction holomorphe sur un voisinage de l'adhérence de $W$. 
Alors la fonction $z=f\circ y$ admet un \dac{} Gevrey d'ordre $\usp$ pour $\e\in S_2$ et $x\in V_1$.
\item
Soit $\f:D(0,x_1)\times D(0,\e_0)\to\C$ une fonction holomorphe bornée telle que $\f(0,0)=0$ et  $\frac{\partial\f}{\partial x}(0,0)=1$. 
Soit $0<\~r\leq x_1$, $\~\mu\in\R$ et 
$\~\a_1,\~\b_1$ vérifiant $\a_1<\~\a_1<\~\b_1<\b_1$ tels que $\f(x,\e)\in V_1(\e)$
quand $\e\in S(\a_2,\b_2,\e_0)$ et $x\in \~V_1(\e)=V(\~\a_1,\~\b_1,\~r,\~\mu\norm\e)$.
Alors la fonction $z:(x,\e)\mapsto y(\f(x,\e),\e)$ admet un \dac{} Gevrey d'ordre $\usp$ 
pour $\e\in S(\a_2,\b_2,\e_0)$ et $x\in \~V_1(\e)$.
\ee
}
\rqs
 1.\ 
 Nous utiliserons le résultat du (c) essentiellement dans deux cas~: le cas où $\f$ est indépendante de $\e$ et le cas du {\sl décalage} $\f(x,\e)=(x+\e T,\e)$ avec $T\in\C$ fixé.
\med\\
2.\ 
La preuve précédente est très indirecte et ne montre pas comment calculer en pratique le \dac{} de la composée, à  gauche ou à  droite, à  partir du  \dac{} $\yc(x,\e)=\ts\sum_{n\geq0}\gk{a_n(x)+g_n\big(\ts\frac x\e\big)}\e^n$ de $y$. Ceci a été fait dans la preuve de la proposition \reff{l2.4} pour les cas (a) et (b), et aussi pour le cas (c) lorsque $\f$ est indépendante de $\e$.
\med\\
3.\ 
En pratique, des valeurs de $\~r,\~\a_1,\~\b_1$  et $\~\mu$ peuvent  être déterminées en fonction de $\f$ dans le cas (c). Par exemple dans la preuve du corollaire \reff{c3.6} nous utiliserons le décalage $\f(x,\e)=x+R\e$ avec un certain $R>0$. Dans ce cas, 
on peut choisir 
$\~r=r-R\e_0$ et, pour tout $\d\in\,\big]0,\tfrac{\b_1-\a_1}2\big[$, $\~\a_1=\a_1+\d$, $\~\b_1=\b_1-\d$ et $\~\mu=\min\{\mu-R,-R/\sin\d\}$. 
Lorsque $\f$ ne dépend que de $x$, la fonction $\psi$ de la preuve ci-dessous est identiquement nulle et on peut prendre $\~\mu$ arbitrairement proche de $\mu$, si on réduit $\~r$.
\med\\
\pr{Preuve}
(a) À l'aide de la proposition \reff{p3.11}, on associe à  $y$ une famille $(y_l^j)$ sur un bon recouvrement cohérent et dont les différences satisfont \rft{expo-l}{expo-y}. 
De manière analogue, on associe à $P$ une famille $(P_l^j)$ sur le même recouvrement
 dont les différences satisfont \rft{expo-l}{expo-y} uniformément par rapport
à $z$.  Pour ceci, il suffit de remplacer les 
coefficients scalaires et les valeurs scalaires des fonctions par des coefficients 
resp.\ valeurs dans l'espace de Banach des fonctions holomorphes bornées dans $D(0,r)$.

On pose $u_l^j(x,\e)=P_l^j\gk{x, y_l^j(x,\e),\e}$. On vérifie que les 
majorations  \rft{expo-l}{expo-y} sont satisfaites pour $u$ et les $u_l^j$, au besoin 
en augmentant la constante $C$. D'après le théorème \reff{t3.1}, les fonctions $u_l^j$ 
admettent un \dac{} commun qui est Gevrey d'ordre $\usp$.
D'après la partie ``réciproque'' de la proposition
\reff{p2.3.1} (a) appliquée à  $u-u_l^j$ sur chacune des intersections, la fonction $u$ admet le même \dac{} Gevrey d'ordre $\usp$. En comparant avec la proposition \reff{l2.4},
on obtient la série formelle du \dac.
\med\\
La démonstration de (b) est complètement analogue.
\med\\
(c) 
Soit $\psi(x,\e)=\De_2\f(x,0,\e)$, où $\De_2\f$ est la différence finie de $\f$ par rapport à la deuxième variable, définie par 
\eq D
{\f(x,\tau)-\f(x,\e)=\De_2\f(x,\e,\tau)(\tau-\e).
}
 Ainsi on a 
$
\f(x,\e)=\f(x,0)+\e\psi(x,\e).
$ 
Quitte à diminuer $\e_0$, 
 la fonction $\psi$ est bornée par une certaine constante $K$. 
Soit d'abord $\bar r>0$, $\bar \e_0>0$ et $\delta>0$ assez petits  et soit $\bar \mu<\mu-K$. 
De même que pour (a), on associe à  $y$ une famille  de fonctions $y_l^j$ définies sur
$V_l^j(\e)=V(\a_l^j,\b_l^j,r_0,\mu|\e|)$, satisfaisant \rft{expo-l}{expo-y}, 
et on pose $z_l^j(x,\e)=y_l^j(\f(x,\e),\e)$ pour $\e\in S(\a_2,\b_2,\bar \e_0)$ et
$x\in\bar V_l^j(\e)=V(\a_l^j+\delta,\b_l^j-\delta,\bar r,\bar \mu|\e|)$. 
D'après les conditions sur $\f$, les fonctions $z_l^j$ sont bien définies et 
 on vérifie 
que les ensembles $S_l$, $V^j$ et $\~V_l^j(\e)$ forment un bon recouvrement cohérent 
 si  $\bar r>0$ et $\bar \e_0>0$ sont suffisamment petits.
Les $z_l^j$ satisfont \rft{expo-l}{expo-y} avec d'autres constantes $\~A,\~B,\~C$,
donc ont un \dac{} commun Gevrey d'ordre $\usp$ d'après le théorème \reff{t3.1}. La fonction $z$ admet alors le même développement pour $\e\in S(\a_2,\b_2,\bar \e_0)$ et
$x\in V(\a_1+\delta,\b_1-\delta,\bar r,\bar \mu|\e|)$, comme en (a).
 À présent, il s'agit de montrer que ce \dac{} est encore valide pour $\e\in S(\a_2,\b_2,\e_0)$ et $u\in\~V_1(\e)$.

La fonction $\Phi:(X,\e)\mapsto\f(\e X,\e)$ satisfait $\Phi(X,\e)=\e X+\e \psi(0,0)+
{\cal O}(\e^2)$ donc, d'après la proposition \reff{p5.6} (b), la fonction $Z:(X,\e)\mapsto z(\e X,\e)$ se prolonge et a un développement Gevrey uniforme pour $\e\in S(\a_2,\b_2,\bar \e_0)$ et 
$x\in V(\~\a_1,\~\b_1,\bar r,\~\mu|\e|)$.
Avec les propositions \reff{matching-gevrey} et \reff{p4.5}, on peut donc prolonger le \dac{} Gevrey
de $z$ pour $\e\in S(\a_2,\b_2,\bar \e_0)$ et 
$u\in V(\~\a_1,\~\b_1,\bar r,\~\mu|\e|)$. De même, avec les propositions \reff{p5.6} (a),
\reff{matching-gevrey} et \reff{p4.5bis}, on prolonge ce \dac{} Gevrey pour 
$\e\in S(\a_2,\b_2,\bar \e_0)$ et $u\in V(\~\a_1,\~\b_1,\~ r,\~\mu|\e|)$.
Enfin, quitte à changer 
les constantes $C, L_1$ et $L_2$ 
 dans \rf{defcombgevrey}, on obtient l'énoncé.
\ep\med
%
%
%
%
%
\sec{5.}{Développements 
combinés et équations \\ différentielles
 singulièrement perturbées}
Notre notion de développement combiné a été motivée par l'étude d'équations différentielles singulièrement perturbées de la forme \rf{fb}, que nous réécrivons ci-dessous par commodité
\eq m{
\eps y'=\Phi(x,y,\eps).
}
Contrairement aux équations de la partie \reff{1.2},  les variables $x$ et $y$ sont à présent dans $\C$ et le petit paramètre $\eps$ est dans un secteur du plan complexe $S=S(-\d,\d,\eps_0)$ contenant une partie de l'axe réel positif, avec $\d,\eps_0>0$ fixés suffisamment petits. \Apriori{} la situation qui nous intéresse est $\eps$ réel positif, mais nous aurons besoin de le considérer dans un secteur.

 On suppose que $\Phi$ est analytique  par rapport à  $x$ et $z$ dans un domaine $\D\subset\C^2$ et Gevrey d'ordre 1 par rapport à  $\eps$  dans $S$~; ceci permet, entre autres, de traiter des équations comportant aussi un paramètre de contrôle et pour lesquelles apparaissent des canards.

On rappelle que l'{\sl ensemble lent} $\L$ est l'ensemble d'équation $\Phi(x,y,0)=0$ et qu'un point $(\~x,\~y)$ de $\L$ est dit  {\sl régulier} si $\frac{\partial\Phi}{\partial y}(\~x,\~y,0)\neq0$.
Au voisinage d'un point régulier, d'après le théorème des fonctions implicites,  $\L$ est le graphe d'une fonction continue $y_0$, dite {\sl lente}, vérifiant $y_0(\~x)=\~y$. 
On note $f(x)=\frac{\partial\Phi}{\partial z}(x,z_0(x),0)$.
Lorsqu'on cherche à  prolonger une telle fonction $y_0$, deux situations peuvent se produire~: ou bien on arrive au bord de $D$, ou bien on arrive à  un point tournant, où $f$ s'annule. \Apriori{} la fonction $f$ a une singularité en un point tournant. Dans ce mémoire, nous nous plaçons dans la siuation très particulière où $f$ est définie et analytique dans un voisinage du point tournant.

Dans la partie \reff{5.1}, nous détaillons l'usage des \dacs{} au voisinage d'un point régulier, où l'on retrouve les développements combinés classiques.
Le cas d'un point tournant est présenté dans la suite, d'abord dans une situation simple, dite  {\sl quasi-linéaire} dans la partie \reff{5.2}, puis dans une situation plus délicate dans la partie \reff{5.3}.
\sub{5.1}{\Dacs\ classiques en un point régulier}
On considère un problème aux valeurs initiales de la forme 
\eq p{
\eps y'=\Phi(x,y,\eps)\,,\qquad y(\~x,\eps)=v(\eps)
}
avec les notations et hypothèses précédentes. On suppose de plus que $v$ est Gevrey d'ordre 1 par rapport à  $\eps\in S$.
On vérifie aisément qu'il existe une unique solution formelle $\yc(x,\eps)=\sum_{n\geq0}a_n(x)\eps^n$ de l'équation \rf m sans condition initiale
avec des coefficients $a_n$ analytiques au voisinage de $\~x$ et avec $a_0=y_0$, la fonction lente, \cf la formule  \rf f ci-dessous.
Il est aussi connu \cit{bfsw,crss} que cette solution formelle satisfait des estimations Gevrey d'ordre 1 en $\eps$ et que, de plus, une solution $y=y(x,\eps)$ de \rf m bornée dans un domaine de la forme $D_1\times S$ avec $D_1=D(\~x,r_0)$ (et $S=S(-\d,\d,\eps_0)$) est asymptotique  Gevrey-1 à $\yc$ dans tout sous-domaine $\~D_1\times\~S$ (\ie $\~D_1=D(\~x,\~r_0)$ et $\~S=S(-\~\d,\~\d,\~\eps_0)$ avec $\~r_0<r_0$, $\~\d<\d$ et $\~\eps_0<\eps_0$).

Par ailleurs, il est également connu \cit{vb,bef} que la solution du problème 
aux valeurs initiales
 \rf p admet un développement combiné de la forme $\sum_{n\geq0}\Big( a_n(x)+g_n\big(\ts\frac {x-\~x}\eps\big)\Big)\eps^n$, où les $a_n$ sont les coefficients de la solution formelle et où les $g_n$ sont des fonctions à  décroissance exponentielle.
Le théorème \reff{t3.6} et le corollaire \reff{c3.6} permettent de retrouver ce résultat et de le généraliser au cadre complexe. Ils donnent en outre des estimations Gevrey de ces développements.
Précisément, nous allons montrer que, si la condition initiale $v$ est suffisamment proche de  $\~y=y_0(\~x)$, alors  la solution du problème \rf p admet un \dac{} Gevrey pour $\eps\in S$ et pour $x$ dans un domaine approprié. Le domaine en $x$ pourra être choisi d'autant plus grand que $v$ sera proche de  $\~y$.
Pour éviter d'avoir un énoncé trop compliqué, et aussi pour que les preuves ne soient pas trop longues, nous avons choisi d'écrire deux énoncés distincts. 

Avant de présenter ces énoncés, nous faisons les réductions suivantes.
Tout d'abord, en changeant $y$ en $y-y_0$ on se ramène à  $y_0\equiv0$. 
Enfin, une translation et une rotation de la variable $x$ permettent de supposer $\~x=0$ et $f(0)<0$. Ceci nous permet de détailler un peu la solution formelle.
Puisque $\Phi(x,0,0)=0$, la fonction $\Phi$ est de la forme $\Phi(x,y,\eps)=y\Psi(x,y)+\eps P(x,y,\eps)$, avec $\Psi$ et $P$ analytiques (et $\Psi(x,0)=f(x)$). Le coefficient du terme d'ordre $n+1$ en $\eps$ dans le développement de Taylor de $\Phi\big(x,\sum_{k\geq1}a_k(x)\eps^k,\eps\big)$
 est donc de la forme $f.a_{n+1}+\phi_n$, où $\phi_n$ ne dépend que des termes $a_0,\ldots,a_n$. Ainsi, la solution formelle est donnée récursivement par 
\eq f{
a_0=0,\qquad a_{n+1}=\frac1f(a'_n-\phi_n).
}
Le premier énoncé concerne des conditions initiales  proches de la courbe lente. Par continuité de la dérivée partielle de $\Phi$, la condition \rf c sera réalisée si $r_0,r_2$ et $\eps_0$ sont assez petits.
\theo{t3.6}{
Soit $\d,\eps_0,r_0,r_2>0$.
On suppose que $\d<\frac\pi4$ et que pour tout $|x|<r_0$, $|y|<r_2$ et $\eps\in S(-\d,\d,\eps_0)$
\eq c{
\Big|\arg\ts\frac{\partial\Phi}{\partial y}(x,y,\eps)-\pi\Big|<\d~\mbox{ et }~\big|\frac{\partial\Phi}{\partial y}(x,y,\eps)\big|\geq\frac12\,|f(0)|.
}
Si $|v|<r_2$ sur $S$, alors la solution $y$ de \rf p admet un développement combiné Gevrey d'ordre 1 en $\eps$, pour $x$ dans le secteur
 $S'=S\big(-\frac\pi2+2\d,\frac\pi2-2\d,r_0\big)$ 
et pour $\eps$ dans $S$. La partie lente de ce développement est la solution formelle \rf f et la partie rapide est constituée de fonctions $g_n$  à  décroissance exponentielle dans le secteur $S''=S\big(-\frac\pi2+\d,\frac\pi2-\d,\infty\big)$.

Précisément, il existe une suite $(g_n)_{n\in\N}, g_n\in\G(S'')$ et des constantes $C_1,C_2>0$ telles que,
 pour tout $x\in D(0,r)\cap S'$, tout $\eps\in S$ et tout entier $N>0$
\eq{8b}{
\Big|y(x,\e)-\sum_{n=0}^{N-1}\Big( a_n(x)+g_n\big(\ts\frac x\eps\big)\Big)\eps^n\Big|
\leq C_1\,C_2^N\,N!\,|\eps|^N.
}
De plus, il existe $C_3,C_4,C_5>0$ telles que, pour tout  $X\in S''$ et tout $n\in\N$, on ait
\begin{eqnarray}
\lb{e6.6}
|g_n(X)|&\leq& C_3\,C_4^n\,
\big(C_5|X|\,\big)^{n+1/2}
\,e^{-C_5\norm X}~\mbox{si}~C_5|X|\geq n\,,\\
|g_n(X)|&\leq &C_3\,C_4^n\,n!~\mbox{sinon}.\nonumber
\end{eqnarray}
}
\rqs 
1.\ L'existence de ces développements combinés en un point régulier est classique 
dans le cadre réel. Les méthodes réelles peuvent s'adapter aisément au cadre complexe.
En revanche, les estimations Gevrey sont nouvelles à  notre connaissance. Nous pensons que ces estimations peuvent être utiles pour des questions conduisant à  utiliser des développements combinés. 
\med\\
2.\ 
Les coefficients $a_n$ sont ceux de la solution formelle, déterminés par \rf{f}~; 
les coefficients
$g_n$ peuvent être obtenus directement à  partir de l'équation {\sl intérieure}, obtenue en posant $x=\eps X,\,y(x)=Y(X)$. 
Rappelons les notations $\Phi(x,y,\eps)=y\Psi(x,y)+\eps P(x,y,\eps)$, avec $\Psi$ et $P$ analytiques et $\Psi(x,0)=f(x)$.
Introduisons aussi $v_n$, donné par $v(\eps)=\sum_{n\geq0}v_n\eps^n$, et $\psi:Y\mapsto\Psi(0,Y)$.
L'équation intérieure est
$$
\frac{dY}{dX}=Y\Psi(\eps X,Y)+\eps P(\eps X,Y,\eps),\qquad Y(0,\eps)=v(\eps).
$$
La solution formelle de ce problème se calcule en résolvant successivement des équations différentielles avec condition initiale~; ces équations sont linéaires non homogènes sauf la première. Pour le premier terme, on obtient
$$
\frac{dY_0}{dX}=Y_0\psi(Y_0),\qquad Y_0(0)=v_0,
$$
ce qui donne implicitement $Y_0$ par  
$\ds\int_{v_0}^{Y_0(X)}\frac{du}{u\psi(u)}=X$, et pour les autres termes
$$
Y'_n=A(X)Y_n+B_n(X)\,,\qquad Y_n(0)=v_n
$$
avec $A=Y_0\psi'(Y_0)+\psi(Y_0)$ et où $B_n$ ne dépend que des termes $Y_1,\ldots,Y_{n-1}$.
On obtient 
$$
Y_n(X)=v_n\,e^{\int_0^XA}+\int_0^Xe^{\int_u^XA}B_n(u)du.
$$
\med

\pr{Preuve}
Soit $N\in\N$ avec $N\geq\frac{2\pi}\d$, et pour $l=0,\ldots,N-1$, soit
$$
S_l=S(\a_l,\b_l,\e_0)~~\mbox{avec}~~\a_l=\frac{2(l-1)\pi}N~\mbox{et}~\b_l=\frac{2(l+1)\pi}N.
$$
De cette manière, $N$ est assez grand pour que $S_0\subset S$.
Soit $v_l:S_l\to\C$ une fonction asymptotique Gevrey-1 à  $\vc=\sum_{n\geq0}v_n\eps^n$, dont l'existence est assurée par le théorème de Borel-Ritt-Gevrey (obtenue par exemple par transformation de Borel de $\vc$ et par transformation de Laplace tronquée, \cf le début de la preuve du lemme \reff{brg1}).
Au besoin en diminuant $\eps_0$, on peut supposer que ces fonctions $v_l$ sont bornées par $r_2$.  Pour fixer les idées on peut choisir pour $v_0$ la restriction de $v$ à  $S_0$, mais ce n'est pas indispensable.
De même, soit $\Phi_l:D(0,r_0)\times D(0,r_2)\times S_l\to\C$ asymptotiques Gevrey-1 en $\eps$ à  la même série que $\Phi$.

Puisque les fonctions $v_0,\ldots,v_{N-1}$ et $v$ sont asymptotiques Gevrey-1 à  une même série, et de même pour $\Phi_0,\ldots,\Phi_{N-1}$ et $\Phi$, il existe $C,A>0$ tels que pour tout $l\in\{0,\ldots,N-1\}$, tout $x\in D(0,r_0)$, tout $y\in D(0,r_2)$ et tout $\eps\in S_{l,l+1}=S_l\cap S_{l+1}$ on ait
\eq w{
|v_{l+1}(\eps)-v_l(\eps)|\leq Ce^{-A/|\eps|}~~\mbox{et}~~
|\Phi_{l+1}(x,y,\eps)-\Phi_l(x,y,\eps)|\leq Ce^{-A/|\eps|}
}
ainsi que 
$$
|v_{l}(\eps)-v(\eps)|\leq Ce^{-A/|\eps|}~~\mbox{et}~~
|\Phi_{l}(x,y,\eps)-\Phi(x,y,\eps)|\leq Ce^{-A/|\eps|}
$$
pour $x\in D(0,r_0)$, $y\in D(0,r_2)$ et $\eps\in S\cap S_{l}$ si cette intersection n'est pas vide.
En vue d'appliquer le théorème \reff{t3.1}, nous allons considérer deux familles $(y_l^1)_{0\leq l<N}$ et $(y_l^2)_{0\leq l<N}$ de solutions d'équations similaires à   \rf m, obtenues en remplaçant $\Phi$ par $\Phi_l$. Ensuite nous vérifierons que les différences sont exponentiellement petites. 

Concernant la première famille, pour chaque $l\in\{0,\ldots,N-1\}$, soit $y_l^1$ la solution de condition initiale $y_l^1(0,\eps)=v_l(\eps)$.
Alors, pour $\eps_0$ suffisamment petit, $y_l^1$ est définie et bornée par $r$ sur $S_l^1\times S_l$, où $S_l^1=S(\a_l^1,\b_l^1,r_0)$, avec
$$
\a_l^1=\frac\pi N\Big(2l-\frac N2+4\Big)~~\mbox{et}~~
\b_l^1=\frac\pi N\Big(2l+\frac N2-4\Big).
$$
\ft
\lem{l3.10}{
On a, pour tout $x$ dans l'intersection $S_{l,\,l+1}^1=S_l^1\cap S_{l+1}^1$ et tout $\eps$ dans $S_{l,\,l+1}=S_l\cap S_{l+1}$
$$
|y_{l+1}^1(x,\eps)-y_l^1(x,\eps)|\leq C(1+r_0)e^{-A/|\eps|},
$$
où $A$ et $C$ sont donnés dans \rf w.

De même, on a, pour tout $x$ dans $S_l^1\cap S'$ et tout $\eps$ dans $S_l\cap S$
$$
|y_{l}^1(x,\eps)-y(x,\eps)|\leq C(1+r_0)e^{-A/|\eps|}.
$$
}
\pr{Preuve}
Posons $z_l=y_{l+1}^1-y_l^1$, notons
$$
b(x,\eps)=\Phi_{l+1}(x,y_l(x,\eps),\eps)-\Phi_{l}(x,y_l(x,\eps),\eps),
$$
rappelons la notation $\De_2$ introduite dans \rf D, définie par 
$$
\Phi_{l+1}(x,z,\eps)-\Phi_{l+1}(x,y,\eps)=\De_2\Phi_{l+1}(x,y,z,\eps)(z-y)
$$
et notons $d(x,\eps)=\De_2\Phi_{l+1}(x,y_l(x,\eps),y_{l+1}(x,\eps),\eps)$.
Puisque 
$$
\De_2\Phi(x,y,z,\eps)=\int_0^1\frac{\partial\Phi}{\partial y}(x,y+t(z-y),\eps)dt,
$$
l'hypothèse \rf c entraîne 
$\big|\arg d(x,\eps)-\pi\big|<\d$.

Alors $z_l$ est solution du problème aux valeurs initiales
$$
\eps z'_l= d(x,\eps)z_l+\eps b(x,\eps),\qquad z_l(0,\eps)=v_{l+1}(\eps)-v_l(\eps).
$$
Étant donné $x\in S_{l,\,l+1}^1$, choisissons pour chemin d'intégration le segment de $0$ à  $x$. La formule de variation de la constante donne,  pour tout $\xi$ sur ce chemin tel que $z_l$ est définie sur $[0,\xi]$
$$
z_l(\xi,\eps)=z_l(0,\eps)\exp\Big(\ts\frac1\eps\int_0^\xi d(u,\eps)du\Big)+\ds
\int_0^\xi b(u)\exp\Big(\ts\frac1\eps\int_u^\xi d(v,\eps)\Big)du.
$$
Pour $\e\in S_{l,\,l+1}$ et  $u\in[0,x]$, donc {\sl a fortiori} pour $u\in[0,\xi]$, on a
$$
\re\big(\ts\frac1\eps d(u,\eps)e^{i\arg u}\big)<0,
$$
 ce qui implique que 
$$
|z_l(\xi,\eps)|\leq |z_l(0,\eps)|+\int_0^{|\xi|}|b(se^{i\arg x})|ds
\leq  C(1+r_0)e^{-A/|\eps|}.
$$
Le principe de majoration \apriori{} entraîne alors
 que $z_l$ est définie sur tout $[0,x]$ et satisfait la même majoration.
 La comparaison entre $y$ et $y_l^1$ est similaire.
\ep

\pr{Suite de la preuve du théorème \reff{t3.6}}
À présent, construisons la deuxième famille. Fixons $\~r_0\in\,]0,r_0[$. Pour chaque $l\in\{0,\ldots,N-1\}$, soit $x_l=-\~r_0e^{2l\pi i/N}$ et soit $y_l^2$ la solution de condition initiale $y_l^2(x_l,\eps)=0$.

Pour $\eps_0$ suffisamment petit, $y_l^2$ est définie sur $D_l\times S_l$, où $D_l$ est l'intersection du disque $D(0,r_0)$ et du secteur de sommet $x_l$ et d'angles $\a_l^1$ et $\b_l^1$,
\cf figure \reff{f6.2}.
\figu{f6.2}{
\vspace{-3mm}
\epsfxsize6cm\epsfbox{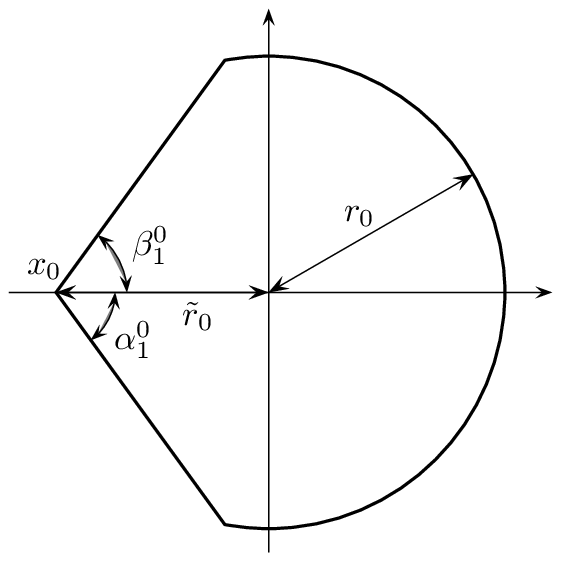}
\vspace{-8mm}
}{En gras, le bord de $D_l$}
Par un argument similaire à  la preuve du lemme \reff{l3.10}, on montre que pour $x\in D_l\cap D_{l+1}$ et $\eps\in S_{l,l+1}$, on a de même
$$
|y_{l+1}^2(x,\eps)-y_l^2(x,\eps)|\leq C(1+2r_0)e^{-A/|\eps|}.
$$
Pour pouvoir appliquer le théorème \reff{t3.1}, il nous reste à  montrer que les solutions $y_l^1$ et $y_l^2$ sont exponentiellement proches dans l'intersection de leur domaine d'existence $D_l\cap S_{l}^1=S_l^1$. Leur différence $w=y_l^2-y_l^1$ satisfait l'équation 
$$
\eps w'=D(x,\eps)w
$$
où $D(x,\eps)=\De_2\Phi(x,y_l^1(x,\eps),y_l^2(x,\eps),\eps)$,
avec la condition initiale $w(0,\eps)=y_l^2(0,\eps)$ $-v_l(\eps)$. On a $|w(0,\eps)|\leq2r_2$~; 
d'après \rf c, on a aussi $|\arg D(x,\eps)-\pi|<\d$ et $|D(x,\eps|\geq\frac12|f(0)|$.
Pour $x\in S_l^1$ tel que $\big|\arg x-\frac{2l\pi}N\big|<\frac\pi3-\d$ et $\eps\in S_{l}$, on obtient $|w(x,\eps)|\leq2r_2e^{-B|x|/|\eps|}$ avec $B=\frac14|f(0)|$.
À présent, les conditions sont réunies pour appliquer le théorème \reff{t3.1}. Chacune des fonctions $y_l^j,\,j=1$ ou $2$, $0\leq l<N$, a donc un développement combiné Gevrey-1.
Puisque la solution $y$ est exponentiellement proche de $y_l^1$ sur $(S'\cap S_l^1)\times (S\cap S_l)$, elle aussi a un développement Gevrey.

De plus, puisque la solution formelle ne contient pas de pôles, c'est qu'elle coïncide avec le développement lent. Autrement dit, tous les coefficients $g_{mn}$ sont nuls. Les fonctions du développement rapide $g_n$ sont donc asymptotiques Gevrey-1 à  la fonction nulle lorsque la variable $X$ tend vers l'infini, donc à  décroissance exponentielle. 

Puisque les fonctions $g_n$ admettent des développements asymptotiques Gevrey 
compatibles au sens de la définition \reff{d2.3.1}, la formule \rf{gnm} donne  \rf{e6.6} avec un choix
adéquat de l'entier $M$ en fonction de $n$ et de $\norm X$.
Précisément, si $C,L_1,L_2>0$ sont des constantes (de \reff{d2.3.1}) telles que pour tout $n,M$ et $X$ (ici $p=1$) 
$$
\norm X^M\big|g_n(X)\big|\leq C L_1^n L_2^M\Gamma(M+n+1),
$$
on obtient \rf{e6.6} en utilisant la majoration $n!\leq n^{n+1/2}e^{1-n}$ et en choisissant $C_3=e\,C,\,C_4=L_1,\,C_5=1/L_2$ et pour l'entier $M$~: $M=0$ si $n\leq \norm X/L_2$ et $M=\big[\norm X/L_2\big]-n$ sinon.
\ep

Le théorème \reff{t3.6} concernait une condition initiale proche de la courbe lente, mais nos résultats de prolongement des \dacs{} permettent de traiter aussi le cas d'une condition initiale loin de la courbe lente. 
\coro{c3.6}{
Avec les notations du théorème \reff{t3.6}, on suppose que $v_0$ est dans le bassin d'attraction de $0$ pour l'équation rapide 
\eq {rr}{
Y'=\Phi(0,Y,0).
}
Alors il existe $\d,\eps_1>0$ tels que la solution $y$ de \rf p admet un développement combiné Gevrey d'ordre 1 en $\eps$, pour $x$ dans le secteur $S(-\d,\d,r_0)$ et pour $\eps$ dans $S=S(-\d,\d,\eps_1)$.

Plus généralement, si $\a<0<\b$ sont tels que $v_0$ est
 dans le bassin d'attraction de $0$ pour l'équation 
$e^{-id}Y'=\Phi(0,Y,0)$ pour tout $d\in[\a,\b]$, alors la solution $y$ de \rf p admet un développement combiné Gevrey d'ordre 1 en $\eps$, pour $x$ dans le secteur $S(\a,\b,r_0)$  et pour $\eps$ dans $S$.
}
\figu{f6.3}{
\vspace{-6mm}
\epsfxsize6cm\epsfbox{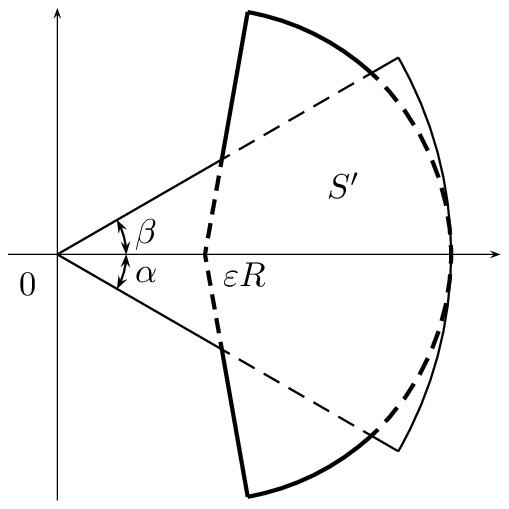}
\vspace{-8mm}
}
{Le secteur $S(\a,\b,r_0)$ et le translaté $S'$ du secteur
$S\big(-\frac\pi2+\d,\frac\pi2-\d,r_0-\eps R\big)$}
\pr{Preuve}
Tout d'abord, si $v_0$ est dans le bassin d'attraction de $0$ pour l'équation rapide \rf{rr} et si $\d$ est suffisamment petit, alors $v_0$ est aussi dans le bassin d'attraction pour l'équation $e^{-id} Y'=\Phi(0,Y,0)$ pour tout $|d|<\d$. Il nous suffit donc de montrer la généralisation. Étant donnés $r_0,r_2$ tels que les hypothèses du théorème \reff{t3.6} soient satisfaites, puisque $v_0$ est dans le bassin d'attraction de $0$, il existe $R>0$ tel que $|w(\eps)|<r_2$ pour tout $\eps\in S$, avec $w(\eps)=y(\eps R,\eps)$. 
Montrons d'abord que $w$ a un développement Gevrey-1. 

Soit $\f_{R}=\f_{R}(v,\eps)$ l'application flot au temps $R$ de l'équation {\sl intérieure}, obtenue en posant $x=\eps X$
\eq M{
\frac{dY}{dX}=\Phi(\eps X, Y,\eps).
}
La fonction $\f_R$ est définie par $\f_{R}(v,\eps)=Y_v(R,\eps)$ si $Y_v$ est la solution de \rf M de condition initiale $Y_v(0,\eps)=v$. D'après le théorème de dépendance des solutions par rapport aux paramètres et aux conditions initiales, $\f_{R}$ est analytique par rapport à  $v$ et $\eps$. Par hypothèse,
 $y(0,\eps)=v(\eps)$ a un développement Gevrey-1. 
Puisque la composée d'une fonction Gevrey-1 et d'une fonction analytique est Gevrey-1 (\cf le début de la section \reff{4.5}), la fonction $w:\eps\mapsto\f_{R}(v(\eps),\eps)$ est donc Gevrey-1.
À présent nous appliquons le théorème \reff{t3.6} au problème
\eq q{
\eps \frac{dz}{du}=\~\Phi(u,z,\eps)\,,\qquad z(0,\eps)=w(\eps)
}
où l'on a posé $\~\Phi(u,z,\eps)=\Phi(u+\eps R,z,\eps)$. La solution $z$, qui est reliée à  $y$ par $z(u,\eps)=y(u+\eps R,\eps)$, a donc un \dac{} 
Gevrey-1, $\sum_{n\geq1}\Big( b_n(u)+h_n\big(\frac u\eps\big)\Big)\eps^n$, pour $u$ dans $S'=S\big(-\frac\pi2+2\d,\frac\pi2-2\d,r_0\big)$ et $\eps$ dans $S=S(-\d,\d,\eps_0)$.
D'après le théorème \reff{t4.7} (\cf aussi la remarque 3 qui suit ce théorème) la solution $y$ de \rf p a aussi un \dac{} pour  $\eps\in S$ et  $x\in V= 
V\big(-\frac\pi2+3\d,\frac\pi2-3\d,\mu|\eps|,r_0-R\eps_0\big)$ avec $\mu=-R/\sin\d$,
si $\eps_0$ est assez petit.
Par ailleurs, d'après les hypothèses, la fonction  $y$ est définie pour $x$ dans le secteur $S(\a,\b,r_0)$.
En utilisant à  nouveau l'équation intérieure \rf M, l'analyticité du 
flot, et le fait que l'image d'une fonction Gevrey par une fonction analytique est Gevrey, on en déduit que $y$ est Gevrey-1 sur  $S(\a,\b,r_0)\setminus V$.
D'après la proposition \reff{p4.5}, $y$ a donc un \dac{} sur  $S(\a,\b,r_0)$.
\ep
\sub{5.2}{\Dacs\ en un point tournant~: le cas quasi-linéaire}
Dans cette partie, $r_0,r_2,\eps_0,\d>0$ sont fixés avec $\d>0$ 
suffisamment petit, $D_1$ et $D_2$ sont les disques 
de centre $0$ et de rayon respectivement $r_0$ et $r_2$, et $\Sigma$ désigne le secteur 
$\Sigma=S(-\d,\d,\eps_0)$.
On considère une équation  de la forme 
\eq1{
\eps y'=px^{p-1}y+\eps P(x,y,\eps),
}
$x,y\in\C$, $\eps$ petit paramètre et $P$ analytique bornée sur $D_1\times 
D_2\times \Sigma$ et Gevrey d'ordre 1 en $\eps$ dans $\Sigma$. Soit $\e=\eps^{1/p}$.
Dans la suite, on utilisera parfois $\eps$ et $\e$ en même temps ; il est sous-entendu qu'on a toujours
la relation $\eps=\e^p$.

L'équation \rf{1} est appelée {\sl quasi-linéaire}, car l'équation intérieure obtenue par
$x=\e X$, $Y(X)=y(\e X)$, \ie $\frac{dY}{dX}=pX^{p-1}Y+\e P(\e X,Y, \e^p)$ se réduit à 
une équation linéaire homogène quand on remplace $\e$ par $0$.

La forme générale de l'équation quasi-linéaire est obtenue en remplaçant $px^{p-1}$ par 
une fonction
$f(x)$ holomorphe dans un voisinage de $0$ et ayant un zéro d'ordre $p-1$ en ce point. Cette forme peut
être réduite à la {\sl forme normale} \rf{1} par un changement de variable 
(voir la remarque 8 dans la partie \reff{6.2.4}). 
Alternativement, on peut
modifier les énoncés et les preuves de cette partie pour les adapter à  cette forme générale.
La description des domaines par rapport à  $x$ étant plus simple pour \rf{1}, nous avons choisi de présenter la théorie pour cette équation dans la suite.

Nous avons choisi une équation simple dépendant de $\eps$, bien qu'on puisse traiter des 
équations plus générales (voir remarque 2 dans \reff{6.2.4} et une généralisation dans la partie
\reff{5.3}). Pour cette équation, deux preuves d'existence de solutions sont possibles. 
 La première, décrite dans les parties \reff{5.2.1} et \reff{4.2a}, procède selon un modèle classique : on détermine d'abord une solution formelle combinée et on montre ensuite l'existence de solutions l'ayant comme \dac. 
La deuxième preuve utilise directement le théorème-clé \reff{t3.1} et montre {\sl l'existence} d'un \dac{} Gevrey sans  donner ses coefficients. 
Bien que la deuxième preuve soit plus courte, 
nous avons choisi de présenter la première car elle permet aussi de traiter des équations qui n'entrent pas dans le cadre de la deuxième preuve 
(par exemple des équations réelles C$^\infty$ seulement) 
et pour lesquelles les \dac{} obtenus ne sont pas nécessairement Gevrey.
\ssub{5.2.1}{Solutions formelles combinées}
Soit $V=V(\a,\b,\infty,\mu)$, où $\mu$ est un nombre réel positif ou négatif, 
 et où $\a$ et $\b$ satisfont $-\frac{3\pi}{2p}<\a<0<\b<\frac{3\pi}{2p}$. 
On rappelle que $\cch(r_0,V)$ désigne l'ensemble des séries formelles combinées associées à  $V$ et au disque $D_1$. Nous introduisons aussi $\~C(V)$, l'ensemble des séries formelles combinées associées à  $V$ et au disque $D_1$ où les coefficients $a_n$ ne sont pas nécessairement bornés sur le disque. 
La motivation pour cette modification réside dans le fait qu'on ne sort pas du cadre en dérivant la partie lente d'un \dac{}~; 
nous n'aurons pas besoin de dériver de partie rapide.
L'espace vectoriel $\~C(V)$ muni de la distance ultramétrique (\cf \rf d) est un espace de Banach. De même on note $\~\H$  l'espace vectoriel des  fonctions holomorphes dans  $D_1$, non nécessairement bornées.
 Voici le résultat principal.
\theo{th4.1}{
Avec les notations précédentes, l'équation \rf1 a une unique solution formelle combinée $\yc$ dans $\e\~C(V)$.
}
Nous verrons dans la partie \reff{4.2} que l'existence de cette solution formelle sera une consé\-quence directe de l'existence de solutions analytiques de \rf1 ayant un développement combiné. 
Cependant il nous a semblé instructif de présenter une preuve indépendante d'existence de cette solution formelle et d'une solution l'ayant comme \dac, qui évite le recours aux techniques développées dans la partie \reff{3.}.

Le principe de la preuve est le suivant. On commence par résoudre l'équation linéaire non homogène

\eq3{
\eps y'=px^{p-1}y+\eps h(x,\e)
}
où $h\in\~C(V)$. On démontre dans le lemme \reff{lm4.4} qu'il existe une unique solution de \rf3 dans $\e\~C(V)$, notée $\Phic(h)$. En notant $\P$ l'opérateur qui, à  une série formelle $\yc$, associe la série formelle obtenue en substituant $\yc$ dans $P$, \ie 
$\P(\yc)(x,\e) = P(x,\yc(x,\e),\e^p)$, on vérifie ensuite  que l'équation 
$$
y=\Phic(\P(y))
$$
satisfait toutes les conditions d'application du théorème du point fixe de Banach dans $\~C(V)$.
\medskip

Par linéarité, l'équation \rf3 se résout en développant $h$ en série formelle combinée et en additionnant les solutions de chacune des équations ne contenant qu'un terme de la série. Il s'agit donc de résoudre \rf3 dans deux cas~: lorsque la fonction $h$ ne dépend que de $x$, \ie $h\in\~\H$, et lorsque $h$ ne dépend que de $X=\frac x\e$, \ie $h(x,\e)=k\big(\frac x\e\big)$ avec $k\in\G(V)$. Nous commençons par le deuxième cas.
\lem{lm4.2}{
Pour tout $k\in\G(V)$, l'équation
\eq{eq2.6}{
\eps y'=px^{p-1}y+\eps k\big(\ts\frac x\e\big)
}
a une unique solution dans $\e\G(V)$.
}
\pr{Preuve}
Le changement d'inconnue $y(x,\e)=\e Y\big(\frac x\e,\e\big)$ aboutit à  l'équation
$$
\frac{dY}{dX}=pX^{p-1}Y+k(X).
$$
Puisque $V$ satisfait $-\frac{3\pi}{2 p}<\a<0<\b<\frac{3\pi}{2 p}$, cette équation a pour unique solution bornée dans $\G(V)$
$$
Y(X)=\exp(X^p)\int_\infty^X\exp(-t^p)k(t)dt
$$
où le chemin d'intégration est dans $V$ et sa partie assez loin de l'origine est 
une demi-droite $\{X_1+t\tq t\in\R^+\}$ avec un certain $X_1\in V$.
\ep

Pour résoudre \rf3 dans le premier cas $h\in\~\H$, nous utiliserons le lemme suivant. Bien qu'il soit relativement classique (des résultats similaires peuvent être trouvés par exemple dans les références \cit{crss,bfsw}), nous joignons une preuve pour la complétude 
 du mémoire.
\lem{lm4.3}{
Pour tout $h\in\~\H$, il existe une unique famille $(\ac_0,...,\ac_{p-2})$ de $p-1$ séries formelles en puissances de $\eps$ sans terme constant (\ie $\ac_l\in\eps\C[\,\![\eps]\,\!]$) telle que l'équation 
\eq{eq2.7}{
\eps y'=px^{p-1}y+\eps h(x)+\sum_{l=0}^{p-2}\ac_lx^l
}
ait une unique solution formelle $\uc$ sans terme constant en puissances de $\eps$, à  coefficients analytiques dans le disque $D_1$~: $\uc(x,\eps)=\sum_{\nu>0}u_\nu(x)\eps^\nu,\;u_\nu\in\~\H$.
}
\pr{Preuve}
Injectons $\uc=\sum_{\nu>0}u_\nu(x)\eps^\nu$ et $\ac_l=\sum_{\nu>0}a_{l\nu}\eps^\nu$ dans \rf{eq2.7}, 
en tenant compte de la contrainte 
 $u_\nu$ sans pôle en $x=0$. En posant $h(x)=\sum_{\nu\geq0}h_\nu x^\nu$, on obtient, en ce qui concerne le terme d'ordre $1$ en $\eps$~:
$$
0=px^{p-1}u_1(x)+h(x)+\ds\sum_{l=0}^{p-2}a_{l1}x^l,
$$
d'où on déduit $a_{l1}=-h_l$ pour $l=0,...,p-2$ et $u_1(x)=-\frac1p\ds\sum_{\nu\geq0}h_{\nu+p-1}x^\nu$, autrement dit, $u_1=-\frac1p\S^{p-1}h$. Concernant le terme d'ordre $n$ en $\eps$ pour $n\geq2$, on obtient 
$$
u'_{n-1}(x)=px^{p-1}u_n(x)+\ds\sum_{l=0}^{p-2}a_{ln}x^l,
$$
La condition $u_n$ sans pôle en $x=0$ est satisfaite si et seulement si $a_{ln}$ est le coefficient du terme de degré $l$ dans le développement en série de $u'_{n-1}$, ce qui donne $u_n=\frac1p\S^{p-1}u'_{n-1}$.
\ep

Revenons à  l'équation \rf3 avec $h\in\~\H$. Le changement d'inconnue $y=\uc+z$, où $\uc$ est donné par le lemme \reff{lm4.3}, aboutit à  l'équation
\eq{eq2.8}{
\eps z'=px^{p-1}z-\sum_{l=0}^{p-2}\ac_lx^l.
}
Par linéarité, il suffit de montrer que pour tout $l\in\{0,...,p-2\}$ il existe une unique solution dans $\e\~C(V)$ de l'équation 
$$
\eps z'=px^{p-1}z-\eps x^l.
$$
Le changement d'inconnue $z(x,\e)=\e^{l+1}Z\big(\frac x\e,\e\big)$ donne
$$
Z'=pX^{p-1}Z-X^l
$$
dont l'unique solution dans $\G(V)$ est donnée par
$$
Z(X)=\exp(X^p)\int_X^\infty t^l\exp(-t^p)dt.
$$
Ainsi, lorsque $h\in\~\H$, l'équation \rf3 a une unique solution dans $\e\~C(V)$. 
En combinant   avec le lemme \reff{lm4.2}, on obtient donc le résultat suivant.
\lem{lm4.4}{
Pour tout $h\in\~C(V)$ il existe une unique solution de \rf3 dans $\e\~C(V)$. 
Elle est notée $\Phic(h)$. L'application $\Phic:\~C(V)\to\e\~C(V),\;h\mapsto\Phic(h)$ ainsi définie est contractante de rapport $\frac12$. 
}
La contractance s'obtient simplement en utilisant la linéarité de l'équation \rf3.
\med\\
{\sl Preuve du théorème \reff{th4.1}} \sep  
 D'après la variante du lemme \reff{lm2.1} pour $\~C(V)$ et le lemme précédent, l'application $\Phic\circ\P:\e\~C(V)\to\e\~C(V)$ est contractante de rapport $\frac12$. Puisque $\e\~C(V)$ est complet, cette  application a donc un unique point fixe dans $\e\~C(V)$.\ep

\rqs 
1.\ 
Une autre façon de présenter le lemme \reff{lm4.4} est de dire que la formule
$$
\Phic h(x,\e)=\ds\int_0^x\exp\Big\{\big(\tfrac x\e\big)^p-\big(\tfrac t\e\big)^p\Big\} h(t,\e)dt
$$
définit un opérateur $\Phic$ de $\~C(V)$ dans $\e\~C(V)$.
\med\\
2.\ 
La preuve du théorème montre que, comme pour le lemme \reff{lm2.1}, l'énoncé reste vrai si on suppose 
seulement que $P(x,y,\eps)$ est une série formelle en $\eps$ et $y$, à  coefficients  dans $\~C(V)$.\med
\med\\
3.\ 
En pratique, on ne calcule pas la solution formelle combinée de \rf1 de la manière suggérée par la
preuve de théorème \reff{th4.1}. Comme on sait qu'elle existe, on peut la calculer 
comme nous l'avons indiqué dans la remarque 
3 qui suit la proposition \reff{combextint}~: 
on commence par déterminer les développements extérieur et intérieur, puis on rejette la partie polaire du développement extérieur pour obtenir le développement lent et la partie polynomiale du développement intérieur pour obtenir le développement rapide.

Le développement extérieur est donné par la solution formelle $\sum_{n\geq1}v_n(x)\eps^n$ de \rf 1. Cette solution est à  coefficients réguliers en dehors du point tournant $0$~; elle est donnée récursivement par
$$
v_0(x)=0, \qquad v_{n+1}(x)=\frac{(v'_n-q_n)(x)}{px^{p-1}}
$$
où $q_n$ est le coefficient (dépendant de $v_1,\ldots,v_{n}$) du terme d'ordre $n$ en $\eps$ obtenu en développant  par Taylor la fonction $\eps\mapsto P\big(x,\sum_{1\leq\nu\leq n}v_\nu(x)\eps^\nu,\eps\big)$.

Pour déterminer le développement intérieur, on pose dans \rf 1 $x=\e X$, $y(x)=Y(X)$ (avec $\e^p=\eps$) et on aboutit à  l'{\sl équation intérieure}
\eq W{
\frac{dY}{dX}=pX^{p-1}Y+\e G(X,Y,\e)
}
avec $G(X,Y,\e)=P(\e X,Y,\e^p)$.
On montre qu'il existe une unique solution formelle $\widehat Y=\sum_{n\geq1}V_n(X)\e^n$ telle que 
toutes les $V_n(X)$ sont à  croissance polynomiale quand $V\ni X\to\infty$. Cette solution formelle est déterminée récursivement par $V_0=0$ et en calculant l'unique solution à  croissance polynomiale sur $V$ de l'équation
\eq k{
\frac{dV_n}{dX}=pX^{p-1}V_n+G_n(X)
} 
où $G_n$ est  le coefficient (dépendant de $V_1,\ldots,V_{n-1}$) du terme d'ordre $n-1$ en $\e$ obtenu en développant  par Taylor la fonction $\e\mapsto G\big(X,\sum_{1\leq\nu<n}V_\nu(X)\e^\nu,\e\big)$.
On trouve $$V_n(X)=\ds\int_{\infty}^X\exp(X^p-s^p)\,G_n(s)ds.$$
\ssub{4.2a}{Solutions analytiques et \dacs{}}
Le but de cette partie est de montrer que l'équation \rf1 admet des solutions ayant des \dacs~; c'est l'objet du théorème \reff{t4.5b}.
Pour des raisons de commodités, nous présentons dans un premier temps un lemme qui sera utilisé à plusieurs reprises.
Rappelons que les nombres $r_0,r_2,\eps_0,\d>0$ ont été fixés au début de \reff{5.2}, avec $\eps_0,\d>0$ suffisamment petit, que la variable $\eps=\e^p$ est dans le secteur $\Sigma=S(-\d,\d,\eps_0)$ et que $D_2$ est le disque $D(0,r_2)$.
On pose $\e_0=\eps_0^{1/p}$. Le premier résultat que  nous montrons est le suivant.
\lem{t4.5a}{On considère l'équation 
\eq{1a}{\e^p y'=p x^{p-1}+ \e^p Q(x,y,\e)}
avec $Q(x,u,\e)$ analytique bornée dans l'ensemble des $(x,u,\e)$ tels que $u\in D_2$, 
$\e\in S_0:=S\big(-\frac\d p,\frac\d p,\e_0\big)$ et  $x\in V_0(\e):=V\big(\a,\b,r_0,\mu|\e|\big)$, 
où $\mu>0$, $\a=-\frac{3\pi}{2p}+\frac{2\d}{p}$ et $\b=\frac{3\pi}{2p}-\frac{2\d}{p}$. 
Soit $r_1\in\,\big]0,r_0\big(\cos(2 \d)\big)^{1/p}\big[\,$.

Alors il existe $\e_1>0$ et une solution $y$ de \rf{1a} définie 
 dans l'ensemble des $\e,x$ tels que 
$\e\in S_1:=S\big(-\frac\d p,\frac\d p,\e_1\big)$ et $x\in V_1(\e):=V\big(\a,\b,r_1,\mu|\e|\big)$.
De plus, $y(x,\e)={\cal O}(\e)$ uniformément dans l'ensemble précédent.
}
\rq Dans cet énoncé, l'équation est quasi-linéaire et les solutions venant d'une montagne se prolongent dans un voisinage de $0$~; les quasi-secteurs peuvent donc être choisis avec $\mu>0$. Dans la partie \reff{5.3} ce ne sera plus le cas et nous aurons besoin de considérer des quasi-secteurs avec $\mu<0$. C'est une des raisons pour lesquelles nous avons défini nos \dacs{} sur des quasi-secteurs avec $\mu$ positif ou négatif, et non sur des secteurs.
\med

\pr{Preuve}
Nous construisons d'abord $y$ pour $\e\in S_0$ assez petit et 
$x$ dans une partie $\W$ de $V_0(\e)$, constituée du secteur $S(\a,\b,r_1)$ et d'un petit ``triangle curviligne'' tangent à  ce secteur. Nous montrerons ensuite que $y$ se prolonge à  $V_1(\e)\setminus\W$ si $\mu$ est assez petit. 

Soit $x_0=r_1\big(\cos(2\d)\big)^{-1/p}$. La condition de l'énoncé sur $r_1$ entraîne que $x_0$ est bien dans $V_0(\e)$ pour tout $\e\in S_0$.

Notons $\W$ la réunion de $S(\a,\b,r_1)$ et de l'intérieur du ``triangle curviligne'' dont l'image par 
\eq{FF}{
F:\C\to\C,\quad x\mapsto x^p
}
est le triangle de sommets $x_0^p, r_1^p\,e^{i2\d}$ et $r_1^p\,e^{-2\d i}$. Ce triangle image est isocèle en $x_0^p$ et a les angles $2\d,\,2\d$ et $\pi-4\d$. 
Le choix de $x_0$ est fait de telle sorte que les parties de $\partial\W$ joignant $x_0$ aux points $r_1\,e^{-2\d i/p}$ et $r_1\,e^{2\d i/p}$ ont pour images par $F$ des segments tangents au cercle de rayon $r_1^p$.
Par construction,  $V_1(\e)\setminus\W$ est une partie du disque $D(0,\mu|\e|)$.
\figu{f6.4}{
\vspace{-2mm}
\epsfxsize4.5cm\epsfbox{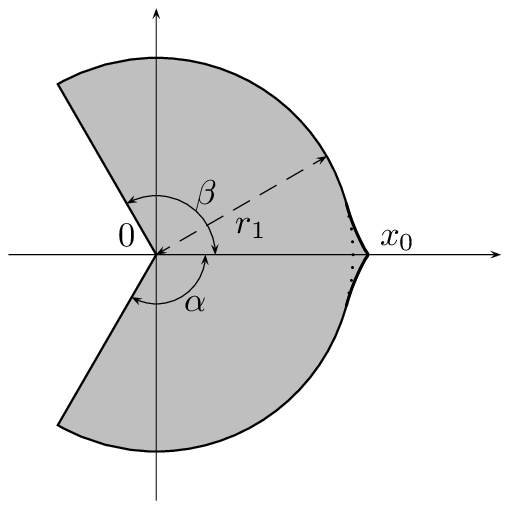}\hspace{1cm}\epsfxsize6cm\epsfbox{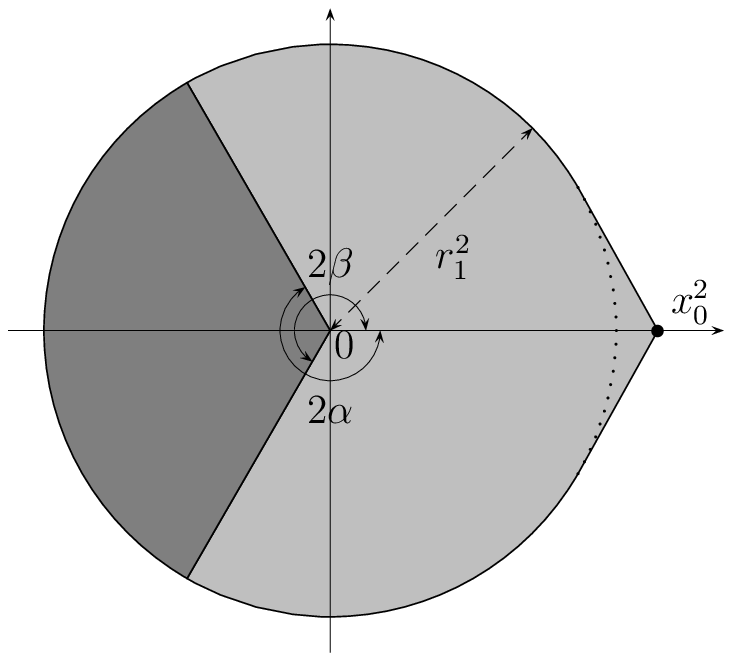}
\vspace{-10mm}
}{À gauche dans le plan des $x$, à droite les images par $F:x\mapsto x^p$ ;
ici $p=2$. En gras, le bord de $\W$. À droite, le secteur en gris foncé 
est couvert par deux feuilles de l'image de $\W$}
La raison pour avoir ajouté ce petit triangle est dans le résultat suivant.
Pour tout $|d|<\d$, la région $\W$ est $\d$-descendante à  partir de $x_0$ par rapport au relief 
\eq{RR}{
R_d: x\mapsto\re(x^pe^{-id})
}
 dans le sens suivant~:
pour tout $x\in\W$, il existe un chemin $\g$ de classe $C^1$ par morceaux, de longueur notée $\ell$,  joignant $x_0$ à  $x$ dans $\W\cup\{0,x_0\}$,  paramétré par son abscisse curviligne, \ie tel que $|\g'(t)|=1$ pour tout $t\in[0,\ell]$,  et tel que pour tout $t\in[0,\ell]$ et tout $d\in\,]-\d,\d[$
\eq t{
\re\big(\g(t)^{p-1}\g'(t)e^{-id}\big)\leq-\s|\g(t)^{p-1}|
}
avec $\s=\sin\d$.
Cela signifie qu'en chaque point de $\g$, l'angle entre la tangente à  $\g$ et la direction de plus grande pente relativement à  $d$ est 
au plus $\frac\pi2-\d$. 
Notons $\g_x$ un tel chemin.

Par exemple,  dans le cas où $|\arg x|\leq\frac\pi p$, on peut choisir $\g_x$ de telle sorte que $\arg\big(x_0^p-\g_x(t)^p\big)$ soit constant et dans le cas où $|\arg x|>\frac\pi p$ on peut choisir pour chemin $\g_x$ la réunion des segments $[x_0,0]$ et $[0,x]$. Autrement dit, on peut choisir $\g_x$ tel que $\g_x([0,\ell])^p=[x_0^p,0]\cup[0,x^p]$ ou $[x_0^p,x^p]$, mais d'autres choix de $\g_x$ sont possibles.
\med

Soit $\E$ l'espace de Banach des fonctions holomorphes $z$ sur $\W\times S_1$ pour 
lesquelles il existe une constante $Z$ telle que $\norm{z(x,\e)}
\leq Z\,\norm\e$ pour tout $(x,\e)$. Pour  $z\in\E$, on définit $\dnorm z$ comme étant la plus petites de ces constantes $Z$.
Avec $\rho$ à  choisir convenablement, soit $\M$ l'ensemble (non vide) 
des éléments $z$ de $\E$ de norme $\dnorm z\leq \rho$.

Pour $z\in\M$, on pose
\eq{TTa}{(Tz)(x,\e)=e^{x^p/\e^p}\int_{\g_x}e^{-\xi^p/\e^p} 
      P(\xi,z(\xi,\e),\e^p)\,d\xi.} 
Un point fixe $y$ de $T$ est  bien une solution de \rf{1} définie sur
$\W\times S_1$  qui satisfait $y(x_0,\e)=0$.

On veut démontrer par le théorème du point fixe de Banach que $T$ admet un point fixe
dans $\M$, si $\e_1$ est assez petit. Pour ceci et les généralisations de la 
partie \reff{5.3}, on montrera plus tard le lemme suivant.\ft
\lem{gamx}{
Soit $\g_x$ un chemin $C^1$ par morceaux joignant $x_0$ à $x$ dans $\W\cup\{0,x_0\}$, paramétré par son abscisse curviligne et vérifiant \rf t.
Pour tout $j\in\{0,1,...,p-1\}$, il existe une constante $L_j>0$ 
telle que 
$$
I_j(x,\e):=\norm{e^{x^p/\e^p}}\int_{\g_x}\norm{e^{-\xi^p/\e^p}}\norm{\xi}^j
  \,\norm{d\xi}  \leq L_j\norm{\e}^{j+1}
$$
pour tout $x\in\W$, $\e\in S_1$.

De plus, si on a $\norm{\g_x(t)}\geq L\norm\e$ pour tout $t$, alors
$$
I_j(x,\e)\leq \frac1{p\sigma} L^{1+j-p}\norm\e^{j+1}
$$
Ici $\int_{\g_x}\norm{e^{-\xi^p/\e^p}}\norm\xi^j\,\norm{d\xi}$ est une notation pour 
$\int_0^\ell \norm{e^{-\g_x(t)^p/\e^p}}\norm{\g_x(t)}^j\,dt$,
où  $\ell$ est la longueur de $\g_x$~; rappelons que $\norm{\g_x'(t)}=1$.
}
\pr{Suite de la preuve du lemme \reff{t4.5a}}
On choisit $\rho=L\sup\norm{Q}$ dans la définition de $\M$. Remarquons d'abord que
$T$ est bien définie pour tout $z\in\M$, si $\e_1$ est assez petit : il suffit
que $\e_1\rho<r_2$.
Alors, on a bien $\norm{(Tz)(x,\e)}\leq L\norm\e\sup\norm Q$ quand
$z\in\M$, $(x,\e)\in\W\times S_1$. Ceci implique $T(\M)\subset\M$.
Notons maintenant 
$$
q(x,z_1,z_2,\eps)=\int_0^1\ts\frac{\partial Q}{\partial z}\big(x,z_1+\tau(z_2-z_1)\big)\,d\tau .
$$
Alors $Q(x,z_2,\eps)-Q(x,z_1,\eps)=q(x,z_1,z_2,\eps)(z_2-z_1)$ pour tout
$x\in D_1$, $z_1,z_2\in D_2$ et $\eps\in\Sigma$.
Quitte à  réduire $r$ un peu, on peut supposer que $q$ est bornée, disons par $K$.
On obtient
$$
\norm{(Tz_2-Tz_1)(x,\e)}\leq \norm{e^{x^p/\e^p}}\int_{\g_x}\norm{e^{-\xi^p/\e^p}}
   \,K\,\norm{z_2(\xi,\e)-z_1(\xi,\e)}\,\norm{d\xi}
$$
et donc
$\dnorm{Tz_2-Tz_1}\leq LK\e_1\,\dnorm{z_2-z_1}$
pour tout $z_1,z_2\in\M$.
Si $\e_1$ satisfait aussi $LK\e_1<1$, alors $T$ est bien une contraction de $\M$ et admet un 
point fixe (unique) dans $\M$.

Ceci montre l'existence de la solution $y$ sur  $\W\times S_1$ et que ses valeurs sont dans $D'(0,\rho|\e|)$. Pour le prolongement
analytique de $y$ à  $V_1(\e)\setminus\W$, introduisons $Y(X,\e)=y(\e X,\e)$.
Elle satisfait l'équation différentielle intérieure 
$$
\frac{dY}{dX}=pX^{p-1}Y + \e Q(\e X,Y,\e^p)
$$
et prend une valeur $\norm {Y(X_0,\eps)}=\O(\norm\e)$ pour un
$X_0>0$ petit arbitraire.  Le théorème de dépendance par rapport aux conditions initiales et 
aux paramètres
implique alors que $Y$ est définie et prend des valeurs $\O(\norm\e)$ sur le
disque $\norm X\leq\mu$. Ceci permet de prolonger $y$ sur l'ensemble des
$(x,\e)$ tels que $\e\in S_1$ et $x\in V_1(\e)$, avec $y(x,\e)=\O(\e)$.
\ep
\\
\pr{Preuve du lemme \reff{gamx}}
Il s'agit de  majorer 
$$
I_j(x,\e)=\norm{e^{x^p/\e^p}}\ds \int_0^\ell \norm{e^{-\g_x(t)^p/\e^p}}
\norm{\g_x(t)}^j\,dt.
$$
On note $d=p\arg\e$ et on introduit d'abord la variable $s=-\re\big(\g_x(t)^pe^{-id}\big)$. D'après \rf t, on a
$$
\frac{ds}{dt}=-p\,\re\big(\g_x(t)^{p-1}\g_x'(t)e^{-id}\big)\geq p\,\sigma\norm{\g_x(t)}^{p-1}
$$ 
avec $\sigma=\sin\d$~; de plus, on a $\norm{\g_x(t)}\geq |s|^{1/p}$. En notant $\nu=\re(x^pe^{-id})/\norm\e^p$, on obtient 
$$
I_j(x,\e)\leq 
e^\nu
  \int_{-\re(x_0^pe^{-id})}^{-\re(x^pe^{-id})}e^{s/\norm\e^p}\,\frac1{p\,\sigma}\,
   \norm s^{\frac{j+1}p-1}\,ds\med\\
  \leq \frac{\norm\e^{j+1}}{p\,\sigma}\, e^{\nu}
   \int_\nu^\infty e^{-u} \norm u^{\frac{j+1}p-1}\,du
$$
On majore enfin $e^{\nu}\ds \int_\nu^\infty e^{-u} \norm u^{\frac{j+1}p-1}\,du=\int_0^\infty e^{-\tau} \norm{\tau+\nu}^{\frac{j+1}p-1}\,d\tau$
par $\Gamma\big(\frac{j+1}p\big)$ si $\nu$ est positif et par 
$1+\ds\int_{\max(0,-\nu-1)}^{-\nu+1}\norm{\tau+\nu}^{\frac{j+1}p-1}\,d\tau$ $\leq 1+
\ts\frac{2p}{j+1}$ si $\nu$ est négatif. 

Pour la deuxième assertion, on utilise $\norm{\g_x(t)}\geq L\norm\e$
à  la place de $\norm{\g_x(t)}\geq s^{1/p}$ et on l'obtient avec
$I_j(x,\e)\leq \frac1{p\sigma} L^{j+1-p}\norm\e^{j+1}
          e^{\nu} \int_\nu^\infty e^{-u} \,du$.
\ep

Maintenant nous sommes en position de montrer le résultat suivant. Les notations  sont celles du début de la partie \reff{5.2}~: 
$\eps=\e^p$, $r_0,r_2,\eps_0,\d>0$, $\Sigma=S(-\d,\d,\eps_0)$, $D_1=D(0,r_0)$ et $D_2=D(0,r_2)$. 
Comme dans le lemme \reff{t4.5a}, on pose $\a=-\frac{3\pi}{2p}+\frac{2\d}{p}$,
$\b=\frac{3\pi}{2p}-\frac{2\d}{p}$ et on fixe $r_1\in\,\big]0,r_0\big(\cos(2 \d)\big)^{1/p}\big[$\,.
\theo{t4.5b}{
On considère l'équation \rf1 avec $P$ analytique bornée dans $D_1\times D_2\times \Sigma$ et admettant un
développement asymptotique uniforme quand $\Sigma\ni\eps\to0$. 

Alors, pour tout $\mu>0$,  il existe $\e_1>0$ et une solution $y(x,\e)$ de \rf1 définie pour 
$\e\in S_1=S\big(-\frac\d p,\frac\d p,\e_1\big)$ et $x\in V(\e)=V\big(\a,\b,r_1,\mu|\e|\big)$.

De plus $y$ a un développement combiné quand $S_1\ni\e\to0$ quand $x\in V(\e)$.
}
\pr{Preuve}
D'après le lemme \reff{t4.5a}, il existe une solution $y(x,\e)={\cal O}(\e)$ de \rf1 quand $\e\in S_1$
et $x\in  \~V(\e):=V\big(\a,\b,\~r_1,\mu|\e|\big)$, si $r_1<\~r_1<r_0\big(\cos(2 \d)\big)^{1/p}$.
D'après le théorème \reff{th4.1}, il existe aussi une solution formelle combinée de \rf1
$\hat y=\sum_{n\geq0}\Big(a_n(x)+g_n\big(\frac x\e\big)\Big)\e^n\in 
\cch\big(r_0,V(\a-\frac\d p,\b+\frac\d p,\infty,\~\mu)\big)$,
avec un $\~\mu>\mu$. 

À présent, fixons $N\in\N$ et choisissons $\bar r_0<r_0$ tel que $\~r_1<\bar r_0\big(\cos(2 \d)\big)^{1/p}$. D'après le lemme \reff{derivasympt} appliqué à la somme partielle
$\hat y_N(x,\e)=\sum_{n=1}^{N-1} (a_n(x)+g_n(\frac x\e))\e^n$, la fonction 
$r_N(x,\e):=px^{p-1}\hat y_N(x,\e)+\e^p P(x,\hat y_N(x,\e),\e^p)-\e^p\hat y_N'(x,\e)$
admet un \dac{} quand $S_1\ni\e\to0$ et $x\in \bar V(\e):=V\big(\a,\b,\bar r_0,\mu|\e|\big)$.
Comme il s'agit d'une somme partielle d'une solution formelle, 
les coefficients de $\e^j$ 
s'annulent pour $j=0,...,N-1$. En particulier $r_N(x,\e)={\cal O}(\e^N)$ uniformément pour $x\in \bar V(\e)$.

On effectue alors le changement de variable 
$y=\hat y_N(x,\e)+\e^{N-p} z$ et on obtient pour $z$ une équation de la forme $\e^p z'=p x^{p-1}z+\e^p Q_N(x,z,\e)$,
 avec $Q_N(x,z,\e)$ bornée quand $\norm z<Z$, $\e\in S_1$ et $x\in \bar V(\e)$ avec un certain $Z>0$. On  applique le lemme \reff{t4.5a} 
à  cette équation~: elle admet
une solution $z_N(x,\e)={\cal O}(\e)$ quand $\e\in S_N:=S(-\frac\d p,\frac\d p,\e_N)$, $x\in \~V(\e)$ pour un certain $\e_N>0$.

On a donc montré que pour tout $N\in\N$, il existe $\e_N$ et une solution $y_N(x,\e)$ de \rf1
définie pour $\e\in S(-\frac\d p,\frac\d p,\e_N)$ et $x\in \~V(\e)$ telle que 
$(y_N-\hat y_N)(x,\e)={\cal O}(\e^{N-p+1})$ 
uniformément sur $\~V(\e)$. On montrera plus bas le lemme suivant.\ft
\lem{expproch}{
Soit $u_1$ et $u_2$ deux solutions de \rf1 vérifiant
$u_1(x,\e),$ $u_2(x,\e)={\cal O}(\e)$ uniformément sur l'ensemble des $(x,\e)$ tels que 
$\e\in S=S\big(-\frac\d p,\frac\d p,\e_N\big),\ x\in \~V(\e)$. Alors $u_1$ et $u_2$ sont exponentiellement proches sur $V(\e)$. 
Précisément, il existe $C>0$ tel que pour tout $\e\in S$ et tout $x\in V(\e)$
$$
|u_1(x,\e)-u_2(x,\e)|\leq Ce^{-\kappa/|\e|^p}.
$$
avec $\kappa=\big(\~r_1^p-r_1^p\big)\cos(3\d)$.
}
\pr{Fin de la preuve du théorème \reff{t4.5b}}
Puisque la solution $y=y_0$ est exponentiellement proche de chacune des solutions
$y_N$, si $\e$ est assez petit, ceci implique  qu'elle admet un \dac{}
pour $\e\in S\big(-\frac\d p,\frac\d p,\e_1\big)$ et $x\in V_0(\e)$, ce qui termine la preuve du théorème \reff{t4.5b}.
Quitte à changer les constantes, le rayon d'un secteur en $\e$ n'importe pas
pour un \dac.\ep

\pr{Preuve du lemme \reff{expproch}} 
Notons $K$ une constante telle que $\norm{u_j(x,\e)}\leq K\norm\e$ pour $j=1,2$ et pour
les $(x,\e)$ de l'hypothèse du 
 lemme.
Posons $z=u_1-u_2$.
Alors $z$ est solution de l'équation
\eq{uu}{
\e^p z'=\big(px^{p-1}+\e^p g(x,\e)\big)z
}
avec $g(x,\e)=\De_2P\big(x,u_1(x,\e),u_2(x,\e),\e^p\big)$, où --- rappelons-le --- $\De_2P$ est défini par
\eq{623}{
P(x,z,\eps)-P(x,y,\eps)=\De_2P(x,y,z,\eps)(z-y).
}
Puisque la fonction $P$ est  bornée sur $D_1\times D_2\times S$, les inégalités de Cauchy montrent que  $\De_2P$ est bornée sur $D_1\times D(0,K\~\e)
\times D(0,K\~\e)\times S$ si $\~\e$ est assez petit, donc $g(x,\e)$ est bornée
pour l'ensemble des $(x,\e)$ tels que $\e\in \~S=S\big(-\frac\d p,\frac\d p,\~\e\big)$ et
$x\in \~V(\e)$. 
\'Etant donnés  $\e\in \~S$ et $x\in V(\e)$, choisissons $\~x=\~r_1e^{i\arg x}$ si $|\arg x|<2\d$ et $\~x=\~r_1$ sinon. 
 Posons $G(x,\e)=\exp\int_{\~x}^xg(u,\e)du$. 
 \Apriori{} $G$ n'est pas continue mais puisque $g$ est bornée, c'est une fonction bornée. 
L'équation \rf{uu} donne, pour $\e\in \~S$ et $x\in \~V(\e)$
$$
z(x,\e)= z(\~x,\e)G(x,\e)\exp\big\{\tfrac1{\e^p}(x^p-\~x^p)\big\}.
$$
Par les choix de $V(\e),\~V(\e)$ et de $\~x$, 
on a pour tout $\e\in\~S$ et tout $x\in V(\e)$,
$$
\re\big(e^{-ip\arg\e}(x^p-\~x^p)\big)<\big(r_1^p-|\~x|^p\big)\cos(3\d)=-\kappa
~\mbox{ si }~|\arg x|<2\d
$$ 
et
\begin{eqnarray*}
\re\big(e^{-ip\arg\e}(x^p-\~x^p)\big)&<&r_1^p\cos\d-|\~x|^p\cos(|\arg x|-\d)\\
&<&
\big(r_1^p-|\~x|^p\big)\cos\d<-\kappa~\mbox{ sinon}.
\end{eqnarray*}
 Il existe donc une constante $C$ telle que
$\norm{z(x,\e)}\leq C \exp(-\kappa/|\e|^p)$ pour ces valeurs de $(x,\e)$. 
Quitte à  augmenter la constante, cette majoration reste valable quand
$\e\in S\big(-\frac\d p,\frac\d p,\e_N\big)$ et $x\in V(\e)$, \ie pour $|\e|\geq\e_1$.
\ep
\ssub{4.2}{Le caractère Gevrey des \dacs{}}
On montre que le \dac{} de la solution de théorème \reff{t4.5b} est 
Gevrey d'ordre $\usp$ au sens de la définition \reff{d2.3.2}. 
Le théorème-clé \reff{t3.1} se révèle particulièrement bien adapté pour obtenir ce résultat.
Le principe de la la preuve est le suivant~: nous construisons d'abord des solutions
pour $\eps$ et $x$ dans un recouvrement cohérent, 
puis nous montrons que ces solutions sont exponentiellement proches les unes  des autres (\cf lemme \reff{l5.8}). Précisément, lorsqu'on change d'un secteur en $\eps$ à  l'autre, les deux solutions sont sur une même 
montagne, donc leur différence est exponentiellement petite de la forme $\exp(-\a/|\eps|)$. 
En revanche, lorsqu'on change de secteur en $x$, les solutions sont définies sur deux montagnes adjacentes et il faut alors descendre dans la vallée les séparant pour qu'elle deviennent exponentiellement proches. 
Leur différence est donc de la forme  $\exp(-\a|x^p/\eps|)$. 
Il se trouve que les estimations obtenues correspondent exactement aux conditions d'application du théorème \reff{t3.1}.
Il nous semble
difficile de démontrer ce résultat directement, \ie dans l'esprit des sous-sections 
précédentes. 
\theo{t4.5}{
On considère l'équation \rf1 avec les hypothèses et notations de théorème \reff{t4.5b}.

Alors, pour chaque $k=0,...,p-1$,  il existe $\mu,\e_1>0$ et une solution $y$ de \rf1 définie pour $\e\in S_1:=S\big(-\frac\d p,\frac\d p,\e_1\big)$ et $x\in V_{k}(\e)=
V\big(\a_k,\b_k,r_1,\mu\norm\e\big)$ avec $\a_k=-\frac{3\pi}{2p}+\frac{2\d}{p}+\frac{2k\pi}{p}$ et
$\b_k=\frac{3\pi}{2p}-\frac{2\d}{p}+\frac{2k\pi}{p}$.

De plus $y$ a un développement combiné Gevrey d'ordre $1/p$.
}
\pr{Preuve}
Dans une première étape, nous construisons une famille de fonctions $(y_l^j)$, solutions d'équations proches de \rf1, où la fonction $P$ est remplacée par des fonctions $P_m$ ayant la même asymptotique Gevrey-1 que $P$ pour $\eps$ dans un recouvrement $(\Sigma_m)_{0\leq m<M}$. Nous vérifions ensuite que les différences de ces fonctions $y^j_l$ satisfont les estimations exponentielles requises pour appliquer le théorème \reff{t3.1}. Au cours de la preuve, nous verrons que les fonctions $P_m$ peuvent être construites comme des fonctions de $x$ et $\eps$, mais que les solutions $y_l^j$ devront être définies comme des fonctions de $x$ et de la variable $\e=\eps^{1/p}$.
\med

Commençons par considérer un recouvrement $(\Sigma_m)_{0\leq m<M}$ du disque épointé $D(0,\eps_1)^*=D(0,\eps_1)\setminus\{0\}$. Soit $M\in\N$ suffisamment grand~; en particulier il nous faut $4\pi\leq M\d$. Pour $0\leq m<M$, soit $\Sigma_m=S\big(\frac{2\pi(m-1)}M,\frac{2\pi(m+1)}M,\eps_1\big)$.
Par une transformation de Borel et Laplace tronquée ({\sl c.f.} la preuve du lemme \reff{brg1}) on construit sur chaque secteur $\Sigma_m$ une fonction $P_m$ bornée et  ayant la même asymptotique Gevrey-1 que $P$.
\med\\
Nous avons besoin de faire un tour complet dans la variable $\e$, ce qui correspond à  $p$ tours dans la variable $\eps$. Ainsi, au recouvrement $(\Sigma_m)_{0\leq m<M}$ de $D(0,\eps_1)^*$ correspond  un recouvrement $(S_l)$ de $D(0,\e_1)^*$, avec $\e_1=\eps_1^{1/p}$, de la façon suivante.
\med\\
Soit $L=pM$. Pour $0\leq l<L$, on pose $\f_l=l\frac{2\pi}L,\,\a_l=\f_l-\frac{2\pi}L(=\f_{l-1}),\,\b_l=\f_l+\frac{2\pi}L(=\f_{l+1})$ et 
$$
S_l=S(\a_l,\b_l,\e_1)=S\big((l-1)\ts\frac{2\pi}L,(l+1)\frac{2\pi}L,\e_1\big).
$$
Ainsi l'image de $S_l$ par l'application $F:\eta\mapsto\eta^p$ est le secteur $\Sigma_{l\!\mod M}$. On étend la famille $(P_m)$ par ``clonage''~: pour $l\in\{M,...,L-1\}$, on pose $P_l=P_{l\!\!\mod M}$.
\med\\
Pour $j\in\{0,\ldots,p-1\}$, soit $V^j=V(\a^j,\b^j,\infty,\mu\big)$ avec
$\a^j=j\frac{2\pi}p-\frac{3\pi}{2 p}+\frac{3\d}{p},\,\b^j=j\frac{2\pi}p+\frac{3\pi}
{2 p} - \frac{3\d}{p}$ et $\mu>0$ assez petit. 
Ainsi, dès que $\d<\frac\pi8$, on a
$$
\b_l-\a_l=\frac{4\pi}L\leq\frac\d p<\frac{\pi-6\d}{2p}=\frac12(\b^j-\a^{j+1}).
$$
À partir de ces secteurs $S_l$ et quasi-secteurs infinis $V^j$, on construit un bon recouvrement cohérent de finesse $\leq\ts\frac\d p$ en considérant le quasi-secteur $V_l^j(\e)=V\big(\a_l^j,\b_l^j,r_1,$ $\mu|\e|\big)$ avec
$$
\a_l^j=\a^j+\f_l+\ts\frac\d p=
\frac{2j\pi}p+\frac{2l\pi}L-\frac{3\pi}{2 p}+\frac{4\d}{p},
\qquad\b_l^j=\b^j+\f_l-\ts\frac\d p=
\frac{2j\pi}p+\frac{2l\pi}L+\frac{3\pi}{2 p}-\frac{4\d}{p}.
$$
Pour simplifier la suite de la présentation, nous ne considérerons les variables $X$ et $x$ que dans des secteurs. Les modifications à apporter pour les quasi-secteurs seront présentées à la fin de la preuve. Précisément on considère $X$ dans le secteur
 $S^j=S(\a^j,\b^j,\infty)$
et $x$ dans le secteur $\~S_l^j=S\big(\~\a_l^j,\~\b_l^j,r_1\big)$ avec
$$
\~\a_l^j=\a_l^j-\ts\frac{2\d}p=\frac{2j\pi}p+\frac{2l\pi}L-\frac{3\pi}{2 p}+\frac{2\d}{p},\qquad
\~\b_l^j=\b_l^j+\ts\frac{2\d}p=
\frac{2j\pi}p+\frac{2l\pi}L+\frac{3\pi}{2 p}-\frac{2\d}{p}.
$$
Le secteur $\~S_l^j$ est bissecté par $2\pi\big(\frac jp+\frac lL\big)$ et de demi-ouverture $\frac{3\pi}{2 p}-\frac{2\d}p$.
On note  $S_{l,l+1}=S_{l}\cap S_{l,l+1}$ et similairement $\~S_l^{j,j+1}$ et $\~S_{l,l+1}^{j}$.
Pour $\e\in S_l$, le secteur $\~S_l^j$ contient essentiellement 
la partie à distance $<r_1$ de 
la $j$-ième montagne du relief $R_d$ donné par \rf{RR} avec $d=\arg\eps=p\,\arg\e$, et 
presque toutes les  deux vallées adjacentes. 
L'intersection  $\~S_l^{j,j+1}$ contient, quant à elle, une grande partie de la vallée entre les montagnes $j$ et $j+1$. En particulier, on a pour tout $x\in\~S_l^{j,j+1}$ et tout $\e\in S_l$, 
\eq{3d}{
\big|(2j+1)\pi-p\,\arg\ts\frac x\e\big|<\frac\pi2-2\d. 
}
Soit $x_l^j=r_1\big(\cos(2\d)\big)^{-1/p}\,\exp\big\{2\pi i\big(\frac jp+\frac lL\big)\big\}$. 
On constate que $x_l^j$ est sur la bissectrice de $\~S_l^j$.
Choisissons pour $y_l^j$ la solution de l'équation
\eq j{
\e^py'=px^{p-1}y+\e^pP_{l}(x,y,\e^p).
}
de condition initiale $y_l^j(x_l^j,\e)=0$.

En ramenant le présent cas à celui traité dans la démonstration du 
lemme \reff{t4.5a} par des 
rotations de $\e$ et de $x$, on obtient que $y_l^j$ est définie et bornée par 
$\rho\eps_1^{1/p}$ avec un certain $\rho<r$ dans $\W_l^j\times S_l$, où  $\W_l^j$ contient le secteur  $\~S_l^j$ et un petit triangle curviligne tangent à   $\~S_l^j$ et ayant un sommet en $x_l^j$.\ft
\lem{l5.8}{
Il existe des constantes $C,A,B>0$ telle que,
pour chaque $(j,l)\in\{1,..,J\}\times\{1,..,L\}$, on ait les majorations
\eq{66}{
\left|y_{l+1}^j(x,\e)-y_{l}^j(x,\e)\right|\leq C\exp\big(-\ts\frac A{|\e|^p}\big)~\mbox{sur}~\~S_{l,l+1}^j\times S_{l,l+1},
}
\eq{77}{
\big|y_{l}^{j+1}(x,\e)-y_{l}^j(x,\e)\big|\leq C\exp\big(-B\big|\ts\frac x\e\big|^p\big)~\mbox{sur}~\~S_l^{j,j+1}\times S_l.
}
}
\pr{Preuve}
Commençons par la majoration \rf{77} et posons $z=y_{l}^{j+1}-y_{l}^j$.
Alors $z$ est solution de l'équation
\eq{zz}{
\e^p z'=\big(px^{p-1}+\e^p g(x,\e)\big)z
}
avec $g(x,\e)=\De_2P_l\big(x,y_{l}^{j}(x,\e),y_{l}^{j+1}(x,\e),\e^p\big)$, $\De_2$ défini dans \rf{623}.
Puisque la fonction $P$ est  bornée sur $D_1\times D_2\times S$, les inégalités de Cauchy montrent que  $\De_2P$ est bornée sur $D_1\times D(0,\rho\eps_1^{1/p})
\times D(0,\rho\eps_1^{1/p})\times S$ si $\eps_1$ est assez petit. 
Posons $G(x,\e)=\exp\int_0^xg(u,\e)du$. C'est une fonction bornée sur $D_1\times S$ par un certain $C_1>0$. L'équation \rf{zz} donne, pour $x\in S_l^{j,j+1}$ et $\e\in S_l$
$$
z(x,\e)= z(0,\e)G(x,\e)\exp\big\{\big(\ts\frac{x}\e\big)^p\big\}.
$$
En utilisant \rf{3d}, on obtient \rf{77} avec $B=\sin(2\d)$ et tout 
$C\geq 2\rho\eps_1^{1/p}C_1$.
\med\\
{\sl Concernant la  majoration \rf{66}}, on pose $w=y_{l+1}^{j}-y_{l}^j$.
Alors $w$ est solution de l'équation
\eq{ww}{
\e^p w'=\big(px^{p-1}+\e^p h(x,\e)\big)w+\e^pQ(x,\e)
}
avec 
$$
h(x,\e)=\De_2P_{l+1}\big(x,y_{l}^{j}(x,\e),y_{l+1}^{j}(x,\e),\e^p\big)
$$
 et 
$$
Q(x,\e)=P_{l+1}\big(x,y_{l}^{j}(x,\e),\e^p\big)-P_{l}\big(x,y_{l}^{j}(x,\e),\e^p\big).
$$
Posons $H(x,\e)=\exp
\int_0^xh(u,\e)du$. C'est une fonction bornée supérieurement par $K$ et inférieure\-ment par $\frac1K$ pour un certain $K>0$.
Soit $\xi\in\partial \W_{l,l+1}^j$ avec $\arg\xi$ bissectant $S_{l,l+1}^j$~; ainsi $\xi$ est de la forme
$\xi=r_3\,\exp\Big(2\pi i\big(\frac jp+\frac{l+1/2}L\big)\Big)$ avec $r_3>r_1$ (on trouve $r_3=r_1\cos\big(2\d-\frac\pi M\big)^{-1/p}$).
\figu{f6.6}{\epsfxsize8cm\epsfbox{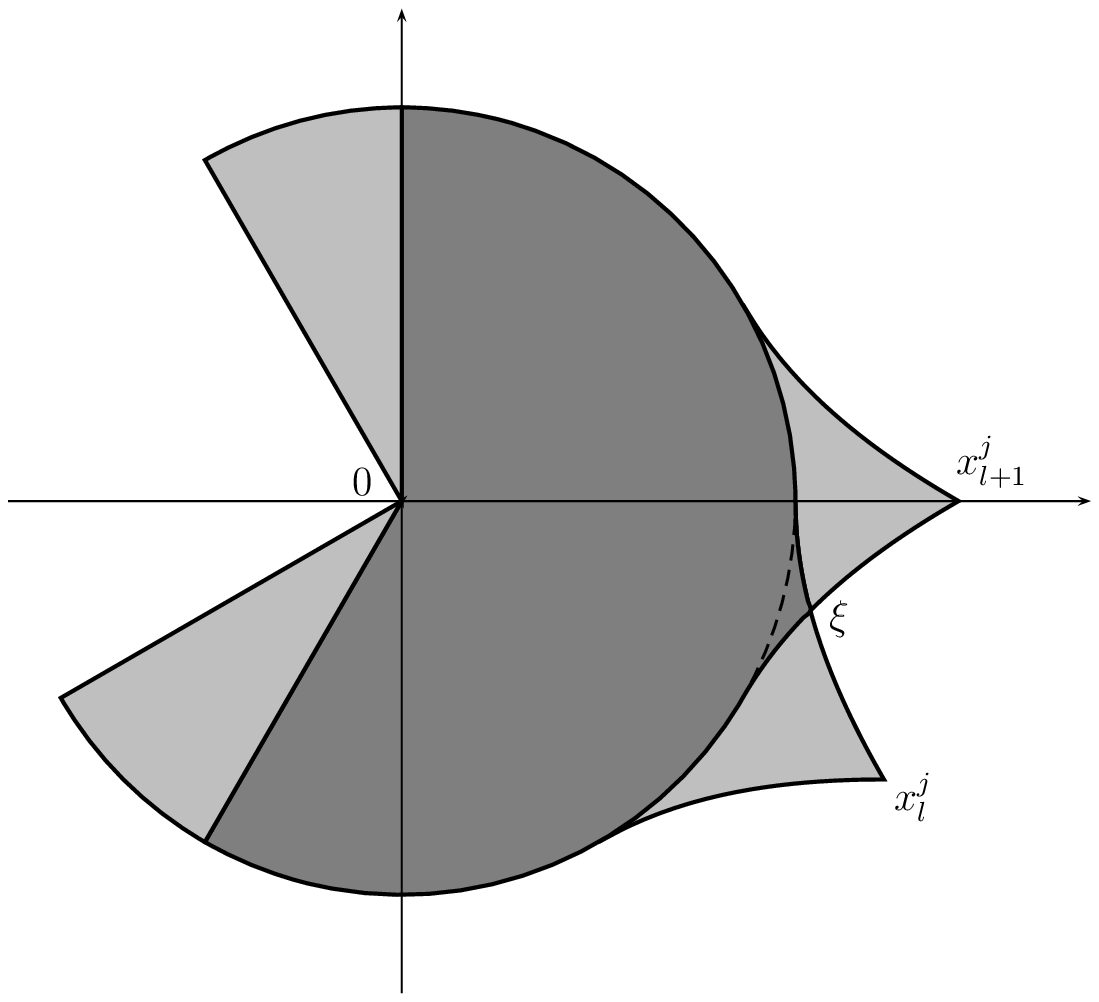}
\vspace{-6mm}
}
{Un dessin caricatural des domaines  $\W_l^j$ et $\W_{l+1}^j$ et du point $\xi$}
La formule de variation de la constante donne, pour tout $x\in  S_{l,l+1}^j$
\eq{vw}{
\begin{array}{rcl}
w(x,\e)&=&w(\xi,\e)\exp\big\{\ts\big(\frac{x}{\e}\big)^p-\big(\frac{\xi}{\e}\big)^p\big\}\dfrac{H(x,\e)}{H(\xi,\e)}\,+\med\\
&&\ds\int_\xi^x\exp\big\{\ts\big(\frac{x}{\e}\big)^p-\big(\frac{s}{\e}\big)^p\big\}\dfrac{H(x,\e)}{H(s,\e)}\,Q(s,\e)ds.
\end{array}
}
où le chemin d'intégration est choisi descendant pour le relief $R_d$ donné par \rf{RR} pour tout $d\in\big]\frac{2l\pi}M,\frac{2(l+1)\pi}M\big[$.

Concernant le premier terme de \rf{vw}, on a $|w(\xi,\e)|<2r_2$ et $\Big|\frac{H(x,\e)}{H(\xi,\e)}\Big|\leq K^2$. De plus, pour tout $x\in S_{l,l+1}^j$ et tout 
 $d\in\big]\frac{2l\pi}M,\frac{2(l+1)\pi}M\big[$, on a $\re\big(x^pe^{-id}\big)-\re\big(\xi^pe^{-id}\big)\leq- A_1$ avec $A_1= r_3\cos\big(\frac\pi M\big)-r_1$, donc 
$\Big|\exp\big\{\ts\big(\frac{x}{\e}\big)^p-\big(\frac{\xi}{\e}\big)^p\big\}\Big|\leq \exp\big(-\ts\frac {A_1}{|\e|^p}\big)$ pour tout $x\in S_{l,l+1}^j$ et tout $\e\in  S_{l,l+1}$.
La condition $A_1>0$ est assurée par le fait que $4\pi\leq M\d$.

Concernant le deuxième terme de \rf{vw}, les fonctions $P_l$ et $P_{l+1}$ sont asymptotiques Gevrey-1 à  la même série, donc il existe $C_2,A_2>0$ tels que pour  tout $x\in S_{l,l+1}^j$ et tout $\e\in  S_{l,l+1}$, $|Q(x,\e)|\leq C_2\exp\big(-\ts\frac {A_2}{|\e|^p}\big)$. Le chemin étant choisi descendant, on a, pour tout $s$ sur ce chemin, tout $\e\in S_{l,l+1}$ et tout $d\in\big]\frac{2l\pi}M,\frac{2(l+1)\pi}M\big[$,
$\re\big\{\ts\big(\frac{x}{\e}\big)^p-\big(\frac{s}{\e}\big)^p\big\}\leq0$. 
On obtient ainsi
$$|w(x,\e)|\leq 2r_2K^2 \exp\big(-\ts\frac {A_1}{|\e|^p}\big)
+(r_1+r_3)K^2C_2 \exp\big(-\ts\frac {A_2}{|\e|^p}\big),
$$
ce qui donne \rf{66} avec $A=\min(A_1,A_2)$ et tout $C\geq(2r_2+(r_1+r_3)C_2)K^2$.
\ep
\\
\pr{Suite de la preuve du théorème \reff{t4.5}}
Quitte à  diminuer $\mu$, les solutions $y_l^j$ se prolongent pour $x$ dans le disque $D(0,\mu|\e|)$, donc pour $x$ dans $V_l^j$.  
Quitte à  augmenter la constante $C$, les inégalités \rf{66} et \rf{77} 
restent valides lorsque $x$ est dans $V_l^j$.
Toutes les conditions sont réunies à  présent pour appliquer le théorème \reff{t3.1}, ce qui entraîne que les fonctions $y_l^j$ ont toutes un \dac\ Gevrey. 
En particulier, pour chaque $k\in\{0,...,p-1\}$, la solution $y=y_0^k$ répond à la question.
\ep
\ssub{6.2.4}{Remarques et extensions}
\noi
1.\ Un corollaire direct du théorème \reff{t4.5} est qu'il existe une 
solution formelle de \rf1, et que cette solution formelle est 
Gevrey d'ordre $1/p$ en $\e$. 
En d'autres termes, cela fournit une autre preuve du théorème \reff{th4.1}.
\med\\
2.\ Le  théorème \reff{t4.5} est valable dans le cas un peu plus général où $P$ est une fonction de la variable $\e$, Gevrey d'ordre $\usp$ en $\e$, au lieu de la variable $\eps$.
\med\\
3.\ Le relief de l'équation \rf1 pour $\arg\e=0$ est constitué d'une succession de $p$ montagnes et $p$ vallées correspondant aux secteurs $S\big(\frac{2k\pi}p-\frac\pi{2p},\frac{2k\pi}p+\frac\pi{2p},r_0\big)$,
respectivement $S\big(\frac{2k\pi}p+\frac\pi{2p},\frac{2k\pi}p+\frac{3\pi}{2p},r_0\big)$, $k=1,...,p$.
Nous avons montré l'existence d'une solution sur une montagne et deux vallées adjacentes ayant un \dac. On peut en déduire que toute solution $y$ de \rf1 de condition initiale $y(x_1)=y_1$ suffisamment petite, en un point $x_1$ sur une montagne, est définie et possède un \dac{} pour $\e$ d'argument petit et pour $x$ dans la partie de la montagne en dessous de $x_1$, dans la majeure partie des deux vallées adjacentes et dans un voisinage de l'ordre de $\e$ du point tournant $0$.

Précisément, étant donné un tel point $x_1\neq0$ avec
$a=\frac{2k\pi}p-\frac\pi{2p}<\arg(x_1)<b=\frac{2k\pi}p+\frac\pi{2p}$, soit 
$\d>0$ tel que $2\d<\min\{\arg(x_1)-a,b-\arg(x_1)\}$ et soit $\W$ l'ensemble des $x\in D(0,x_0)$ accessibles depuis $x_1$ par un chemin $\d$-descendant pour le relief $R_d$ donné par \rf{RR}, pour tout $d\in\,]-\d,\d[$.
La même preuve que la première étape de la preuve du théorème \reff{t4.5} montre que $y$ est définie pour $\e\in S(-\d,\d,\e_0)$ et $x\in\W\cup D(0,K|\e|)$ pour $\e_0$ et $K$ assez petits, et un argument analogue à  celui pour la preuve de \rf{77} montre que $y$ est exponentiellement proche de la solution, notée ici $\~y$, donnée par le théorème \reff{t4.5} en tout point $x$ de $\W$ suffisamment loin de $x_1$, avec un coefficient dans l'exponentielle donné par $\min_{|d|\leq\d}\{R_d(x_1)-R_d(x)\}$. Avec la proposition 
 \reff{p2.3.1} (a) ``réciproque'', il est alors immédiat que $y$ admet le même \dac{} que $\~y$.
\med\\
4.\ La région $\W$ précédente ne contient que des parties d'une montagne et de deux vallées adjacentes, mais en réalité la solution $y$ se prolonge aussi à  de grandes parties de toutes les autres vallées. En effet, 
une modification de la preuve du théorème \reff{t4.5b} montre que 
$y$ se prolonge dans chacun des secteurs $S\big(\frac{2j\pi}p+\frac\pi{2p}+\frac{2\d}p,\frac{2j\pi}p+\frac{3\pi}{2p}-\frac{2\d}p,r\big)$. Il se pose la question si ces 
prolongements admettent des \dacs{}~; on montre ci-dessous que ceci est vrai et que les 
coefficients rapides de ces \dacs{} sont les prolongements analytiques de ceux du 
théorème \reff{t4.5b} dans les autres vallées.
\med\\
5.\ Plus généralement, soit maintenant $\~y$ une solution de \rf1 de condition initiale
$\~y(0,\e)=u_0(\e)$, où $u_0:S=S(-\gamma,\gamma,\e_0)\to\C$, $\gamma>0$ assez petit%
\footnote{\ 
D'autres arguments de $\e$ peuvent être  ramenés à ce cas par des rotations
$x=\xi e^{i\varphi}$. Si on veut un secteur de plus grande ouverture en $\e$,
on le recouvre par de petits secteurs en $\e$ et on passe à l'intersection des 
secteurs en $x$ correspondants. On peut considérer une condition initiale en un
point $x=L\e$, $L\in\C$ au lieu de $x=0$ ; dans ce cas on doit d'abord étudier l'existence et 
l'asymptotique de la solution correspondante en $x=0$ à l'aide de l'équation intérieure.
},
satisfait 
$\sup_S\norm{u_0}\leq\~r<r$ et admet un développement 
asymptotique $\hat u_0(\e)$ Gevrey d'ordre $\usp$ quand $\e\to0$. 
Comme dans la remarque précédente
 et comme dans la preuve du théorème \reff{t4.5b}, on montre que 
$\~y(x,\e)$ se prolonge dans $D(0,K\norm\e)\cup \bigcup_{j}S_j$ avec $K>0$
et les secteurs $S_j$ précédents, quand $\e\in S(\a,\b,\e_1)$, si $\e_1>0$ est assez petit.
Cette solution $\~y$ diffère dans $S_j$ des deux solutions données par
théorème \reff{t4.5b} venant des montagnes adjacentes par une quantité de l'ordre
${\cal O}(e^{-B\norm x^p/\norm\e^p})$, mais ceci ne suffit pas pour montrer que $\~y$
admet un \dac{}. En effet la réciproque de la proposition \reff{p2.3.1} (b) est fausse.
Nous sommes amenés à 
utiliser de nouveau le théorème \reff{t3.1}, ce que nous  faisons ci-dessous
de manière concise selon le modèle de la preuve du théorème \reff{t4.5}.
\\
D'abord on utilise les secteurs $S_l,\ l=1,...,L$ en $\e$ de cette preuve et on construit
sur chaque $S_l$ des fonctions $P_l$ et $u_{0l}$ avec même asymptotique Gevrey que
$P$, resp. $u_0$.  On note $\~y_l$ la solution de $\e^p y'=px^{p-1}y+\e^p P_l(x,y,\e^p)$ de 
condition initiale $\~y_l(0,\e)=u_{0l}(\e)$ quand $\e\in S_l$.
On réduit les $V^j$ presque aux montagnes : $V^j=V(\a^j,\b^j,\infty,\mu)$ avec
$\a^j=j\frac{2\pi}p-\frac\pi{2p}-\frac{4\d}{p},\,\b^j=j\frac{2\pi}p+\frac\pi {2p}+
  \frac{4\d}{p}$ et $\mu>0$ assez petit. On complète en un recouvrement de $\C$ par
les quasi-secteurs $\~V^j=V(\~\a^j,\~\b^j,\infty,\mu)$, avec
$\~\a^j=j\frac{2\pi}p+\frac\pi{2p}+\frac{\d}{p}$ et $\~\b^j=(j+1)\frac{2\pi}p-\frac\pi {2p}-\frac{\d}{p}$.
Ainsi $\~V^j$ a une intersection non vide avec $V^j$ et $V^{j+1}$.
Comme dans la preuve du théorème \reff{t4.5}, on associe les quasi-secteurs $V_l^j$ et $\~V_l^j$ à
ces $S_l$, $V^j, \~V^j$ et on obtient ainsi un bon recouvrement cohérent.

De manière analogue à $\~y$, la solution $\~y_l$ se prolonge analytiquement dans 
l'union des $\~V_l^j$, $j=1,\dots p$~; on note $\~y_l^j$ sa restriction à $\~V_l^j$. Dans $V_l^j$,
on considère les (restrictions des) solutions $y_l^j$ de la démonstration du théorème
\reff{t4.5}. On montre de manière analogue à cette démonstration que des majorations 
analogues au lemme \reff{l5.8} sont satisfaites et donc le théorème \reff{t3.1} peut être appliqué. On en déduit que 
$\~y(x,\e)\sim_\usp\hat{\~y}(x,\e):=
\sum_{n\geq0}\Big(a_n(x)+\~g_n\big(\tfrac x\e\big)\Big)\e^n$ 
admet un \dac{}
Gevrey d'ordre $\usp$ dans chacun des $S_j$. Ici $a_0(x)=0$ et $\hat{\~y}(x,\e)$ est la solution 
formelle combinée de \rf1 de condition initiale $\hat{\~y}(0,\e)=\hat u_0(\e)$, 
où $\hat u_0(\e)$
est la série associée à $u_0(\e)$. En particulier, $\~g_0(X)$ est la solution de
$Y'=pX^{p-1}Y$ avec $Y(0)=u_0(0)$, \ie $\~g_0(X)=u_0(0)e^{X^p}$.
Pour le prolongement d'une solution $y$ dans les vallées discuté au point 4, on obtient
ainsi qu'ils y admettent le même \dac, \ie que leurs coefficients sont les prolongements
analytiques des coefficients du \dac{} de $y$ dans $\W$.
\med\\
6.\ Il n'était pas dans notre intention de présenter les meilleures 
constantes $A$ et $B$ dans la
démonstration du lemme \reff{l5.8}. La constante
$B=\sin(2\d)$ peut être 
largement améliorée si on considère les solutions sur des quasi-secteurs 
plus petits~; si on réduit l'angle d'ouverture à $2\pi/p+\d/p$, on peut obtenir
$B=\cos\d$ et donc une constante aussi proche de l'optimum 1 qu'on désire.
La constante $A_2$ dépend directement du type Gevrey de la fonction $P$ et ne peut pas être améliorée. La constante $A_1$ peut être choisie en fonction du point $x$. Selon les données de la 
démonstration, ce choix est d'autant plus mauvais que $r_1$ est proche de $r_0$, mais pour $r_1$ fixé il est possible de choisir les points $x_l^j$ de module proche de $r_0$, ce qui permet de récupérer une constante $A_1$ aussi proche que possible de $r_0-r_1$.
La constante optimale $A$ dépend du point $x$ et de propriétés globales de 
l'équation différentielle en dehors du cadre local que nous nous sommes fixés dans ce mémoire.
\med\\
7.\ Nous avons choisi de présenter la théorie pour un {\sl petit} secteur 
$S(-\delta,\delta,\eps_0)$ en $\eps$, car c'est le plus proche de certaines 
applications, dans lesquelles on n'a besoin que de valeurs positives de $\eps$.
Bien sûr, le résultat reste valable pour tout autre secteur en $\eps$ de petite 
ouverture car il suffit de changer $\arg\eps$ et $\arg x$. 
De plus, il est possible de déduire du théorème \reff{t4.5} des résultats sur de grands secteurs.

Par exemple, si $P(x,y,\eps)$ est réel pour $x,y,\eps$ réels, on peut montrer,
pour tout $\delta>0$ petit,
l'existence de la solution de \rf{1} avec condition initiale $0$ en $x=r_0$
pour $x\in[0,r_0]$ et 
$\eps\in S\big(-\frac\pi2+\delta,\frac\pi2-\delta,\eps_1\big)$, si $\eps_1>0$ est assez petit.
Choisissons $r_1\in\,]0,r_0[$~;
puisque cette solution est  exponentiellement proche sur $[0,r_1]$ d'une
solution du théorème pour un petit sous-secteur $\norm{\arg \eps-\varphi}<\delta$,
quelque soit $\varphi\in\,\big]\!-\frac\pi2+2\delta,\frac\pi2-2\delta\big[$, elle admet aussi
un \dac{} pour $x\in[0,r_1]$, $\eps\in S\big(-\frac\pi2+\delta,\frac\pi2-\delta,\eps_1\big)$.
Ceci sera utile pour l'étude des {\sl canards globaux} de la partie \reff{6.1}.
\med\\
8.\ 
La forme générale d'une équation dite {\sl quasi-linéaire} est
\eq v{
\eps v'=f(x)v+\eps P(x,v,\eps)
}
avec $P$ analytique dans $D\times D(0,\rho)\times S$ pour tout $\rho>0$, si $\eps_0$ est assez petit.
Nous pouvons lui appliquer le théorème \reff{4.5} en effectuant un changement de la variable $x$ et en utilisant le théorème \reff{t4.7}. En effet, la fonction 
$F(x)=\int_0^x f(s)\,ds$ satisfait $F(x)=x^{p}(1+\O(x)),\,x\to0$, donc il existe une 
fonction analytique $h$ avec $h(0)=0$ et $h'(0)=1$, telle que $F(x)=h(x)^p$. 
Soit $\f=h^{-1}$ le difféomorphisme réciproque~; ainsi on a
$F(\f(t))=t^p$. Le changement de variable $x=\f(t)$ aboutit alors à  l'équation (en posant $w(t,\eps)=v(\f(t),\eps)$ )
\eq T{
\eps \frac{dw}{dt}=pt^{p-1}w+\eps\~P(t,w,\eps)
}
avec $\~P(t,w,\eps)=\f'(t)P(\f(t),w,\eps)$.
D'après le théorème \reff{t4.5} et la remarque 2, 
une solution 
$w$ de \rf T de condition initiale $w(t_1)=w_1$ suffisamment petite en un point $t_1$ sur une montagne, est définie et possède un \dac{} pour $\e$ d'argument petit et pour $t$ dans la partie de la montagne en dessous de $t_1$, dans la majeure partie des deux vallées adjacentes et dans un voisinage de l'ordre de $|\e|$ de $t=0$.
Le théorème \reff{t4.7} (c) permet alors d'obtenir le résultat suivant, donné sans preuve.
\propo{p5.9}{   
Soit $j\in\{1,...,p\}$ et  $x_1$ un point de la $j$-ième montagne $\M_j$.
Soit $\d>0$ arbitrairement petit (en particulier tel que $|\arg F(x_1)|<\frac\pi2-2\d$).
On note $\De(x_1)$ la partie de $\V_{j-1}\cup \M_j\cup \V_j$ délimitée par
courbes d'équation $|\arg F|=\frac{3\pi}2-2\d$ et $|\arg(F-F(x_1))|=\frac{\pi}2+2\d$.
Soit $v$ une solution de \rf v de condition initiale $v(x_1)=v_1$ suffisamment petite.
Alors, pour $\a,\e_0>0$ suffisamment petits, $v$ est définie pour
$\e\in S= S(-\d,\d,\e_0)$ et $x\in \De(x_1)\cup D(0,\a|\e|)$.
De plus $v$ possède un \dac{} Gevrey pour $\e\in S$ et
$x\in \De(x_1)\cup D(0,\a|\e|),\,|x-x_1|>\d$.\vspace{-0.7cm}
}
\figu{f6.7}{\epsfxsize5cm\epsfbox{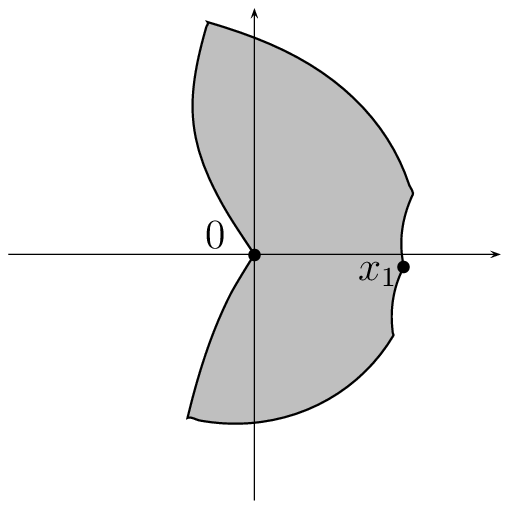}\hspace*{1cm}
   \epsfxsize5cm\epsfbox{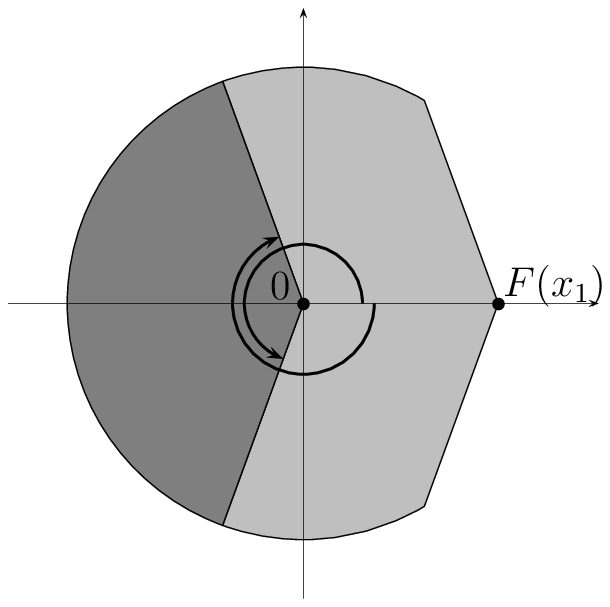}\vspace{-0.6cm}}
{Le domaine  $\De(x_1)$ et son
image par $F$~; cette image est un secteur d'ouverture $3\pi-4\d$ privé
d'un secteur de sommet $F(x_1)$ et d'ouverture $\pi+4\d$. Ici 
$F(x)=x^2+\frac{1+0,7i}4\,x^3$}
De même que dans les remarques 4 et 5, $v$ se prolonge aux autres vallées~; elle 
admet aussi un \dac{} dans ces vallées, qui est le même que celui que nous venons de 
calculer.
\sub{5.3}{\dac{} en un point tournant~: généralisations}
Nous voulons appliquer les \dacs{} à  des équations plus générales que \rf1.
Nous reprenons l'équation \rf{m}
$$
\eps z'=\Phi(x,z,\eps)
$$
avec $\Phi$ analytique par rapport à  $x$ et $z$ dans un domaine $\D\subset\C^2$ et Gevrey d'ordre 1 par rapport à  $\eps$ dans le secteur $S=S(-\d,\d,\e_0)$
 (la variable $y$ a été changée en $z$ pour des raisons de commodités).
On suppose que l'ensemble lent $\L$ contient le graphe d'une fonction lente $z_0$, analytique dans un domaine simplement connexe $D$. Ainsi on a pour tout $x\in D$, $(x,y_0(x))\in\D$ et  $\Phi(x,z_0(x),0)=0$. 
On rappelle la notation $f(x)=\frac{\partial\Phi}{\partial z}(x,z_0(x),0)$.

L'hypothèse importante que nous faisons à  présent est que {\sl $D$ contient un point tournant}  $x^*$, autrement dit que $f$ s'annule en  $x^*$.
C'est une hypothèse restrictive car généralement la fonction lente $z_0$ présente une ramification en un point tournant~; d'autres singularités de $z_0$ sont aussi possibles, par exemple des pôles. Ici nous supposons que $z_0$ reste régulière.

Quitte à  effectuer une translation et/ou une homothétie sur la variable $x$, on peut supposer sans perte de généralité que $x^*=0$ 
et que $f(x)=px^{p-1}(1+\O(x))$ lorsque $x\to0$, pour un certain $p\in\N,\,p\geq2$.
On note $F$ la primitive de $f$ s'annulant en $0$ et $R$ la partie réelle de $F$.

Le changement d'inconnue $z=z_0+y$ aboutit à  l'équation
\eq h{
\eps y'=\big(f(x)+\eps g(x,\eps)\big)y+\eps h(x,\eps)+ y^2P(x,y,\eps)
}
avec $h(x,\eps)=\frac1\eps\big(\Phi(x,z_0(x),\eps)-\Phi(x,z_0(x),0)\big)- z_0'(x)$
où $g$ et $P$ sont donnés par 
$$
\De_2\Phi(x,z_0(x),z_0(x)+y,\eps)=f(x)+\eps g(x,\eps)+ yP(x,y,\eps)~;
$$
rappelons la notation $\De_2\Phi$~:
$$
\Phi(x,z_2,\eps)-\Phi(x,z_1,\eps)=\De_2\Phi(x,z_1,z_2,\eps)\,(z_2-z_1).
$$
Le changement d'inconnue $y=\eps u$ aboutit à  l'équation
\eq u{
\eps u'=\big(f(x)+\eps g(x,\eps)\big)u+ h(x,\eps)+ \eps u^2P(x,\eps u,\eps)
}
La nouvelle fonction lente $u_0$ est obtenue en posant $\eps=0$ dans \rf u, \ie $u_0=-h_0/f$, où l'on a posé  $h_0:x\mapsto h(x,0)$.
En général, cette fonction lente a donc un pôle en $x=0$, mais dans le cas particulier où   $h_0$ a un zéro d'ordre au moins égal à  celui de $f$ en $0$, $u_0$ est  à  nouveau régulière en $x=0$. Le changement d'inconnue $u=u_0+v$ aboutit alors à  une équation quasi-linéaire \rf{v}. 
\med

Pour le cas où $h_0/f$ admet un pôle en $x=0,$ nous présentons le théorème suivant.
L'équation intérieure réduite correspondant à \rf{h} ne sera plus linéaire. 
Le prix à  payer pour cette généralisation est que la solution
n'est \apriori{} pas définie dans un disque de rayon proportionnel à  $\e$ autour de $0$.
C'est essentiellement pour cette raison que nous avons introduit les quasi-secteurs avec $\mu<0$.

Comme avant, nous supposons d'abord que $f(x)=px^{p-1}$ et nous généraliserons par la suite.
Nous reprenons les notations du début de la partie \reff{5.2} et des lemme \reff{t4.5a} et théorème \reff{t4.5b}~: $\eps=\e^p$, $r_0,r_2,\eps_0,\d>0$, $\Sigma=S(-\d,\d,\eps_0)$,  $D_1=D(0,r_0)$, $D_2=D(0,r_2)$, $\a=-\frac{3\pi}{2p}+\frac{2\d}{p}$,
$\b=\frac{3\pi}{2p}-\frac{2\d}{p}$, $r_1\in\,\big]0,r_0\cos(2 \d)^{1/p}\big[$\,.
\theo{t5.3}{
On considère l'équation 
\eq{hh}{\eps y'=px^{p-1}y+\eps h(x,\eps)+ y\,P(x,y,\eps)}
avec $h$ et $P$ analytiques bornées dans $D_1\times \Sigma$, resp.\ 
$D_1\times D_2\times \Sigma$ et admettant chacune un
développement asymptotique uniforme Gevrey d'ordre 1
quand $\Sigma\ni\eps\to0$. 

On suppose qu'il existe $r\in\{1,...,p-1\}$ tel que,
d'une part $h(x,0)=\O(x^{r-1}),\,x\to0$ et d'autre part on a un développement 
$$
P(x,y,0)=\sum_{k\geq0,l\geq1,\,k+rl\geq p-1}p_{kl}x^ky^l.
$$
Alors il existe $\mu\in\R$, $\e_1>0$ et une solution $y(x,\e)={\cal O}(\e^r)$ de \rf{hh} 
définie pour 
$\e\in S_1:=S\big(-\frac\d p,\frac\d p,\e_1\big)$ et $x\in V(\e)=V\big(\a,\b,r_1,\mu|\e|\big)$.

De plus, $y$ a un \dac{} Gevrey d'ordre $\usp$ quand $S_1\ni\e\to0$ et $x\in V(\e)$.
}
\rqs 
1.\ 
La condition sur $P$ est équivalente à  la condition que $P(x,0,0)$ $=0$ et
que la $\l$-valuation de
la série de $P(\l X,\l^rY,0)$ est au moins $p-1$. Une autre condition équivalente est 
qu'on peut écrire 
\eq{yPc}{
P(x,y,\eps)=\sum_{\ell=1}^{q-1}x^{p-1-r\ell}P_\ell(x)y^{\ell+1}+
      y^{q+1}Q_0(x,y)+\eps yQ_1(x,y,\eps)
}
avec des fonction holomorphes $P_\ell(x)$,  $Q_0(x,y), Q_1(x,y,\eps)$ sur 
$D_1$, resp. $D_1\times D_2$ et $D_1\times D_2\times\Sigma$, où $q$ désigne le plus petit entier tel que $qr\geq p-1$.

L'équation intérieure, obtenue
en posant $x=\e X$, $y=\e^r Y$, est 
$$
Y'=pX^{p-1}Y+\~P(X,Y,\e)
$$
avec
$$
\~P(X,Y,\e)= \e^{1-r}h(\e X,\e)+ \e^{1-p}Y\,P(\e X,\e^r Y,\e^p).
$$ 
Dans la limite $\e\to0$, on obtient donc
l'équation non linéaire $Y'=pX^{p-1}Y+c X^{r-1}+Y\,Q(X,Y)$,
 où $c$ est le coefficient de $x^{r-1}$ dans $h(x,0)$  et
$Q(X,Y)=\sum_{k+rl=p-1}p_{kl}X^kY^l$. 
\med\\
2.\ 
L'exemple ci-dessous montre que la condition imposée sur $P$ est nécessaire et naturelle.
Il s'agit de l'équation
\eq{e1}{
\eps y'=4x^3y-4\eps-xy^2.
}
On a donc $p=4$, $r=1$ et $p_{1\,1}\neq0$. Nous montrons à la fin de cette partie \reff{5.}
que cette équation ne peut pas avoir de solution ayant un \dac.
\med\\
\pr{Preuve du théorème \reff{t5.3}}
Elle est une modification de celle du théorème \reff{t4.5b} -- nous ne présentons que la preuve utilisant
le théorème-clé \reff{t3.1}. 
Nous montrons d'abord l'existence d'une solution $y$, ensuite celle d'une famille 
$y^j_l$ sur un
bon recouvrement et enfin que leurs différences sont exponentiellement petites.
Nous concluons en utilisant le théorème \reff{t3.1}.
Notons $S_1= S\big(-\frac\d p,\frac\d p,\e_1\big)$ avec un petit $\e_1>0$ à  déterminer.
Le domaine en $x$ est différent de celui de lemme \reff{t4.5a}, car il ne contient
pas \apriori{} $x=0$.

Soit $x_0=r_1\cos(2\d)^{-1/p}$.
Pour $m>0$, notons $\W(m)$ l'union du quasi-secteur $V(\a,\b,r_1,-m)$
et de l'intérieur du triangle curviligne $T(r_1,\d)$ de la preuve de \reff{t4.5b},%
\footnote{\ 
\ie 
l'ensemble dont l'image par  $F$ est le triangle
de sommets $x_0^p, r_1^p\,e^{2\d i}$ et $r_1^p\,e^{-2\d i}$.
} 
privée du triangle
curviligne dont l'image par l'application $F:x\mapsto x^p$ est $\Delta\big(0, m^pe^{\pi-2\d},$ $\frac{m^p}
{\sin\d}\,e^{i(\frac{3\pi}2-3\d)}\big)$, ainsi que du triangle curviligne
dont l'image par ${\cal P}$ est $\Delta\big(0,$ $m^pe^{-\pi+2\d},$ $\frac{m^p}
{\sin\d}\,e^{i(-\frac{3\pi}2-3\d)}\big)$~; voir figure \reff{fig6.8}. 
\figu{fig6.8}{
\epsfysize4.8cm\epsfbox{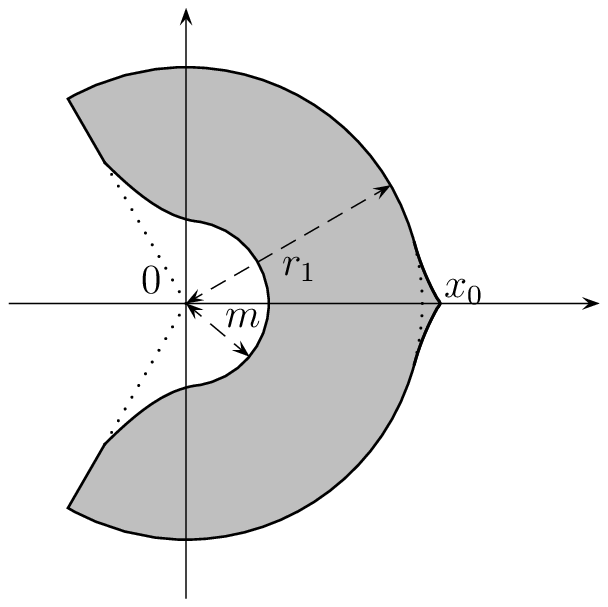}\hspace*{1.6cm}
                  \epsfxsize5.4cm\epsfbox{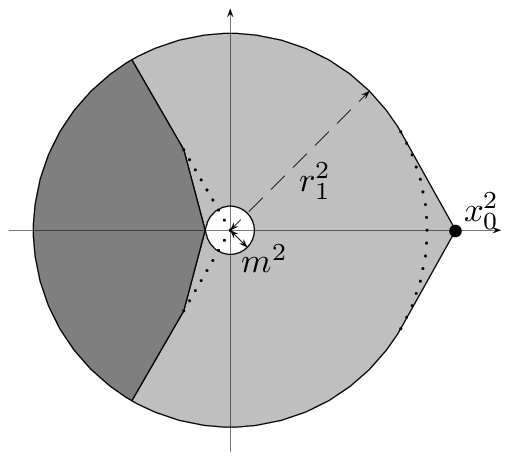}
}
{
Le domaine $\W(m)$ et son image par $F:x\mapsto x^p$~; ici $p=2$. On pourra comparer avec la figure \reff{f6.4}
}
Comme dans la preuve de lemme \reff{t4.5a}, le choix du domaine $\W(m)$ est motivé
par le fait qu'il contient $V(\a,\b,r_1,-m (\sin\d)^{-1/p})$  et 
qu'il est $\d$-descendant à partir
de $x_0$ par rapport au relief $R_d$ donné par \rf{RR} pour tout $\norm d<\d$.
La deuxième propriété peut être exprimée de la façon suivante :
pour tout $x\in\W(m)$ il existe un chemin $\g:[0,\ell]\to \W(m)\cup\{x_0\}$,
$\norm{\g'(t)}=1$ pour tout $t$, de $x_0$ à $x$ tel que pour tout $t\in[0,\ell]$
et tout $d\in]-\d,\d[$ l'inégalité \rf{t}, \ie
$$\re\big(\g(t)^{p-1}\g'(t)e^{-id}\big)\leq-\s|\g(t)^{p-1}|$$
est satisfaite avec $\sigma=\sin\d$. On note $\g_x$ un tel chemin ; pour lui,
on peut donc appliquer le lemme \reff{gamx}.

Notons $D\geq1$ un nombre réel à choisir convenablement plus tard. 
Soit $\E$ l'espace de Banach des fonctions $z$ holomorphes sur l'ensemble des
$(x,\e)$ avec $\e\in S_1$, $x\in\W(D|\e|)$ pour lesquelles il existe une constante
$Z$ telle que $\norm{z(x,\e)}
\leq Z\,\norm\e^r$ pour tout $(x,\e)$. Pour  $z\in\E$, on définit $\dnorm z$ comme étant la plus petite de ces constantes $Z$.
Avec $\rho$ à  choisir convenablement, soit $\M$ l'ensemble (non vide) 
des éléments $z$ de $\E$ de norme $\dnorm z\leq \rho$.

Pour $z\in\M$, on pose
\eq{TTc}{
\begin{array}{rl}
(Tz)(x,\e)=\tfrac1{\e^p}e^{x^p/\e^p}\ds\int_{\g_x}e^{-\xi^p/\e^p} 
      \big(\!\!\!\!&\e^ph(\xi,\e^p)+\med\\
      &z(\xi,\e)P(\xi,z(\xi,\e),\e^p)\;\big)d\xi.
\end{array}
} 
Un point fixe $y$ de $T$ est  bien l'unique solution de \rf{hh} 
satisfaisant $\lim_{x\to x_0}y(x,\e)=0$ pour tout $\e$.
L'écriture \rf{yPc} implique l'existence de constantes $C_0,...,C_q$ telles que
$$
\norm{\e^ph(\xi,\e^p)+z(\xi,\e)P(\xi,z(\xi,\e),\e^p)}\leq
  C_0 \norm \xi^{r-1}\norm\e^p + \hspace{3cm}$$
  $$ \hspace{3cm}
  \sum_{\ell=1}^{q-1} C_\ell\rho^{\ell+1}\norm \xi^{p-1-r\ell}\norm\e^{r(\ell+1)}+
  C_q\norm\e^{p-1+r}(\rho+\rho^{q+1}),
$$
quand $z\in\M$. Puisque $\norm{\g_x(x)}\geq D\norm\e$, on en déduit à  l'aide 
du lemme \reff{gamx}.2 que
$$
\norm{(Tz)(x,\e)}\leq \frac {K(\rho)}{p\sigma D}\norm\e^{r}
$$
avec $K(\rho)=C_0+ \sum_{\ell=1}^{q-1} C_\ell\rho^{\ell+1}+
C_q(\rho+\rho^{q+1})$, si $D\geq1$.
Nous choisissons donc $\rho>0$ arbitraire et $D\geq1$ assez grand pour que
$K(\rho)\leq p\sigma D\rho$. Notre majoration assure alors que 
$T:\M\to\M$.

À présent nous montrons que $T$ est une contraction, sous une condition additionnelle pour $D$.
Notons $q(x,z_1,z_2,\eps)$ la fonction holomorphe telle que
$z_2P(x,z_2,\eps)-z_1P(x,z_1,\eps)=q(x,z_1,z_2,\eps)(z_2-z_1)$. Alors
$q(x,0,0,0)=0$ pour tout $x$ et on peut écrire $q(x,z_1,z_2,0)$ de manière analogue à 
$P$ :
\eq{sum-qc}{q(x,z_1,z_2,0)=\sum_{k+\ell_1+\ell_1\geq p-1} q_{k,\ell_1,\ell_2} 
  x^kz_1^{\ell_1}z_2^{\ell_2}.}
Ceci implique une majoration 
\eq{maj-qc}{\norm{q(x,z_1(x,\e),z_2(x,\e),\e^p)}\leq
  \sum_{\ell=1}^{q-1} \~C_\ell\rho^{\ell}\norm x^{p-1-r\ell}\norm\e^{r\ell}+
  \~C_q\norm\e^{p-1}(1+\rho^{q})}
avec certaines constantes $\~C_\ell$, si $z_1,z_2\in\M$.
Comme précédemment, ceci implique
$$\dnorm{Tz_2-Tz_1}\leq\frac{\~K(\rho)}{\sigma D}\dnorm{z_2-z_1}$$
avec $\~K(\rho)=\sum_{\ell=1}^{q-1} \~C_\ell\rho^{\ell}+\~C_q(1+\rho^{q})$.
Si $D$ satisfait de plus 
$\~K(\rho)<\sigma D$, alors $T$ est une contraction sur $\M$.
On obtient ainsi l'existence de la solution $y$ énoncée dans le théorème
sur l'ensemble des $(x,\e)$ avec $\e\in S_1$, $x\in\W(D|\e|)$  et donc aussi pour
$x\in V(\a,\b,r_1,\mu\norm\e)$ avec $\mu=-D/\sin\d$.

Ensuite nous utilisons les mêmes nombres $L,M$, les mêmes angles $\f_l,\a_l,\b_l,
\a^j,$ $\b^j,\a_l^j,\b_l^j,\~\a_l^j,\~\b_l^j$ et les mêmes secteurs $S_l$, $l=1,...,L$,
$j=0,...,p-1$ que dans la démonstration de théorème \reff{t4.5}. Par une transformation
de Borel et par des transformations de Laplace tronquées, nous construisons comme
dans ladite démonstration des fonctions $h_l(x,\e)$, $P_l(x,\e)$ ayant le même 
développement asymptotique Gevrey que les données $h$ et $P$.

Au lieu des secteurs en $x$ de la preuve antérieure, nous utilisons
des quasi-secteurs $V^j=V(\a^j,\b^j,\infty,\mu)$ et $\~V_l^j(\e)=V\big(\~\a_l^j,\~\b_l^j,r_1, 
-D|\e|(\sin\d)^{1/p}\big)$ ainsi que des ensembles $\W_l^j(D|\e|)$ correspondants (\cf figure \reff{fig6.8}),
avec $D\geq1$ à déterminer.
Pour chacune des équations différentielles
$\eps y'=px^{p-1}y+\eps h_l(x,\e)+yP_l(x,y,\eps)$ --- satisfaisant les même conditions que
\rf{hh} dans l'hypothèse du théorème --- et tout $j\in\{0,...,p-1\}$,
nous construisons sa solution 
$y_l^j$ sur l'ensemble des $(x,\e)$ avec $\e\in S_l,\ x\in \W_l^j(D|\e|)$
de la même façon que $y$ ci-dessus.
Le nombre $D$ est choisi tel que ceci soit possible pour tout $(j,l)$.

Montrons à  présent que leurs différences sont exponentiellement petites de façon 
analogue au lemme \reff{l5.8}. Il s'agit de montrer qu'il existe des constantes
$C,A,B>0$ telles que pour chaque $(j,l)$ on ait
\eq{66c}{
\left|y_{l+1}^j(x,\e)-y_{l}^j(x,\e)\right|
\leq C\exp\big(-\ts\frac A{|\e|^p}\big)~
\mbox{quand}~\e\in S_{l,l+1},\ x\in \~V^j_{l,l+1}(\e),
}
\eq{77c}{
\big|y_{l}^{j+1}(x,\e)-y_{l}^j(x,\e)\big|\leq C\exp\big(-B\big|\ts\frac x\e\big|^p\big)
~\mbox{quand}~\e\in S_l,\ x\in \~V_l^{j,j+1}(\e),}
où on a utilisé comme précédemment $\~V_{l,l+1}^j$ et $\~V_l^{j,j+1}$ pour désigner les intersections correspondantes.

Commençons par la majoration \rf{77c} et posons $z=y_{l}^{j+1}-y_{l}^j$.
Alors $z$ est solution de l'équation
\eq{zza}{
\e^p z'=\big(px^{p-1}+ g(x,\e)\big)z
}
avec $g(x,\e)=q_l(x,y_l^{j}(x,\e),y_l^{j+1}(x,\e),\e^p)$, où $q_l$ est l'analogue de $q$, satisfaisant $z_2P_l(x,z_2,\eps)-z_1P_l(x,z_1,\eps)=q_l(x,z_1,z_2,\eps)(z_2-z_1)$.
De manière analogue à  \rf{maj-qc}, on obtient
\eq{maxg}{
\norm{g(x,\e)}\leq \sum_{\ell=1}^{q-1} \bar C_\ell\rho^{\ell}\norm x^{p-1-r\ell}
  \norm\e^{r\ell}+\bar C_q\norm\e^{p-1}(1+\rho^{q})
    =:\norm\e^{p-1} h\big(\big|\tfrac x\e\big|\big)
}    
avec certaines constantes $\bar C_\ell$ et donc
avec un certain polynôme $h$ de degré au plus  $p-2$.
Ceci implique 
$$
\norm{z(x,\e)}\leq K\exp\big\{\re\big(\tfrac x\e\big)^p+
   H\big(\big|\tfrac x\e\big|\big)\big\}
$$ 
où $H$ est  la primitive de $h$ s'annulant en $0$ et où $K>0$ est une constante.
Comme le degré de $H$ est au maximum $p-1$, ceci implique \rf{77c}, si
$\big|\frac x\e\big|$ est assez grand, autrement dit si $\mu$ est assez grand.

Concernant la  majoration \rf{66c}, on procède de manière analogue à  
la démonstration de \rf{66} et l'on pose $w=y_{l+1}^{j}-y_{l}^j$.
Alors $w$ est solution de l'équation
\eq{wwa}{
\e^p w'=\big(px^{p-1}+ \~g(x,\e)\big)w+ Q(x,\e)
}
avec 
$$
\~g(x,\e)=q_{l+1}\big(x,y_{l}^{j}(x,\e),y_{l+1}^{j}(x,\e),\e^p\big)
$$
 et 
$$
Q(x,\e)=\e^p(h_{l+1}-h_l)(x,\e)+
   yP_{l+1}\big(x,y_{l}^{j}(x,\e),\e^p\big)-yP_{l}\big(x,y_{l}^{j}(x,\e),\e^p\big).
$$
Posons $\~G(x,\e)=\exp\big(\e^{-p}\int_0^x\~g(u,\e)du\big)$. De manière analogue à  \rf{maxg},
on montre l'existence de constantes $L$ et $M$ telles que
$\frac1L\exp\big(-M\norm\e^{1-p}\big)\leq\norm{\~G(x,\e)}\leq L\exp\big(M\norm\e^{1-p}\big)$.
Soit de nouveau $\xi\in\partial \W_{l,l+1}^j(D|\e|)$ 
avec une distance maximale de $0$ et avec
$\arg\xi$ bissectant $V_{l,l+1}^j(\e)$~; ainsi $\xi$ est de la forme
$\xi=r_2\,\exp\big\{2\pi i\big(\frac jp+\frac{l+1/2}L\big)\big\}$ avec $r_2>r_1$ (on trouve $r_2=r_1\cos\big(2\d-\frac\pi M\big)^{-1/p}$).
La formule de variation de la constante donne, pour tout $x\in  V_{l,l+1}^j(\e)$
\begin{eqnarray}\lb{vwa}
w(x,\e)&=&w(\xi,\e)\exp\big\{\ts\big(\frac{x}{\e}\big)^p-\big(\frac{\xi}{\e}\big)^p\big\}
\
\dfrac{\~G(x,\e)}{\~G(\xi,\e)}\
+\\\nonumber&&\e^{-p}\ds\int_\xi^x\exp\big\{\ts\big(\frac{x}{\e}\big)^p-\big(\frac{s}{\e}\big)^p\big\}\dfrac{\~G(x,\e)}{\~G(s,\e)}\,Q(s,\e)ds.
\end{eqnarray}
où le chemin d'intégration est choisi descendant pour le relief $R_d$ donné par \rf{RR} pour tout $d\in\,\big]\frac{2l\pi}M,\frac{2(l+1)\pi}M\big[\,$.

Concernant le premier terme de \rf{vwa}, on a 
$$
|w(\xi,\e)|<2r~\mbox{ et }~
\Big|\frac{\~G(x,\e)}{\~G(\xi,\e)}\Big|\leq L^2\exp(2M\norm\e^{1-p}).
$$ 
De plus, pour tout $x\in V_{l,l+1}^j(\e)$ et tout 
 $d\in\big]\frac{2l\pi}M,\frac{2(l+1)\pi}M\big[$, on a 
 $$
 \re\big(x^pe^{-id}\big)-\re\big(\xi^pe^{-id}\big)\leq- A_1
 $$
  avec $A_1= r_2\cos\big(\frac\pi M\big)-r_1$, donc 
$$
\Big|\exp\big\{\ts\big(\frac{x}{\e}\big)^p-\big(\frac{\xi}{\e}\big)^p\big\}\Big|\leq \exp\big(-\ts\frac {A_1}{|\e|^p}\big)
$$
 pour tout $\e\in S_{l,l+1}$ et tout 
$x\in  V_{l,l+1}^j(\e)$.
Ceci implique que le premier terme est exponentiellement petit.

Concernant le deuxième terme de \rf{vwa}, les fonctions $h_l$ et $h_{l+1}$, respectivement
$P_l$ et $P_{l+1}$, sont asymptotiques Gevrey-1 à  la même série, donc il existe $C_2,A_2>0$ tels que pour  tout $\e\in S_{l,l+1}$ et tout $x\in  V_{l,l+1}^j(\e)$, 
$|Q(x,\e)|\leq C_2\exp\big(-\ts\frac {A_2}{|\e|^p}\big)$. Le chemin étant choisi descendant pour tout $d\in\,\big]\frac{2l\pi}M,\frac{2(l+1)\pi}M\big[\,$, on a, pour tout $s$ sur ce chemin et tout $\e\in S_{l,l+1}$,
$\re\big\{\ts\big(\frac{x}{\e}\big)^p-\big(\frac{s}{\e}\big)^p\big\}\leq0$. 
On obtient que le deuxième terme est majoré par 
$(r_1+r_2)C_2 L^2 \exp\big(- {A_2}{|\e|^{-p}+2M\norm\e^{1-p}}\big),
$
et donc lui aussi exponentiellement petit si $\e_1$ est assez petit.
Ceci démontre enfin \rf{66c}.

On termine cette preuve comme celle du théorème \reff{t4.5} 
en appliquant le théorème \reff{t3.1}.\ep

Revenons à l'équation générale \rf{h} et rappelons que $F$ est la primitive s'annulant en $0$ du
coefficient $f$ de $y$ et que $R$ est sa partie réelle.
Le relief associé, \ie le graphe de $R:\C\simeq\R^2\to\R$, consiste en une succession de $p$ montagnes $\M_j,\ j=0,...,p-1$ où $R>0$ et $p$ vallées $\V_j$ où $R<0$, délimitées par les séparatrices de col $R=0$. On numérote ces montagnes et vallées de façon à  ce qu'elles soient connexes, que $\M_0$ contienne un bout de l'axe réel positif, et qu'elles alternent $\M_j,\V_j,\M_{j+1}$ (modulo $p$). Ainsi au voisinage de $x=0$, $\M_j$ est ``tangente'' au secteur $S\big(\frac{2j\pi}p-\frac\pi{2p},\frac{2j\pi}p+\frac\pi{2p},\infty\big)$,
et $\V_j$ au secteur $S\big(\frac{2j\pi}p+\frac\pi{2p},\frac{2j\pi}p+\frac{3\pi}{2p},\infty\big)$.
Nous montrons comme corollaire du théorème \reff{t5.3} pour $j=1,...,p$
qu'il existe une solution de \rf{h} admettant un \dac{} Gevrey pour tout quasi-secteur inclus dans $\M_j\cup\V_j\cup\V_{j-1}$, \ie pour chaque montagne et 
les deux vallées adjacentes. 
Les notations  $r_0,r_1,r_2,\eps_0,\d,D_1,D_2$ et $\Sigma$ sont les mêmes que dans le théorème \reff{t5.3}. Les notations $\a$ et $\b$ ont été changées par commodité.
\coro{co616}{
On considère l'équation 
\eq{hf}{\eps y'=f(x)y+\eps h(x,\eps)+ y\,P(x,y,\eps)}
avec $f$ analytique dans $D_1$, $f(x)=px^{p-1}+{\cal O}(x^p)$ quand $x\to0$ et 
avec $h$ et $P$ analytiques bornées dans $D_1\times \Sigma$, resp.\ 
$D_1\times D_2\times \Sigma$ et admettant chacune un
développement asymptotique uniforme Gevrey d'ordre 1
quand $\Sigma\ni\eps\to0$. 

On suppose qu'il existe $r\in\{1,...,p-1\}$ tel que,
d'une part $h(x,0)=\O(x^{r-1}),\,x\to0$ et d'autre part on a un développement 
$$
P(x,y,0)=\sum_{k\geq0,l\geq1,\,k+rl\geq p-1}p_{kl}x^ky^l.
$$

Enfin, on suppose que $\a<\b$, $j$ et $r_3>0$ sont tels que
$S(\a,\b,r_3)\subset (\V_{j-1}\cup \M_j\cup \V_j)\cap D(0,r_1)$, et que $\d<\frac p6(\b-\a)$.

Alors, il existe $\mu\in\R$, $\e_1>0$ et une solution $y(x,\e)={\cal O}(\e^r)$ de \rf{hf} 
définie pour 
$\e\in S_1:=S(-\frac\d p,\frac\d p,\e_1)$ et $x\in V(\e)=V\big(\a+\frac{3\d}p,\b-\frac{3\d}p,r_3-\d,\mu|\e|\big)$.

De plus $y$ a un \dac{} Gevrey d'ordre $1/p$ quand $S_1\ni\e\to0$ et $x\in V(\e)$.
}
\pr{Preuve}
Il existe une fonction $\f$ avec $\f(u)=u\,e^{2\pi i j/p}+{\cal O}(u^2)$ telle que
$F(\f(u))=\f(u^{p})$. 
Pour $\rho>0$ assez petit, $\f$ est un difféomorphisme de $D(0,\rho)$ sur son image.
Le changement de variable $x=\f(u)$ transforme alors
\rf{hf} en \rf{hh} (avec la variable indépendante notée $u$). 
On applique le théorème \reff{t5.3}
et on obtient l'existence d'une solution $z(u,\e)$ admettant un \dac{} quand
$\e\in S_1:=S\big(-\frac\d p,\frac\d p,\e_1\big)$ et 
$u\in V\big(\~\a,\~\b,r_1,\~\mu|\e|\big)$ où $\~\a=-\frac{3\pi}{2p}+\frac{2\d}{p}$ et 
$\~\b=\frac{3\pi}{2p}-\frac{2\d}{p}$ avec certains $\e_1>0$, $\~\mu\in\R$.

On applique le théorème
\reff{t4.7} (c) au changement de variable $u=\f^{-1}(x)$.
On utilise $\rho>0$ assez petit pour que l'image de 
$V\big(\a+\frac{3\d}p,\b-\frac{3\d}p,\frac\rho2,\mu|\e|\big)$ par $\f^{-1}$ soit contenue dans 
$V\big(\~\a,\~\b,\rho,\~\mu|\e|\big)$, pour un certain $\mu\in\R$. 
Ceci est possible car la condition de l'énoncé implique $\frac{2j\pi}p-\frac{3\pi}{2p}\leq\a<\b\leq\frac{2j\pi}p+\frac{3\pi}{2p}$.
On obtient ainsi un \dac{} dans $V\big(\a+\frac{3\d}p,\b-\frac{3\d}p,\frac\rho2,\mu|\e|\big)$. Pour prolonger ce \dac{} aux $x\in V(\e)$ tels que $\frac\rho2<|x|<r_3-\d$, on utilise la proposition \reff{p4.5bis}.
\ep
\med

On détermine la série formelle combinée du \dac{} d'une manière 
analogue à la fin de \reff{5.2.1}. 
Pour la partie lente, on détermine comme avant la partie
non polaire en $x=0$ de la solution formelle extérieure de \rf{hf}.

Pour déterminer la partie rapide, on passe à l'équation intérieure 
en posant $x=\e X$ et $y=\e^rY$, comme on l'avait fait pour l'équation \rf{hh} 
dans la remarque qui suit le théorème \reff{t5.3}. Avec les notations
$f(x)=px^{p-1}+x^pf_1(x)$ et $h(x,\eps)=cx^{r-1}+x^rh_1(x)+\eps k(x,\eps)$, on a
\eq Y{
\begin{array}{rcl}
\dfrac{dY}{dX}&=&\big(pX^{p-1}+\e X^pf_1(\e X)\big)Y+cX^{r-1}+\med\\
&&\,\e M(X,\e)+\e^{1-p}\,Y\,
   P(\e X,\e^r\,Y,\e^p)
\end{array}
}
avec $M(X,\e)= X^rh_1(\e X)+\e^{p-r}k(\e X,\e^p)$.
La condition du théorème sur $P$ signifie que le dernier terme de \rf Y,  
$\e^{1-p}\,YP(\e X,\e^r\,Y,\e^p)$, est borné lorsque $X$ et $Y$ restent dans des 
compacts, uniformément par rapport à  $\e$. 
D'après cette hypothèse, en notant 
\eq{QQ}{Y\,Q(X,Y)=\ds\sum_{k+rl=p-1}p_{k,l}X^kY^{l}}
la partie quasi-homogène de plus bas degré de $P$, on obtient que ce terme est de 
la forme
$$
\e^{1-p}\,YP(\e X,\e^r\,Y,\e^p)=Y^2Q(X,Y)+\e YP_1(X,Y,\e).
$$
L'équation intérieure a pour limite lorsque $\e\to0$
\eq 0{
\frac{dY}{dX}=pX^{p-1}Y+cX^{r-1}+Y^2Q(X,Y).}
Dans la théorie des points singuliers irréguliers, on montre que \rf0 admet une solution
unique $Y_0(X)\sim -\frac cpX^{r-p}$ quand $\norm X\to\infty$ dans un quasi-secteur
$V=V\big(-\frac{3\pi}{2p}+\gamma,\frac{3\pi}{2p}-\gamma,\infty,\mu\big)$ avec $\gamma>0$ arbitraire et $\mu<0$, 
$|\mu|$ assez grand. 
La fonction 
$Y_0$ admet un développement Gevrey $\usp$
$$
Y_0(X)\sim \sum_{l\geq1}d_lX^{r-pl}\mbox{ quand }V\ni X\to\infty.
$$
La solution intérieure formelle complète (sur $V$) est déterminée en injectant 
$\hat Y(X,\e)=\sum_{n=0}^\infty Y_n(X)\e^n$ dans \rf{Y}. Ainsi, on détermine $Y_n(X)$
pour $n\geq1$ comme solution à croissance polynomiale d'une équation linéaire
non homogène $\frac{dY_n}{dX}=pX^{p-1}Y_n+h_n(X)$, où $h_n$ contient des termes
de $f_1,M$ et $P$ ainsi que $Y_0,...,Y_{n-1}$. La partie rapide du \dac{} du
corollaire \reff{co616} sur $V$, \ie pour $j=1$, est 
$\sum_{n=0}^\infty g_n(X)\e^{n+r}$ avec les parties
non polynomiales $g_n$ des $Y_n$. Pour déterminer les parties rapides des \dac{} pour
d'autres valeurs de $j$, on effectue une rotation $X\leftarrow Xe^{2\pi i/p}$.

Le théorème \reff{t5.3} et le corollaire \reff{co616} ont le défaut de ne pas contenir
d'information sur la quantité $\mu$ déterminant la distance du domaine de validité
du \dac{} et d'un secteur. L'énoncé suivant permet d'avoir cette information à partir des informations sur la solution $Y_0$ de \rf{0}.
\coro{co617} {
Sous les conditions du corollaire \reff{co616}, on suppose que la
solution $Y_0$ de l'équation intérieure réduite \rf0  satisfaisant 
$Y_0(X)\sim-\frac cpX^{r-p}$ quand $V\ni X\to\infty$
peut être prolongée sur un voisinage de l'adhérence de 
$V\big(\~\a-\frac{2j\pi i}p-\frac\d p,\~\b-\frac{2j\pi i}p+\frac\d p,\infty,\~\mu\big)$ 
avec certains $\~\a,\~\b$ satisfaisant $\a<\~\a<\~\b<\b$ et $\~\mu\in\R$. Alors 
les solutions $y$  du corollaire \reff{co616} peuvent être prolongées et admettent des \dacs{} Gevrey 
pour l'ensemble des $(x,\e)$ tels que $\e\in S\big(-\frac\d p,\frac\d p,\e_1\big)$
et $x\in V(\~\a,\~\b,\~r_1,\~\mu\norm\e)$.}
Nous proposons un cas particulier séparément.
\coro{co618}{
Sous les conditions du corollaire \reff{co616}, on suppose que
$p_{kl}=0$ si $k+rl=p-1$. Alors, pour tout $\~\mu\in\R$, il existe
$\e_1>0$ et une solution $y(x,\e)$ de \rf{hf} 
définie pour 
$\e\in S_1:=S(-\frac\d p,\frac\d p,\e_1)$ et $x\in V(\e)=V\big(\a,\b,\~r_1,\~\mu\norm\e\big)$.
\\ 
De plus $y$ a un \dac{} Gevrey d'ordre $1/p$ quand $S_1\ni\e\to0$ et $x\in V(\e)$.
}
\pr{Preuve du corollaire \reff{co617}}
D'après le corollaire \reff{co616}, les solutions
$y$ existent et admettent des \dacs{} Gevrey  quand $\e\in S\big(-\frac\d p,\frac\d p,\e_1\big)$
et $x\in V(\a,\b,\~r_1,\mu\norm\e)$ avec un certain $\mu\in\R$, éventuellement
négatif et tel que $\norm\mu$ est grand. Fixons une telle solution $y$. D'après la proposition
\reff{matching-gevrey}, $Y(X,\e):=y(\e X,\e)$ admet un développement asymptotique
uniforme Gevrey d'ordre $\usp$ de la forme $Y(X,\e)\sim_{\usp} \sum_{n\geq0}Z_n(X)\e^n$
sur des compacts de $V\big(\a-\frac\d p,\b+\frac\d p,\infty,\mu\big)$. Ici en particulier
$Z_0(X)=Y_0\big(Xe^{-\frac{2j\pi i}p}\big)$.

Choisissons $X_0\in V\big(\~\a-\frac\d p,\~\b+\frac\d p,L,\mu\big)$ avec $L>-\mu$. Alors en 
particulier $Y_0(\e):=Y(X_0,\e)$ admet un développement asymptotique
Gevrey d'ordre $\usp$.

 À présent, considérons  la solution $\~Y$ de l'équation intérieure \rf Y
avec la condition initiale $\~Y(X_0,\e)=Y_0(\e)$. 
D'une part, elle coïncide avec $Y$  sur des compacts de 
$V\big(\a-\frac\d p,\b+\frac\d p,\infty,\mu\big)$. 
D'autre part, elle se réduit à la solution $Z_0$ de \rf{0} lorsque $\e$ tend vers $0$
puisque $Y_0(0)=Z_0(X_0)$. Or d'après l'hypothèse du corollaire, $Y_0$ peut être prolongée
et donc aussi $Z_0$
peut être prolongée analytiquement sur un voisinage de l'adhérence de
$$
\W=\ts V\big(\~\a-\frac\d p,\~\b+\frac\d p,\infty,\~\mu\big)
  \setminus V\big(\a-\frac\d p,\b+\frac\d p,\infty,-L\big).
  $$
Comme \rf{Y} est régulièrement perturbée (sur des compacts), le théorème de dépendance par rapport aux paramètres et aux conditions initiales entraîne
que $\~Y(X,\e)$ est holomorphe pour $X$ dans l'adhérence de $\W$ quand
$\e\in S\big(-\frac\d p,\frac\d p,\~\e_1\big)$, si $\~\e_1>0$ est assez petit.
De plus, on obtient que $\~Y$ admet un développement asymptotique uniforme Gevrey
d'ordre $\usp$ sur cette adhérence. 

Or on avait vu que $Y$ et $\~Y$ coïncident sur un ouvert de $\C^2$, donc l'une
est un prolongement de l'autre. Les hypothèses de la proposition \reff{p4.5}
relative au 
prolongement intérieur des \dac{} sont donc satisfaites et on peut conclure.\ep  

\pr{Preuve du corollaire \reff{co618}}
Dans ce cas, l'équation \rf{0} est linéaire et la condition du corollaire  \reff{co617}
est donc trivialement satisfaite.\ep
\med

\rqs
1.\ 
Dans ses travaux \cit{ma,ma1,ma2}, É.\ Matzinger  démontre l'existence de solutions avec des développements extérieur et intérieur pour une 
équation différentielle plus générale, qui en particulier peut contenir des pôles
en $x=0$,  \cf par exemple le théorème 1 de \cit{ma2}. D'après notre proposition \reff{p3.10bis}, ceci implique aussi l'existence de \dacs{} pour ces solutions. Il se pose naturellement la question du caractère Gevrey de ces
\dacs{}, et si on peut montrer ce caractère de manière analogue à la preuve de notre théorème \reff{t5.3}.\med\\
2.\ On peut traiter de la même manière que \rf{hh} l'équation suivante, dans laquelle
le côté droit admet un développement en puissances de $\e$ au lieu de $\eps$.
Il s'agit de
\eq{hheta}{\e^p y' = f(x)y+ h(x,\e)+yP(x,y,\e),}
où $f(x)=px^{p-1}+O(x^p)$ comme avant, où $h$ et $P$ ont des développements asymptotiques Gevrey-$\usp$ et où les coefficients $h_j(x)$ et $P_j(x,y)$
 satisfont, avec un certain $r\in\{1,...p-1\}$, les conditions 
$$
h_j(x)=\sum_{l\geq 0}h_{jl}x^{l+p+r-1-j}
$$ 
$P_0(x,0,0)=0$ et pour $j=0,...,p-2$
$$
P_j(x,y,0)=\sum_{k\geq0,l\geq0,\,k+rl\geq p-1-j}p_{jkl}x^ky^l.
$$
Si $\norm{h_{00}}$ est assez petit, alors l'équation
$$
p\lambda+h_{00}+\sum_{k+rl= p-1} p_{0,k,l}\lambda^l=0
$$ 
admet une solution simple $\lambda$ de petite valeur absolue.
Il existe alors une unique transformation $y=\sum_{j=0}^{r-1}a_j(x)\e^j+z$ avec 
$a_0(x)=\lambda x^r+{\cal O}(x^{r+1})$ réduisant
\rf{hheta} à une équation de la même forme, satisfaisant les mêmes conditions
et  satisfaisant de plus $h_0=...=h_{r-1}=0$.

Pour cette dernière équation, on montre 
l'existence d'une solution $z(x,\e)={\cal O}(\e^r)$ 
pour $\e\in S_1$ et $x\in V(\e)$ avec $S_1,V(\e)$
comme dans le corollaire \reff{co616}~; ceci entraîne de suite le même énoncé pour $y$. 
Notons toutefois que la solution formelle dans le cas présent peut être très différente de celle
du corollaire \reff{co616}~; par exemple {\sl chaque} puissance de $\e$ peut contenir un facteur \og lent\fg.
\med\\
3.\ Toutes les équations de cette partie \reff{5.} 
 pourraient dépendre analytiquement
(ou seulement continûment) des paramètres additionnels. La preuve avec le théorème
du point fixe et (dans le théorème \reff{t3.1}) avec des formules intégrales
montre alors que les solutions dépendent analytiquement (resp.\ continûment)
de ces paramètres et que
leurs \dac{} Gevrey sont uniformes par rapport à ces paramètres.
\med\\
4.\ 
Nous reprenons ici plus en détails l'étude de l'exemple \rf{e1}.
Sur l'axe réel positif, on peut voir qu'il existe une unique solution $y^+$ tendant vers $0$ lorsque $x\to+\infty$. Cependant cette solution ne peut pas se prolonger jusqu'à  des $x$ de l'ordre de $\e=\eps^{1/4}$. On peut même voir qu'elle présente des singularités sur l'axe réel pour des $x$ de l'ordre de $\mu=\eps^{1/5}$. En effet, le changement de variables et inconnues $x=\mu X, y=\mu^2Y,\eps=\mu^5$ conduit à  l'équation 
\eq{mu}{
\mu \frac{dY}{dX}=4X^3Y-4-XY^2.
}
Il s'agit d'une équation singulièrement perturbée (avec $\mu$ tenant le rôle de $\eps$).
La courbe lente $X\mapsto Y_0(X)$ qui nous concerne est la branche asymptote à  l'axe $0X$ en $+\infty$ de la courbe algébrique $4X^3Y-4-XY^2=0$~; elle n'est pas définie si $0\leq X<1$. La solution $Y^+$ correspondant à  $y^+$, \ie $Y^+(X,\mu)=\mu^{-2}y^+(\mu X,\mu^5)$, satisfait $Y^+(X,\mu)\to Y_0(X)$ lorsque $\mu\to0$, pour tout $X\in\,]1,+\infty[$.
Pour tout $\d>0$, il ne peut pas exister de solution de \rf{mu} bornée sur $[1-\d,1]$ puisque le membre de droite de l'équation, $4X^3Y-4-XY^2$, est majoré par une constante négative sur $[1-\d,1-\d/2]$ pour tout $Y\in\R$. En utilisant que \rf{mu} est une équation de Riccati et en effectuant le changement d'inconnue $Y=1/Z$, on peut aussi montrer que $Y^+$ a des singularités polaires dans tout intervalle $[1-\d,1]$, si $\mu$ est assez petit.

On  constate aussi que la solution formelle de \rf{e1} présente des pôles d'ordre trop élevé par rapport à  l'entier $p=4$.
En effet, la recherche d'une solution formelle $\hat y(x,\eps)=\sum_{n\geq1}y_n(x)\eps^n$ conduit à  la récurrence
$$
y_1(x)=\frac1{x^3},\qquad y_n(x)=\frac1{4x^3}\Big(y'_{n-1}(x)+ x\sum_{k=1}^{n-1}y_k(x)y_{n-k}(x)\Big).
$$
Il s'ensuit que $y_n$ a un pôle d'ordre $5n-2$ en $x=0$. Précisément $y_n$ est de la forme $y_n(x)=x^{-5n+2}(a_n+xP_n(x))$ où $P_n$ est un polynôme et où les nombres $a_n$ satisfont $a_1=1$ et pour $n\geq2,\ a_n=\ds\sum_{k=1}^{n-1}a_ka_{n-k}$ donc sont strictement positifs. Ceci est incompatible avec l'existence d'un \dac, \cf la remarque 1 après la proposition \reff{combextint}.

Par ailleurs, il est encore possible d'appliquer notre théorie des \dacs{} sur l'équation \rf{mu} au point tournant $X=1,Y=2$, mais ceci ne sera pas développé dans le présent mémoire.
%
%
%
%
%
%
%
\sec{6.}{Applications}
Nous présentons trois situations dans lesquelles l'usage des \dacs{} se montre approprié. Notre première application concerne un problème de canard en un point tournant multiple pour une équation analytique, non seulement par rapport à la variable $x$, mais aussi par rapport au petit paramètre $\eps$. 
Il s'agit d'une part de donner une condition ``formelle'' nécessaire et suffisante pouvant être testée sur les coefficients de l'équation pour qu'il existe une solution proche de la courbe lente sur tout un intervalle ouvert contenant le point tournant, ce que nous appelons {\sl canard local}, 
et d'autre part de montrer, sous l'hypothèse d'analyticité en $\eps$, 
que lorsqu'il existe un tel canard local, alors il existe 
aussi un canard global. Autrement dit, il n'y a pas de phénomène de butée pour ce type d'équation. Cette absence de butée avait déjà été démontrée dans \cit{fs1} dans le cas d'un point tournant simple, puis par Peter De Maesschalck dans le cas d'un point tournant multiple. Nous ajoutons ici la condition formelle.  Nous présentons ce résultat d'abord pour l'équation quasi-linéaire de la partie \reff{5.2}, puis pour la généralisation de \reff{5.3}.

La deuxième application concerne des canards dits {\sl non lisses} ou {\sl angulaires}, car ils longent une courbe lente non dérivable. Ces canards avaient déjà été étudiés  par Marc Diener, Emmanuel Isambert et Véronique Gautheron \cit{di,gi,i}, et nous retrouvons une partie de leurs résultats. La théorie des \dacs{} apporte deux nouveautés~: d'une part elle permet de donner une approximation des solutions canards uniforme, d'autre part elle fournit des estimations Gevrey des développements asymptotiques.
Nous en profitons aussi pour présenter un exemple de \dac{} convergent.

Enfin, en dernière application, nous résolvons un problème de résonance au sens d'Ackerberg-O'Malley. \`A l'origine, le problème  de la résonance exposé dans \cit{aom} est un problème aux limites pour une équation linéaire du second ordre.
 Dans ce mémoire, nous ne décrivons pas ce problème original mais un problème voisin, sans conditions aux limites. Nous ne décrivons pas non plus sa relation avec le problème de surstabilité pour l'équation de Riccati associée~; nous renvoyons le lecteur à  \cit{fs1} et à  la littérature citée là.
\sub{6.1}{Canards en un point tournant multiple}
On considère l'équation
\eq{yy}{
\eps y'=f(x)y+\eps P(x,y,\eps)
}
où $f$ est analytique dans un voisinage complexe d'un intervalle réel $[a,b]$ avec $a<0<b$, $f$ réelle sur $\R$, $xf(x)>0$ si $x\neq0$, et $P$ analytique au voisinage de $[a,b]\times\{0\}\times\{0\}\subset\C^3$, $P(x,y,\eps)$ réel si $x,y,\eps$ le sont. On suppose de plus que $x=0$ est un point tournant {\sl multiple}, \ie $f(x)=\l x^{p-1}\big(1+\O(x)\big)$ si $x\to0$ avec $p$ pair, $p\geq4$ et $\l>0$. Par commodité, on se ramène à $\l=p$. Ceci peut être fait une homothétie en $x$ ou $\eps$.

Un {\sl canard local} est une solution de \rf{yy} bornée sur un intervalle ouvert contenant $0$ (``bornée'' sous-entend uniformément par rapport à  $\eps$).

Un {\sl canard global} est une solution de \rf{yy} bornée sur tout $[a,b]$.

Nous sommes dans les conditions d'application de la proposition \reff{p5.9}. Chacune des deux montagnes de l'axe réel contient une solution qui admet un \dac{} Gevrey. Notons $\M^-$ la montagne à  l'ouest contenant $[a,0[$ et $\M^+$ pour celle à  l'est contenant $]0,b]$~; notons $y^-$ la solution sur $\M^-$ et $ \ds\sum_{n\geq1}\gk{a_n(x)+g^-_n\big(\ts\frac x\e\big)}\e^n$ son \dac. De même sur $\M^+$, avec les notations similaires~:
$$
y^+(x,\e)\sim_\usp\ds\sum_{n\geq1}\gk{a_n(x)+g^+_n\big(\ts\frac x\e\big)}\e^n.
$$
\theo{t6.1}{
Les assertions suivantes sont équivalentes.
\be[\rm(a)]
\item
Il existe un canard local.
\item
Il existe un canard global.
\item
Pour tout $n\in\N$, on a $g_n^-\equiv g_n^+$.
\item
Pour tout $n\in\N$, on a $g_n^-(0) = g_n^+(0)$.
\ee
}
\rqs
1.\ 
Dans \cit{fs1}, nous avions établi un lien entre les solutions formelles, les solutions surstables et ce que nous avions appelé les {\sl canards-$\cc^\infty$}, c'est-à-dire dont toutes les dérivées sont bornées sur un intervalle $]-\d,\d[$ uniformément par rapport à  $\eps$. En utilisant ces solutions formelles et surstables, nous montrons dans \cit{fs1} l'équivalence entre  l'existence d'un canard-$\cc^\infty$ local et d'un canard-$\cc^\infty$ global.

Concernant les canards non nécessairement $\cc^\infty$, la première preuve de l'équivalence entre l'existence d'un canard local et d'un canard global est due à Peter De Maesschalck, \cf \cit{m}, à  l'aide d'éclatements 
ainsi que des résultats Gevrey de \cit{dm}.
Ici  nous ajoutons une condition  ``formelle'', 
portant sur les \dacs{}~: il existe un canard si et seulement si les fonctions de la partie rapide associées aux deux montagnes de l'axe réel coïncident.
\med\\
2.\ 
 Pour $m\geq1$, appelons {\sl canard-$\cc^m$} une solution de \rf{yy} dont toutes les dérivées jusqu'à  l'ordre $m$ inclus sont bornées sur un intervalle $]-\d,\d[$ (uniformément par rapport à  $\eps$). Alors nous pouvons préciser l'énoncé précédent : il existe un canard-$\cc^m$ (avec $m\geq1$) si, et seulement si, d'une part un des énoncés équivalents (a)--(d) est vérifié et d'autre part
$g_n^+\equiv g_n^-\equiv0$ pour $n=1,...,m-1$.
\med\\
3.\ La preuve montre qu'il y a équivalence entre l'existence d'un canard local et d'une solution bornée pour $x$ dans un domaine contenant non seulement $[a,b]$ mais aussi des secteurs autour de $[a,b]$ (correspondant aux montagnes $\M^-$ et $\M^+$ et à  une partie de leurs vallées adjacentes)  et pour $\e$ dans un secteur.
\med\\
4.\ 
En pratique, la condition $g_n^-\equiv g_n^+$ peut se vérifier sur les développements formels de solutions de l'équation intérieure. 
On procède de manière analogue à  la remarque à  la fin de \reff{5.2.1}.
L'équation intérieure pour $x=\e X\,, y(x)=Y(X),\,\eps=\e^p$ est de la forme
\eq Z{
\frac{dY}{dX}=pX^{p-1}Y+\e G(X,Y,\e)
}
avec $G(X,Y,\e)=X^pf_1(\e X)Y+P(\e X,Y,\e^p)$, $f_1$ donnée par $f(x)=px^{p-1}+x^pf_1(x)$.
Chacune des deux solutions $Y^+$ et $Y^-$ correspondant à  $y^+$ et $y^-$ a un développement en puissances de $\e$ de la forme 
$\sum_{n\geq1}Y_n^\pm(X)\e^n$
où les $Y_n^\pm$ sont donnés récursivement par 
$$
\begin{array}{rclrcl}
Y^\pm_0\!\!&=\!\!&0,&\quad 
Y_n^+(X)\!\!&=\!\!&\ds\int_{+\infty}^X\exp(X^p-s^p)\,G_n^+(s)ds,\med\\
&&&Y_n^-(X)\!\!&=\!\!&\ds\int_{-\infty}^X\exp(X^p-s^p)\,G_n^-(s)ds,
\end{array}
$$
où $G_n^\pm$ est  le coefficient (dépendant de $Y_1^\pm,\ldots,Y_{n-1}^\pm$) du terme d'ordre $n-1$ en $\e$ obtenu en développant $G\big(X,\sum_{1\leq\nu<n}Y^\pm_\nu(X)\e^\nu,\e\big)$  par la formule de Taylor.
Par récurrence, si pour $k<n$ on a $Y_k^+\equiv Y_k^-$, alors $G_n^+\equiv G_n^-$, et la condition $Y_n^+ \equiv Y_n^-$ est équivalente à $\ds\int_{-\infty}^{+\infty}\exp(-s^p)\,G_n^+(s)ds=0$ (et donc à $Y_n^+(0)=Y_n^-(0)$). 
Comme les $Y_n^\pm$ et les $g_n^\pm$ diffèrent par des polynômes et que ces polynômes sont les mêmes
pour les cas ``$+$'' et ``$-$'', cette dernière condition est équivalente à (d). Ceci
démontre déjà l'équivalence de (c) et (d).
\med\\
5.\ La condition (d) est une suite de conditions polynomiales en les 
coefficients $(c_k)_{k\geq p}$, $(q_{jkl})_{(j,k,l)\in\N^3}$ des développements
de Taylor  
$$f(x)=px^{p-1}+\sum_{k\geq p}c_kx^k\mbox{ et }P(x,y,\eps)=\sum_{(j,k,l)\in\N^3}
q_{jkl}x^jy^k\eps^l.$$ 
Pour voir ceci, on montre par récurrence que les fonctions 
$Y_n^-(X)$ de la remarque précédente sont des polynômes en 
$(c_k)_{k\geq p}$, $(q_{jkl})_{(j,k,l)\in\N^3}$   dont les coefficients sont dans 
l'ensemble $\cal F^-$ des fonctions contenant $1$ et $X$, stable par sommes, différences,
produits et par l'opérateur ${\cal T}^-$ défini par
$${\cal T}^-(H)(X)=e^{X^p}\int_{-\infty}^{X}e^{-T^p}H(T)\,dT.$$
On remarque que toutes ces fonctions sont à croissance au plus polynomiale 
quand $X\to -\infty$ et donc l'intégrale converge toujours.

De manière analogue, les fonctions $Y_n^+(X)$ sont des polynômes en 
$(c_k)_{k\geq p}$, $(q_{jkl})_{(j,k,l)\in\N^3}$ dont les coefficients sont dans 
l'ensemble ${\cal F}^+$ analogue à ${\cal F}^-$ pour $+\infty$ à la place de
$-\infty$. La condition (d), équivalente aux conditions $Y_n^+(0)=Y_n^-(0)$,
$n=0,1,...$ d'après la remarque 4, est donc une suite de conditions polynomiales en
$(c_k)_{k\geq p}$, $(q_{jkl})_{(j,k,l)\in\N^3}$ dont les coefficients sont certains
nombres connus \apriori, mais définis de manière \og transcendente \fg.

Par ailleurs, si $f$ et $P$ sont des polynômes  en toutes leurs variables, on a un nombre fini de $c_k$ et $q_{jkl}$. D'après le {\em Nullstellensatz} de Hilbert, un nombre fini de ces conditions polynomiales suffit.
\med\\
\pr{Preuve}
Nous montrons que (a) implique (c) et que (c) implique (b). L'équivalence de (c) et (d) a été établie dans la remarque 4 et l'implication de (b) à  (a) est triviale.

\bul
Supposons (a) et appelons $y$ ce canard local. Comme il doit être exponentiellement 
proche de $y^-$ en descendant un peu le relief à partir de $a$, 
il a le même \dac{} que $y^-$ pour $\e$ dans un intervalle $I=]0,\e_0[$ et pour $x$ dans un quasi-secteur $V^-(\e)=V(\pi-\d,\pi+\d,\mu|\e|)$ autour de $\R^-$, pour certains $\d,\e_0,\mu>0$. De même, il admet aussi le même \dac{} que $y^+$ pour $\e$ dans  $I$ et pour $x$ dans un quasi-secteur $V^+(\e)=V(-\d,\d,|a|-\d,\mu|\e|)$ (avec les mêmes $\d,\e_0,\mu$, quitte à  les diminuer). On a donc, pour tout $\e\in I$, tout $x\in V^-(\e)\cap V^+(\e)$ et tout $N\in\N$
$$
\sum_{n=0}^{N-1}\gk{a_n(x)+g^-_n\big(\ts\frac x\e\big)}\e^n=\sum_{n=0}^{N-1}\gk{a_n(x)+g^+_n\big(\ts\frac x\e\big)}\e^n+\O(\e^N)
$$
d'où on déduit que $g_n^-(X)=g^+_n(X)$ pour tout $n\in\N$ et tout $|X|<\mu$, donc $g_n^-$ et $g^+_n$ sont des prolongements analytiques l'une de l'autre.

\bul
Pour montrer  l'implication de (c) à  (b), nous utilisons un argument analogue à  celui dans \cit{fs1}. On considère $y^-$ la solution de condition initiale $y^-(\~a,\e)=0$ et $y^+$ la solution de condition initiale $y^+(\~b,\e)=0$, où $\~a<a<b<\~b$ sont tels que, d'une part $f$ et $P$ sont analytiques pour $x$ dans un voisinage complexe $U$ de $[\~a,\~b]$ et d'autre part $F(\~a)\neq F(\~b)$, où $F(x)=\int_0^xf$. Pour fixer les idées, on suppose $F(\~a)<F(\~b)$. 
Étant donné $\d>0$ petit, notons $D^-(\d,d)$ le domaine contenant $]\~a,0]$ dont l'image par $F$ est la réunion du triangle de sommets 
$F(\~a),iF(\~a)\tan\d,-iF(\~a)\tan\d$ et du disque centré en $0$ de rayon $d$.
De même, soit  $D^+(\d,d)$ le domaine contenant $[0,\~b[$ dont l'image par $F$ est la réunion du triangle de sommets 
$F(\~b),iF(\~b)\tan\d,-iF(\~b)\tan\d$ et du disque centré en $0$ de rayon $d$.
Si on choisit $\d$  assez petit, pour $\~\d>\d$, $\~\d$ arbitrairement proche de $\d$ et pour  $\e_0$ assez petit, les solutions $y^-$, resp. $y^+$,  sont définies et ont un \dac{} Gevrey pour $\e$ dans le secteur $S=S\big(-\frac\pi{2p}+\frac{\~\d}p,
 \frac\pi{2p}-\frac{\~\d}p,\e_0)$ et pour $x$ dans 
$D^-(\d,\d |\e|^p)$, resp. $D^+(\d,\d|\e|^p)$. 
Ceci est démontré en appliquant la remarque 7 de la partie \reff{6.2.4}.
\figu{f7.1}{
\vspace{-1.2cm}
\epsfxsize5cm\epsfbox{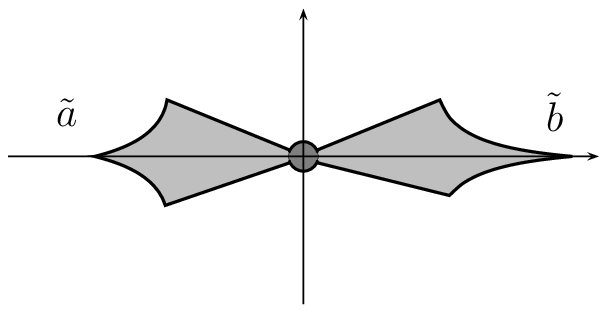}\hspace*{1cm}
\epsfxsize5cm\epsfbox{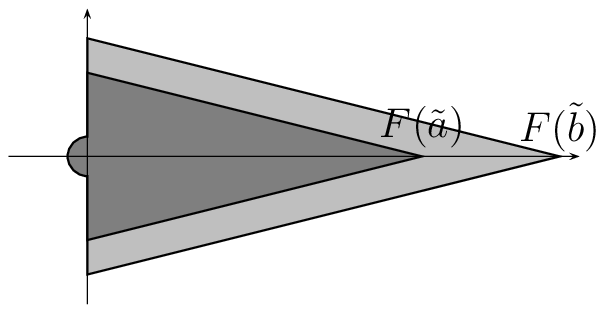}
\vspace{-1.8cm}
}
{Les domaines $D^-(r)$, $D^+(r)$ et leurs images par $F(x)=x^4$. Pour une meilleure
visibilité, les longueurs pour les images par $F$ ne correspondent pas à celles 
des originales}

Considérons la fonction $\f:S\to\C,\,\e\mapsto y^+(0,\e)-y^-(0,\e)$. 
Puisque les solutions $y^\pm$ ont des \dacs{} Gevrey avec les mêmes coefficients $g_n$, 
leur différence au point 0 admet un développement asymptotique Gevrey nul. Comme il est bien connu, ceci implique $\f(\e)=\O(\exp(-r(\~\delta)/|\e|^p))$ pour un certain 
$r=r(\~\delta)>0$ dépendant de $\~\delta$.
On montrera plus tard, pour $\~\delta,\e_0>0$ assez petits et pour tout $\d\in\,]0,\~\d[$, que
\eq{expo71}{\f(\e)=y^+(0,\e)-y^-(0,\e)={\cal O}(-F(\~a)\sin\d/\norm\e^p)}
quand $\e\in S=S\big(-\frac\pi{2p}+\frac{\~\d}p,
 \frac\pi{2p}-\frac{\~\d}p,\e_0)$.
Pour $p|\arg\e|=\frac\pi2-\~\d$, la majoration \rf{expo71} implique
$\f(\e)=\O\Big(\exp\big(-q\frac{ F(\~a)}{\e^p}\big)\Big)$ avec 
$q=\frac{\sin\d }{\sin\~\d}$ arbitrairement proche de $1$. 
D'après le théorème de Phagmén-Lindelöf, ceci reste valide pour $\arg\e=0$. Le lemme de Gronwall montre alors que la solution $y^+$ est définie et bornée pour des valeurs de $x$ arbitrairement proches de $\~a$, donc sur $[a,b]$ tout entier.

Pour la démonstration de \rf{expo71}, puisqu'on a supposé que $f$ ne s'annule qu'en $0$, il existe $g$ analytique réelle telle que $F(g(t))=t^p$. Par le changement de variable $x=g(t)$ (et en conservant la notation $x$ au lieu de $t$) on se ramène au cas où
$F(x)=x^p$ et donc $f(x)=px^{p-1}$.
\`A présent on utilise l'existence de $\mu>0$, et pour tout $\g>0$ d'un entier 
$L>\frac\pi{4p\g}$ et de $\e_0,\rho>0$ assez petits, tels qu'il existe des solutions 
$z_\ell^\pm(x,\e),\ \ell=-L,...,L$ de \rf{yy} définies, holomorphes et bornées quand
$\e\in S_\ell=S(\ell\frac\pi{2pL}-\g,\ell\frac\pi{2pL}+\g,\e_0)$ et 
$x\in V_\ell^\pm=\pm V(-\frac{3\pi}{2p}+\ell\frac\pi{2pL} + 2\g,
     \frac{3\pi}{2p}+\ell\frac\pi{2pL}-2\g,\rho,\mu\norm\e )$ 
qui admettent des
\dac{} Gevrey d'ordre $\usp$ 
$$z_\ell^\pm\sim_{\usp}\sum_{n=0}^{N-1}\gk{a_n(x)+g^\pm_n\big(\ts\frac x\e\big)}\e^n$$
Ceci est une conséquence du théorème \reff{t4.5} et de rotations 
$\e=e^{i\psi}\~\e$, $x=\pm e^{i\psi}\~x$ avec $\psi=\ell\frac\pi{2pL}$ ; le fait que
les fonctions $a_n$ et $g_n^\pm$ sont indépendantes de $\ell$ est dû à l'unicité
de la solution formelle de \rf{yy}, \ie\ au théorème \reff{th4.1}.

Or d'après notre hypothèse (c), pour tout $\ell$,
les fonctions $z^+_\ell(0,\e)$ et $z_\ell^-(0,\e)$
ont le même développement asymptotique qui est aussi Gevrey d'ordre $\usp$ en $\e$ et sont donc exponentiellement proches l'une de l'autre. Il existe donc $s>0$ tel que
\eq{expoz}{z^+_\ell(0,\e)-z_\ell^-(0,\e)={\cal O}(e^{-s/\norm\e^p})}
quand $\ell\in\{-L,...,L\}$ et $\e\in S_\ell$. Quitte à réduire $\rho$, on peut supposer
que $\rho^p<s$. 

Considérons maintenant $\~\d>0$ tel que $F(\~\a)\sin\~\d<\rho^p$ et $\d\in\,]0,\~\d[$ arbitraire.
Soit $\ell\in\{-L,...,L\}$ et $\e\in S_\ell\cap S$, $\arg\e=\psi$. Alors $y^+(x,\e)$ et 
$z_\ell^+(x,\e)$ sont holomorphes bornées sur un voisinage du segment 
$x\in [0,Te^{i\psi}]$ où $T^p=F(\~\a)\sin\d$.
Comme leur différence $D=y^+-z_\ell^+$ satisfait
$\e^p D'=\big(px^{p-1}+\e^p \Delta_2P(x,z_\ell^+(x,\e),y^+(x,\e),\e^p)\big)D,$
on en déduit que
\eq{expoyz}{
\begin{array}{rcl}
y^+(0,\e)-z_\ell^+(0,\e)&=&{\cal O}\gk{\exp\big(-(Te^{i\psi})^p/\e^p\big)}\med\\
& =&{\cal O}\gk{e^{-T^p/\norm\e^p}}\mbox{ 
quand }\e\in S_\ell\cap S.
\end{array}}
De manière analogue, on montre que $y^-(0,\e)-z_\ell^-(0,\e)=
  {\cal O}\gk{e^{-T^p/\norm\e^p}}$ quand $\e\in S_\ell\cap S.$ Avec \rf{expoyz}
et \rf{expoz}, ceci implique \rf{expo71} quand $\e\in S_\ell\cap S.$
Comme l'union des $S_\ell$ contient $S$, la démonstration de \rf{expo71} est donc 
complète.
\ep

Le théorème \reff{t6.1} peut être généralisé pour l'équation \rf{hf} du corollaire
\reff{co616}, \ie pour
$$\eps y'=f(x)y+\eps h(x,\eps)+ y\,P(x,y,\eps)$$
avec $f$ analytique dans $D_1$, $f(x)=px^{p-1}+{\cal O}(x^p)$ quand $x\to0$, $p$ pair et 
avec $h$ et $P$ analytiques bornées dans $D_1\times \Sigma$, resp.\ 
$D_1\times D_2\times \Sigma$ et admettant chacune un
développement asymptotique uniforme Gevrey d'ordre 1
quand $\Sigma\ni\eps\to0$, où $D_1$ est un voisinage de l'intervalle $[a,b]$, 
$D_2=D(0,r_2)$ 
 et $\Sigma=S(-\d,\d,\eps_0)$, et où $-a,b,r_2,\eps_0,\d>0$, 
$\d$ assez petit. On suppose de plus que
les valeurs des fonctions $f,h,P$ sont réelles quand leurs arguments sont réels et
qu'il existe $r\in\{1,...,p-1\}$ tel que,
d'une part $h(x,0)=\O(x^{r-1}),\,x\to0$ et d'autre part  
$$
P(x,y,0)=\sum_{k\geq0,l\geq1,\,k+rl\geq p-1}p_{kl}x^ky^l.
$$

D'après le corollaire \reff{co616}, 
l'équation admet des solutions $y^\pm(x,\e)$ holomorphes dans l'ensemble
des $(x,\e)$ avec $\norm{\arg\e}<\d$, $\norm\e<\e_0$, 
$x \in V^+(\e)=V(-\b,\b,b,\mu|\e|)$,
respectivement $x \in V^-(\e)=V(\pi-\b,\pi+\b,|a|,\mu|\e|)$
avec certains $\d,\e_0,\b>0$ et un certain $\mu<0$ pour $\d$ assez petit.
Ces solutions admettent des \dacs{} 
\eq{dacgen}{y^\pm(x,\e)\sim_\usp\sum_{n=r}^{p-1}g^\pm_n\big(\tfrac x\e\big)\e^n+
   \sum_{n=p}^{\infty}\Big(a_n(x)+g^\pm_n\big(\tfrac x\e\big)\Big)\e^n.}

Contrairement à l'équation \rf{yy},
les domaines d'existence des solutions ont \apriori{} une intersection vide.
Ceci nécessite des modifications de l'énoncé du théorème et de sa preuve. 
\theo{t6.2}{
Sous les conditions et avec les notations précédentes,
les énoncés suivants sont équivalents.
\be[\rm(a)]\item
 Il existe un canard local $y$ avec $y(x,\e)={\cal O}(\e^r)$  uniformément dans un voisinage de $x=0$.
\item
 Il existe un canard global $y$ avec $y(x,\e)={\cal O}(\e^r)$ uniformément
sur $[a,b]$.
\item
Les fonctions $g^\pm_n(X)$ peuvent être prolongées 
analytiquement en un voisinage de $\R$ et ces prolongements coïncident, 
\ie $g_n^+\equiv g_n^-$ pour tout $n$.
\ee
}
\pr{Preuve}
La preuve que (c) implique (a) est presque identique à 
celle du théorème \reff{t6.1}. La seule observation à ajouter est la suivante~: puisque $g_r^+=g_r^-$ sont analytiques sur $\R$, 
on peut appliquer le corollaire \reff{co617} et on obtient que $y^+(x,\e)$ et son \dac{} 
se prolongent dans un  certain
quasi-secteur $V(-\~\b,\~\b,b,\~\mu|\e|)$ avec $\~\b,\~\mu>0$~; 
pour $y^-$ et les $z_\ell^\pm$ de la preuve de théorème \reff{t6.1}, l'énoncé analogue est vrai.
Le reste de la preuve est inchangé.

La preuve que (c) est une conséquence de (a) nécessite
l'utilisation de l'équation intérieure. Le théorème sur les \dacs{} de \rf{hf} dit 
seulement que le canard local $y(x,\e)$ admet un \dac{} sur $[a+\d,-L\e]$, 
$L=\norm{\mu}$, (et même dans  $V^+(\e)\cap D(0,\norm a-\d)$) et un autre \dac{}
sur $[L\e,b-\d]$. Pour la solution $Y(X,\e)=\e^{-r}y(\e X,\e)$ de l'équation 
intérieure, ceci implique que ses valeurs $U^\pm(\e)=Y(\pm L,\e)$ admettent des
développements asymptotiques Gevrey d'ordre $\usp$ quand $\e\to 0.$  Par hypothèse,
 la fonction
$Y$  est bornée sur $[-L,L]\times]0,\e]$. Comme les 
valeurs de $Y'(X,\e)$ sont exprimées par l'équation intérieure
$$
Y'=pX^{p-1}Y+\~P(X,Y,\e)$$\mbox{ avec } $$\~P(X,Y,\e)=\e^{1-r}h(\e X,\e)+
      \e^{1-p}\, Y P(\e X,\e^r Y,\e^p),
$$
et donc par une équation différentielle régulièrement perturbée,
ces valeurs sont
aussi bornées. Le théorème d'Arzela-Ascoli montre que toute suite $(\e_k)_{k\in\N}$ de
nom\-bres strictement positifs tendant vers $0$ admet une sous-suite 
$(\e_{k_l})_{l\in\N}$ telle que la suite $\big(Y(X,\e_{k_l})\big)_{l\in\N}$ converge uniformément sur $[-L,L]$ ;
la limite est nécessairement la solution $G$ de l'équation intérieure réduite 
\eq{intred}
{Y'=pX^{p-1}Y+cX^{r-1}+Y\,Q(X,Y), \ \ Q(X,Y)=\sum_{\stackrel{k\geq0,l\geq1}{k+rl=p-1}} p_{kl}X^kY^l,
}
avec la condition initiale $G(L)=\lim_{\e \to0} U^+(\e)$. On obtient donc ici
que les premiers termes $g_r^\pm(X)$ des \dacs{} peuvent être prolongés sur $\R$ et 
coïncident. Comme ci-dessus, on conclut que les \dacs{} peuvent être prolongés sur
des quasi-secteurs avec un certain $\~\mu>0$. 
Le reste de la preuve est comme celle du théorème \reff{t6.1}. \ep

\rq Les conditions ne sont plus équivalentes à une suite de conditions
polynomiales comme dans la remarque 5.\ après le théorème \reff{t6.1} ; elles sont 
transcendentes en $c$ et les $p_{kl}$ de \rf{intred}. 
\sub{6.2}{Canards non lisses}
{\noi\bf Equations de type \og Union Jack \fg} \sep
\med\\
On considère une équation différentielle de la forme
\eq{uj}{\eps y'=y(y-x)(y+x)+P(x,y,\eps)+\eps c,}
où $P$ est analytique sur $D_1\times D(0,r)\times D(0,\eps_1)$, $D_1$
un voisinage d'un intervalle $[a,b]$, $a<0<b$,  $r>\max(\norm a,b)$ et
$c\in\C$ un paramètre additionnel.  
\Apriori{} les exemples étudiés concernent des fonctions $P$ réelles pour des arguments réels de $x,y,\eps$, mais pour des fonctions $P$ à valeurs complexes cela ne changer rien.
On fait l'hypothèse 
que $P(0,0,\eps)=\O(\eps^2)$ et 
que la valuation  homogène de
$P(x,y,0)$  est au moins 4, \ie il existe des $p_{kl}\in\C$
tels que
\eq{hypuj}{
P(x,y,0)=\sum_{k+l\geq4}p_{kl}x^ky^l,\mbox{\ \ 
pour $\norm x,\norm y$ assez petits.}
}
L'ensemble lent de \rf{uj}, d'équation $y(y-x)(y+x)+P(x,y,0)=0$, peut être désingularisé
par un éclatement $y=xz$. On obtient l'équation $z(z-1)(z+1)+xQ(x,z)=0$ avec
$Q(x,z)=x^{-4}P(x,xz,0)$ analytique dans $D_1\times D(0,1+\d),\ \d>0$,
à laquelle on peut appliquer le théorème des fonctions implicites aux points 
$(0,0)$ et $(0,\pm1)$.
\vspace{-0.4cm}
\figu{f7.2}{\epsfxsize4.4cm\epsfbox{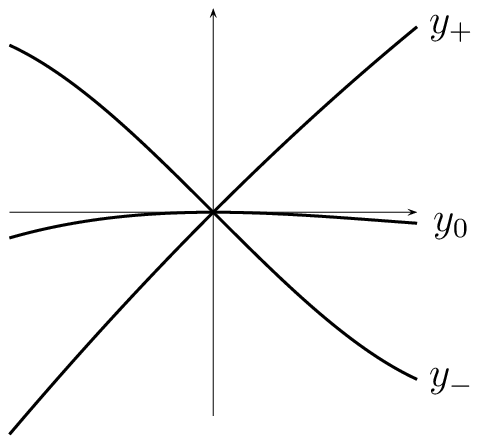}\vspace{-0.7cm}}
{Les branches $y_-$, $y_+$ et $y_0$}
On a ainsi localement trois courbes lentes analytiques $y_0(x)={\cal O}(x^2)$
et $y_\pm(x)=\pm x+ {\cal O}(x^2)$.
L'équation \rf{uj} est dite {\sl de type Union Jack}
car l'équation \og modèle \fg\ $\eps y'=y(y-x)(y+x)+\eps c$ a pour ensemble lent la réunion des trois droites $y=0$ et $y=\pm x$, et ressemble donc au drapeau du Royaume-Uni. 
Les trois courbes lentes $y_0$, $y_+$ et $y_-$ du cas général forment ainsi un \og Union Jack modifié\fg.

Par ailleurs, les droites $y=0$ et $y=\pm x$ sont aussi des solutions particulières de l'équation modèle pour certaines valeurs de $c$~: la solution $y\equiv0$ pour $c=0$, $y\equiv x$ pour  $c=1$ et $y\equiv-x$ pour $c=-1$. 
Il en est de même pour l'équation intérieure réduite de \rf{uj}, obtenue en
posant $x=\e X$, $y=\e Y$, $\eps=\e^3$ et en faisant tendre $\e$ vers $0$, \ie
\eq{UJ}{Y'=Y(Y-X)(Y+X)+c.}
\figu{f7.3}{\epsfxsize3.7cm\epsfbox{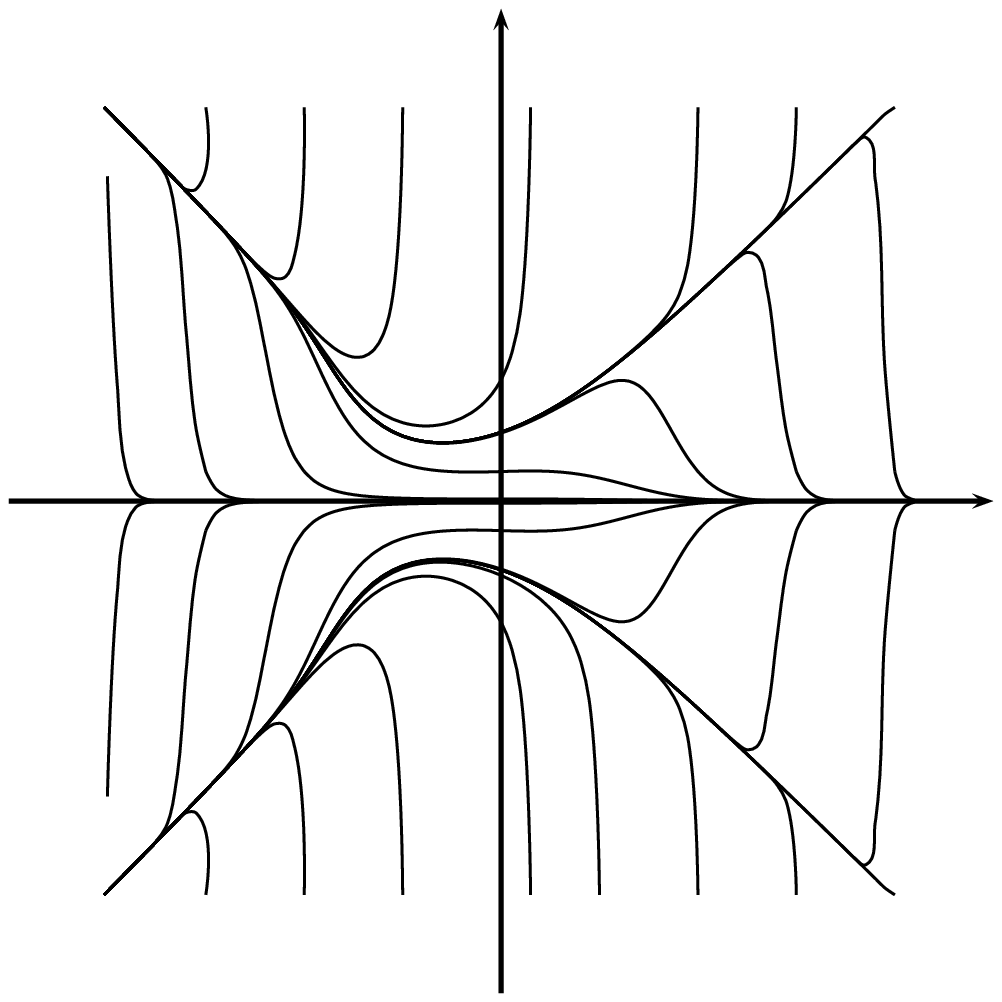}\hspace*{0.1cm}
            \epsfxsize3.7cm\epsfbox{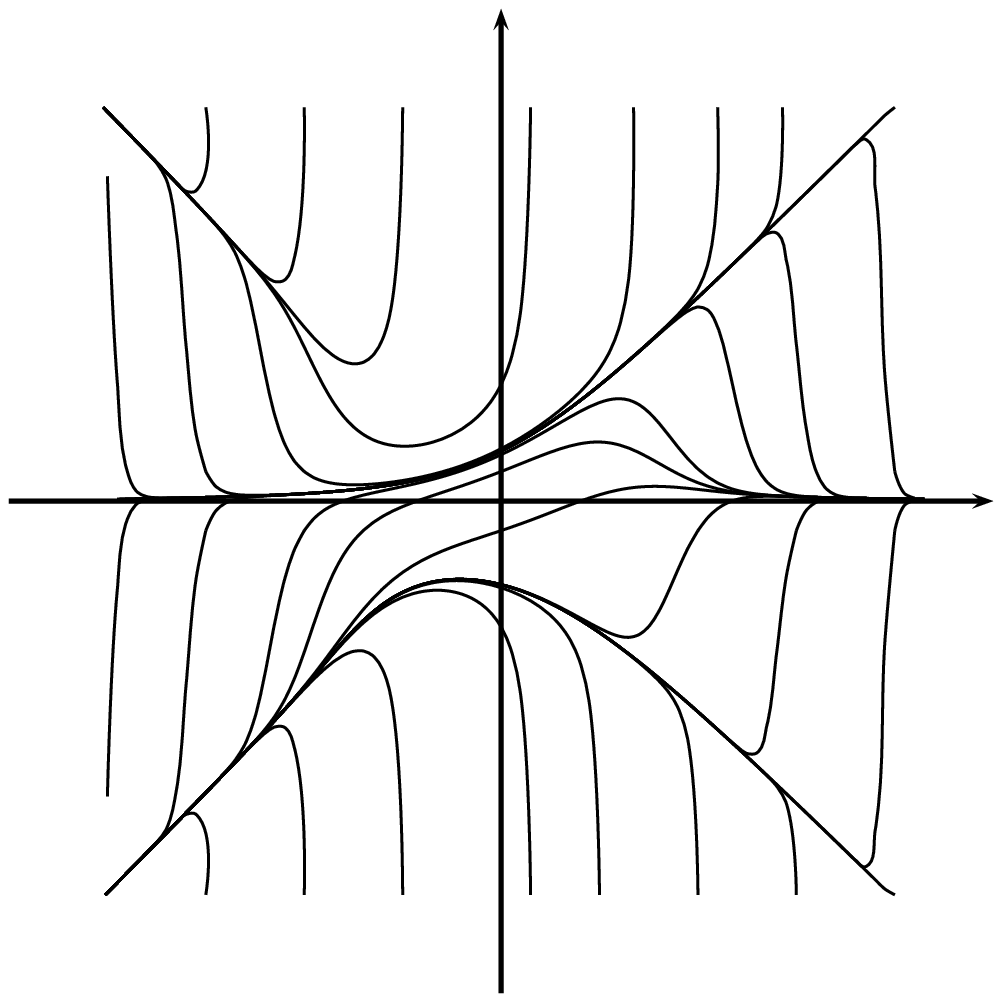}
             \hspace*{0.1cm}\epsfxsize3.7cm\epsfbox{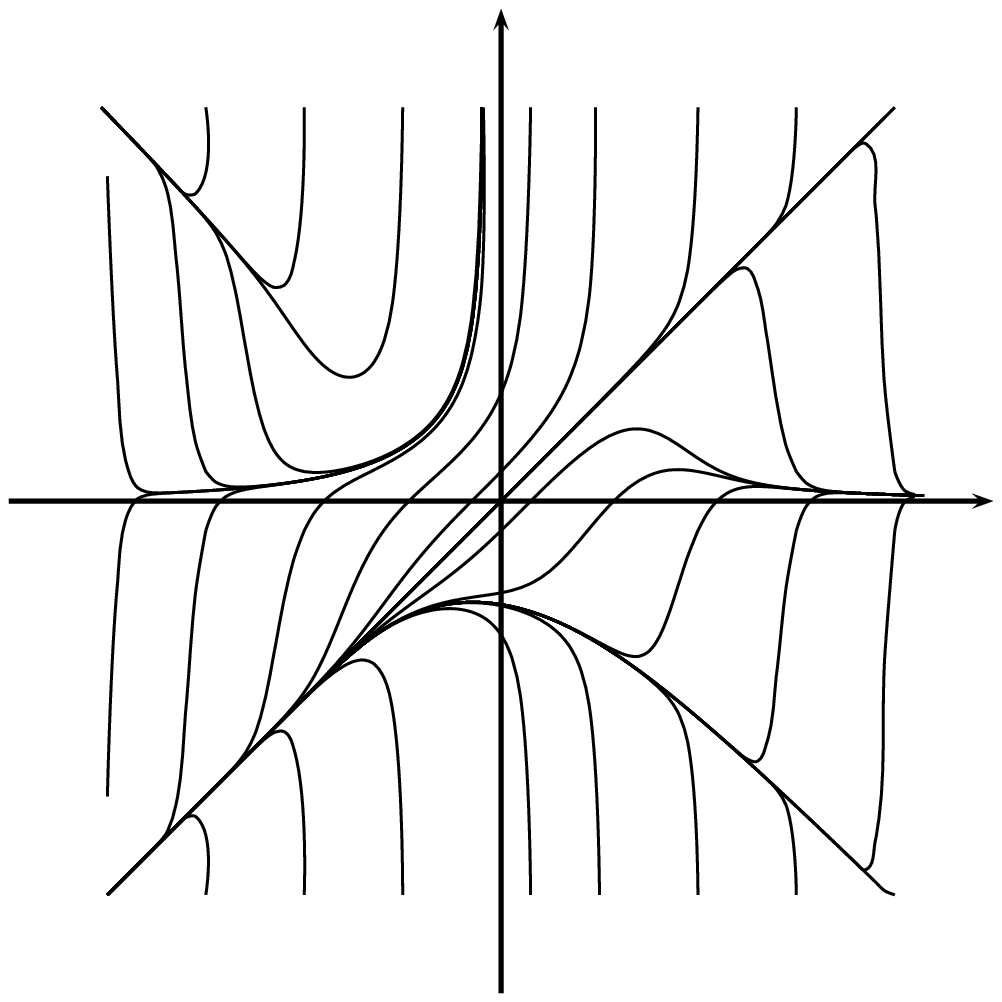}}
{Trois portraits de phase de l'équation \rf{UJ}, pour $c=0$, $c=0.3621759411$ et $c=1$}
Quelque soit la valeur de $c$, l'équation \rf{UJ} admet une unique solution
$Y_g(X,c)$ telle que $Y_g(X,c)\to0$ quand $X\to-\infty$~; en effet, cette équation est de la forme
$Y'=-X^2 Y + {\cal O}(1+\norm {X})$ quand $Y$ reste bornée. Pour une raison analogue,
il existe, pour toute valeur de $c$, deux solutions uniques $Y_d^\pm(X,c)$ 
telles que
$Y_d^\pm(X,c)\sim \pm X$ quand $X\to+\infty$. On vérifie que $Y_0(X,c)={\cal O}(X^{-2})$
quand $X\to -\infty$ et $Y^\pm_d(X,c)=\pm X + {\cal O}(X^{-2})$ quand $X\to\infty$ ;
ceci est le cas uniformément pour des compacts en $c$.

Il a été démontré dans \cit{di} que \rf{UJ} admet une valeur unique $c=c_0\in]0,1[$, telle que
$Y_g$ et $Y^+_d$ coïncident, \ie $Y_g(X,c_0)\equiv Y_d^+(X,c_0)\equiv: Y_0(X)$. 
Cette valeur de $c$ est une valeur
à long canard pour l'équation correspondante $\eps y'=y(y-x)(y+x)+\eps c$ car
la solution $y(x,\eps)=Y_0(\tfrac x\e)$ est attractive quand $x\leq\d<0$ et répulsive
quand $x\geq\d>0$ ; il s'agit d'un canard {\sl non lisse} car la limite uniforme
$z(x)$ de $y(x,\eps)$ quand $\eps\to 0$ est $z(x)=0$ quand $x\leq0$ et 
$z(x)=x$ quand $x>0$ : la courbe lente associée n'est pas dérivable.
Par ailleurs, la valeur $c=-c_0$ est aussi une valeur à canard non lisse~: par symétrie les solutions $Y_g$ et $Y^-_d$ coïncident pour cette valeur de $c$. 

Le résultat qui suit répond à la question naturelle si ce phénomène persiste pour l'équation complète \rf{uj} et si on peut décrire les valeurs à canards et les solutions 
canards correspondantes. Cette question a été résolue
par M.~Diener  \cit{di} et E.~Isambert  \cit{i} à l'exception 
d'approximations uniformes des solutions canards et du caractère Gevrey des développements asymptotiques. Notre théorie des \dacs{} se révèle particulièrement bien adaptée à ce contexte.
\theo{t-uj}{
Avec les hypothèses et notations précédentes, 
on suppose que la courbe lente $y_0(x)$ peut être prolongée analytiquement sur $[a,\d]$  
et la courbe lente $y_+(x)$ sur $[-\d,b]$ avec un certain $\d>0$. 
On suppose que le prolongement, encore noté $y_0(x)$, est attractif, \ie 
$3y_0(x)^2-x^2+\frac{\partial P}{\partial y}(x,y_0(x),0)<0$ sur $[a,0[$, 
tandis que le prolongement encore noté $y_+(x)$ est répulsif sur $]0,b]$.

Alors l'équation \rf{uj} admet une valeur à {\sl canard
non lisse} $c=c(\e)$ et  une solution canard $y(x,\e)$ correspondante 
telles que $y(x,\e) - z(x)={\cal O}(\e)$, où $z(x)=y_0(x)$ quand $x\leq0$ et 
$z(x)=y_+(x)$ quand $x>0$.

De plus, la fonction $c=c(\e)$ admet un développement asymptotique Gevrey d'ordre
${\frac13}$ de la forme 
$$
c(\e)\sim_{\frac13}\sum_{n=0}^\infty c_n \e^n
$$ 
avec pour premier terme la valeur $c_0$ introduite précédemment.
De même, la fonction $y$ admet des \dacs{} Gevrey d'ordre ${\frac13}$
\eq{uj-dac-l}{
y(x,\e)\sim_{\frac13} y_0(x)+\sum_{n=1}^\infty \Big(a_{gn}(x)+b_{gn}\big(\tfrac x\e\big)\Big)
   \e^n 
}
quand $\e\to0$ et $x\in[a,L\e]$ pour $L>0$ arbitraire et
\eq{uj-dac-r}{y(x,\e)\sim_{\frac13} y_+(x)+ \sum_{n=1}^\infty \Big(a_{dn}(x)+b_{dn}\big(\tfrac x\e\big)\Big)
    \e^n }
quand $\e\to0$ et $x\in[-L\e,b]$ pour $L>0$ arbitraire où les
$a_{gn}$ sont analytiques sur un voisinage (complexe)
 de $[a,\d]$, les $a_{dn}$ sur un voisinage
de $[-\d,b]$ et où les $b_{gn}$ et les $b_{dn}$ sont analytiques sur un voisinage de
$\R$. Les fonctions $b_{gn}$ admettent des développements Gevrey d'ordre ${\frac13}$ compatibles
au sens de \rf{gnm} quand $X\to-\infty$, les $b_{dn}$ quand $X\to+\infty$.

Enfin, l'énoncé analogue est vrai pour $y_-(x)$ à la place de $y_+(x)$.   }

\rqs 
1.\ 
La série asymptotique $\hat c(\e)$ de $c(\e)$ et 
les séries formelles combinées du théorème 
se calculent comme avant à partir des développements extérieurs et intérieurs.  
Néanmoins, ici il faut commencer par le développement intérieur pour pouvoir déterminer $\hat c(\e)$~; ceci a été fait par E.~Isambert dans \cit{i}.  
Remarquons que les deux développements intérieurs doivent coïncider, \ie
$y(\e X,\e)\sim\sum_{n=1}^\infty u_n(X)\e^n=:\e\hat Y(X,\e)$ avec une solution formelle
de l'équation intérieure $\hat Y(X,\e)$ avec $c=\hat c(\e)$
dont les  coefficients $u_n(X)$ 
ont une croissance polynomiale quand $X\to+\infty$ et $X\to-\infty$. On montre
(comme dans \cit{i}) que ceci détermine uniquement 
les valeurs des $c_n$ et les fonctions $u_n(X)$. Par exemple, on a
la relation $b_{g1}(X)=u_1(X)=Y_0(X)=X+b_{d1}(X)$.

Les parties lentes des \dacs~\rf{uj-dac-l} et \rf{uj-dac-r} sont ensuite déterminées de 
façon usuelle comme étant les parties non polaires des solutions formelles extérieures de
\rf{uj} avec $c=\hat c(\e)$ pour les courbes lentes $y_0(x)$, respectivement $y_+(x)$.
Ceci implique par exemple $a_{gn}=0$ et $a_{dn}=0$ pour $n=1$ et $2$.
\med\\
2.\ 
Comme cela est classique dans les problèmes de canards, il n'y a pas unicité des valeurs à canard $c(\e)$ ni des solutions canards correspondant à une valeur à canard~: 
une modification de $c(\e)$ ou de la condition 
initiale par un terme exponentiellement petit ${\cal O}(e^{-K/\norm\e^p})$ avec $K>0$
suffisamment grand ne change pas la conclusion  du théorème. 
 Réciproquement, 
  deux valeurs à canard sont exponentiellement proches. 
Ceci peut se démontrer 
 avec la  variation de la constante, \cf par exemple \cit{bcdd, 
 bfsw, crss}, ou encore \cit{fs3} et les références qui y sont incluses.
\med\\
\pr{Preuve} Sur un voisinage de $0$, disons $\norm x < \g$, on fait le changement de variable
$y=y_0(x)+z$. L'équation obtenue s'écrit d'abord sous la forme
\begin{eqnarray*}
\eps z'&=&Q(x,z,\eps)+\eps c\\
&:=&(z+y_0(x))(z+y_0(x)-x)(z+y_0(x)+x)+\\&&P(x,z+y_0(x),\eps)
  -\eps y_0'(x) + \eps c.
  \end{eqnarray*}
Par construction, on a $Q(x,0,0)\equiv0$ ; la fonction $\~P:(x,z,0)\mapsto Q(x,z,0)-z(z-x)(z+x)$ satisfait une propriété analogue à \rf{hypuj}. 
L'équation précédente peut donc être écrite
\eq{ujmod}{\eps z'= r(x)z + \eps (c+s(x,\eps)) + z R(x,z,\eps) ,}
où on a décomposé $Q(x,z,\eps)= r(x)z +\eps s(x,\eps)+ z R(x,z,\eps)$, \ie  $r(x)=\frac{\partial Q}{\partial z}(x,0,0)=
-x^2+{\cal O}(x^3)$ et 
$s(x,\eps)=\frac1\eps Q(x,0,\eps)$
 satisfait $s(0,0)=0$.
D'après notre hypothèse sur $P$ et l'observation précédente pour $\~P$, 
la fonction $R$ peut donc être développée pour  $\eps=0$
$$R(x,z,0)=z^2+\sum_{k\geq0,l\geq1,k+l\geq3}R_{kl}x^kz^l ,$$
quand $\norm x<\g$, $\norm z$ petite.

Après l'homothétie $x\to-\sqrt[3]3x$,
l'équation \rf{ujmod} entre donc dans le cadre du corollaire \reff{co616}
avec $p=3$ et $r=1$.
De plus, son équation intérieure réduite est $Z'=Z(Z-X)(Z+X)+c$. 
La solution de cette dernière équation dont le  comportement asymptotique est de la forme $ \mbox{const.}\,X^{-2}$ 
quand $X$ tend vers $-\infty$ est $Y_g$. 
Quand $c=c_0$, celle-ci coïncide avec $Y_d^+$~;  en particulier, elle peut être
prolongée analytiquement sur $\R$. Il existe donc un voisinage $\norm{c-c_0}<\rho$ tel
que $Y_g$ peut être prolongée sur $[-\infty,M]$ avec un certain $M>0$.
Maintenant on applique le corollaire \reff{co617}. On obtient que \rf{ujmod}
admet une solution $z=z(x,c,\e)={\cal O}(\e)$ ayant un \dac{} quand $\e\to0$,
$x\in V(-\b,\b,\g,L\e)$ et $\norm{c-c_0}<\rho$ avec un certain $L>0$. 
Or il est bien connu que
notre hypothèse d'attractivité de $y_0(x)$ sur $[a,0[$ entraîne que la solution de
\rf{ujmod} avec condition initiale $0$ en un point un peu avant $a$ admet
un développement asymptotique Gevrey d'ordre 1 en $\eps$ sans terme constant
sur l'intervalle $[a,-\g/2]$, disons. Elle est donc exponentiellement proche
de la solution $z=z(x,c,\e)$ construite ci-dessus, uniformément pour $x\in[-\g,-\g/2]$
et $\norm{c-c_0}<\rho$,
Ainsi, on obtient l'existence d'une solution $y_g(x,c,\e)$
de \rf{uj} analytique pour $\e\in S(-\a,\a,\e_1)$, $x\in V(\pi-\b,\pi+\b,\norm a,L\e)$ 
et $\norm{c-c_0}<\rho$, pour certains $L,\e_1,\a,\b>0$,
ayant un \dac{}  Gevrey d'ordre ${\frac13}$
\eq{dac-l}{
y_g(x,c,\e)\sim_{\frac13} y_0(x)+
    \sum_{n=1}^\infty \Big(a_{gn}(x,c)+b_{gn}\big(\tfrac x\e,c\big)\Big)\e^n ,
}
où $b_{g1}(X,c)=Y_g(X,c)$.

Comme l'équation obtenue de \rf{uj} par $y=y_+(x)+z$
entre aussi dans le cadre du corollaire \reff{co617}, on obtient en utilisant ici
la répulsivité de $y_+$ sur $]0,b]$ l'existence d'une
solution $y_d(x,c,\e)$ de \rf{uj} analytique  pour $\e\in S(-\a,\a,\e_1)$, 
$x\in V(-\b,\b,b,$ $L\e)$ et $\norm{c-c_0}<\rho$ admettant un \dac{} Gevrey d'ordre 
${\frac13}$
\eq{dac-r}{y_d(x,c,\e)\sim_{\frac13} y_+(x)+
    \sum_{n=1}^\infty \Big(a_{dn}(x,c)+b_{dn}\big(\tfrac x\e,c\big)\Big)\e^n ,}
où $b_{d1}(X,c)+X=Y_d^+(X,c)$.

On obtient des valeurs à canards non lisses du théorème en résolvant l'équation
\eq{eqc}{
y_g(0,c,\e)=y_d(0,c,\e).
} 
Il faut donc montrer que la solution $c=c(\e)$ existe
et qu'elle a  les propriétés énoncées, ainsi que les fonctions $y_{g|d}(x,c(\e),\e)$.

On applique le théorème des fonctions implicites à l'équation \rf{eqc} modifiée en 
\eq{eqcmod}{
f(c,\e)=0,\mbox{ où }f(c,\e)=\tfrac1\e (y_g(0,c,\e)- y_d(0,c,\e)).
} 
On a d'abord l'égalité
$\lim_{\e\to0}\tfrac1\e y_g(0,c_0,\e)=Y_0(0)=\lim_{\e\to0}\tfrac1\e y_d(0,c_0,\e)$
et donc $\lim_{\e\to0}f(c_0,\e)=0$.
On montrera que
\eq{neqc}{\lim_{\e\to0}\tfrac1\e\frac{\partial y_g}{\partial c}(0,c_0,\e)\neq 
   \lim_{\e\to0}\tfrac1\e\frac{\partial y_d}{\partial c}(0,c_0,\e)}
et donc $\lim_{\e\to0}\frac{\partial f}{\partial c}(c_0,\e)\neq0$.
Les conditions du théorème  des fonctions implicites sont donc satisfaites et on obtient
l'existence d'une solution $c=c(\e)$ de \rf{eqc} avec $c(0)=c_0$. \Apriori{},
ce théorème dit seulement que $c$ est une fonction continue, mais la formule
$$c(\e)=\tfrac1{2\pi i}\int_{\norm {x-c_0}=\rho/2} \frac{x\,
   \tfrac{\partial f}{\partial c}(x,\e)}{f(x,\e)}\,dx$$
en combinaison avec la compatibilité des développements Gevrey avec les 
opérations élémentaires montre qu'elle admet un développement Gevrey d'ordre 
${\frac13}$ quand $\e\to0$ dans un secteur $S(-\a,\a,\e_2)$ avec un certain $\e_2>0$.
À l'aide du théorème \reff{t4.7} (a), on obtient
que les compositions $(x,\e)\mapsto y_g(x,c(\e),\e)$ et 
$(x,\e)\mapsto y_d(x,c(\e),\e)$ admettent des \dacs{} Gevrey d'ordre ${\frac13}$. 
En tant que
solutions de la même équation \rf{uj} avec la même condition initiale
$y_g(0,c(\e),\e)=y_d(0,c(\e),\e)$, elles coïncident. Ceci démontre les énoncés du 
théorème, en particulier \rf{uj-dac-l} et \rf{uj-dac-r} dont les coefficients peuvent 
être obtenus en développant $a_{g|d,n}(x,\hat c(\e))=a_{g|d,n}(x,c_0)+...$ 
respectivement $b_{g|d,n}(X,\hat c(\e))$ par la formule de Taylor dans \rf{dac-l} et \rf{dac-r}.

Pour la démonstration de \rf{neqc}, on utilise que les \dacs{} de $y_{g|d}$ sont 
uniformes par rapport à $c$ dans un voisinage complexe de $c_0$. On peut donc
obtenir les dérivées partielles par rapport à $c$ en utilisant la formule de Cauchy ;
par conséquent les dérivées partielles ont aussi des \dacs{} et ces \dacs{} sont obtenus
en dérivant ceux de $y_{g|d}$ terme à terme. Une comparaison avec le développement
extérieur (voir section \reff{2.3}) montre que $a_{g1}=0$ et $a_{d1}=0$, donc
$\frac{\partial y_g}{\partial c}(0,c_0,\e)=\e Z_g(0,c_0)+{\cal O}(\e^2)$ avec
$Z_g=\tfrac{\partial Y_g}{\partial c}$ et 
$\frac{\partial y_d}{\partial c}(0,c_0,\e)=\e Z_d(0,c_0)+{\cal O}(\e^2)$ avec
$Z_d=\tfrac{\partial Y_d^+}{\partial c}$. 
Or les fonctions $X\mapsto Z_{g|d}(X,c_0)$ sont des solutions de l'équation \rf{UJ} dérivée
par rapport à $c$ prise en $c=c_0$ et pour $Y=Y_g(X,c_0)=Y_0(X)$, 
resp.\ $Y=Y_d^+(X,c_0)=Y_0(X)$. 
Autrement dit, ce sont des solutions de
$$
Z'=(3Y_0(X)^2-X^2)Z+1.
$$
Précisément, $Z_g(X,c_0)$ en est la solution tendant vers
$0$ quand $X\to-\infty$ et $Z_d(X,c_0)$ celle tendant vers 0 quand $X\to+\infty$.
Notons $I_0(X)$ la primitive de $3Y_0(X)^2$ s'annulant en $X=0$ et
$J_0(X)$ celle de $3Y_0(X)^2-3X^2$ s'annulant en $X=0$.
Puisque  $3Y_0(X)^2 ={\cal O}(X^{-1})$
 quand $X\to-\infty$ et $3Y_0(X)^2-3X^2 ={\cal O}(X^{-1})$
 quand $X\to+\infty$, ces primitives ont une croissance au plus logarithmique en $-\infty$, resp. $+\infty$. La formule de variation de la constante donne alors
$$Z_g(0,c_0)=\ds\int_{-\infty}^0 \exp\big(X^3/3-I_0(X)\big)\,dX>0$$ et
$$Z_d(0,c_0)=-\ds\int_0^\infty \exp\big(-2 X^3/3 - J_0(X)\big)\,dX<0.$$
Ceci démontre \rf{neqc} et la preuve du théorème est complète.\ep

Bien entendu, la théorie des \dacs{} n'est pas indispensable pour démontrer l'existence
de valeurs à canards pour \rf{uj} ou des équations similaires, 
ni le fait que le côté droit de l'équation différentielle soit analytique. En effet,
on pourrait utiliser le théorème du point fixe comme dans la preuve du théorème
\reff{t5.3} pour montrer l'existence de $y_g(x,c,\e)$ et de $y_d(x,c,\e)$ pour 
$x\in[a,-L\e]$,  resp.\ $x\in[L\e,b]$, et ensuite les prolonger sur $[a,0]$ resp.\ $[0,b]$ 
en utilisant
l'équation différentielle intérieure. Une valeur à canard $c=c(\e)$ possible est 
alors la solution de l'équation $y_g(0,c,\e)=y_d(0,c,\e)$ dont on assure l'existence
par le théorème des fonctions implicites. 
Le fait crucial qu'une certaine dérivée partielle
ne s'annule pas est alors démontré --- comme ci-dessus --- en se ramenant à
une certaine équation différentielle linéaire pour le paramètre $\e=0$.
De ce point de vue, le problème de l'existence de valeurs à canards traité
dans cette partie est plus simple que, par exemple, celui de conditions nécessaires
et suffisantes pour l'existence de canards d'équations différentielles 
{\sl sans paramètres} traité dans la partie \reff{6.1} ou celui de la résonance
d'Ackerberg-O'Malley traité dans la partie  \reff{6.3}.
Toutefois, la théorie des \dacs{} apporte deux nouvelles propriétés des valeurs à 
canards non lisses : d'une part elle permet de donner une approximation uniforme des solutions canard , précisément notre \dac, et d'autre part les développements  asymptotiques de la fonction $c=c(\e)$ et des solutions canard sont  Gevrey. Ceci peut servir par exemple à démontrer qu'une sommation ``au plus petit terme'' dans l'esprit de ce qui est fait dans \cit{fs} fournit une valeur à canard.
\bigskip

{\noi\bf Equations non analytiques} \sep
\med

Les canards du théorème \reff{t-uj} sont appelés non lisses car la courbe lente associée
n'est pas lisse. En revanche, pour chaque $\eps$, les solutions canards de l'équation analytique \rf{uj} sont
bien sûr analytiques en $x$. Dans la suite nous voulons indiquer comment la théorie des \dacs{}
s'applique aussi à des problèmes de canards pour des équations différentielles non lisses.

Il s'agit des équations du type {\sl canard angulaire}
\eq{c-ang}{
\eps y'=\gk{y-f(x)}\gk{y-g(x)}+\eps c,
}
où $f(x)=\a x+ \O(x^2)$, $g(x)=\b x + \O(x^2)$, $\a\neq\b$ et où les restrictions
de $f$ et $g$ à des intervalles $[0,b]$, resp.\ $[a,0]$ sont des fonctions réelles 
analytiques, mais où $f$ et $g$ ne sont pas supposées $C^\infty$ sur $[a,b]$. 
\figu{f7.4}{
\epsfxsize6cm\epsfbox{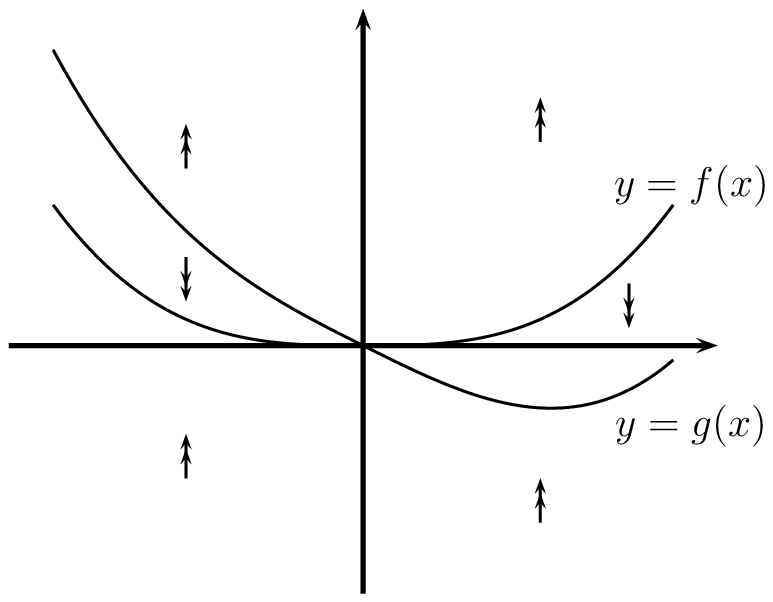}
\vspace{-8mm}
}{Les courbes lentes de \rf{c-ang} et 
l'orientation du champ}
On suppose que la courbe lente $y=f(x)$ est attractive pour $x\in[a,0[$ et
répulsive pour $]0,b]$, \ie  $x(f(x)-g(x))\geq 0$ sur $[a,b]$~;   
en particulier on a $\a>\b$. L'exemple classique de \cit{i} correspond à $f(x)=\norm x^3/3$ et $g(x)=-x+\norm x^3/3$. Dans un autre exemple de \cit{gi,i}, on a
$f(x)=x(1+\norm x)$ et $g=-f$.

Traitons d'abord les solutions sur l'intervalle $[0,b]$.
On écrit l'équation différentielle sur un voisinage complexe de $[0,b]$
avec les prolongements analytiques $f_+,g_+$ 
des restrictions $f\mid_{[0,b]}$ et $g\mid_{[0,b]}$
\eq{c-ang+}{
\eps y'=\gk{y-f_+(x)}\gk{y-g_+(x)}+\eps c.
}
Le changement de variable $y=f_+(x)+z$ mène alors
à l'équation
$$
\eps z'=(f_+(x)-g_+(x))z+z^2+\eps(c-f_+'(x))
$$
qui entre, après une homothétie sur la variable $x$, 
dans le cadre du corollaire \reff{co616} avec $p=2$, $r=1$. Pour un voisinage
$\norm c< \rho$, $\rho$ assez petit,
on peut comme avant utiliser le corollaire \reff{co617} ; la répulsivité permet comme
dans la preuve ci-dessus d'obtenir une solution ayant un \dac{} 
jusqu'à $b$. On obtient qu'il
existe une solution holomorphe $y=y_d(x,c,\e)$ de \rf{c-ang+}
pour $\e\in S(-\a,\a,\e_1)$, $x\in V(-\b,\b,b,L\e)$ et $\norm c< \rho$
avec certains $\a,\e_1,\b,L,\rho>0$ assez petits
et qu'elle admet un \dac{} Gevrey d'ordre $\frac12$
\eq{dac+}{
y_d(x,c,\e)\sim_{\tfrac12}\hat y_d(x,c,\e):=f_+(x)+\sum_{n=1}^\infty 
  (a_{dn}(x,c)+b_{dn}(\tfrac x\e,c))\e^n
}
avec $a_{d1}=0$ et $b_{d1}(X,c)=U_d(X,c)$, où $U_d$ est la solution de l'équation 
intérieure réduite 
\eq{er}{
U'=(\a-\b)XU+U^2+c
}
tendant vers 0 quand $X\to+\infty$. 
On vérifie, comme à la fin de la preuve du
théorème \reff{t-uj}, que $\frac{\partial U_d}{\partial c}(0,0)<0$.

Pour l'intervalle $[a,0]$, l'équation \rf{c-ang} se simplifie aussi en une équation
analytique
\eq{c-ang-}{
\eps y'=\gk{y-f_-(x)}\gk{y-g_-(x)}+\eps c 
}
avec les prolongements analytiques $f_-$ et $g_-$ 
des restrictions $f\mid_{[a,0]}$ et $g\mid_{[a,0]}$ sur un voisinage complexe
de $[a,0]$. On obtient ici l'existence d'une solution holomorphe $y=y_g(x,c,\e)$ de \rf{c-ang-}
pour $\e\in S(-\a,\a,\e_1)$, $x\in V(\pi-\b,\pi+\b,\norm a,L\e)$ et $\norm c< \rho$
avec certains $\a,\e_1,\b,L,\rho>0$ petits
et ayant un \dac{} Gevrey d'ordre $\frac12$
\eq{dac-}{
y_g(x,c,\e)\sim_{\tfrac12}\hat y_g(x,c,\e):=f_-(x)+\sum_{n=1}^\infty 
  (a_{gn}(x)+b_{gn}(\tfrac x\e))\e^n
}
avec $a_{g1}=0$ et $b_{g1}(X,c)=U_g(X,c)$, où $U_g$ est la solution de la même équation 
intérieure réduite \rf{er} tendant vers 0 quand $X\to-\infty$. 
On vérifie que $\frac{\partial U_g}{\partial c}(0,0)>0$.

En appliquant le théorème des fonctions implicites à l'équation
$\frac1\e y_g(0,c,\e)$ $=\frac1\e y_d(0,c,\e)$ au voisinage de $c=0$
comme dans la preuve du théorème \reff{t-uj},
on déduit de nouveau l'existence de valeurs à canards non lisses $c=c(\e)$ ayant
un développement asymptotique $\hat c(\e)$ Gevrey d'ordre $\frac12$ et l'existence
de \dacs{} Gevrey 
pour la solution canard définie par $y(x,\e)=y_g(x,c(\e),\e)$ quand $x\leq0$
et par $y(x,\e)=y_d(x,c(\e),\e)$ quand $x\geq0$ sur les intervalles $[a,0]$ 
resp.\ $[0,b]$.
\med

\Apriori{}, $\hat c(\e)=\sum_{n=1}^\infty c_n \e^n$ est une série formelle en puissances
de $\e$~;
 la théorie des \dacs{} nous apprend déjà que cette série est Gevrey. En complément à cette étude, nous redémontrons et améliorons ci-dessous un
 résultat de \cit{i}.
C'est l'occasion de présenter un \dac{} convergent, le seul dans ce mémoire. Ce \dac~est en fait obtenu par un développement de Taylor d'une fonction spéciale.
\propo{isam+} {
Dans le cas du canard angulaire classique, \ie \rf{c-ang} avec 
$f(x)=\norm x^3/3$ et $g(x)=-x+\norm x^3/3$, la série formelle associée
aux valeurs à canard 
\eq{form-c}{
c(\e)=\sum_{m=1}^\infty c_{4m} \e^{4m}
} 
ne contient que des puissances multiples de 4 et elle converge.
}
\rqs 1.\ E.~Isambert \cit{i} a démontré la forme \rf{form-c} de la série formelle associée
aux valeurs à canards. Nous montrons, de plus, sa convergence.\med\\
2.\ Pour cette équation modèle, on peut aussi exprimer les solutions de façon plus simple
et on obtient des \dacs{} très particuliers comme nous le verrons dans la preuve
et à la fin de cette partie.\med\\
\pr{Preuve}Dans la suite, nous utilisons seulement $\eps$, car la plupart 
des fonctions qui entrent en jeu sont des fonctions de cette variable, et non de $\e=\sqrt\eps$. Il s'agit donc de montrer que la valeur à canard $c=c(\eps)$ est holomorphe dans un voisinage de $0$ et paire.

Considérons d'abord l'équation sur $[0,+\infty[$~:
\eq{c-ang-cl+}{
\eps y'=\gk{y-\tfrac {x^3}3}\gk{y+x-\tfrac{x^3}3}+\eps c.
}
Le changement de variables $y=\frac{x^3}3+z$ la ramène à une équation  simple
$$
\eps z'=xz+z^2+\eps(c-x^2),
$$
qu'on peut encore simplifier. Notons $d=d(\eps)$ la solution de 
$d+d^2 = \eps $ holomorphe dans un voisinage de $\eps=0$ et telle que $d(0)=0.$ Alors
le changement de variables $z=d(\eps)x+u$ mène à l'équation
$$
\eps u'=\big(1+2 d(\eps)\big)xu + u^2 + \eps\big(c-d(\eps)\big).
$$
L'homothétie $x=t/\g(\eps)$, $u=\g(\eps)v$ avec la fonction holomorphe
$\g(\eps)$ vérifiant $\g(0)=1$ et $\g(\eps)^2=1+2d(\eps)$ aboutit à l'équation
$$
\eps \frac{dv}{dt} = t v + v^2 + \eps C(\eps),\ \ \mbox{où}\ \ 
             C(\eps)=(c-d(\eps))/\g(\eps)^{2}.
$$
Pour cette équation, on peut éliminer le petit paramètre $\eps$, sauf dans
l'argument de $C$. En effet, le changement de variables $t=\sqrt\eps\, T$, $v=\sqrt\eps\, V$ mène à 
\eq{reduit}{
\frac{dV}{dT}=TV+V^2+D
}
avec $D=C(\eps)$,
qui est l'équation intérieure réduite hormis le fait que $D$ n'est pas indépendant de $\eps$~; c'est une fonction analytique  $D=C(\eps)$, avec $C(0)=c$.
Pour $D$ arbitraire dans $\C$, notons maintenant $V_d(T,D)$ la solution de \rf{reduit} tendant 
vers $0$ quand $T$ tend vers $+\infty$. Elle peut être exprimée par des solutions
de l'équation de Weber, mais cela n'apporte rien ici. Constatons simplement qu'elle est une fonction
entière de $D$, méromorphe de $T$, qu'elle admet un développement asymptotique
Gevrey d'ordre $\frac12$ uniforme par rapport à $D$ dans tout compact, de la forme
$V_d(T,D)\sim_{\frac12}\sum_{m=0}^\infty W_m(D) T^{-2m-1}$ avec $W_0(D)=-D$ et enfin
que $V_d(T,0)=0$ et $\frac{\partial V_d}{\partial D}(0,0)<0$.

En résumé, si l'on tient compte de tous les changements de variables,  pour 
$\norm c$ petit, \rf{c-ang-cl+} admet une unique solution
$y_d$ telle que $y_d(x,c,\eps)-\tfrac{x^3}3-d(\eps)x$ 
tend vers 0 quand $x$ tend vers $+\infty$. 
Il s'agit de la fonction
$$
y_d(x,c,\eps)=\tfrac{x^3}3+d(\eps)x+\sqrt\eps\g(\eps) 
V_d\gk{\g(\eps)\tfrac x{\sqrt\eps}  ,   \tfrac{c-d(\eps)}{\g(\eps)^{2}}}.
$$
Cette solution est holomorphe pour $\eps\in S(-\a,\a,\eps_1)$, 
$x\in V(-\b,\b,\infty,L\sqrt\eps)$ et 
$\norm c < \rho$ avec
$\a,\b,\eps_1,L,\rho>0$ peut-être petits. 

Sur $]-\infty,0]$, l'équation du canard angulaire classique est
\eq{c-ang-cl-}{
\eps y'=\gk{y+\tfrac {x^3}3}\gk{y+x+\tfrac{x^3}3}+\eps c.
}
De manière 
analogue, et en utilisant que \rf{reduit} ne change pas sous la transformation
$T\to-T,$ $V\to-V$, on obtient que la fonction
$$
y_g(x,c,\eps)=-\tfrac{x^3}3+d(-\eps)x-\sqrt\eps\g(-\eps) V_d\gk{-\g(-\eps)
   \tfrac x{\sqrt{\eps}}  ,
    \tfrac{c-d(-\eps)}{\g(-\eps)^{2}}}
$$
est son unique solution telle que $y_g(x,c,\eps)+\tfrac{x^3}3-d(-\eps)x$ 
tende vers 0 quand $x$ tend vers $-\infty$. Elle est holomorphe 
pour $\eps\in S(-\a,\a,\eps_1)$, $x\in V(\pi-\b,\pi+\b,\infty,L\sqrt\eps)$ et 
$\norm c < \rho$ avec certains $\a,\b,\eps_1,L,\rho>0$.

L'équation déterminant la valeur à canard non lisse $c=c(\eps)$ est donc
\eq{eqc-cl}{
 \g(\eps) V_d\gk{0,\tfrac{c-d(\eps)}{\g(\eps)^{2}}}=-
  \g(-\eps) V_d\gk{0,\tfrac{c-d(-\eps)}{\g(-\eps)^{2}}}. 
  }
À cette équation, on peut appliquer le théorème des fonctions implicites
pour les fonctions holomorphes. On obtient d'abord que $c=c(\eps)$ est une fonction
holomorphe de $\eps$ dans un voisinage de $0$ et ensuite qu'elle est
paire en vertu de la symétrie de l'équation \rf{eqc-cl}. 
Ceci démontre l'énoncé.\ep

\rq 
Non seulement la valeur à canard $c=c(\eps)$, mais aussi les solutions canards $y_g$ et $y_d$, ont des développements convergents en puissances de $\eps$. 
Pour $x\geq0$, la solution canard est donnée par 
$$y(x,\eps)=y_d(x,c(\eps),\eps)=\tfrac{x^3}3+d(\eps)x+
   \sqrt\eps\g(\eps) V_d\gk{\g(\eps)\tfrac x{\sqrt\eps}  ,
    \tfrac{c(\eps)-d(\eps)}{\g(\eps)^{2}}}\ \ $$
avec les fonctions $y_d$ et $V_d$ de la preuve.
Or la fonction $(X,\eps)\mapsto \g(\eps)V_d\big(\g(\eps)X,$ $C(\eps)\big)$ avec 
$C(\eps)=(c(\eps)-d(\eps))/{\g(\eps)^{2}}$ est holomorphe bornée dans
$V(-\a,\a,$ $\infty,L)\times D(0,\eps_1)$, si $\a,L,\eps_1>0$ sont assez petits.
Comme $c(0)=d(0)=0$, on en déduit que la série de Taylor
$$
\g(\eps)V_d(\g(\eps)X,C(\eps))=\sum_{n=1}^\infty W_n(X)\eps^n
$$
converge uniformément sur $V(-\a,\a,\infty,L)$ et sur tout compact de $D(0,\eps_1)$.
Par conséquent, le \dac{} de $y$ est une série convergente en $\e=\sqrt\eps$.
Ce \dac~est constitué d'une partie lente ne contenant que des puissances paires
de $\e$, avec pour terme principal  $x^3/3$ et les autres termes des multiples scalaires de $x$, et d'une partie rapide ne contenant que des puissances impaires
de $\e$ et dont les coefficients sont des fonctions $W_n(X)$ ayant des développements
asymptotiques Gevrey d'ordre $\frac12$ quand $X$ tend vers $+\infty$.
Comme pour la fonction $V_d$, ces derniers développements ne contiennent 
que des puissances impaires de $X^{-1}$.
Les propriétés de la solution canard pour $x\leq0$ sont analogues.
\sub{6.3}{Résonance d'Ackerberg-O'Malley}
On considère l'équation linéaire d'ordre 2
\eq{eao}{
\eps z''-f(x,\eps)z'+g(x,\eps)z=0
}
où $f$ et $g$ sont analytiques dans un voisinage 
complexe
$\U$ de $(0,0)$. Dans  \cit{fs1}, on avait étudié, entre autres, des {\sl solutions $\cc^\infty$-résonantes locales} de \rf{eao}. Ce sont des solutions $z=z(x,\eps)$ définies pour $0<\eps\leq\eps_0$ et $-\delta\leq x\leq \delta$ avec certains $\eps_0,\delta>0$, 
qui tendent vers une solution non triviale de l'équation différentielle réduite
\eq{resred}{-f(x,0)z'+g(x,0)z=0}
quand $\eps\to 0$, uniformément sur $]\!-\d,\d[$, et
telles que toutes les dérivées $z^{(m)}$ sont bornées sur $]\!-\d,\d[$ (uniformément quand $\eps$ tend vers 0).
Nous étudierons aussi des {\sl solutions résonantes locales}, 
\ie des solutions qui tendent vers une solution non triviale de \rf{resred} quand
$\eps\to 0$ uniformément sur $]\!-\d,\d[$, 
mais dont les dérivées ne sont plus nécessairement bornées uniformément par rapport à  
$\eps$. 
Dans cette partie, on fait l'hypothèse
\eq H{
f(x,0)=\a x^{p-1}+\O(x^p)~\mbox{ \sl et }~ 
g(x,0)=\b x^{p-2}+\O(x^{p-1})
}
avec $p$ entier naturel pair, $\a,\b\in\R$ et $\a>0$. 
Le résultat que nous présentons (théorème \reff{reso}) donne des conditions nécessaires et suffisantes, portant sur les solutions formelles de l'équation \rf{eao} et de l'équation intérieure associée, pour l'existence de solutions
résonantes et $\cc^\infty$-résonantes locales. 
Il convient donc d'introduire cette équation intérieure~; on l'obtient avec
le changement de variables $x=\e X$, $Z(X)=z(\e X)$ et $\eps=\e^p$
\eq{eaoint}{\frac{d^2Z}{dX^2}-\~f(X,\e)\frac{dZ}{dX}+\~g(X,\e)Z=0,}
où $\~f(X,\e)=\e^{1-p}f(\e X,\e^p)$ et
$\~g( X,\e)=\e^{2-p}g(\e X,\e^p)$. 
L'hypothèse \rf H entraîne
que 
$$
\~f(X,\e)=\sum_{n=0}^\infty p_n(X)\e^n,\ 
\~g(X,\e)=\sum_{n=0}^\infty q_n(X)\e^n\ ,
$$
séries convergentes pour $X$ dans un compact de $\R$ et pour $\norm\e$ assez petit, où $p_n,q_n$ sont des polynômes, dont les premiers sont $p_0(X)=\a X^{p-1}$ et $q_0(X)=\b X^{p-2}$.
Avant d'énoncer le résultat, nous décrivons les solutions formelles, dites {\sl extérieures}, resp. {\sl intérieures}, de \rf{eao} et \rf{eaoint}. Les solutions
formelles extérieures de \rf{eao} sont de la forme 
$\hat z(x,\eps)=\sum_{n=0}^\infty z_n(x)\eps^n$ ; leurs coefficients satisfont
la récurrence
$$f(x,0)z_n'-g(x,0)z_n= h_n(x),$$
où $h_0\equiv0$ et où  $h_n(x)$ est le coefficients de
$\eps^n$ dans
 le développement de Taylor de 
$$
\eps \hat z\,''(x,\eps)-\big(f(x,\eps)-f(x,0)\big)\hat z\,'(x,\eps)+
\big(g(x,\eps)-g(x,0)\big)\hat z(x,\eps).
$$
Les fonctions $h_n$ dépendent de $f,g$, de $z_0,...,z_{n-1}$ et 
de leurs dérivées. 
On voit facilement que ces solutions formelles
sont uniques à  un facteur constant%
\footnote{\ 
Ici et dans la suite, le mot ``constant'' signifie constant par rapport à $x$.
} près. Précisément, si $\hat z_0=\sum_{n\geq0} z_n(x)\eps^n$ est une solution formelle non triviale de \rf{eao}, alors les solutions formelles de \rf{eao} sont les séries formelles de la forme $\hat c(\eps)\,\hat z_0$ où
$\hat c(\eps)=\sum_{n=0}^\infty c_n\eps^n$.
Par ailleurs, les coefficients $z_n$
peuvent avoir des singularités en $x=0$, mais sont prolongeables le long de tout chemin
restant  dans le voisinage $\U$ et évitant $0$.

Les solutions formelles intérieures de \rf{eaoint} sont de la forme
$\hat Z(X,\eps)=\sum_{n=0}^\infty U_n(X)\e^n$~; leurs coefficients satisfont
la récurrence
\eq{rec}{
\dfrac{d^2U_n}{dX^2}-\a X^{p-1}\dfrac{dU_n}{dX}+ \b X^{p-2} U_n=H_n(X)\ ,
}
où $H_0\equiv0$ et $H_n(x)$ est le coefficient de $\e^n$
dans  le développement de Taylor de 
$$
\big(\~f(X,\e)-\a X^{p-1}\big)\frac{d\hat Z}{dX}(X,\e)-
    \big(\~g(X,\e)-\b X^{p-2})\big)\hat Z(X,\e).
$$
Les fonctions $H_n$ dépendent de $\~f,\~g$, de $U_0,...,U_{n-1}$ et 
de leurs dérivées. 
Ces solutions formelles intérieures ne sont pas uniques, même
à  un facteur $\hat c(\e)$ constant près, mais dépendent de deux paramètres  $\hat c_1(\e)$ et  $\hat c_2(\e)$. Par contre, les fonctions $U_n$ sont entières car les coefficients de $\~f,\~g$ sont des polynômes et la récurrence est linéaire. Si on demande de plus que les coefficients $U_i$
soient des fonctions à  croissance polynomiale quand $X$ tend vers $+\infty$ (resp.\ $-\infty$),
on obtient des solutions formelles intérieures notées 
$\hat Z^\pm(X,\eps)=\sum_{n=0}^\infty U^\pm_n(X)\e^n$ qui sont uniques à  un facteur
constant près et qui joueront un rôle important dans la suite.

Signalons que les deux types de solutions formelles peuvent 
contenir des termes logarithmiques~:
les solutions formelles extérieures en la singularité $x=0$, les solutions intérieures
dans le comportement asymptotique de leurs coefficients quand $X\to\pm\infty$. Ceci
sera source de complications dans la preuve.
À présent nous sommes en mesure d'énoncer le résultat.
\theo{reso}{
\be[\rm(a)]\item
 Sous l'hypothèse \rf H, l'équation différentielle \rf{eao} admet
des solutions résonantes locales si et seulement si les deux conditions suivantes sont réalisées~:

\bul le quotient $D=\b/\a$
est un entier positif ou nul congru à  0 ou 1 modulo $p$ et

\bul il existe une solution formelle
non triviale $\hat Z(X,\eps)=Z_0(X)+\sum_{n=1}^\infty$ $U_n(X)\e^n$ de \rf{eaoint}
dont les coefficients sont à  croissance polynomiale à  la fois quand $X\to+\infty$ et quand $X\to-\infty$.
\item
Sous l'hypothèse \rf H, l'équation \rf{eao} admet
des solutions $\cc^\infty$-résonantes locales si et seulement si~:

\bul $D$ est un entier positif ou nul congru à  0 ou 1 modulo $p$ et 

\bul il existe une solution formelle non triviale
$\hat z(x,\eps)=\sum_{n=0}^\infty z_n(x)\eps^n$ de \rf{eao}
dont les coefficients sont analytiques au voisinage de $0$.
\med\\
Ceci est le cas si et seulement si 
$\e^{-D}\hat z(\e X,\e^p)=\sum_{n=0}^\infty Z_n(X)\e^n$, 
où $Z_n$ sont des polynômes de degré inférieur ou égal à  $n+D$.
\ee
Dans les deux parties de l'énoncé, $Z_0$ est un polynôme de degré $D$ exactement~; c'est le même si on le choisit unitaire.
}
\rqs 
1.\ 
Pour le lien avec le problème original de la résonance de \cit{aom} et avec la surstabilité, on pourra consulter \cit{fs1}.
\med\\
2.\ 
Dans \cit{fs1}, nous avions  conjecturé la  partie (b) et nous l'avions 
démontrée dans deux cas~: 
d'une part dans le cas $p=2$ et d'autre part dans le cas $p>2$ et $\b=0$, 
\ie $g(x,0)={\cal O}(x^{p-1})$. Dans \cit{m1}, Peter De Maesschalck présente des résultats analogues.  Nous les décrivons après  le théorème \reff{resog}.
\med\\
3.\ 
Dans le cas $\b=0$, l'équation \rf{eao} réduite $f(x,0)z'=g(x,0)z$
n'a pas de point singulier en $x=0$ et donc $Z_0\equiv1$. Par conséquent, le passage à  l'équation de Riccati du début de la preuve ci-dessous n'introduit pas de pôle si $\b=0$. 
Par contre, si $\b\neq0$, les solutions $y^\pm$ de cette équation de Riccati peuvent avoir des
pôles dans $[-L\e,L\e]$ (où $L$ est introduit au début de la preuve), ce qui explique pourquoi 
 nous revenons aux solutions de l'équation
linéaire d'ordre deux correspondante.\med\\
4.\ Quand on fixe $\a,\b$ tels que $D=\b/\a$ est un entier positif congru à 0 ou 1
mod $p$, alors la condition du théorème \reff{reso} (a) est une suite de conditions
polynomiales dans les
coefficients $a_{nm}$ et $b_{nm}$ des séries 
$$f(x,\eps)=\a x^{p-1}+\sum_{m\geq p} a_{0m}x^m+
  \sum_{n\geq1}\sum_{m\geq0}a_{nm}x^m\eps^n$$ respectivement
$$g(x,\eps)=\b x^{p-2}+\sum_{m\geq p-1} b_{0m}x^m+
          \sum_{n\geq1}\sum_{m\geq0}b_{nm}x^m\eps^n.$$
Ceci est démontré de manière analogue à la remarque 5 après le théorème \reff{t6.1}. 
\med\\
\pr{Preuve}
La première étape, suivant une idée de Jean-Louis Callot, est de passer à  l'équation de Riccati correspondante. Ceci se fait en posant $y=\eps z'/z$ et aboutit à 
\eq{ric}{
\eps y'=f(x,\eps)y-\eps g(x,\eps)-y^2.
}
D'après notre hypothèse, le corollaire \reff{co616} peut être appliqué avec $r=p-1$.
En tenant compte de la remarque qui suit ce corollaire, on obtient qu'il existe 
$\e_0,x_0,L>0$ tels que \rf{ric} admet deux solutions $y^\pm(x,\e)$, définies pour 
$\e=\eps^{1/p}\in\,]0,\eps_0^{1/p}]$ et $\pm x \in [L\e,x_0]$, et 
que ces solutions admettent des \dacs{} Gevrey
\eq{dacric}{
y^\pm(x,\e)\sim_\usp g^\pm_{p-1}\big(\tfrac x\e\big)\e^{p-1}+
\sum_{n=p}^\infty \gk{a_n(x)+g^\pm_{n}\big(\tfrac x\e\big)}\e^n
}
quand $\e$ tend vers 0, uniformément sur $\pm[L\e,x_0]$. Ici et dans la suite, nous
combinons deux énoncés pour $x$ positifs resp.\ $x$ négatifs en utilisant le symbole
$\pm$~; les fonctions $a_n$ sont réelles analytiques sur $[-x_0,x_0]$ et les
fonctions $g^{\pm}_n$ sont réelles analytiques sur $\pm[L,+\infty[$. La définition 
de \dac{} Gevrey est l'usuelle définition  \reff{d2.3.2}, 
en remplaçant les secteurs et quasi-secteurs 
par des intervalles.
Remarquons que, d'après \rf{ric}, les dérivées de $y^\pm$ admettent aussi des \dacs{}
Gevrey d'ordre $\usp$,
à  cause de la compatibilité des \dacs{} Gevrey avec les opérations élémentaires ;
ces \dacs{} des dérivées sont donc obtenus en dérivant terme à  terme d'après
la remarque après le lemme \reff{derivasympt}.
Par ailleurs, d'après la remarque 3 qui suit la proposition
\reff{combextint}, les fonctions $a_n$ dans \rf{dacric} sont les mêmes pour 
$y^+$ et pour $y^-$.

En injectant ces \dacs{} dans \rf{ric}, on obtient que $g^\pm_{p-1}$ est l'unique solution de l'équation de Riccati réduite
\eq{redric}{\frac{dY}{dX}=\a X^{p-1}Y-\b X^{p-2}-Y^2}
qui tend vers 0 quand $X$ tend vers $\pm\infty$. On a donc 
$g^\pm_{p-1}(X)=\frac{D}X+\O(X^{-2}),\ X\to\pm\infty$ (avec $D=\b/\a$). Nous utiliserons aussi que 
$g_{p-1}^\pm(X)={Z_0^\pm}'(X)/Z_0^\pm(X)$ où $Z_0^\pm$ est l'unique solution de
l'équation linéaire réduite associée à  \rf{eaoint}
\eq{red}{
\dfrac{d^2Z}{dX^2}-\a X^{p-1}\dfrac{dZ}{dX}+ \b X^{p-2} Z=0
}
telle que $Z_0^\pm(X)\sim_{\usp} (\pm X)^{D}\bigg(1+\ds\sum_{m=1}^\infty
     {d_m}{ X^{-mp}}\bigg)$ quand $X$ tend vers $\pm \infty$.
Ceci facilitera le passage du \dac{} \rf{dacric} à  l'équation linéaire.

Ce passage se fait par $z=\exp\big(\tfrac1\eps\int y\big)$ et 
nécessite donc l'intégration d'un \dac. En utilisant les développements
$g_n^\pm(X)\sim \sum_{m=1}^\infty g_{nm}X^{-m}$ et le fait que $g_{nm}$ est nul lorsque $n+m\not\equiv0\mod p$, on obtient d'après la proposition \reff{p4}
$$
\frac1\eps\ds\int_{\pm r}^x y^\pm(\xi,\e)\,d\xi \sim 
\hat Y^\pm(x,\e)-\hat Y^\pm(\pm r,\e)
$$
avec 
$$
\hat Y^\pm(x,\e)= \log\Big(\e^D Z^\pm_0\big(\tfrac x\e\big)\Big)+
  A_0(x) + \hat R(\eps)\log({x^p+\eps})+
  \ds\sum_{n=1}^\infty \gk{ A_n(x) + 
            G_{n}^\pm\big(\tfrac x\e\big)}\e^{n}
$$
où  $ A_n(x)=\ds\int_{0}^x a_{n+p}(\xi)\,d\xi$, 
$\hat R(\eps)=\frac1p\ds\sum_{l=2}^\infty g_{lp-1,1}\,\eps^{l-1}$
et
$$
G_n^\pm(X)=\ds\int^X_{\pm\infty} 
 \gk{g_{n+p-1}^\pm(T)-g_{n+p-1,1}T^{p-1}(T^p+1)^{-1}}\,dT.
$$
Comme nous l'avions fait dans la preuve de la proposition \reff{p4}, on identifie ici $\hat Y^\pm(\pm r,$ $\e)$ avec la série formelle en puissances de $\e$
obtenue en développant $\log(r^p+\eps)$ par  la formule de Taylor et en utilisant les 
développements asymptotiques à  l'infini de $\log\big(\e^D Z^\pm_0(\tfrac {\pm r}\e)\big)$
 et de $G_n^\pm\big(\frac {\pm r}\e\big)$.

En utilisant la compatibilité des \dacs{} avec la composition 
(ici avec l'exponentielle) et en multipliant par une fonction de $\e$ seulement
ayant l'asymptotique
$\e^{-D}\exp\big(\hat Y^\pm(\pm r,$ $\e)\big)$, on en déduit l'existence de deux solutions ayant une forme de \dac{} généralisé \rf{dacz} ci-dessous. Nous décrivons ces solutions dans un énoncé séparé car il peut être utile pour l'étude d'équations linéaires du deuxième ordre, indépendamment du problème de la résonance.
\ft
\propo{p4.7}{
Sous l'hypothèse \rf H, il existe $x_0,L>0$, une série formelle 
$\hat R=\hat R(\eps)$ Gevrey d'ordre 1 sans terme constant,
des fonctions $B_n$ réelles analytiques sur $[-x_0,x_0]$ avec $B_0(0)=1$ 
et $H_n^\pm$ réelles analytiques sur $\pm[L,\infty[$ 
et  deux solutions $z^\pm$ de l'équation \rf{eao}, telles que 
\eq{dacz}{
\begin{array}{rl}
z^\pm(x,\e)\!\!\sim_{\usp}& Z_0^\pm\big(\tfrac x\e\big)
 \,e^{\hat R(\eps)\log(x^p+\eps)}\bigg(B_0(x)+
 \\&\ds
 \sum_{n=1}^\infty \gk{B_n(x)+
H_n^\pm\big(\tfrac x\e\big)}\e^n\bigg)
\end{array}
}
quand $\e\to0$, uniformément sur $[L\e,x_0]$, resp. $[-x_0,-L\e]$, où
$Z_0^\pm(X)$ sont les uniques solutions de \rf{red} avec 
$Z_0^\pm(X)\sim (\pm X)^{D}\big(1+\O(X^{-1})\big)$ quand $X\to\pm\infty$.

Précisément, la relation \rf{dacz} signifie que, pour toute fonction $R:\;]0,\eps_0]\to\R$
ayant $\hat R(\eps)$ comme développement asymptotique Gevrey d'ordre 1 (par rapport à $\eps$),
la fonction 
$$
e^{-R(\eps)\log(x^p+\eps)}z^\pm(x,\e)/Z_0^\pm\big(\tfrac x\e\big)\sim_{\usp}
  B_0(x)+\sum_{n=1}^\infty \gk{B_n(x)+H_n^\pm\big(\tfrac x\e\big)}\e^n
$$
admet un \dac{} Gevrey d'ordre $\usp$.
}
\rqs 
1.\  
En termes des fonctions $A_n$ et $G_n^\pm$ précédentes,
les fonctions $B_n$ et $H_n^\pm$ sont définies par $B_0(x)=\exp(A_0(x))$ et par
\eq{bnhn}{
\begin{array}{rcl}\ds
B_0(x)\!\!\!\!&+&\!\!\ds\sum_{n=1}^\infty \Big(B_n(x)+
H_n^\pm\big(\tfrac x\e\big)\Big)\e^n
\\
&=&\!\! B_0(x)\,\exp\bigg( \ds\sum_{n=1}^\infty \gk{ A_n(x) 
+ G_{n}^\pm\big(\tfrac x\e\big)}\e^{n}\bigg).
\end{array}
}
2.\ 
L'unicité des \dacs{} à  gauche et à  droite pour l'équation de Riccati implique que 
les développements 
\rf{dacz} à  gauche et à  droite de solutions de \rf{eao} sont uniques 
à  une constante multiplicative $\hat c(\e)$ près.
\med\\
3.\  
De même que pour les fonctions $a_n$ dans \rf{dacric} concernant  $y^\pm$  et les fonctions $A_n$ pour $\hat Y^\pm$ ainsi que l'asymptotique des fonctions $H_n^\pm$ quand 
$X\to\pm\infty$, les fonctions $B_n(x)$ sont les mêmes pour $z^+$ et $z^-$.
\med\\
4.\ 
Pour $x$ fixé, $\e^D Z_0^{\pm}\big(\tfrac x\e\big)$ et les deux autres facteurs des
développements \rf{dacz} de $z^\pm$ admettent un développement asymptotique
en puissances de $\e^p=\eps$ quand $\e\to0$. C'est la raison pour laquelle nous avons
choisi $\ell(x,\e)=\log(x^p+\e^p)$ en appliquant la proposition \reff{p4}.
\med\\
5.\
Si $\re(D)>0$, alors $\e^D z^\pm(x,\e)$ tend vers $x^D B_0(x)$ uniformément
sur $\pm[L\e,x_0]$ quand $\e\to0$. En effet, comme produit de deux fonctions
admettant des \dacs{}, $e^{-R(\eps)\log(x^p+\eps)}\gk{\frac x\e}^{-D}z^\pm(x,\e)$
admet pour \dac{} Gevrey d'ordre $\usp$ la série formelle obtenue en développant
$$
\big(\tfrac x\e\big)^{-D}Z_0^\pm\big(\tfrac x\e\big) \bigg(B_0(x)+\sum_{n=1}^\infty \gk{B_n(x)+
H_n^\pm\big(\tfrac x\e\big)}\e^n\bigg)
=B_0(x)+G_0^\pm\big(\tfrac x\e\big)+{\cal O}(\e)
$$
avec une certaine fonction $G_0^\pm(X)$ ayant un développement asymptotique 
sans terme constant quand $X\to\pm\infty$.
En multipliant avec $x^{D}$ et en utilisant  $x^DG_0^{\pm}(\frac x \e)=
{\cal O}(\e^{\min(1,\re D)})$, on obtient le résultat voulu.
\med\\
6.\ 
Une comparaison avec les développements présentés dans la suite montre que les produits
$Z_0^\pm(X)H_n^\pm(X)$ sont analytiques dans un voisinage de l'axe réel~; puisque 
$Z_0^\pm$ peut avoir des zéros, cela n'empêche pas que les $H_n^\pm$ puissent avoir 
des singularités, même sur l'axe réel, mais ces singularités sont nécessairement des pôles. 
En revanche, ces produits $Z_0^\pm(X)H_n^\pm(X)$ 
 peuvent avoir d'autres singularités en dehors de l'axe réel. Ceci 
provient du fait qu'on a mis en facteur une série contenant 
$\log(X^p+1)$ qui a des points de ramification en les zéros de $X^p+1$.
\med\\
7.\ 
Comme avant, on peut aussi dériver ces développements terme à  terme.
\med

\pr{Suite de la preuve du théorème \reff{reso}}
Sur des intervalles de la forme $\pm[L,M]$ avec un certain $M>0$, 
on a donc $Z^\pm(X,\e):=z^\pm(\e X,\e)
=Z_0^\pm(X)+{\cal O}(\e)$. Or $Z_0^\pm$ est solution de \rf{eaoint}
qui est régulièrement perturbée et se réduit à  \rf{red} quand $\e=0$. 
Le prolongement
de $z^\pm$ sur $[-M\e,M\e]$ satisfait donc $z^\pm(x,\e)=
Z_0^\pm\big(\tfrac x\e\big)+{\cal O}(\e)$ d'après le théorème de dépendance
analytique d'équations différentielles par rapport aux paramètres.
L'existence d'une solution résonante locale de \rf{eao} nécessite donc d'abord que
$Z_0^+$ soit proportionnelle à $Z_0^-$. Nous montrerons plus bas le lemme qui suit.
\ft
\lem{z0}{
La fonction $Z_0^+$ est proportionnelle à $Z_0^-$ 
si et seulement si $D=\b/\a$ est un entier positif 
congru à  0 ou 1 modulo $p$. De plus, dans ce cas, $Z_0:= (\pm1)^D Z_0^\pm$ est un
polynôme de degré $D$.
}
\fp
La condition sur $D$ est supposée désormais.
Alors, en utilisant la proposition \reff{matching-gevrey}, la relation
\rf{dacz} entraîne les développements intérieurs suivants pour $\pm X\in[L,M]$
\eq{devint}{
\begin{array}{rl}
z^\pm(\e X,\e)e^{-\hat R(\eps)\log\eps} \sim_{\usp}\!\!&\ds
e^{\hat R(\eps)\log(X^p+1)}\bigg(Z_0(X)+
\\&\ds
\sum_{n=1}^\infty\gk{Z_n(X)+K_{n}^\pm(X)}\e^n\bigg)
\end{array}
}
avec des polynômes $Z_n$ de degré au plus $n+D$ et des fonctions $K_n^\pm\in\G_L^\pm$, où
$\G_L^\pm$ est l'ensemble des fonctions réelles analytiques sur $[L,+\infty[$ ayant un 
développement asymptotique sans terme constant quand $X$ tend vers $+\infty$.
 
 Précisément, notons $P_n$ la partie polynomiale du produit $Z_0H_n^\pm$~; son degré est au plus $D-1$ puisque $H_n^\pm(X)=\O\big(X^{-1}\big)$ lorsque $X\to\pm\infty$.
On a alors $K_n^\pm=Z_0H_n^\pm-P_n$.
Le polynôme $Z_n$, quant à  lui, est obtenu en développant les fonctions $X\mapsto B_k(\e X)\e^k,\,k\leq n$ par la formule de Taylor et en ne retenant que les termes en $\e^n$, puis en multipliant par $Z_0$, et enfin en ajoutant $P_n$. En d'autres termes
$$
Z_n(X)=P_n(X)+Z_0(X)\,\sum_{k=0}^n\ts\frac{B_{n-k}^{(k)}(0)}{k!}X^k.
$$
En développant l'exponentielle du côté droit de \rf{devint}, on obtient donc des fonctions
$U_n^\pm\in\G_L^\pm[X,\log(X^p+1)]$
(\ie polynomiales en $X$ et $\log(X^p+1)$, avec des coefficients dans $\G_L^\pm$) 
telles que
\eq{devint2}{z^\pm(\e X,\e)e^{-\hat R(\eps)\log\eps} \sim_{\usp}
 Z_0(X)+ \sum_{n=1}^\infty U_{n}^\pm(X)\e^n}
quand $\e\to0$ uniformément pour $\pm X\in [L,M]$. Or la série formelle dans \rf{devint2}
est solution formelle de \rf{eaoint},  donc ses coefficients peuvent être prolongés en des fonction entières.
Par construction, la croissance de $U^\pm(X)$ quand $X$ tend vers $\infty$
est polynomiale. Puisque \rf{eaoint} est régulièrement perturbée (et linéaire), 
on peut appliquer le théorème de dépendance analytique et la proposition \reff{p4.5}. On
obtient que \rf{devint2} reste valable pour $X\in[-M,M]$, ainsi que \rf{devint}.
Comme auparavant, on peut dériver \rf{devint} et \rf{devint2} 
terme à  terme par rapport à  $X$.
\med

Supposons d'abord que \rf{eao} admette une solution résonante locale, notée $z(x,\e)$, 
définie et bornée quand $-\delta\leq x\leq\delta$ et $0<\e\leq\e_0$.
Alors, quitte à diminuer $\d$, la fonction $y:(x,\e)\mapsto\eps z'(x,\e)/z(x,\e)$ est une solution de \rf{ric}
de conditions initiales en $\pm\delta$ bornées. 
Elle admet donc des \dacs{} \rf{dacric} quand
$0<\e\leq\e_0$ et $\pm x\in[L\e,\delta]$ avec un certain $L>0$ et des $\e_0,\delta$ diminués.
À un facteur constant près, 
les fonctions $\exp\!\big\{\tfrac1\eps\int_{\pm\delta}^xy(\xi)\,d\xi\big\}$ 
admettent donc les \dacs{} généralisés \rf{dacz} et aussi les développements
intérieurs \rf{devint2}. Par construction, ces fonctions sont proportionnelles à la solution $z$ donnée, donc chacune d'elle est proportionnelle à l'autre. 
Le quotient des 
deux solutions doit donc avoir un développement asymptotique en série de puissances 
de $\e$, que nous notons $\hat c(\e)$. 
Par conséquent la série formelle $Z_0(X)+\sum_{n=1}^\infty U^+_n(X)\e^n$ 
doit être égale à  $\big((-1)^D Z_0(X)+\sum_{n=1}^\infty U^-_n(X)\e^n\big)\hat c(\e)$. 
Ceci implique
que les coefficients des deux séries ont une croissance polynomiale dans les deux
directions $+\infty$ et $-\infty$, ce qui montre la nécessité pour la première partie du théorème.
\med

Pour montrer qu'il s'agit d'une condition suffisante, on utilise l'existence
des deux solutions satisfaisant \rf{dacz} et donc aussi \rf{devint} et \rf{devint2} pour $X\in[-M,M]$.
L'existence d'une solution formelle de \rf{eaoint} dont les coefficients sont à  croissance polynomiale quand $X$ tend vers $+\infty$ et $-\infty$ entraîne que la série formelle $Z_0(X)+\sum_{n=1}^\infty U^+_n(X)\e^n$ 
doit être proportionnelle à $Z_0(X)+\sum_{n=1}^\infty (-1)^D U^-_n(X)\e^n$. Or les coefficients lents $B_n(x)$ sont les mêmes dans \rf{dacz}, donc le quotient doit être égal à  1.
À ce point, on utilise  que nos \dacs{} sont Gevrey; en particulier \rf{devint2}.
De plus, ce sont {\sl les mêmes développements} pour $z^+$ et $z^-$ et en particulier
ces deux fonctions et leurs dérivées
ont le même développement asymptotique Gevrey en $x=0$. 
Les différences $z^+(0,\e) - z^{-}(0,\e)$ et ${z^+}'(0,\e) - {z^{-}}'(0,\e)$ sont 
donc exponentiellement petites
(d'ordre $p$ en $\e$, \ie d'ordre 1 en $\eps$). De manière analogue à la 
preuve du théorème \reff{t6.1}, on en déduit par le lemme de Gronwall 
que $z^+$ et $z^-$ sont exponentiellement proches
sur un intervalle $[-\d,\d]$ avec $\d>0$ indépendant de $\e$. 
Les prolongements des solutions $\e^D z^{\pm}(x,\e)$,
sur un intervalle indépendant de $\e$ contenant $0$ 
dans son intérieur, tendent donc vers $x^DB_0(x)$ 
uniformément en $\e$ d'après la remarque 5.
Ce sont donc des solutions résonantes locales.
\med

Supposons maintenant que \rf{eao} admette une solution $\cc^\infty$-résonante locale.
Alors les fonctions $z^+$ et $z^-$ doivent avoir des dérivées bornées (uniformément par rapport à  $\eps$) dans un voisinage de $0$ indépendant de $\eps$. Les fonctions $Z^\pm:(X,\e)\mapsto z^\pm(\e X,\e)$ doivent donc satisfaire ${Z^\pm}^{(m)}(X,\e)={\cal O}(\e^m)$
uniformément sur $\pm[L,M]$ pour tout $m\in\N$.
Si l'un des coefficients de $\hat R$ ou l'une des fonctions $K_n^\pm$
était non nul, et
donc si l'un des $U_n^\pm$ n'était pas un polynôme en $X$, 
alors en dérivant \rf{devint2} terme à  terme, 
on obtiendrait que ${Z^\pm}^{(m)}(X,\e)$ ne pourrait pas être ${\cal O}(\e^m)$ pour $m$ grand. On doit donc avoir $\hat R(\eps)=0$ et $K^\pm_{n}=0$ pour tout $n$. 
Avec cette information, \rf{dacz} entraîne un développement pour $z^\pm$ de la forme
\eq{daczmod}{\e^D\,z^\pm(x,\e)\sim \sum_{n=0}^\infty C_n(x)\e^n := 
    \e^D Z_0\big(\tfrac x\e\big)\sum_{n=0}^\infty B_n(x)+ 
     \sum_{n=1}^\infty\e^DP_n\big(\tfrac x\e\big)\e^n\ \ ;} 
ici $\e^D Z_0\big(\tfrac x\e\big)$ et les $\e^DP_n\big(\tfrac x\e\big)$ sont des polynômes en $x$ et 
$\e$ et $C_0(x)=x^D B_0(x)\not\equiv 0$.

Dans le cas où il existe une solution locale $\cc^\infty$-résonante, 
on a $\e^Dz^\pm(x,\e)\sim 
\sum_{n=0}^\infty C_n(x)\e^n$ et, en particulier, le côté droit est une solution
formelle non triviale de \rf{eao} avec des coefficients analytiques en $x=0$.
Quitte à  multiplier cette solution par une série $\hat c(\eta)$, on peut imposer par exemple $C_0(0)=1, C_1(0)=\dots=C_{p-1}(0)=0$, et on vérifie alors que $C_n=0$ si $n\not\equiv0\mod p$~; c'est donc une solution formelle en puissances de $\eps=\e^p$.
\med

Pour voir qu'il s'agit d'une condition suffisante ici aussi, on utilise à  nouveau
l'existence
des deux solutions satisfaisant \rf{dacz}.
L'existence d'une solution formelle non triviale de \rf{eao} sans singularité 
en $x=0$ implique en remplaçant $x=\e X$, que \rf{eaoint} admet une solution formelle
à  coefficients polynomiaux en $X$.
Puisque les solutions formelles de \rf{eaoint}
$Z_0(X)+\sum_{n=1}^\infty U_n^\pm(X)\e^n$, avec 
$U_n^\pm(X)$ à  croissance polynomiale quand $X$ tend vers l'infini, sont uniques
à  un facteur constant près, les fonctions $U_n^\pm$ sont des polynômes.
Une comparaison avec \rf{devint} montre que $\hat R(\eps)=0$ et que $K_n^\pm\equiv0$ pour 
tout $n$. En utilisant leur définition (en dessous de \rf{devint}), on obtient donc que les fonctions
$z^\pm$ satisfont \rf{daczmod}. Ce développement asymptotique peut être  considéré comme un \dac{} sans partie rapide, et ceci
pour $X\in[L\e,r]$ resp.\  $X\in[-r,-L\e]$ avec un certain $L>0$.
On peut étendre la validité à  $\pm[0,r]$ comme précédemment en utilisant le fait 
que $X\mapsto z^\pm(X\e,\e)$ satisfait une équation régulièrement perturbée sur 
$\pm[0,L]$.
Enfin, on utilise de nouveau que \rf{dacz} et donc
\rf{daczmod} sont  Gevrey. 
Puisque c'est le même développement pour $z^+$ et $z^-$,
on en déduit comme précédemment par Gronwall que les fonctions  $z^\pm$
sont des solutions résonantes locales.
En dérivant terme à  terme le \dac{} \rf{daczmod},
on peut faire le même raisonnement pour toutes leurs dérivées~;
il s'agit donc de solutions $\cc^\infty$-résonantes locales.\ep

\rq
Une preuve de la suffisance pour la deuxième partie comme dans \cit{crss} 
 est possible (voir aussi \cit{fs1}). On montre
d'abord le caractère Gevrey de la solution formelle $\hat z$ 
(à  un facteur $\hat c(\eps)$ près)~; ensuite on construit une quasi-solution de \rf{eao}
et enfin on  montre que la solution avec la même condition initiale en $0$ est une solution
$\cc^\infty$-résonante locale.
\med\\
Nous terminons avec la\!\!\!\!\pr{preuve du lemme \reff{z0}}
Supposons que $Z_0^+=c\,Z_0^-$ 
pour un certain $c\in\R$. Comme $p$ est pair, le changement
de variable $X\mapsto-X$ laisse l'équation \rf{red} inchangée. 
D'après leur définition par l'asymptotique,
on a donc $Z_0^+(-X)=Z_0^-(X)$. Par conséquent on a aussi $Z_0^-=c\,Z_0^+$, donc
$c=\pm1$. 
Dans le cas $c=1$, la fonction $Z_0:=Z_0^-=Z_0^+$ est paire, donc $Z_0'(0)=0$.
Dans le cas $c=-1$, la fonction $Z_0:=Z^-_0=-Z^+_0$ est impaire, donc $Z_0(0)=0$.
D'après la théorie générale des points singuliers irréguliers  des équations linéaires d'ordre 2, 
on a $Z_0(X)=X^D(1+{\cal O}(X^{-1}))$ dans le secteur $\norm{\arg(X)}<\frac{3\pi}{2p}$ 
(correspondant à  une montagne et aux deux vallées adjacentes).

L'équation  \rf{red} admet d'autres symétries. 
Posons $\rho=\exp(2\pi i/p)$~; le changement de variable $X\mapsto\rho X$ laisse \rf{red} 
inchangée.
La fonction $\~Z_0$ définie par $\~Z_0(X)=Z_0(\rho X)$ satisfait donc \rf{red} et
$\~Z_0(X)= \exp(2D \pi i/p) X^D (1+ {\cal O}(X^{-1}))$ quand $\norm{\arg X + 
  \frac{2\pi}p}<\frac{3\pi}{2p}$.

Dans le cas $c=1$, les valeurs initiales en $X=0$ impliquent que $\~Z_0=Z_0$ et donc
l'unicité du comportement asymptotique dans l'intersection $\{X\,\mid\,
\arg X\in[-\frac{3\pi}{2p},-\frac\pi{2p}]\}$ 
des deux secteurs mentionnés ci-dessus
implique que $\exp(2D \pi i/p)=1$ et donc $D$ doit être entier et multiple de $p$.
En utilisant $Z_0(\rho^kX)=Z_0(X)$, $k=0,1,...,p-1$, 
on obtient $Z_0(X)=X^D(1+{\cal O}(X^{-1}))$ quand $\norm X\to+\infty$ quelque soit
$\arg X$. Ceci implique  que $D$ est positif et que $Z_0(X)$ est un polynôme de degré $D$.

Dans le cas $c=-1$, on obtient $\~Z_0=\rho Z_0$ et donc 
$\exp(2D \pi i/p)=\rho=\exp(2\pi i/p)$. Ceci implique que $D-1$ est un multiple de $p$ ;
comme avant on déduit que $D$ est positif et que $Z_0(X)$ est un polynôme de degré $D$.
Ceci démontre la nécessité de la condition du lemme. 
On démontre facilement que \rf{red} admet une solution polynomiale dans les deux cas, qui est forcément
proportionnelle à $Z_0^\pm$~; ainsi la condition est suffisante.
\ep

Les résultats précédents sont de nature locale. 
Comme dans la partie \reff{6.1}, les conditions du théorème \reff{reso}
donnent aussi des résultats globaux.
\theo{resog}{
Soit $a<0<b$ et $f,g$ analytiques dans un voisinage complexe $\cal U$ de 
$[a,b]$.
On fait l'hypothèse \rf{H} et on suppose de plus que 
$f(x,0)$ est réel sur l'axe réel et $xf(x,0)>0$ quand  $x\in[a,b]\setminus\{0\}$.

Alors les conditions du théorème \reff{reso} (a) sont équivalentes à l'existence d'une
{\em solution résonante  globale} de \rf{eao}, \ie une solution qui tend 
vers une solution non triviale de 
l'équation réduite \rf{resred} uniformément  sur $[a,b]$ quand $\eps\to0$.

De même, les conditions du théorème \reff{reso} (b) sont équivalentes à l'existence d'une
{\em solution $\cc^\infty$-résonante  globale} de \rf{eao}, \ie une solution qui tend 
vers une solution non triviale de 
l'équation réduite \rf{resred} uniformément  sur $[a,b]$ et dont toutes les dérivées 
sont bornées uniformément sur $[a,b]$  quand $\eps\to0$.
}
\rq 
 Dans \cit{m1}, Peter De Maesschalck démontre que la condition du  théorème \reff{reso} (b) est suffisante pour l'existence d'une solution $\cc^\infty$-résonante globale, \cf théorème 2. Il démontre aussi, \cf théorème 5 du même article, que l'existence d'une solution résonante locale entraîne celle d'une solution résonante globale.
Par contre, la nécessité de la condition dans  \reff{reso} (b) et la nécessité et suffisance de la
condition dans  \reff{reso} (a) pour l'existence d'une solution résonante locale ou 
globale sont nouvelles.
\med\\
\pr{Preuve} La nécessité des conditions du théorème \reff{reso} (a) pour l'existence 
d'une solution résonante globale est évidente, car elles sont déjà nécessaires
pour l'existence d'une solution résonante locale. 

La preuve de la suffisance est analogue à celle de l'implication (c)$\Rightarrow$(b) 
du théorème \reff{t6.1}.
Quitte a faire un changement de variable $x=h(t)$ tel que $F_0(h(t))=t^p$ pour
$F_0(x)=\int_0^x f(t,0)\,dt$, 
on peut supposer que $f(x,0)=p x^{p-1}$. L'hypothèse sur 
$f$ assure qu'on peut faire ce changement de variable globalement sur
un certain voisinage de $[a,b]$. 
Quitte à faire un changement de variable $z=e^{G(x,\eps)}\~z$ avec 
$G'(x,\eps)=\frac1\eps\big(f(x,\eps)-px^{p-1}\big)$, on peut supposer que
$f(x,\eps)=px^{p-1}$.
Sans perte en généralité, nous pouvons supposer l'existence de $\~a,\~b\in{\cal U}$  
avec $\~a<a<b<\~b$ et supposer $\~a^p<\~b^p$.

Il convient d'introduire le {\em wronskien} de deux solutions $z_1,z_2$ de \rf{eao}.
On pose 
$$[z_1,z_2](\eps)=e^{-x^p/\eps}(z_1z_2'-z_1'z_2)(x,\eps).$$ 
D'après le théorème
de Liouville, il s'agit d'une fonction de $\eps$ seulement. Le wronskien est
une fonction bilinéaire et antisymétrique. On vérifie, pour 
plusieurs solutions $z_i$ de \rf{eao}, les formules
\begin{eqnarray}\lb{wronskia} 
[z_1,z_2]z_3+[z_2,z_3]z_1+[z_3,z_1]z_2&=&0 ,\\{}\lb{wronskib}
[z_1,z_2][z_3,z_4]-[z_1,z_3][z_2,z_4]+[z_1,z_4][z_2,z_3]&=&0.
\end{eqnarray}

Supposons maintenant les hypothèses du théorème satisfaites et considérons les solutions $z_a,z_b$ de \rf{eao} déterminées par
$$
z_a(\~a,\eps)=1,\,z_a'(\~a,\eps)=0,\ \ \ z_b(\~b,\eps)=1,\,
z_b'(\~b,\eps)=0.
$$
Il est connu que ces solutions ont des développements asymptotiques extérieurs
sur $[\~a+d,-d]$, respectivement $[d,\~b-d]$, où $d>0$ est arbitrairement petit et qu'elles
sont prolongeables sur tout ${\cal U}$, mais les développements des prolongements 
peuvent être différents 
ou ne pas exister.

Nous allons montrer que $[z_a,z_b](\eps)$ est exponentiellement petite~;  précisément nous montrerons qu'il existe $C>a^p$ telle que
\eq{zab}{
[z_a,z_b](\eps)={\cal O}\big(e^{-C/\eps}\big)\mbox{  quand  }\eps\to0. 
}
Ensuite, nous montrerons que cette propriété implique que $z_b$ est une solution résonante
globale.

Pour cela, nous introduisons d'autres solutions et nous procédons de 
manière analogue à la preuve de \rf{expo71}.
D'abord on utilise l'existence de $\mu\in\R$, et pour tout $\g>0$ d'un entier 
$L>\frac\pi{4p\g}$ et de $\e_0,\rho>0$ assez petits, tels qu'il existe des solutions 
$z_\ell^\pm(x,\e),\ \ell=-L,...,L$ de \rf{eao}  holomorphes quand
$\e\in S_\ell=S\big(\ell\frac\pi{2pL}-\g,\ell\frac\pi{2pL}+\g,\e_0\big)$ et 
$x\in V_\ell^\pm=\pm V\big(-\frac{3\pi}{2p}+\ell\frac\pi{2pL} + 2\g,
     \frac{3\pi}{2p}+\ell\frac\pi{2pL}-2\g,\rho,\mu\norm\e \big)$ 
ayant des \dacs{} généralisés
\eq{dacz-glob}{
z^\pm_\ell(x,\e)\sim_{\usp} Z_0^\pm\big(\tfrac x\e\big)
 \,e^{\hat R(\eps)\log(x^p+\eps)}\bigg(B_0(x)+\sum_{n=1}^\infty \gk{B_n(x)+
H_n^\pm\big(\tfrac x\e\big)}\e^n\bigg)
}
avec les fonctions $Z_0^\pm$, $B_n$ et $H_n^\pm$ de la proposition \reff{p4.7} ; les $B_n$ sont
holomorphes bornées dans $D(0,\rho)$, les $H_n^\pm$ dans  
$\pm V\big(-\frac{3\pi}{2p} + 2\g, \frac{3\pi}{2p}-2\g,\infty,\mu \big)$ et admettent
des développements asymptotiques sans termes constants quand $X\to\pm\infty$.
Ceci est démontré de manière analogue à la proposition \reff{p4.7} en passant par
l'équation de Riccati \rf{ric} et en utilisant les rotations 
$\e=e^{i\psi}\~\e$, $x=\pm e^{i\psi}\~x$ avec $\psi=\ell\frac\pi{2pL}$ ; le fait que
les fonctions $a_n$ et $g_n^\pm$ soient indépendantes de $\ell$ est dû à l'unicité
de la solution formelle de \rf{ric}, \ie au théorème \reff{th4.1}, et à la construction
des solutions de la proposition \reff{p4.7} à partir des solutions de \rf{ric}.
Notons que \rf{dacz-glob} entraîne, pour $\e\to0$ 
\eq{extint-z}{
\begin{array}{l}
z^\pm_\ell(x,\e)=\e^{-D}\gk{x^DB_0(x)+\O(\e)}\mbox{ pour }x\in\pm V_\ell(\e),\norm x>d,\med\\
               z^\pm_\ell(x,\e)= Z_0^\pm\big(\tfrac x\e\big)+\O(\e)\mbox{ pour }
                \norm x\leq K\norm\e  .
\end{array}}
pour tout $K>0$. Comme dans la preuve du théorème \reff{reso}, la dernière relation est une conséquence de \rf{dacz-glob} et du fait que l'équation intérieure  est régulièrement perturbée.

Comme nous l'avons vu dans la preuve du théorème \reff{reso}, la condition (a) est 
équivalente aux conditions $Z_0^+=Z_0^-=:Z_0$ et $H_n^+=H_n^-=:H_n$ pour tout $n$.
Puisque le développement asymptotique est Gevrey d'ordre $\usp$ en $\e$, ceci entraîne que les wronskiens
$[z_\ell^+,z_\ell^-]$ sont exponentiellement petits. Il existe donc $s>0$ tel que
\eq{expoz-glob}{
[z_\ell^+,z_\ell^-](\e)={\cal O}\big(e^{-s/\norm\e^p}\big)
}
quand $\ell\in\{-L,...,L\}$ et $\e\in S_\ell$. Quitte à réduire $\rho$, on peut supposer
que $\rho^p<s$.

On a aussi besoin de solutions avec un comportement asymptotique complémentaire
aux fonctions $z_\ell^\pm$. Pour les construire, on effectue le
changement de variable $z=e^{x^p/\eps}w$ qui transforme \rf{eao} en
\eq{eao2}{
\eps w'' + p x^{p-1}w'+\~g(x,\eps)w=0\mbox{\ \  où\ \ }
                 \~g(x,\eps)=g(x,\eps)+p(p-1)x^{p-2}\ .
                 }
Cette équation est très similaire à \rf{eao}, mais le signe du coefficient de la dérivée
a changé. On la traite de manière analogue à \rf{eao} 
-- la seule différence est que les montagnes et les vallées
de l'équation de Riccati 
\eq{ric2}{\eps u' = - px^{p-1}u - \~g(x,\eps) - u^2}
satisfaite par $u=\eps w'/w$ ont été interchangées par rapport à celles  de \rf{ric}. 
On obtient ainsi, avec $\mu$, $\g$, $L$, $\e_0$ et $\rho$ comme avant, des solutions 
$v_\ell^\pm(x,\e),\ \ell=-L,...,L$ de \rf{eao} définies, holomorphes et bornées quand
$\e\in S_\ell=S\big(\ell\frac\pi{2pL}-\g,\ell\frac\pi{2pL}+\g,\e_0\big)$ et 
$x\in \~V_\ell^\pm=\pm V\big(-\frac{\pi}{2p}+\ell\frac\pi{2pL} + 2\g,
     \frac{5\pi}{2p}+\ell\frac\pi{2pL}-2\g,\rho,\mu\norm\e \big)$ 
ayant des \dacs{} généralisés
\eq{dacz-glob2}{
v^\pm_\ell(x,\e)\sim_{\usp} 
\~Z_0^\pm\big(\tfrac x\e\big)
 \,e^{-\hat R(\eps)\log(x^p+\eps)}\bigg(\~B_0(x)+\sum_{n=1}^\infty \gk{\~B_n(x)+
\~H_n^\pm\big(\tfrac x\e\big)}\e^n\bigg)
}
où $\~Z_0^\pm$ sont les uniques solutions de \rf{red} avec
$\~Z_0^\pm(X)\sim e^{X^p}X^{-D-p+1}\big(1+\O(X^{-1})\big)$ quand $\arg(\pm X)=\frac\pi p$,
$\norm X\to\infty$, les $\~B_n$ sont analytiques sur $D(0,\rho)$, $\~B_0$ ne 
s'annule pas, $\~B_0(0)=1$ et enfin les $\~H_n$ sont holomorphes dans
$\pm V\big(-\frac{\pi}{2p} + 2\g, \frac{5\pi}{2p}-2\g,\infty,\mu \big)$ et admettent
des développements asymptotiques sans termes constants quand $X\to\pm\infty$.
Notons encore que \rf{dacz-glob2} entraîne, comme pour les $z_\ell^\pm$, quand $\e\to0$ 
\eq{extint-v}{
\begin{array}{l}
v_\ell^\pm(x,\e)=e^{x^p/\e^p}\gk{(\tfrac x\e)^{-D-p+1}\~B_0(x)+\O(\e)}\mbox{ pour }
       x\in \~V_\ell^\pm(\e), \norm x\geq d,\\
v_\ell^\pm(x,\e)=\~Z_0^\pm\big(\tfrac x\e\big)+\O(\e)\mbox{ pour }\norm x \leq K\norm\e
\end{array} 
}
ainsi que $v_0^\pm(x,\e)=\O\big(e^{x^p/\e^p}\big)$
quand $\e>0$ et $\pm x\in [0,\rho[$ ; le fait que l'axe réel négatif descende
dans une vallée  de \rf{ric2} implique que ceci reste vrai pour $v_\ell^-$ et tout
l'intervalle $[a,0]$.

En utilisant les développements extérieurs de $z_{a|b}$, les estimations
\rf{extint-z} et \rf{extint-v}, et la définition des wronskiens, on obtient 
\eq{nonul}{[z_a,v_\ell^-](\e)\asymp \e^{D-1},\ [z_b,v_\ell^+](\e)\asymp \e^{D-1},\ 
   [z_\ell^\pm,v_\ell^\pm](\e)\asymp \e^{-1}}
ainsi que $[v_\ell^+,v_\ell^-](\e)=\O(\e^{-1})$, où le symbole $f(\e)\asymp g(\e)$ signifie ici que le quotient $f(\e)/g(\e)$ tend vers une limite
non nulle quand $\e\to0$.
 
La démonstration de \rf{zab} utilise le théorème de Phragmén-Lindelöf et
suit celle de \rf{expo71}. Étant donné $\d,d>0$ petit, 
notons $D^-(\d,d)$ le domaine contenant $]\~a,-d]$ dont l'image par $F(x)=x^p$ est 
le triangle de sommets $F(\~a),iF(\~a)\tan\d,-iF(\~a)\tan\d$ privé
du disque centré en $0$ de rayon $d$.
De même, soit  $D^+(\d,d)$ le domaine contenant $[d,\~b[$ dont l'image par $F$ est 
le triangle de sommets $F(\~b),iF(\~b)\tan\d,-iF(\~b)\tan\d$ privé du disque centré en 
$0$ de rayon $d$.
Si on choisit $\d$  assez petit, pour $\~\d>\d$, $\~\d$ arbitrairement proche de $\d$ et pour  $\e_0,d$ assez petit, les solutions $z_a(x,\e)$, resp.\ $z_b(x,\e)$,  
sont holomorphes bornées dans le secteur $S=S\big(-\frac\pi{2p}+\frac{\~\d}p,
 \frac\pi{2p}-\frac{\~\d}p,\e_0\big)$ et pour $x$ dans 
$D^-(\d,d)$, resp. $D^+(\d,d)$. 
Ceci est classique.
\figu{f7.5}{
\vspace{-12mm}
\epsfxsize5cm\epsfbox{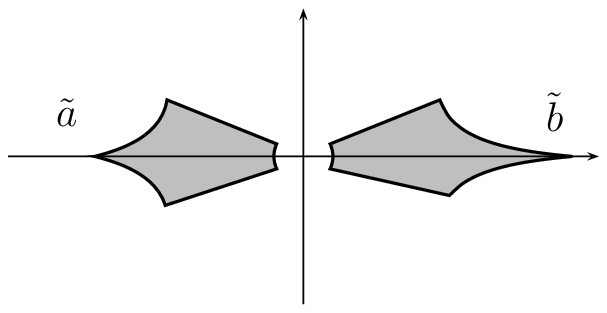}\hspace*{1cm}
\epsfxsize5cm\epsfbox{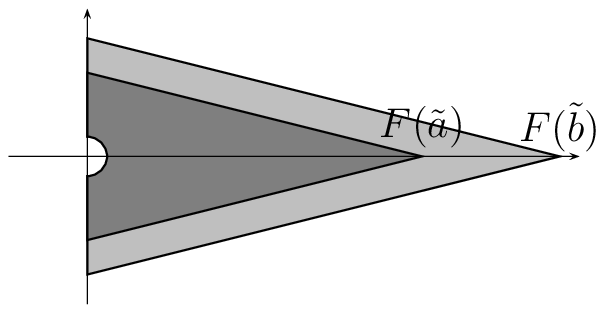}
\vspace{-18mm}
}
{Les domaines $D^-(\d,d)$, $D^+(\d,d)$ et leurs images par $F(x)=x^4$. Ici encore l'échelle n'est pas respectée entre les deux dessins pour une meilleure visibilité
}

Considérons maintenant $\~\d>0$ tel que $F(\~a)\sin\~\d<\rho^p$ et $\d\in\,]0,\~\d[$ arbitraire.
Soit $\ell\in\{-L,...,L\}$ et $\e\in S_\ell\cap S$, $\arg\e=\psi$. Alors $z_b(x,\e)$ et 
$z_\ell^+(x,\e)$ sont holomorphes bornées sur un voisinage du point 
$x = Te^{i\psi}$ où $T^p=F(\~a)\sin\d$.
On en déduit avec \rf{extint-z} que
\eq{expob+}{
[z_b,z_\ell^+](\e)={\cal O}\gk{\norm\e^{-D}\exp\big(-(Te^{i\psi})^p/\e^p\big)}
 ={\cal O}\gk{\norm\e^{-D} e^{-T^p/\norm\e^p}}
}
quand $\e\in S_\ell\cap S$.
De manière analogue, on montre que $[z_a,z_\ell^-](\e)=
  {\cal O}\gk{\norm\e^{-D}e^{-T^p/\norm\e^p}}$ quand $\e\in S_\ell\cap S.$ 

La formule \rf{wronskib} appliquée à $z_a,z_\ell^-,z_\ell^+$ et $v_\ell^-$
s'écrit
$$
[z_a,z_\ell^-][z_\ell^+,v_\ell^-]-[z_a,z_\ell^+][z_\ell^-,v_\ell^-]+
   [z_a,v_\ell^-][z_\ell^+,z_\ell^-]=0\ .	
$$
Avec \rf{expoz-glob}, la majoration de $[z_a,z_\ell^-]$ ci-dessus et \rf{nonul}, ceci implique que 
le wronskien $[z_a,z_\ell^+]$ satisfait  $[z_a,z_\ell^+](\e)={\cal O}\gk{\norm\e^{-D}e^{-T^p/\norm\e^p}}$.
En appliquant \rf{wronskib} à $z_a,z_\ell^-,v_\ell^+$ et $v_\ell^-$, on obtient 
$[z_a,v_\ell^+](\e)\asymp\e^{D-1}$ ; de manière analogue, on obtient
\eq{asympb-}{[z_b,v_\ell^-](\e)\asymp\e^{D-1}.}

En appliquant \rf{wronskib} cette fois à $z_a,z_\ell^+,z_b$ et $v_\ell^+$
et en utilisant \rf{expob+}, on obtient
$[z_a,z_b](\e)={\cal O}\gk{e^{-T^p/\norm\e^p}}$ quand $\e\in S_\ell\cap S$.
Comme l'union des $S_\ell$ couvre $S$, ceci donne
$$
[z_a,z_b](\e)={\cal O}\big(e^{-T^p/\norm\e^p}\big)
$$
pour $\e\in S$.

Pour $p|\arg\e|=\frac\pi2-\~\d$, cette majoration  implique
$[z_a,z_b](\e)=\O\Big(\exp\big(-q\frac{ F(\~a)}{\e^p}\big)\Big)$ avec 
$q=\frac{\sin\d }{\sin\~\d}$ arbitrairement proche de $1$. 
D'après le théorème de Phragmén-Lindelöf, ceci reste valide pour $\arg\e=0$. 
Ceci démontre enfin \rf{zab}.

La formule \rf{wronskia} appliquée à $z_b,z_0^+$ et $v_0^+$ donne
$$
[z_b,z_0^+]v_0^++[z_0^+,v_0^+]z_b+[v_0^+,z_b]z_0^+=0.
$$
Avec la majoration \rf{expob+}, les estimations \rf{nonul} et
avec $v_0^+(x,\e)=\O\big(e^{x^p/\e^p}\big)$ pour $x\in [0,d]$, $\e>0$,
ceci montre que $z_b$ est exponentiellement proche d'une solution proportionnelle
à $\e^Dz_0^+$ sur
$[0,d]$ et tend donc vers une limite aussi sur cet intervalle, si $d<T$.
Le coefficient de proportionnalité, $\frac{[z_b,v_0^+]}{[z_0^+,v_0^+]}\,\e^{-D}$, dépend de $\e$ mais tend vers une limite non nulle lorsque $\e$ tend vers $0$.
En utilisant \rf{asympb-}, on traite de la même manière $z_b$,  $z_0^-$ et $v_0^-$ sur $[-d,0].$ 
Enfin \rf{wronskia} appliquée à $z_b,z_a$ et $v_0^-$ donne
$$
[z_b,z_a]v_0^-+[z_a,v_0^-]z_b+[v_0^-,z_b]z_a=0.
$$
La majoration \rf{zab} entraîne alors, avec \rf{nonul} et $v_0^-(x,\e)=\O\gk{e^{x^p/\e^p}}$ sur $[a,0]$,
que $z_b$
reste exponentiellement proche d'une solution proportionnelle à $z_a$ sur  $[a,-d]$ aussi~; 
c'est donc une solution résonante globale.

Si la condition du théorème \reff{reso} (b) est satisfaite, alors celle de 
\reff{reso} (a) aussi. La solution $z_b$ est donc une solution résonante globale
et il ne reste plus qu'à montrer que toutes ses dérivées restent bornées.
D'après la preuve du théorème \reff{reso} (b), les dérivées de
$\e^Dz_\ell^\pm(x,\e)$ restent bornées sur $\pm[0,d]$~; puisque $z_b$
est exponentiellement proche d'elles, ceci reste vrai pour $z_b$.
Enfin à cause de l'équation différentielle, toutes les dérivées de $z_a$
restent bornées sur $[a,-d]$~; puisque $z_b$ est exponentiellement proche
d'elle, il suit que $z_b$ est une solution $\cc^\infty$-résonante globale
sur tout $[a,b]$.
\ep
%
%
%
\sec{8.}{Remarques historiques}
La littérature concernant le matching, \ie {\em the method of matched asymptotic 
expansions} est abondante. C'est le sujet principal du livre général de 
W.~Eckhaus \cit{e} et des chapitres VII et VIII du livre de  W.~Wasow \cit{wa1}.
Dans ce dernier livre, la méthode est présentée pour des systèmes linéaires d'équations 
différentielles singulièrement perturbées.

Le livre de Wasow présente aussi une méthode alternative au matching et aux \dacs{}
pour les systèmes de dimensions 2 : celle de la simplification uniforme en une
équation de fonctions spéciales. Cette méthode est aussi le sujet du mémoire
de Y.~Sibuya \cit{Ured}. Sous certaines conditions, l'équation  $\eps y'=A(x,\eps)y $
en $\C^2$ avec $A(x,0)=\left(\begin{array}{cc}0&1\\x^2/4&0\end{array}\right)$
peut être réduite à $\eps z'= \left(\begin{array}{cc}0&1\\\frac{x^2}4+\eps c(\eps)&0
\end{array}\right)z$ par une transformation $y=T(x,\eps)z$,
où $T(x,\eps)$ admet un développement asymptotique quand $\eps\to0$ uniformément
dans un voisinage de $x=0$. Comme l'équation réduite est équivalente à
$\eps^2u''=\gk{\frac{x^2}4+\eps c(\eps)}u$ et peut donc être ramenée à 
$\frac{d^2Z}{d\,X^2}=\gk{\frac{X^2}4+c(\eps)}Z$ par $Z(X,\eps^{1/2})=
z(\frac x{\eps^{1/2}},\eps)$, on pourrait voir cette simplifications uniforme comme
une séparation de l'asymptotique en une partie lente $T(x,\eps)$ et une
partie rapide venant de l'équation réduite. Il serait intéressant de comparer cette 
approche en détails avec le \dacs{}.

Pour décrire des couches limites, des \dacs{} dans le cas régulier (comme dans le 
paragraphe \reff{5.1}, mais sans étudier leur caractère Gevrey) étaient étudiés
dans plusieurs travaux. Nous mentionnons Vasil'eva-Butuzov \cit{vb} et 
Benoît-El Hamidi-Fruchard \cit{bef}.

Nous avons trouvé peu d'articles contenant des analogues de \dacs{} au voisinage
de points tournants avec une importante exception : la thèse de Th.~Forget \cit{fo}.
Il utilise une variante des \dacs{} pour obtenir des approximations uniformes
de solutions (et valeurs) canards d'équations analogues à \rf{1} ; ceci a été mentionné
à plusieurs endroits dans notre mémoire (parties \reff{e1.7}, \reff{1.3} 
et \reff{6.1}).

Parmi l'abondante littérature concernant les \dacs{} classiques et le matching, outre les travaux de W.~Eckhaus \cit{e} et W.~Wasow \cit{wa,wa1} déjà cités, on peut mentionner les travaux de L.~A.~Skinner \cit{s1,s2,s3} et de L.~E.~Fraenkel \cit{f}. 
Chacun de ces auteurs a mis en place un formalisme de ``composite asymptotic expansions''. Par exemple L.~E.~Fraenkel définit des opérateurs $E_p$ et $H_p$ qui, dans notre contexte, consistent essentiellement à  prendre les $p$ premiers termes du développement extérieur, resp. intérieur, d'une fonction $f=f(x,\eps)$. Il montre ensuite que la fonction $(E_p+H_p-E_pH_p)f$ est une approximation à  l'ordre $p$ de $f$ en $\eps$ uniformément pour $x\in[0,r]$. Si l'on répartit judicieusement les termes de $E_pH_pf$ pour en regrouper certains avec $E_pf$ et d'autres avec $H_pf$, le résultat semble être les premiers termes d'un \dac{} pour $f$. 
Toutefois L.~E.~Fraenkel ne fait pas une étude systématique de ces opérateurs. Il démontre cette approximation uniforme en faisant l'hypothèse \apriori{} que $f$ admet une approximation intérieure et une approximation extérieure et que les régions de validité s'intersectent (``overlapping assumption''). Il ne montre pas la stabilité de ces développement par produit, dérivation, intégration ou composition à  gauche ou à  droite par les fonctions d'une variable. Par ailleurs, seules des sommes finies sont prises en compte. Enfin l'aspect Gevrey, indispensable pour les applications à  la perturbation singulière, est absent. 

Mentionnons encore quelques travaux liées à notre mémoire.
Le lien avec les travaux d'Éric Matzinger \cit{ma,ma1,ma2} a été détaillé 
dans la partie \reff{5.3}. 
Le lien entre un \dac{} sur une couronne en $X=x/\e$ et l'asymptotique monomiale 
de \cit{cms} est détaillé dans la première remarque après
la définition \reff{defcomb}. Le lien avec les travaux de Peter De Maesschalck a été présenté dans la remarque 1 après le théorème \reff{t6.1} et dans la remarque après le théorème \reff{resog}.

\noi Adresses des auteurs :\med\\ 
Augustin Fruchard\\
Laboratoire de Mathématiques, Informatique et Applications\\
Faculté des Sciences et Techniques,
Université de Haute Alsace\\
4, rue des Frères Lumière,
F-68093 Mulhouse cedex,
France\med\\
Courriel : augustin.fruchard@uha.fr
\bigskip\\ Reinhard Schäfke\\ Institut de Recherche Mathématique Avancée\\ U.F.R.\ de Mathématiques et Informatique \\ Université Louis Pasteur et C.N.R.S.\\ 7, rue René Descartes, 67084 Strasbourg cedex, France\med\\ Courriel : schaefke@unistra.fr
\end{document}